\documentclass[11pt,fullpage]{article}

\usepackage{amsmath}
\usepackage{amssymb}
\usepackage{amsfonts}
\usepackage{xr}

\numberwithin{equation}{section}
\newtheorem{thm}{Theorem}[section] 
\newtheorem{prp}[thm]{Proposition}
\newtheorem{lmm}[thm]{Lemma}   
\newtheorem{crl}[thm]{Corollary} 
\newtheorem{dfn}[thm]{Definition}

\def\e_ref#1{(\ref{#1})}
\def\under#1{\underline{#1}}
\def\ov#1{\overline{#1}}
\def\lra{\longrightarrow}
\def\Lra{\Longrightarrow}

\def\l{\left}
\def\r{\right}
\def\lan{\langle}
\def\ran{\rangle}
\def\lr#1{\lan{#1}\ran}
\def\llan{\langle\langle}
\def\rran{\rangle\rangle}

\def\al{\alpha}
\def\be{\beta}
\def\de{\delta}
\def\ep{\epsilon}
\def\ga{\gamma}
\def\io{\iota}
\def\ka{\kappa}
\def\la{\lambda}
\def\na{\nabla}
\def\om{\omega}
\def\si{\sigma}
\def\th{\theta}
\def\ve{\varepsilon}
\def\ups{\upsilon}
\def\vp{\varpi}
\def\vt{\vartheta}

\def\Ga{\Gamma}
\def\La{\Lambda}
\def\Om{\Omega}
\def\Si{\Sigma}

\def\rk{\text{rk}}
\def\ev{\text{ev}}
\def\coker{\text{coker}}
\def\codim{\text{codim}}
\def\RT{\text{RT}}

\def\P{\Bbb{P}^n}
\def\PP{\Bbb{P}^2}
\def\PPP{\Bbb{P}^3}
\def\i{\infty}
\def\eset{\emptyset}
\def\w{\wedge}

\begin{document}

\title{Enumeration of Genus-Two Curves with a Fixed Complex Structure
in $\PP$ and~$\PPP$}
\author{Aleksey Zinger
\thanks{Partially supported by NSF Graduate Research Fellowship
and NSF grant DMS-9803166}}
\date{\today}
\maketitle

\begin{abstract}
\noindent
We express the genus-two fixed-complex-structure enumerative 
invariants of $\PP$ and $\PPP$  in terms of the genus-zero  
enumerative invariants. 
The approach is to relate each genus-two fixed-complex-structure 
{\it enumerative} invariant to the corresponding 
{\it symplectic} invariant.
\end{abstract}

\thispagestyle{empty}

\tableofcontents

\section{\bf Introduction} 

\subsection{Background and Results}
\label{back}

\noindent
Let $(\Si,j_{\Si})$ be a nonsingular Riemann surface of genus $g\!\ge\! 2$,
and let $d,n$ be positive integers with $d\!\ge\! 1$ and $n\!\ge\! 2$. 
Denote by ${\cal H}_{\Si,d}(\P)$ the set of simple holomorphic maps 
from $\Si$ to $\P$ of degree~$d$. 
Let 
$\mu\!=\!(\mu_1,\ldots,\mu_N)$
be an $N$-tuple of proper complex submanifolds of $\Bbb{P}^n$ such~that 
\begin{equation}\label{c_codim}
\sum_{l=1}^{l=N}\hbox{codim}_{\Bbb{C}}\mu_l=
d(n+1)-n(g-1)+N.
\end{equation}
If these submanifolds are in general position, the cardinality of the set
\begin{equation}
\label{intro_1}
{\cal H}_{\Si,d}(\mu)=
\{(y_1,\ldots,y_N;u)\!: u\!\in\!{\cal H}_{\Si,d}(\P);~
y_l\!\in\!\Si,~u(y_l)\!\in\!\mu_l~\forall l\!=\!1,\ldots,N\}
\end{equation}
is finite and depends only on the homology classes of $\mu_1,\ldots,\mu_N$.
The group $\hbox{Aut}(\Si)$ of holomorphic automorphisms of $\Si$ acts 
freely on ${\cal H}_{\Si,d}(\mu)$.
For this reason, algebraic geometers prefer to consider 
the ratio of the cardinality of the set ${\cal H}_{\Si,d}(\mu)$
and the order of the group $\hbox{Aut}(\Si)$.
The cardinality of the set ${\cal H}_{\Si,d}(\mu)$
depends only on the genus of $\Si$.
Furthermore, for a dense open subset of complex structures on $\Si$,
the group of holomorphic automorphisms of $\Si$ has the same order.
If $j_{\Si}$ lies in this open subset, we denote 
the above ratio by $n_{g,d}(\mu)$. 
This number is precisely the number of 
irreducible, nodal degree-$d$ genus-$g$ curves in $\P$ 
with a  fixed generic complex structure on the normalization
and  passing through the constraints $\mu_1,\ldots,\mu_N$.\\

\noindent
For $g\!=\!0,1$, one can define the numbers $n_{g,d}(\mu)$
for constraints of appropriate total codimension by
counting the number of equivalence classes under the action of 
the now infinite group $\hbox{Aut}(\Si)$ on the set
${\cal H}_{\Si,d}(\mu)$ defined as in \e_ref{intro_1} above.
It is shown in \cite{RT} that 
$$n_{0,d}(\mu)=\RT_{0,d}(\mu_1,\mu_2,\mu_3;\mu_4,\ldots,\mu_N),$$
where $\RT_{0,d}(\cdot;\cdot)$ denotes the symplectic invariant of $\P$
as defined in~\cite{RT}.
For $g\!=\!1$, in \cite{I} the~difference 
$$\RT_{g,d}(\mu_1;\mu_2,\ldots,\mu_N)-2n_{g,d}(\mu)$$
is expressed as an intersection number on a blowup of the space
of degree-$d$ $(N+1)$-marked rational curves passing through
the constraints $\mu_1,\ldots,\mu_N$.
This number is shown to be computable, and 
explicit formulas are given in the $n\!=\!2,3$ cases.
On the other hand, the symplectic invariant is easily computable
from the two composition laws of~\cite{RT}. 
A completely algebraic approach for the $n\!=\!2$, $g\!=\!1$ case
is given in \cite{P1}.
Using this approach, \cite{KQR} express $n_{2,d}$ in the $n\!=\!2$ case
in terms of the numbers $n_{0,d'}$ with $d'\le d$.\\

\noindent
In this paper, we extend the approach of \cite{I} to compute
the difference 
$$\RT_{2,d}(\cdot;\mu_1,\ldots,\mu_N)-2n_{2,d}(\mu)$$
in the $n\!=\!2,3$ cases.
The reason for the factor of two above is that the automorphism
group of a generic genus-two Riemann surface has order two.
The following two theorems are the main results of this paper.
The two tables list some low-degree genus-two numbers.
Evidence in support of the two formulas is described in 
Subsection~\ref{final_form}, where more
low-degree numbers for $\PPP$ are also given.

\begin{thm}
\label{p2_final} 
Let $n_{2,d}$ denote the number of genus-two degree-$d$
curves that pass through $3d\!-\!2$ points 
in general position in $\PP$ and have a fixed generic complex
structure. With $n_d\!=\!n_{0,d}$,
$$ n_{2,d}=3(d^2-1)n_d+
\frac{1}{2}\sum\limits_{d_1+d_2=d}
\Big(d_1^2d_2^2+28-16\frac{9d_1d_2-1}{3d-2}\Big)
\binom{3d-2}{3d_1-1}d_1d_2n_{d_1}n_{d_2}.$$
\end{thm}

\begin{center}
\begin{tabular}{||c|c|c|c|c|c|c|c||}
\hline\hline
$d$&        1& 2& 3&    4&       5&        6& 7\\
\hline
$n_{2,d}$&  0& 0& 0& 14,400&  6,350,400& 3,931,128,000&
               3,718,909,209,600\\
\hline\hline
\end{tabular}
\end{center}

\vspace{.1in}

\begin{thm}
\label{p3_fin}
If $d$ is a positive integer and $\mu$ is a tuple  of $p$ points
and $q$ lines in general position in $\PPP$ with 
\hbox{$2p\!+\!q\!=\!4d\!-\!3$},
$$2n_{2,d}(\mu)=RT_{2,d}(\cdot;\mu)-CR(\mu),$$
where $CR(\mu)$ is the sum of the intersection numbers
of explicit tautological classes in the space of stable rational
maps into~$\PPP$.
\end{thm}

\begin{center}
\begin{tabular}{||c|c|c|c|c|c|c||}
\hline\hline
degree&         \multicolumn{3}{c|}{4}& \multicolumn{2}{c|}{5}& 6\\
\hline
(p,q)&   (3,7)& (2,9)& (1,11)& (8,1)& (0,17)& (10,1)\\
\hline
$n_{2,d}(\mu)$&   14,400& 307,200& 4,748,160& 9,600&7,494,574,433,280&
      1,301,760\\
\hline\hline
\end{tabular}
\end{center}

\vspace{.1in}

\noindent
A formula for $CR(\mu)$ is given in Theorem~\ref{p3_final}.
Intersection numbers of tautological classes
are shown to be computable in~\cite{P2}.
In fact, we give a method of computing these numbers
along the lines of that in~\cite{I}, which is slightly different
from the method of~\cite{P2}; see Subsection~\ref{chern_class}.\\

\noindent
The numbers we obtain in the $n\!=\!2$ case are different from the numbers
given in~\cite{KQR}.
However, our numbers can be recovered via the approach of~\cite{KQR}.
In particular,
$$n_{2,d}=6\big(n_{2,d}^{KQR}+\tau_d),$$
where $\tau_d$ is the number of degree-$d$ tacnodal
rational curves passing through $(3d\!-\!2)$ points in general
position in~$\PP$.
The factor of six is a minor omission on the authors' part.
The contribution of~$6\tau_d$ arises from a three-component strata
\cite{KQR} rule out by Remark~3.12,
which is stated without a proof.
Details can be found in~\cite{Z2}.
\\

\noindent
This paper combines the topological tools of Section~\ref{top_sect}
with the explicit analytic structure theorems of~\cite{Z}.
Together these give a general framework that will
hopefully provide a way of computing 
positive-genus enumerative invariants
from the symplectic ones in any homogeneous Kahler manifold.
In fact, the methods of this paper should also apply, 
with very little change, at least
up to genus-seven in~$\PP$, to the $g\!=\!3$ case~$\PPP$, 
and to the $g\!=\!2$ case in~$\Bbb{P}^4$.
Along the way, we enumerate cuspidal rational curves in~$\PP$ and 
two-component rational curves connected at a tacnode in~$\Bbb{P}^3$;
see Lemmas~\ref{n2cusps} and~\ref{n3tangents}.
Analogous methods can be used to enumerate rational
curves with higher-order singularities in any projective space
or, perhaps, even in any homogeneous Kahler manifold.\\

\noindent
The author is grateful to T.~Mrowka for pointing out the paper~\cite{I}
and many useful discussions, and  
G.~Tian  for first introducing him to Gromov's symplectic invariants.
The author also thanks R.~Vakil for sharing some of his
expertise in enumerative algebraic geometry, and
A.~J.~de~Jong and J.~Starr for help with understanding~\cite{KQR}.

\subsection{Summary}
\label{summary}

\noindent
If $\nu\!\in\!\Ga(\Si\!\times\!\P;
            \La^{0,1}\pi_{\Si}^*T^*\Si\otimes\pi_{\P}^*T\P)$, 
let ${\cal M}_{\Si,\nu,d}$ denote the set of all smooth maps~$u$ 
from $\Si$ to $\P$ of degree $d$ such that
$\bar{\partial}u|_z\!=\!\nu|_{(z,u(z))}$ for all $z\!\in\!\Si$.
If $\mu$ is as above, put
$${\cal M}_{\Si,\nu,d}(\mu)=
\big\{(y_1,\ldots,y_N;u)\!: u\!\in\!{\cal M}_{\Si,\nu,d};~
y_i\!\in\!\Si,~ u(y_i)\!\in\!\mu_i ~\forall i\!=\!1,\ldots,N\big\}.$$
For a generic $\nu$, ${\cal M}_{\Si,\nu,d}$ is a smooth
finite-dimensional oriented manifold, and 
${\cal M}_{\Si,\nu,d}(\mu)$
is a zero-dimensional finite submanifold of 
${\cal M}_{\Si,\nu,d}\!\times\!\Si^N$, 
whose cardinality (with sign) depends only the homology
classes of $\mu_1,\ldots,\mu_N$; see~\cite{RT}.
The symplectic invariant $\RT_{g,d}(;\mu)$ is the signed
cardinality of the set~${\cal M}_{\Si,\nu,d}(\mu)$.\\

\noindent
If $\|\nu_i\|_{C^0}\!\lra\! 0$ and 
   $(\under{y}_i,u_i)\!\in\! {\cal M}_{\Si,\nu_i,d}(\mu)$, 
then a subsequence of $\{(\under{y}_i,u_i)\}_{i=1}^{\i}$ must converge 
in the Gromov topology to one of the following:\\
(1) an element of ${\cal H}_{\Si,d}(\mu)$;\\
(2) $(\Si_T,\under{y},u)$, where $\Si_T$ is a bubble tree of 
$S^2$'s attached to $\Si$ with marked points $y_1,\ldots,y_N$,
and \hbox{$u\!:\Si_T\!\lra\!\P$} is a holomorphic map such that
$u(y_i)\!\in\!\mu_i$ for $i\!=\!1,\ldots,N$, and\\
(2a) $u|\Si_T$ is simple and the tree contains at least one $S^2$;\\
(2b) $u|\Si_T$ is multi-covered;\\
(2c) $u|\Si_T$ is constant and the tree contains at least one $S^2$.\\
By Proposition~\ref{empty_spaces}, the case (2a)
does not occur if the constraints are in general position.
Furthermore, if $g\!=\!2$, (2b) cannot occur either
if $n\!=\!2,3$ or if $n\!=\!4$ and $d\!\neq\! 2$.
It~is well-known that $n_{2,2}(\mu)\!=\!0$, 
and thus the case $n\!=\!4$ and $g\!=\!d\!=\!2$ presents no interest.
Our approach will be to take $t$ very small and to count the number
of elements of ${\cal M}_{\Si,t\nu,d}(\mu)$ that lie near 
the maps of type~(2c).
The rest of the elements of ${\cal M}_{\Si,t\nu,d}(\mu)$
must lie near the space ${\cal H}_{\Si,d}(\mu)$.
By Proposition~\ref{gl-holom_prp} in~\cite{Z} and Corollary~\ref{reg_crl},
there is a one-to-one correspondence between
the elements of ${\cal H}_{\Si,d}(\mu)$
and the nearby elements of ${\cal M}_{\Si,t\nu,d}(\mu)$,
at least if $d\!\ge\!3$.
If $d\!=\!1,2$, ${\cal H}_{\Si,d}(\mu)\!=\!\eset$;
see the proof of Proposition~\ref{empty_spaces}.
Thus, we are able to compute the cardinality of 
${\cal H}_{\Si,d}(\mu)$ by computing the total number
of elements of ${\cal M}_{\Si,t\nu,d}(\mu)$ that lie
near the maps of type~(2c).\\

\noindent
In Subsection~\ref{notation_sec}, we summarize our notation
for spaces of bubble maps and vector bundle over them. 
For details, the reader is referred to~\cite{Z}. 
In Section~\ref{analysis_sec}, 
we describe an obstruction-bundle setup and state Theorem~\ref{si_str}, 
which relates the elements of ${\cal M}_{\Si,d,t\nu}(\mu)$ lying near 
the maps of type~(2c) to the zero set of a map between two bundles.
We also describe the local structure of certain spaces of stable 
rational maps.
These spaces are very familiar in algebraic geometry,
but for our computations in Section~\ref{comp_sect}
we need the analytic estimate of Theorem~\ref{str_global}.\\

\noindent
In Section~\ref{resolvent_sec}, 
we use the topological tools of Subsection~\ref{top_sec1}
to show that the number of zeros of the maps of  Theorem~\ref{si_str}
is the same the number of zeros of explicit affine maps between 
vector bundles 
over  cartesian products of spaces of rational maps with~$\Si^k$.
The results of this simplification are summarized in 
Subsection~\ref{summary4}.
In Section~\ref{comp_sect}, we relate the zeros of
these affine maps to the intersection numbers of spaces of 
stable rational maps into~$\P$.
We use Theorem~\ref{str_global} and Section~\ref{top_sect}
for {\it local} excess-intersection type of computations.
If $n\!=\!2$, most excess contributions arise from points,
and we could instead use Theorem~\ref{str_global} to construct
blowups around such points.
If $n\!=\!3$, Theorem~\ref{str_global} does not provide 
enough justification for such a blowup construction,
since many excess contributions arise from positive-dimensional
compact manifolds.
We conclude our computations with
the very explicit formula of Theorem~\ref{p2_final}
in the $n\!=\!2$ case and a somewhat less explicit one
of Theorem~\ref{p3_final} in the $n\!=\!3$ case.

\subsection{Notation}
\label{notation_sec}

\noindent
In this subsection, we give a brief description
of the most important notation used in this paper.
See Section~\ref{gl-bubble_sect} in~\cite{Z} for more details.\\

\noindent
Let $q_N,q_S\!: \Bbb{C}\!\lra\! S^2\!\subset\!\Bbb{R}^3$ 
be the stereographic
projections mapping the origin in $\Bbb{C}$ to the north and 
south poles, respectively. Explicitly,
\begin{equation}
\label{stereo_pr}
q_N(z)=\Big(\frac{2z}{1+|z|^2},\frac{1-|z|^2}{1+|z|^2}\Big)
\in\Bbb{C}\!\times\!\Bbb{R},\quad
q_S(z)=\Big(\frac{2z}{1+|z|^2},\frac{-1+|z|^2}{1+|z|^2}\Big).
\end{equation}
We denote the south pole of $S^2$, i.e.~the point 
$(0,0,-1)\!\in\!\Bbb{R}^3$, by $\i$. Let
\begin{equation}
\label{special_vect}
e_{\i}=(0,0,1)=dq_S\Big|_0
   \Big(\frac{\partial}{\partial s}\Big)\in T_{\i}S^2,
\end{equation}
where we write $z=s+it\in\Bbb{C}$.
We identify $\Bbb{C}$ with  $S^2-\{\i\}$ via the map $q_N$.\\

\noindent
If $N$ is any nonnegative integer, let $[N]\!=\!\{1,\ldots,N\}$.
If $I_1$ and $I_2$ are two sets, we denote the disjoint
union of $I_1$ and $I_2$ by $I_1\!+\!I_2$.

\begin{dfn}
\label{rooted_tree}
A finite partially ordered set $I$ is a \under{linearly ordered set}
if for all \hbox{$i_1,i_2,h\!\in\! I$} such that $i_1,i_2\!<\!h$, 
either $i_1\!\le\! i_2$ \hbox{or $i_2\!\le\! i_1$.}\\
A linearly ordered set $I$ is a \under{rooted tree} if
$I$ has a unique minimal element, 
i.e.~there exists \hbox{$\hat{0}\!\in\! I$} such that $\hat{0}\!\le\! i$ 
for {all $i\!\in\! I$}.
\end{dfn}

\noindent
If $I$ is a linearly ordered set, let $\hat{I}$ be
the subset of the non-minimal elements of~$I$.
For every $h\!\in\!\hat{I}$,  denote by $\io_h\!\in\!I$
the largest element of $I$ which is smaller than~$h$.
We call $\io\!:\hat{I}\!\lra\!I$ the attaching map of~$I$.
Suppose $I\!=\!\bigsqcup\limits_{k\in K}\!I_k$
is the splitting of $I$ into rooted trees 
such that $k$ is the minimal element of~$I_k$.
If $\hat{1}\!\not\in\!I$, we define
the linearly ordered set $I\!+_k\!\hat{1}$ to
be the set $I\!+\!\hat{1}$ with all partial-order relations of~$I$
along with the relations
$$k<\hat{1},\qquad \hat{1}<h\hbox{~~if~~}h\!\in\!\hat{I}_k.$$
If $I$ is a rooted tree, we write $I\!+\!\hat{1}$ 
\hbox{for $I\!+_k\!\hat{1}$}.
\\

\noindent
If $S=\Si$ or $S=S^2$ and $M$ is a finite set, 
a {\it $\P$-valued bubble map with $M$-marked 
points} is a tuple
$$b=\big(S,M,I;x,(j,y),u\big),$$
where $I$ is a linearly ordered set, and
$$x\!:\hat{I}\!\lra\!S\cup S^2,~~j\!:M\!\lra\! I,~~ 
y\!:M\!\lra\!S\cup S^2,\hbox{~~and~~} 
u\!:I\!\lra\!C^{\i}(S;\P)\cup C^{\i}(S^2;\P)$$
are maps such that
$$x_h\in
\begin{cases}
S^2-\{\i\},&\hbox{if~}\io_h\!\in\!\hat{I};\\
S,&\hbox{if~}\io_h\!\not\in\!\hat{I},
\end{cases}\qquad
y_l\in
\begin{cases}
S^2-\{\i\},&\hbox{if~}j_l\!\in\!\hat{I};\\
S,&\hbox{if~}j_l\!\not\in\!\hat{I},\qquad
\end{cases}
u_i\in
\begin{cases}
C^{\i}(S^2;\P),&\hbox{if~}i\!\in\!\hat{I};\\
C^{\i}(S;\P),&\hbox{if~}i\!\not\in\!\hat{I},
\end{cases}$$  
and  $u_h(\i)\!=\!u_{\io_h}(x_h)$ for all 
$h\!\in\!\hat{I}$.
We associate such a tuple with Riemann surface
$$\Si_b=
\Big(\bigsqcup_{i\in I}\Si_{b,i}\Big)\Big/\!\sim,
\hbox{~~where}\qquad 
\Si_{b,i}=
\begin{cases}
\{i\}\times S^2,&\hbox{if~}i\!\in\!\hat{I};\\
\{i\}\times S,&\hbox{if~}i\!\not\in\!\hat{I},
\end{cases}
\quad\hbox{and}\quad
(h,\i)\sim (\io_h,x_h)
~~\forall h\!\in\!\hat{I},$$
with marked points $(j_l,y_l)\!\in\!\Si_{b,j_l}$,
and continuous map $u_b\!:\Si_b\!\lra\!\P$,
given by $u_b|\Si_{b,i}\!=\!u_i$ for \hbox{all $i\!\in\! I$}.
We require that all the singular points of $\Si_b$,
i.e.~$(\io_h,x_h)\!\in\!\Si_{b,\io_h}$ for $h\!\in\!\hat{I}$,
and all the marked points be distinct.
In addition, if $\Si_{b,i}\!=\!S^2$ and 
$u_{i*}[S^2]\!=\!0\!\in\! H_2(\P;\Bbb{Z})$,
then $\Si_{b,i}$ must contain at least two singular and/or
marked points of~$\Si_b$ other than~$(i,\i)$.
Two bubble maps $b$ and $b'$ are {\it equivalent} if there
exists a homeomorphism $\phi\!:\Si_b\!\lra\!\Si_{b'}$
such that $u_b\!=\!u_{b'}\!\circ\phi$,
$\phi(j_l,y_l)\!=\!(j_l',y_l')$ for \hbox{all $l\!\in\! M$},
$\phi|_{\Si_{b,i}}$ is holomorphic for \hbox{all $i\!\in\!I$},
and $\phi|_{\Si_{b,i}}\!=\!Id$ if $S\!=\!\Si$ 
\hbox{and $i\!\in\! I\!-\!\hat{I}$.}
\\

\noindent
The general structure of bubble maps is described
by tuples ${\cal T}\!=\!(S,M,I;j,d)$,
with \hbox{$d_i\!\in\!\Bbb{Z}$} 
describing the degree of the map $u_b$ on~$\Si_{b,i}$.
We call such tuples {\it bubble types}.
Bubble type ${\cal T}$ is {\it simple} if  $I$ is a rooted tree;
${\cal T}$ is {\it is basic} if $\hat{I}\!=\!\eset$;
${\cal T}$ is {\it semiprimitive} if 
$\io_h\!\not\in\!\hat{I}$ for all $h\!\in\!\hat{I}$.
We call semiprimitive bubble type ${\cal T}$
{\it primitive} if $j_l\!\in\!\hat{I}$ for \hbox{all $j_l\!\in\! M$}.
The above equivalence relation on the set of bubble maps
induces an equivalence relation on the set of bubble types.
For each $h,i\!\in\! I$, let
\begin{gather*}
D_i{\cal T}=\{h\!\in\!\hat{I}\!:i\!<\!h\},\quad
\bar{D}_i{\cal T}=D_i{\cal T}\cup\{i\},\quad
H_i{\cal T}=\{h\!\in\!\hat{I}\!:\io_h\!=\!i\},\quad
M_i{\cal T}=\{l\!\in\! M\!:j_l\!=\!i\},\\
\chi_{\cal T}h=\begin{cases}
0,&\hbox{if~} \forall i\!\in\! 
             I\hbox{~s.t.~}h\!\in\!\bar{D}_i{\cal T}, d_i\!=\!0;\\
1,&\hbox{if~} d_h\!\neq\! 0,\hbox{~but~}
\forall i\!\in\! I\hbox{~s.t.~}h\!\in\! D_i{\cal T}, d_i\!=\!0;\\
2,&\hbox{otherwise}.
\end{cases}
\end{gather*}
Let ${\cal H}_{\cal T}$ denote the space of all holomorphic
bubble maps with structure~${\cal T}$.
\\

\noindent
If $S\!=\!\Si$, we denote by ${\cal M}_{\cal T}$ the set of
equivalence classes of bubble maps in~${\cal H}_{\cal T}$.
Then there exists ${\cal M}_{\cal T}^{(0)}\!\subset\!{\cal H}_{\cal T}$
such that ${\cal M}_{\cal T}$ is the quotient of ${\cal M}_{\cal T}^{(0)}$ 
by an~$(S^1)^{\hat{I}}\!\times\!\hbox{Aut}({\cal T})$-action.
If the group $\hbox{Aut}({\cal T})$ is trivial,
corresponding to this action, we obtain 
$|\hat{I}|$ line orbi-bundles 
$\{L_h{\cal T}\!\lra\!{\cal M}_{\cal T}\!:h\!\in\!\hat{I}\}$.
If $\hbox{Aut}({\cal T})$ is not trivial,
the fiber product and the direct sum of 
the bundles $L_h{\cal T}$ taken over any orbit
of the $\hbox{Aut}({\cal T})$ are well-defined.
The bundle of gluing parameters in this case~is
$$F{\cal T}=\bigoplus_{h\in\hat{I}}F_h{\cal T},
\hbox{~~where}\qquad
F_{h,[b]}{\cal T}=\begin{cases}
L_{h,[b]}{\cal T}\otimes L_{\io_h,[b]}^*{\cal T},&
\hbox{if~}\io_h\!\in\!\hat{I};\\
L_{h,[b]}{\cal T}\otimes T_{x_h}\Si,&
\hbox{if~}\io_h\!\not\in\!\hat{I}.
\end{cases}$$
Let
$F^{\eset}{\cal T}\!=\!\{\ups=(\ups_h)_{h\in{\hat{I}}}\!\in\! F{\cal T}\!: 
\ups_h\!\neq\!0~\forall h\!\in\!\hat{I}\}$.
Each line orbi-bundle $F_h{\cal T}\!\lra\!{\cal M}_{\cal T}$ is 
the quotient of 
a line bundle  $F_h^{(0)}{\cal T}\!\lra\!{\cal M}_{\cal T}^{(0)}$ by a 
\hbox{$G_{\cal T}\!
  \equiv\!(S^1)^{\hat{I}}\!\times\!\hbox{Aut}({\cal T})$}-action.
We denote by $F^{(\eset)}{\cal T}$ the preimage of
$F^{\eset}{\cal T}$ in 
$F^{(0)}{\cal T}\!\equiv\!\bigoplus\limits_{h\in\hat{I}}F_h^{(0)}{\cal T}$.
The bundles $F^{\eset}{\cal T}$, $F^{(\eset)}{\cal T}$, 
and $F_h^{(0)}{\cal T}$ are defined even if 
the automorphism group of ${\cal T}$ is nontrivial.
\\

\noindent
For each bubble type ${\cal T}=(S^2,M,I;j,d)$,
let
$${\cal U}_{\cal T}=\big\{[b]\!:
b\!=\!\big(S^2,M,I;x,(j,y),u\big)\!\in\!{\cal H}_{\cal T},~
u_{i_1}(\i)=u_{i_2}(\i)~\forall i_1,i_2\!\in\! I\!-\!\hat{I}
\big\}.$$
Similarly to the $S\!=\!\Si$ case above,
${\cal U}_{\cal T}$ is the quotient of a subset ${\cal B}_{\cal T}$
of ${\cal H}_{\cal T}$ by a
\hbox{$\tilde{G}_{\cal T}\!\equiv\!
         (S^1)^I\!\times\!\hbox{Aut}({\cal T})$}-action.
Denote by ${\cal U}_{\cal T}^{(0)}$ the quotient of ${\cal B}_{\cal T}$
by 
$$G_{\cal T}\!\equiv\!(S^1)^{\hat{I}}\!\times\!\hbox{Aut}({\cal T})
\subset\tilde{G}_{\cal T}.$$
Then ${\cal U}_{\cal T}$ is the quotient of
${\cal U}_{\cal T}^{(0)}$ the residual 
$G_{\cal T}^*\!\equiv\!(S^1)^{I-\hat{I}}\!\subset\!\tilde{G}_{\cal T}$
action.
Corresponding to these quotients, we obtain line orbi-bundles
$\{L_h{\cal T}\!\!\lra\!{\cal U}_{\cal T}^{(0)}\!\!:h\!\in\!\hat{I}\}$
and $\{L_i{\cal T}\!\!\lra\!{\cal U}_{\cal T}\!\!:i\!\in\! I\}$
if the automorphism group of ${\cal T}$ is trivial.
Let 
\begin{gather*}
F{\cal T}=\bigoplus_{h\in\hat{I}}F_h{\cal T}\lra{\cal U}_{\cal T}^{(0)},
\hbox{~~where}\quad
F_{h,[b]}{\cal T}=\begin{cases}
L_{h,[b]}{\cal T}\otimes L_{\io_h,[b]}^*{\cal T},&
\hbox{if~}\io_h\!\in\!\hat{I};\\
L_{h,[b]}{\cal T},&
\hbox{if~}\io_h\!\not\in\!\hat{I};
\end{cases}\\
{\cal FT}=\bigoplus_{h\in\hat{I}}{\cal F}_h{\cal T}\lra{\cal U}_{\cal T},
\hbox{~~where}\quad
{\cal F}_{h,[b]}{\cal T}=L_{h,[b]}{\cal T}\otimes L_{\io_h,[b]}^*{\cal T}.
\end{gather*}
The orbi-bundles  $F_h{\cal T}$ and ${\cal F}_i{\cal T}$ are
quotients of line bundles over ${\cal B}_{\cal T}$
similarly to the \hbox{$S\!=\!\Si$ case}.\\

\noindent
Gromov topology on the space of equivalence classes of bubble maps
induces a partial ordering on the set of bubble types
and their equivalence classes such that the spaces
$$\bar{\cal M}_{\cal T}=
\bigcup_{{\cal T}'\le{\cal T}}{\cal M}_{{\cal T}'},\quad
\bar{\cal U}_{\cal T}^{(0)}=
\bigcup_{{\cal T}'\le{\cal T}}{\cal U}_{{\cal T}'}^{(0)}\quad
\bar{\cal U}_{\cal T}=
\bigcup_{{\cal T}'\le{\cal T}}{\cal U}_{{\cal T}'}$$
are compact and Hausdorff.
The $G_{\cal T}^*$-action on ${\cal U}_{\cal T}^{(0)}$ extends
to an action on $\bar{\cal U}_{\cal T}^{(0)}$,
and thus line orbi-bundles $L_{{\cal T},i}\lra{\cal U}_{\cal T}$
with $i\!\in\! I\!-\!\hat{I}$ extend over~$\bar{\cal U}_{\cal T}$.
The evaluation maps
$$\ev_l\!:{\cal H}_{\cal T}\lra\P,\quad
\ev_l\big((S,M,I;x,(j,y),u)\big)=
u_{j_l}(y_l),$$
descend to all the quotients and induce continuous
maps on $\bar{\cal M}_{\cal T}$, $\bar{\cal U}_{\cal T}$,
and~$\bar{\cal U}_{\cal T}^{(0)}$.
If $\mu\!=\!\mu_M$ is an $M$-tuple of submanifolds of~$\P$,
let
$${\cal M}_{\cal T}(\mu)=
\{b\!\in\!{\cal M}_{\cal T}\!:\ev_l(b)\!\in\!\mu_l~\forall 
l\!\in\! M\}$$
and define spaces ${\cal U}_{\cal T}(\mu)$,
$\bar{\cal U}_{\cal T}(\mu)$, etc.~in a similar way.
If $S=S^2$, we define another evaluation~map,
$$\ev\!: 
{\cal B}_{\cal T}\lra\P\quad\hbox{by}\quad
\ev\big((S^2,M,I;x,(j,y),u)\big)=u_{\hat{0}}(\i),$$
where $\hat{0}$ is any minimal element of $I$.
This map descends to ${\cal U}_{\cal T}^{(0)}$ and~${\cal U}_{\cal T}$.
If $\mu=\mu_{\hat{0}+M}$ is a tuple of constraints, let
$${\cal U}_{\cal T}(\mu_{\hat{0}};\mu_M)=
\{b\!\in\!{\cal U}_{\cal T}(\mu_{M})\!:\ev(b)\!\in\!\mu_{\hat{0}}\}$$
and define ${\cal U}_{\cal T}^{(0)}(\mu_{\hat{0}};\mu_M)$, etc.~similarly.
If $S=\Si$, ${\cal T}$ is a simple bubble type,
and $d_{\hat{0}}=0$, define
$$\ev\!: 
{\cal H}_{\cal T}\lra\P\quad\hbox{by}\quad
\ev\big((\Si,M,I;x,(j,y),u)\big)=u_{\hat{0}}(\Si).$$
This map is well-defined, since $u_{\hat{0}}$ is a degree-zero
holomorphic map and thus is constant.\\

\noindent
If ${\cal T}$ is any bubble type, 
let $\lr{\cal T}$ be the basic bubble such 
that \hbox{${\cal T}\!\le\!\lr{\cal T}$}.
If ${\cal T}$ is a simple bubble type,
let $\bar{\cal T}$ be the bubble type obtained from ${\cal T}$
by dropping the minimal element $\hat{0}$ from the indexing set~$I$
and the subset $M_{\hat{0}}{\cal T}$ from~$M$.
Note that if ${\cal T}$ is primitive, 
$\bar{\cal T}$ is basic.\\

\noindent
Finally, if $X$ is any space, $F\!\lra\! X$ a normed vector bundle,
and $\de\!: X\!\lra\!\Bbb{R}$ is any function, let
$$F_{\de}=\big\{(b,v)\!\in\! F\!: |v|_b<\de(b)\big\}.$$
Similarly, if $\Om$ is a subset of $F$, let 
$\Om_{\de}=F_{\de}\cap\Om$.
If $\ups=(b,v)\!\in\! F$, denote by $b_{\ups}$ the image of $\ups$
under the bundle projection map, i.e. $b$ in this case.

\section{Analysis}
\label{analysis_sec}

\subsection{The Basic Setup}
\label{basic_setup}

\noindent
In this section, we focus on bubble types
${\cal T}\!=\!\big(S,M,I;j,d\big)$ 
such that either $S\!=\!S^2$ or $d_{\hat{0}}\!=\!0$.
In the first case, we describe a small neighborhood of 
${\cal U}_{\cal T}(\mu)$ in $\bar{\cal U}_{\lr{\cal T}}(\mu)$
and the behavior of sections of certain bundles over
$\bar{\cal U}_{\lr{\cal T}}(\mu)$ near ${\cal U}_{\cal T}(\mu)$;
see Theorem~\ref{str_global}.
This theorem is deduced from 
Theorem~\ref{gl-str_global} in~\cite{Z}.
If ${\cal T}$ is a simple bubble type, $S\!=\!\Si$, and $d_{\hat{0}}\!=\!0$,
we describe the elements of ${\cal M}_{\Si,t\nu,d}(\mu)$ lying near 
${\cal M}_{\cal T}(\mu)$ as the zero set of a map
defined on an open subset of the bundle~$F{\cal T}$;
see Theorem~\ref{si_str}.
The map takes values in a bundle over ${\cal M}_{\cal T}(\mu)$,
which is the analogue of Taubes's obstruction bundle  of~\cite{T}
in this setting.
Theorem~\ref{si_str} is a consequence of Theorem~\ref{gl-si_str}
in~\cite{Z}, which requires us to make two major choices.
This is done in the next two subsections.\\

\noindent
If ${\cal T}\!=\!\big(S,M,I;j,d\big)$ and $S\!=\!S^2$,
by Corollaries~\ref{reg_crl2} and~\ref{reg_crl},
${\cal T}$ is a $(\P,J)$-regular bubble type in the sense
of Definition~\ref{gl-reg_dfn} in~\cite{Z}.
This regularity property implies that\\
(R1) ${\cal H}_{\cal T}$ is a smooth manifold;\\
(R2) for any $b=\big(S,M,I;x,(j,y),d\big)
 \!\in\!{\cal H}_{\cal T}$, a neighborhood of 
$b$ in ${\cal H}_{\cal T}$, is modeled on 
$$\ker\big(D_b\!:\Ga(b)\!\lra\!\Ga^{0,1}(b)\big)\oplus
\bigoplus_{h\in I}T_{x_h}\Si_{b,\io_h}\oplus
\bigoplus_{l\in M}T_{y_l}\Si_{b,j_l}.$$
(R3) $D_b\!:\Ga(b)\!\lra\!\Ga^{0,1}(b)$ is surjective for all
 $b\!\in\!{\cal H}_{\cal T}$.\\
Here $\Ga^{0,1}(b)$ denotes the space of 
$u_b^*T\P$-valued $(0,1)$-forms on the components of~$\Si_b$,
while $\Ga(b)$ is the set of vector fields $\xi$ on the components
of~$\Si_b$ that agree at the nodes and such that 
\hbox{$\xi(i_1,\i)\!=\!\xi(i_2,\i)$} for 
\hbox{all $i_1,i_2\!\in\! I\!-\!\hat{I}$}.
The operator $D_b$ is the linearization of the $\bar{\partial}$-operator
with respect to a connection in~$T\P$.
Along ${\cal H}_{\Si}$, it is independent of the choice of
the connection.
On the other hand, if ${\cal T}$ is a simple bubble type,
$S\!=\!\Si$, and $d_{\hat{0}}\!=\!0$, by the same two corollaries,
${\cal T}$ is a $(\P,J)$-semiregular bubble type in the sense
of Definition~\ref{gl-semireg_dfn} in~\cite{Z}.
This means that (R1) and (R2) are satisfied,
with $\Ga(b)$ defined as above but omitting the last condition.
Property (R3) is not satisfied, and in fact by the two corollaries,
$$\coker~D_b\approx {\cal H}_{\Si}^{0,1}\otimes T_{\ev(b)}\P
\qquad\forall b\!\in\!{\cal H}_{\cal T},$$
where ${\cal H}_{\Si}^{0,1}$ is the space of 
harmonic $(0,1)$-forms on~$\Si$.
This cokernel bundle descends to 
a bundle $\Ga_-^{0,1}\!\lra\!{\cal M}_{\cal T}$,
which will be our obstruction bundle.\\

\noindent
If $S\!=\!\Si$, for the gluing construction in~\cite{Z},
we choose a smooth family 
\hbox{$\{g_{b,\hat{0}}\!:b\!\in\!{\cal H}_{\cal T}\}$} 
of metrics on~$\Si$ such that for all 
$$b=\big(\Si,M,I;x,(j,y),u\big)\in{\cal H}_{\cal T},$$
the metric $g_{b,\hat{0}}$ is flat on a neighborhood of $x_h$
in $\Si$ for all $h\!\in\!\hat{I}$ such \hbox{that $\io_h\!=\!\hat{0}$}.
This family of metrics, in fact, depends only 
on the sets~$\{x_h\!:\io_h\!=\!\hat{0}\}$.
Along with the standard metric on~$S^2$, 
the metric $g_{b,\hat{0}}$ induces a Riemannian metric
$g_b\!=\!(g_{b,i})_{i\in I}$ on $\Si_b\!=\!\bigcup\limits_{i\in I}\Si_{b,i}$.
If $S\!=\!S^2$, we take $g_{b,i}$ to be the standard metric
on $\Si_{b,i}\!=\!S^2$ for \hbox{all $i\!\in\! I$}.
With notation as above, if $x_h,z\!\in\!\Si_{b,\hat{0}}\!=\!\Si$,
let $r_{b,h}(z)\!=\!d_{g_{b,\hat{0}}}(x_h,z)$.
If $x_h,z\!\in\!\Si_{b,i}\!=\!S^2$ and $z\!\neq\!\i$, 
\hbox{let  $r_{b,h}(z)\!=\!|z\!-\!x_h|$}.\\

\noindent
For each $\ups\!=\!(b,v_h)_{h\in\hat{I}}\!\in\! F^{(0)}{\cal T}$
sufficiently small, in~\cite{Z}
we then define a complex curve $\Si_{\ups}$,
smooth maps $q_{\ups}\!:\Si_{\ups}\!\lra\!\Si_b$
and  $q_{\ups,i}\!:\Si_{\ups,\hat{0}}\!\lra\!\Si_b$ for $i\!\in\! I$,
and Riemannian metric $g_{\ups}$ on $\Si$ on $\Si_{\ups}$
such~that\\
(G1) the linearly ordered set corresponding to $\Si_{\ups}$
is $I(\ups)\!\equiv\! I\!-\!\{h\!\in\!\hat{I}\!:v_h\!\neq\!0\}$;\\
(G2) the map $q_{\ups}|_{\Si_{\ups,\hat{0}}}$ factors 
through each of the maps $q_{\ups,i}$;\\
(G3) $q_{\ups}\!: (\Si_{\ups},g_{\ups})\!\lra\!(\Si_b,g_b)$ is an isometry
(and thus holomorphic) outside of the annuli
\begin{equation}\label{ann_def1}\begin{split}
&A_{\ups,h}^+=q_{\ups,\io_h}^{-1}\Big(\big\{z\!\in\!\Si_{b,\io_h}\!: 
|v_h|^{\frac{1}{2}}\le r_{b,h}(z)\le 2|v_h|^{\frac{1}{2}}\big\}\Big);\\
&A_{\ups,h}^-=q_{\ups,\io_h}^{-1}\Big(\big\{z\!\in\!\Si_{b,\io_h}\!: 
\frac{1}{2}|v_h|^{\frac{1}{2}}\le r_{b,h}(z)\le |v_h|^{\frac{1}{2}}
\big\}\Big).
\end{split}\end{equation}
(G4) $q_{\ups,\io_h}\!: (A_{\ups,h}^{\pm},g_{\ups})\!\lra\!
 \big(q_{\ups,\io_h}(A_{\ups,h}^{\pm}),g_b\big)$
is an isometry.\\
The map $q_{\ups}$ collapses disjoint circles on $\Si_{\ups}$
and identifies the resulting surfaces with~$S^2$
in a manner encoded by~$\ups$.
Alternatively, $(\Si_{\ups},g_{\ups})$ can be viewed
as the surface obtained by smoothing (some~of) the nodes of~$\Si_b$.
The maps $q_{\ups}$ and $q_{\ups,i}$ are constructed
explicitly by fixing a smooth bump function 
\hbox{$\be\!:\Bbb{R}\!\lra\![0,1]$} such that
\begin{equation}\label{cutoff_fun}
\be(t)=\begin{cases}
0,&\hbox{if~}t\le 1;\\
1,&\hbox{if~}t\ge 2,
\end{cases}
\qquad\hbox{and}\qquad
\be'(t)>0~~\hbox{if~}t\!\in\!(1,2).
\end{equation}
If $r>0$, let $\be_r\!\in\! C^{\i}(\Bbb{R};\Bbb{R})$ be given by 
$\be_r(t)\!=\!\be(r^{-\frac{1}{2}}t)$.
Note that
\begin{equation}\label{cutoff_fun2}
\hbox{supp}(\be_r)=[r^{\frac{1}{2}},2r^{\frac{1}{2}}],~~
\|\be_r'\|_{C^0}\le C_{\be}r^{-\frac{1}{2}},\hbox{~~and~~}
\|\be_r''\|_{C^0}\le C_{\be}r^{-1}.
\end{equation}
These cutoff functions will not appear in the main statements of this paper,
but they do show up in the proofs of 
Lemma~\ref{flat_metrics},
Theorem~\ref{str_global}, and  Proposition~\ref{dbar_gen2prp}.
Having constructed the maps $q_{\ups}$, we let
$b(\ups)\!=\!(\Si_{\ups},u_{\ups})\!=\!(\Si_{\ups},u_b\circ q_{\ups})$.
The marked points on $\Si_{\ups}$ are the preimages of the marked
points of $\Si_b$ under the map~$q_{\ups}$.\\

\noindent
We also need to choose a smooth family 
$\{g_{\P,b}\!:b\!\in\!{\cal M}_{\cal T}^{(0)}\}$ of metrics on $\P$
invariant under the equivalence relation 
on ${\cal M}_{\cal T}^{(0)}$ if $S\!=\!\Si$ and on ${\cal B}_{\cal T}$
\hbox{if $S\!=\!S^2$}.
While taking $g_{\P,b}$ to be the standard metric on~$\P$
may be the canonical choice, for computational reasons
it is more convenient to take $g_{\P,b}\!=\!g_{\P,\ev(b)}$,
where $\{g_{\P,q}\!: q\!\in\!\P\}$ is the family of metrics
of Lemma~\ref{flat_metrics}.

\begin{lmm}
\label{flat_metrics}
There exist $r_{\P}\!>\!0$ and a smooth family of Kahler metrics
$\{g_{\P,q}\!: q\!\in\!\P\}$ on~$\P$ with the following property.
If $B_q(q',r)\!\subset\!\P$ denotes the $g_{\P,q}$-geodesic ball about~$q'$
of radius~$r$,
the triple $(B_q(q,r_{\P}),J,g_{\P,q})$ is isomorphic
to a ball in $\Bbb{C}^n$ for \hbox{all $q\!\in\!\P$}.
\end{lmm}

\noindent
{\it Proof:} On the open set 
$U_0\!=\!\{[X_0\!:\!\ldots\!:\!X_n]\!\in\!\P\!: X_0\!\neq\! 0\}$,
the Fubini-Study symplectic form is given~by
\begin{equation}\label{flat_metrics_e1}
\om_{\P}=\frac{i}{2\pi}\partial\bar{\partial}\ln(1+f_0),
\hbox{~~where~~}
f_0([X_0:\ldots:X_n])=\sum_{k\in[n]}|X_k/X_0|^2;
\end{equation}
see \cite[p31]{GH}. Let $q=[1:0:\ldots:0]$. Set
\begin{equation}\label{flat_metrics_e2}
\om_{\P,q,\ep}=\frac{i}{2\pi}\partial\bar{\partial}
\big\{f_0+(\be_{\ep^2}\circ f_0)
\big(\ln(1+f_0)-f_0\big)\big\}.
\end{equation}
Note that $\om_{\P,q,\ep}$ agrees with $\om_{\P}$ outside of the set
$\{f_0\le 2\ep\}$ and with the standard symplectic form $\om_{\Bbb{C}^n}$
on $\{f_0\le\ep\}$.
Here we view $\om_{\Bbb{C}^n}$ as a form on $U_0$ via the coordinates
$z_{0,i}\!=\!X_i/X_0$, $i\!\in\![n]$.
In particular, $\om_{\P,q,\ep}$ is globally defined, and 
the corresponding Riemannian metric on $\{f_0\le\ep\}$ is flat.
Furthermore,
\begin{equation}\label{flat_metrics_e3}\begin{split}
\om_{\P,q,\ep}=&\big\{(1-\be_{\ep^2}\circ f_0)\om_{\Bbb{C}^n}
+(\be_{\ep^2}\circ f_0)\om_{\P}\big\}\\
&+\frac{i}{2\pi}\Big\{
\big(\partial(\be_{\ep^2}\circ f_0)\big)\big(\bar{\partial}\tilde{f}_0\big)
-\big(\bar{\partial}(\be_{\ep^2}\circ f_0)\big)
\big(\partial\tilde{f}_0\big)
+\big(\partial\bar{\partial}(\be_{\ep^2}\circ f_0)\big)\tilde{f}_0\Big\},
\end{split}\end{equation}
where $\tilde{f}_0=\ln(1+f_0)-f_0$.
On the set $\{f_0\le 2\ep\}$ with $\ep\le\frac{1}{2}$,
\begin{equation}\label{flat_metrics_e4}
\big\|\tilde{f}_0\big\|_{C^0}\le C\ep^2\hbox{~~and~~}
\big\|d\tilde{f}_0\big\|_{C^0}\le C\ep^{\frac{3}{2}},
\end{equation}
where $\big\|d\tilde{f}_0\big\|_{C^0}$ denotes the $C^0$-norm with
respect to the standard metric on $\Bbb{C}^n$.
Furthermore, by~\e_ref{cutoff_fun2},
\begin{equation}\label{flat_metrics_e5}
\big\|d(\be_{\ep^2}\circ f_0)\big\|_{C^0}\le C\ep^{-1}\ep^{\frac{1}{2}},
\quad
\big\|\na^2(\be_{\ep^2}\circ f_0)\big\|_{C^0}\le 
C\big(\ep^{-2}\ep^{\frac{1}{2}}\ep^{\frac{1}{2}}+\ep^{-1}\big),
\end{equation}
where again all the norms are computed with respect to the standard metric
on $\Bbb{C}^n$. Equations \e_ref{flat_metrics_e4} and \e_ref{flat_metrics_e5}
imply that the term on the second line of \e_ref{flat_metrics_e3}
tends to $0$ as $\ep$ goes $0$.
Thus by \e_ref{flat_metrics_e3}, 
we can choose $\ep>0$ such that $\om_{\P,q}\!\equiv\!\om_{\P,q,\ep}$
is a symplectic form on all of~$\P$.
Note that $\om_{\P,q}$ is invariant under the action of 
the stabilizer of $q$ in $SU_{n+1}$, which is the subgroup
$$\hbox{Stab}_p(SU_{n+1})=\l\{
\l(\begin{array}{cc}
\ov{\det(h)}& 0\\
0& h
\end{array}\r)\!: h\!\in\! U_n\r\}\subset SU_{n+1}.$$
We can define a smooth family of symplectic Kahler forms on $\P$ by 
$$\om_{\P,g\cdot q}=g^*\om_{\P,q},\qquad g\in SU_{n+1}.$$
The above invariance property of $\om_{\P,q}$ insures that 
$\om_{\P,g\cdot q}$ depends only on~$g\cdot q$.
We can now take $g_{\P,g\cdot q}$ to be the metric corresponding
to the symplectic form $\om_{\P,g\cdot q}$ and the standard complex
structure~$J$ on~$\P$.\\

\noindent
We denote by $\exp_b$ and $\Pi_{b,X}$ for $X\!\in\!T\P$
the $g_{\P,b}$-exponential map and $g_{\P,b}$-parallel transport
along the $g_{\P,b}$-geodesic for~$X$, respectively.
If $\ups\!\in\! F^{(0)}{\cal T}$, let
$$g_{\P,\ups}=g_{\P,b_{\ups}},\quad
\exp_{\ups}=\exp_{b_{\ups}},\quad
\Pi_{\ups,X}=\Pi_{b_{\ups},X}.$$
If $\ups\!\in\! F^{(0)}$ is sufficiently small,
we define $L^2$-norms inner-products on 
$$\Ga(\ups)\equiv\Ga\big(b(\ups)\big)
\qquad\hbox{and}\qquad
\Ga^{0,1}(\ups)\equiv\Ga^{0,1}\big(b(\ups)\big)$$
via the metrics $g_{\P,\ups}$ and $g_{\ups}$ in the usual way.
Denote by $D_{\ups}$ the linearization of 
the $\bar{\partial}$-operator with respect to 
the metric~$g_{\P,\ups}$ on~$\P$ and by $D_{\ups}^*$
its formal adjoint with respect to 
the above $(L^2,\ups)$ inner-product.
We fix $p\!>\!2$ and denote by 
$\|\cdot\|_{\ups,p,1}$ and $\|\cdot\|_{\ups,p}$ 
the modified Sobolev $(L^p_1,g_{\P,\ups},g_{\ups})$
and $(L^p,g_{\P,\ups},g_{\ups})$ norms of~\cite{LT}
on $\Ga(\ups)$ and $\Ga^{0,1}(\ups)$, respectively.
Let $L^p_1(\ups)$ and $L^p(\ups)$ be the corresponding completions.
A description of the modified Sobolev norms in the notation
of this paper can be found in~\cite{Z}.
They are needed only for certain technical aspects of this paper.

\subsection{Obstruction Bundle}
\label{obs_sec}

\noindent
In this subsection, in the case $S\!=\!\Si$,
we choose an obstruction bundle over~$F^{(\eset)}{\cal T}_{\de}$ 
in the sense of  Definition~\ref{gl-obs_setup_dfn2} in~\cite{Z}
with $\de\!\in\! C^{\i}({\cal M}_{\cal T};\Bbb{R}^+)$
sufficiently small.\\

\noindent
Let $\de_{\cal T}\!\in\! C^{\i}({\cal M}_{\cal T};\Bbb{R}^+)$ such that 
$$4\de_{\cal T}(b)\|du_i\|_{b,C^0}<r_{\P}
\qquad\forall
b=\big(\Si,M,I;x,(j,y),u\big)\!\in\!{\cal M}_{\cal T},~
i\!\in\!I.$$
We assume that the above function $\de$ is such that
$8\de^{\frac{1}{2}}\!\!<\!\de_{\cal T}$.
For all $\ups\!\in\! F^{(\eset)}{\cal T}_{\de}$ and
\hbox{$X\psi\!\in\! T_{\ev(b_{\ups})}\P\!\otimes\!{\cal H}_{\Si}^{0,1}$}\!, 
define $R_{\ups}X\psi\!\in\!\Ga^{0,1}(u_{\ups})$ as follows.
If $z\!\in\!\Si_{\ups}\!=\!\Si$ is such that 
$q_{\ups}(z)\!\in\!\Si_{b_{\ups},h}$ 
for some $h\!\in\!\hat{I}$ with $\chi_{\cal T}h\!=\!1$ and 
$\big|q_S^{-1}(q_{\ups}(z))\big|\le 2\de_{\cal T}(b_{\ups})$,
by our assumption on~$\de_{\cal T}$,
we can define \hbox{$\bar{u}_{\ups}(z)\!\in\! T_{\ev(b_{\ups})}\P$ by}
$$\exp_{\ups,\ev(b_{\ups})}\bar{u}_{\ups}(z)=u_{\ups}(z), \quad
|\bar{u}_{\ups}(z)|< r_{\P}.$$
Given $z\!\in\!\Si$, let $h_z\!\in\! I$ be such that 
$q_{\ups}(z)\!\in\!\Si_{b_{\ups},h_z}$.
If $w\!\in\! T_z\Si$, put
$$R_{\ups}X\psi|_zw=
\begin{cases}
0, & \hbox{if~}\chi_{\cal T}h_z=2;\\
\be\big(\de_{\cal T}(b_{\ups})|q_{\ups}z|\big)
(\psi|_zw)\Pi_{\ups,\bar{u}_{\ups}(z)}X,&  
                \hbox{if~}\chi_{\cal T}h_z=1;\\
(\psi|_zw)X,& \hbox{if~}\chi_{\cal T}h_z=0.
\end{cases}$$
Let $\Ga_-^{0,1}(\ups)$ be the image of 
$T_{\ev(b_{\ups})}\P\otimes{\cal H}_{\Si}^{0,1}$ under the map~$R_{\ups}$.
Denote by $\pi_{\ups,-}^{0,1}$ the $(L^2,\ups)$-orthogonal 
projection of $L^p(\ups)$ onto~$\Ga_-^{0,1}(\ups)$.\\

\noindent
The spaces $\Ga_-^{0,1}(\ups)$ form our obstruction 
bundle over~$F^{(\eset)}{\cal T}$.
We need to show that these spaces satisfy the requirements of
Definition~\ref{gl-obs_setup_dfn2}.
First, the rate of change of $\pi_{\ups,-}^{0,1}$
with respect to changes in $\ups$ should be controlled by
a function of~$b_{\ups}$ only.
The proof of this is similar to the proof of 
the second statement of (5) of Lemma~\ref{gl-inj_l1} in~\cite{Z}.
The next lemma implies that the remaining conditions are also
satisfied.
For any $h\!\in\!\hat{I}$, put 
$$|\ups|_h=\prod_{i\in\hat{I},h\in\bar{D}_i{\cal T}}|v_i|.$$

\begin{lmm}
\label{approx_coker}
For any $\ups\!\in\! F^{(\eset)}{\cal T}_{\de}$
and $X\psi\!\in\! T_{\ev(b_{\ups})}\P\otimes{\cal H}_{\Si}^{0,1}$,
$D^*_{\ups}R_{\ups}X\psi$ vanishes outside of the~annuli
$$\tilde{A}_{\ups,h}\equiv q_{\ups}^{-1}\big(
\big\{(h,z)\!\in\!\Si_{b_{\ups},h}\!: 
\de_{\cal T}(b_{\ups})\le|q_S^{-1}(z)|\le 
2\de_{\cal T}(b_{\ups})\big\}\big)$$
with $h\!\in\!\hat{I}$ such that $\chi_{\cal T}h\!=\!1$.
Furthermore,
there exists $C\!\in\! C^{\i}({\cal M}_{\cal T};\Bbb{R}^+)$ 
such that\\
(1) $\|D_{\ups}^*R_{\ups}X\psi\|_{\ups,C^0}\le  C(b_{\ups})
\Big(\sum\limits_{\chi_{\cal T}h=1}|\ups|_h\Big)
|X|_{\ups}\|\psi\|_2$;\\
(2) $\big(1-C(b_{\ups})^{-1}|\ups|^{\frac{2}{\tilde{p}}}\big)
\|X\psi\|_{\ups,\tilde{p}}
\le\|R_{\ups}X\psi\|_{\ups,\tilde{p}}
\le \big(1+C(b_{\ups})^{-1}|\ups|^{\frac{2}{\tilde{p}}}\big)
\|X\psi\|_{\ups,\tilde{p}}$, where $\tilde{p}=2,p$.
\end{lmm}

\noindent
{\it Proof:} 
The first statement and estimate (2) are immediate from 
the definition of $R_{\ups}X\psi$ and of the norms; see~\cite{Z}.
Let $(s,t)$ be the  conformal coordinates on $\tilde{A}_{\ups,h}$ given by
$q_{\ups}(s,t)=s+it\in\Bbb{C}$.
Write $g_{\ups}=\th^{-2}(s,t)(ds^2+dt^2)$. Then
\begin{equation}\label{coker_e1}
\th=\frac{1}{2}\l(1+s^2+t^2\r).
\end{equation}
Put
\begin{equation}\label{coker_e2}
\xi(s,t)=\big\{R_{\ups}X\psi\big\}_{(s,t)}\partial_s=
\be\big(\de_{\cal T}(b_{\ups})\sqrt{s^2+t^2}\big)
\big(\psi|_{(s,t)}\partial_s\big)\Pi_{\ups,\bar{u}_{\ups}(s,t)}X.
\end{equation}
Then by \cite[p29]{MS},
\begin{equation}\label{coker_e3}
D^*_{\ups}R_{\ups}X\psi|_z=\th^2\Big(
-\frac{D}{ds}\xi+J\frac{D}{dt}\xi\Big),
\end{equation}
where $\frac{D}{ds}$ and $\frac{D}{dt}$ denote covariant differentiation
with respect to the metric $g_{\P,\ups}$ on~$\P$.
Since this metric is flat on the support of~$\xi$ and 
$\psi\!\in\!{\cal H}_{\Si}^{0,1}$, 
equations \e_ref{coker_e1}-\e_ref{coker_e3} give
\begin{equation}\label{coker_e4}\begin{split}
D^*_{\ups}R_{\ups}X\psi|_z=\frac{\big(1+s^2+t^2\big)^2}{4}
\Big\{\be'|_{\de_{\cal T}(b_{\ups})\sqrt{s^2+t^2}}\de_{\cal T}(b_{\ups})
\frac{-s+it}{\sqrt{s^2+t^2}}\Big\}
\big(\psi|_{(s,t)}\partial_s\big)\Pi_{\ups,\bar{u}_{\ups}(s,t)}X.&
\end{split}\end{equation}
Since the right hand-side of \e_ref{coker_e4} vanishes unless
$\de_{\cal T}(b_{\ups})^{-1}\!\le\!\sqrt{s^2+t^2}
\!\le\! 2\de_{\cal T}(b_{\ups})^{-1}$,
it follows~that
\begin{equation}\label{coker_e5}
\big|D^*_{\ups}R_{\ups}X\psi\big|_{\ups,z}\le 
C(b_{\ups})\big|\psi|_{(s,t)}\partial_s\big||X|_{\ups}
\le C'(b_{\ups})|\ups|_h\|\psi\|_2|X|
\end{equation}
Claim (1) follows from \e_ref{coker_e5}.\\

\noindent
Let 
$\tilde{R}_{\ups}\!:{\cal H}_{\Si}^{0,1}\otimes T_{\ev(b_{\ups})}\P
\!\lra\!\Ga_-(\ups)$ be the adjoint of~$R_{\ups}^{-1}$,
i.e.
\begin{equation}\label{new_isom}
\big\lan\big\lan \tilde{R}_{\ups}X\psi,R_{\ups}X'\psi'
\big\ran\big\ran_{\ups,2}=
\big\lan\big\lan X\psi,X'\psi'\big\ran\big\ran_{b_{\ups},2}
=\lan X,X'\ran_{b_{\ups}}\lan\psi,\psi'\ran_2
\qquad
\end{equation}
for all $X,X'\!\in\!T_{\ev(b_{\ups})}\P$ and 
$\psi,\psi'\!\in\!{\cal H}_{\Si}^{0,1}$.
By Lemma~\ref{approx_coker}, 
$\|\tilde{R}_{\ups}-R_{\ups}\|_2\le C(b_{\ups})|\ups|$.

\subsection{Tangent-Bundle Model}
\label{tangent_sec}

\noindent
We now describe our choice for a tangent-bundle model,
which is the subject of Definition~\ref{gl-obs_setup_dfn2} in~\cite{Z}.\\

\noindent
For any $\ups\!\in\! F^{(0)}{\cal T}$ sufficiently small and
$\xi\!\in\!\Ga(b_{\ups})$,
define \hbox{$R_{\ups}\xi\!\in\! L^p_1(\ups)$~by}
$\{R_{\ups}\xi\}(z)=\xi\big(q_{\ups}(z)\big)$.
Let $\Ga_-(\ups)$ be the image of $\ker(D_{b_{\ups}})$ 
under the map~$R_{\ups}$. 
Denote by $\Ga_+(\ups)$ its $(L^2,\ups)$-orthogonal complement in 
$L^p_1(\ups)$.
Let $\pi_{\ups,\pm}$ be the 
$(L^2,\ups)$-orthogonal projection onto~$\Ga_{\pm}(\ups)$.\\

\noindent
If $x\!\in\!\Si$, let
${\cal H}_{\Si}^-(x)\!=\!\{\psi\!\in\!{\cal H}_{\Si}^{0,1}\!: 
\psi|_x\!=\!0\}$.
This is a codimension-one subspace of 
${\cal H}_{\Si}^{0,1}$ for all $x\!\in\!\Si$; see~\cite{GH}.
Denote by ${\cal H}_{\Si}^+(x)$ its $L^2$-orthogonal complement.
The space ${\cal H}_{\Si}^+(x)$ is independent of
the choice of a Kahler metric on~$(\Si,j_{\Si})$.
For any $h\!\in\!\hat{I}$, we put
$\tilde{x}_h(\ups)\!=\!q_{\ups,\io_h}^{-1}(\io_h,x_h)$.
Fix $h^*\!\in\!\hat{I}$ such that $\chi_{\cal T}h^*\!=\!1$.
Let
$$\bar{\Ga}_-(\ups)=D_{\ups}^*R_{\ups}
\big( {\cal H}_{\Si}^+(\tilde{x}_{h^*}(\ups))\otimes 
T_{\ev(b_{\ups})}\P\big).$$
Denote by $\bar{\Ga}_+(\ups)$ the $(L^2,\ups)$-orthogonal complement
of $\bar{\Ga}_-(\ups)$ in $L^p_1(\ups)$
and by $\bar{\pi}_{\ups,\pm}$ 
the $(L^2,g_{\ups})$-orthogonal projections
onto $\bar{\Ga}_{\pm}(\ups)$.
Let $\tilde{\Ga}_-(\ups)$ be the $(L^2,\ups)$-orthogonal complement
of $\tilde{\Ga}_+(\ups)\!\equiv\!\bar{\pi}_{\ups,+}\big(\Ga_+(\ups)\big)$
in $L^p_1(u_{\ups})$.\\

\noindent
The spaces  $\tilde{\Ga}_-(\ups)$ will be our tangent-space model.
We need to check that the requirements of 
Definition~\ref{gl-obs_setup_dfn1} are satisfied.
Let 
$$\{h\!\in\!\hat{I}\!:\chi_{\cal T}(h)\!=\!1\}
=\{h_1=h^*,h_2,\ldots,h_m\}.$$
If $z\!\in\!\Si_{b,h_r}$ is such that $|q_S^{-1}(z)|\le 2\de_{\cal}(b)$,
define $\bar{u}_{h_r}(z)\in T_{\ev(b)}\P$ by
$$\exp_{b,\ev(b)}\bar{u}_{h_r}(z)=u_{h_r}(z),\quad
|\bar{u}_{h_r}(z)|_b<r_{\P}.$$
If $X\!\in\! T_{\ev(b)}\P$, define $R_{b,h_r}X\!\in\!\Ga(u_{h_r})$ by
$$R_{b,h_r}X(z)=
\begin{cases}
0,&\hbox{if~}|z|\ge 2\de_{\cal T}(b)^{-1};\\
\be'\big|_{\de_{\cal T}(b)|z|}
\frac{(1+|z|^2)^2z}{|z|}
\Pi_{b,\bar{u}_{h_r}(z)}X,&\hbox{otherwise}.
\end{cases}$$
Since $R_{b,h_r}X$ vanishes at all the nodes of $\Si_b$
by assumption on $\de_{\cal T}$,
we can extend $R_{b,h_r}X$ by zero to an element of $\Ga(b)$.
If $c\!=\!c_{[m]}\!\in\!\Bbb{C}^{[m]}$ is different from zero, let
$$\bar{\Ga}_-(b;c)=\Big\{ \sum_{r\in[m]}c_rR_{b,h_r}X\!: 
 X\!\in\! T_{\ev(b)}\P\Big\}.$$
Denote by $\bar{\Ga}_+(b;c)$ the $(L^2,b)$-orthogonal complement
of $\bar{\Ga}_-(b;c)$ in $\Ga(b)$.
Let $\bar{\pi}_{(b;c),\pm}$ be the corresponding $(L^2,b)$-orthogonal
projection maps.
Let $\tilde{\Ga}_+(b;c)=\bar{\pi}_{(b,c),+}\big(\Ga_+(b)\big)$
and let $\tilde{\Ga}_-(b;c)$ be its $(b,L^2)$-orthogonal complement.

\begin{lmm}
\label{inverse_p2_l1}
There exist $\de,C\!\in\! C^{\i}({\cal M}_{\cal T}^{(0)};\Bbb{R}^+)$
such that for all $\ups\!\in\! F^{(\eset)}{\cal T}_{\de}$ and 
$\xi\!\in\!\bar{\Ga}_-(\ups)$,
$$\|\xi\|_{\ups,p,1}\le C(b_{\ups})\|\xi\|_{\ups,2}.$$
In addition, $\dim_{\Bbb{C}}\bar{\Ga}_-(\ups)\!=\!
     \dim_{\Bbb{C}}\bar{\Ga}_-(b_{\ups};c)\!=\!n$
                             for any nonzero $c\!\in\!\Bbb{C}^m$.
Furthermore, if 
\hbox{$\ups_k\!\!\lra\! b\!\in\!{\cal M}_{\cal T}^{(0)}$}
and $\xi_k\!\in\!\bar{\Ga}_-(\ups)$ is such that $\|\xi_k\|_{\ups_k,2}=1$,
then there exists a nonzero $c\!\in\!\Bbb{C}^m$ and 
$\xi\!\in\!\bar{\Ga}_-(b;c)$
with $\|\xi\|_{b,2}\!=\!1$ such that
a subsequence of $\{\xi_k\}$ $C^0$-converges to~$\xi$.
\end{lmm}

\noindent
{\it Remark:} The last statement means that 
a subsequence of $\{\xi_k\}$ $C^0$-converges to~$\xi$
on compact subsets of $\Si_b^*$ and the norms $\|\xi_k\|_{\ups_k,p,1}$
are uniformly bounded; see Definition~\ref{gl-C0conv_dfn} in~\cite{Z}.\\

\noindent
{\it Proof:} (1) Let $\psi$ be a generator of 
              ${\cal H}^{0,1}_{\Si,+}\big(\tilde{x}_{h_1}(\ups)\big)$.
If $X\!\in\! T_{\ev(b_{\ups})}\P$ and $r\!\in\![m]$, define
$R_{\ups,h_r}X\!\in\!\Ga(u_{\ups})$ as follows.
If $q_{\ups}(z)\!\in\!\Si_{b_{\ups},h_r}$, let
\begin{equation*}\begin{split}
R_{\ups,h_r}X(z)=
\Big(\sum_{r\in[m]}
\big|\psi_{\tilde{x}_r(\ups)}d(q_{\ups,h_r}^{-1}\circ q_N)\partial_s\big|
\Big)^{-1}
\frac{(1+|q_{\ups}z|^2)^2q_{\ups}z}{|q_{\ups}z|}&\\
\times\be'\big|_{\de_{\cal T}(b_{\ups})|q_{\ups}z|}
\big(\psi_zd(q_{\ups,h_r}^{-1}\circ q_N)\partial_s\big)
\Pi_{b_{\ups},\bar{u}_{\ups}(z)}X.&
\end{split}\end{equation*}
Note that the sum is not zero, since 
$\psi|_{\tilde{x}_{h_1}(\ups)}\!\neq\! 0$. 
If $q_{\ups}(z)\!\not\in\!\Si_{b_{\ups},h_r}$, 
we let $R_{\ups,h_r}X(z)\!=\!0$.
Since the modified Sobolev norms are equivalent to
the standard ones away from the thin necks 
of~$(\Si_{\ups},g_{\ups})$,
\begin{equation}\label{inverse_p2_l1_e2}\begin{split}
\big\|R_{\ups,h_r}X\big\|_{\ups,p,1}
&\le C(b_{\ups})\Big(\sum_{r\in[m]}
\big|\psi_{\tilde{x}_r(\ups)}d(q_{\ups,h_r}^{-1}\circ q_N)\partial_s\big|
\Big)^{-1}
\big|\psi_zd(q_{\ups,h_r}^{-1}\circ q_N)\partial_s\big||X|_{\ups}\\
&\le C'(b_{\ups})\|R_{\ups,h_r}X\|_{\ups,2}.
\end{split}\end{equation}
By Lemma~\ref{approx_coker}, if $\xi\in\bar{\Ga}_-(\ups)$,
$$\xi=R_{\ups}X\equiv
\sum\limits_{r\in[m]}R_{\ups,h_r}X,$$
for some $X\!\in\! T_{\ev(b_{\ups})}\P$.
Thus, the first two statements of the lemma follow
from \e_ref{inverse_p2_l1_e2}.\\
(2) If $\ups_k\!\lra\! b$ and 
$\xi_k\!=\!R_{\ups}X_k\!\in\!\bar{\Ga}_-(\ups_k)$
is such that $\|\xi_k\|_{\ups_k,2}\!=\!1$, then it is immediate from (1)
that a subsequence of $\xi_k$ $C^0$-converges to 
$\sum\limits_{r\in[m]}c_rR_{b,h_r}X$, where
\begin{equation}\label{inverse_p2_l1_e4}
X=\lim_{k\lra\i} X_k,\quad
c_r=\lim_{k\lra\i}\Big(\sum_{r\in[m]}
\big|\psi_{\tilde{x}_r(\ups)}d(q_{\ups,h_r}^{-1}\circ q_N)
 \partial_s\big|\Big)^{-1}
\big(\psi_{\tilde{x}_r(\ups)}d(q_{\ups,h_r}^{-1}\circ q_N)\partial_s\big).
\end{equation}
The two limits in \e_ref{inverse_p2_l1_e4} exist after passing to 
a subsequence of the original sequence.
This proves the last statement of the lemma.

\begin{lmm}
\label{inverse_p2_l2}
There exist $\de,C\!\in\! C^{\i}({\cal M}_{\cal T}^{(0)};\Bbb{R}^+)$
such that for all $\ups\!\in\! F^{(\eset)}{\cal T}_{\de}$ and 
$\xi\!\in\!\tilde{\Ga}_-(\ups)$,
$$\|\xi\|_{\ups,p,1}\le C(b_{\ups})\|\xi\|_{\ups,2}.$$
\end{lmm}

\noindent
{\it Proof:}
Let $\Ga_{-+}(\ups)$ be the $(\ups,L^2)$-orthogonal complement
of $\pi_{\ups,-}\big(\bar{\Ga}_-(\ups)\big)$ in $\Ga_-(\ups)$.
Then
$$\tilde{\Ga}_-(\ups)=\Ga_{-+}(\ups)\oplus\bar{\Ga}_-(\ups).$$
Since this decomposition is $(L^2,\ups)$-orthogonal,
we can assume that either $\xi\!\in\!\Ga_{-+}(\ups)$ or 
$\xi\!\in\!\bar{\Ga}_-(\ups)$.
In the first case, the statement is obvious, since 
$\Ga_{-+}(\ups)\!\subset\!\Ga_-(\ups)$.
The second case is proved in Lemma~\ref{inverse_p2_l1}.

\begin{crl}
\label{inverse_p2_c3}
Suppose  $\ups_k\!\in\! F^{(0)}{\cal T}_{\de}$ and
$\ups_k\!\lra\! b\!\in\! {\cal M}_{\cal T}^{(0)}$.
If $\{\xi_{\ups_k,l}\}$ is an $(L^2,\ups)$-orthonormal basis
for $\tilde{\Ga}_-(\ups_k)$, then there exists
a nonzero $c\!\in\!\Bbb{C}^m$ and  
an $(L^2,b)$-orthonormal basis $\{\xi_{b,l}\}$ for $\tilde{\Ga}_-(b;c)$ 
such that after passing to a subsequence
$\xi_{\ups_k,l}$ $C^0$-converges to $\xi_{b,l}$ for all~$l$. 
\end{crl}

\noindent
{\it Proof:} 
If $\xi_{k,l}\!\in\!\bar{\Ga}_-(\ups_k)$, by Lemma~\ref{inverse_p2_l1}
a subsequence of  $\{\xi_{k,l}\}$ $C^0$-converges 
to an element of $\xi_l\!\in\!\bar{\Ga}_-(b;c)$
for some nonzero $c\!\in\!\Bbb{C}^n$ dependent on the sequence~$\{\ups_k\}$.
Furthermore, orthonormal pairs of such elements 
$C^0$-converge to an orthonormal pair in $\bar{\Ga}_-(b)$.
If $\xi_{k,l}\!\in\!\bar{\Ga}_{-+}(\ups_k)\!\subset\!\Ga_-(\ups_k)$,
then by definition of $\Ga_-(\ups_k)$, 
a subsequence of $\{\xi_{k,l}\}$ 
$C^0$-converge to an element $\xi_l\!\in\!\Ga_-(b)$,
which must be orthogonal to $\bar{\Ga}_-(b;c)$;
see~Lemma~\ref{gl-C0L2} in~\cite{Z}.
Thus, a subsequence of $\big\{\{\xi_{k,l}\}\big\}$ 
$C^0$-converges to an orthonormal set
of vectors in $\tilde{\Ga}_-(b)$, which implies that
$\dim_{\Bbb{C}}\tilde{\Ga}_-(b;c)\ge\dim_{\Bbb{C}}\tilde{\Ga}_-(\ups_k)$.
However,
\begin{alignat*}{2}
&\dim_{\Bbb{C}}\tilde{\Ga}_-(b;c)&&
=\dim_{\Bbb{C}}\Ga_{-+}(b;c)+\dim_{\Bbb{C}}\bar{\Ga}_-(b;c)\\
&&&=\dim_{\Bbb{C}}\Ga_-(b)+
\l(\dim_{\Bbb{C}}\bar{\Ga}_-(b;c)-
    \dim_{\Bbb{C}}\pi_{b,-}\bar{\Ga}_-(b;c)\r);\\
&\dim_{\Bbb{C}}\tilde{\Ga}_-(\ups_k)&&
=\dim_{\Bbb{C}}\Ga_{-+}(\ups_k)+\dim_{\Bbb{C}}\bar{\Ga}_-(\ups_k)\\
&&&=\dim_{\Bbb{C}}\Ga_-(\ups_k)+
\l(\dim_{\Bbb{C}}\bar{\Ga}_-(\ups_k)-
                         \dim_{\Bbb{C}}\pi_{\ups_k,-}\bar{\Ga}_-(\ups_k)\r),
\end{alignat*}
where $\Ga_{-+}(b;c)$ denotes the $(L^2,b)$-complement of
$\pi_{b,-}\bar{\Ga}_-(b;c)$ in $\Ga_-(b)$.
Since $\Ga_-(\ups_k)$ and $\Ga_-(b)$ have the same dimension,
in order to conclude the proof, it is sufficient to show that
$$\pi_{b,-}: \bar{\Ga}_-(b;c)\lra \Ga_-(b;c)$$
is an isomorphism; Lemma~\ref{inverse_p2_l4}.

\begin{lmm}
\label{inverse_p2_l4}
There exists $C\!\in\! C^{\i}({\cal M}_{\cal T}^{(0)};\Bbb{R}^+)$
such that for all $b\!\in\! {\cal M}_{\cal T}^{(0)}$\!,
nonzero $c\!\in\!\Bbb{C}^m$\!, \hbox{and $\xi\!\in\!\bar{\Ga}_-(b;c)$}
$$\|\xi\|_{b,2}\le C(b_{\om})\|\pi_{b,-}\xi\|_{b,2}.$$
\end{lmm}

\noindent
{\it Proof:} Let $X\!\in\! T_{\ev(b)}\P$.
Outside of $\i\!\in\!\Si_{b,h_r}$, by the same computation as 
in the proof of Lemma~\ref{approx_coker},
$R_{b,h_r}X\!=\!D_{b,u_{h_r}}^*R_{b,h_r}X$, where we define 
$R_{b,h_r}X\!\in\!\Ga^{0,1}(u_{h_r})$ as follows.
Let $d\bar{z}$ denote the usual $(0,1)$-form on $\Bbb{C}$.
Then
$$\de_{\cal T}(b) R_bX\big|_x=
4\be\big(\de_{\cal T}(b)|q_N^{-1}(x)|\big)
\big(q_N^{-1*}d\bar{z}\big) \Pi_{b,\bar{u}_{h_r}(x)}X.$$
Thus, if $\xi\!=\!\xi_{\hat{I}}\!\in\!\ker D_b$ and 
$2\de<\de_{\cal T}(b)$,
by integration by parts,
\begin{equation}\label{inverse_p2_l2_e2}\begin{split}
\llan\xi,R_{b,h_r}X\rran_b
&=\llan\xi_{h_r},D_{b,u_{h_r}}^*R_{b,h_r}X\rran_b\\
&=2i\de_{\cal T}(b)^{-1}\!\!\!
\int\limits_{|q_N^{-1}(x)|=\de^{-1}}\!\!\!
\lan\xi_{h_r}(x),\Pi_{b,\bar{u}_{h_r}(x)}X\ran_b q_N^{-1*}dz,
\end{split}\end{equation}
since $D_{b,u_{h_r}}\xi_{h_r}\!=\!0$.
Using the change of variables with $x=q_N(w^{-1})$, we obtain
\begin{equation}\label{inverse_p2_l2_e3}\begin{split}
\int\limits_{|q_N^{-1}(x)|=\de^{-1}}&
\lr{\xi_{h_r}(x),\Pi_{b,\bar{u}_{h_r}(x)}X}_b q_N^{-1*}dz
=-\int\limits_{|w|=\de}
\lan\xi_{h_r}|_{q_N(w^{-1})},\Pi_{b,\bar{u}_{h_r}(q_N(w^{-1}))}X
                   \ran_b \frac{dw}{w^2}\\
&= -2\pi i\frac{d}{dw}
\lr{\xi_{h^*}|_{q_N(w^{-1})},
      \Pi_{b,\bar{u}_{h_r}(q_N(w^{-1}))}X}_b
\Big|_{w=0}\\
&=-2\pi i\frac{d}{d\bar{z}}
\lr{\xi_{h_r}|_{q_S(z)},
      \Pi_{b,\bar{u}_{h_r}(q_S(z))} X}_b\Big|_{z=0}
=-2\pi i\Big\lan\frac{D}{ds}(\xi_{h_r}\circ q_S)\Big|_{z=0},X\Big\ran,
\end{split}\end{equation}
since $D_{b,u_{h_r}}\xi_{h_r}\!=\!0$.
It follows from \e_ref{inverse_p2_l2_e2} and \e_ref{inverse_p2_l2_e3}
that for any $\xi\!=\!\xi_{[M]}\!\in\!\ker(D_b)$,
\begin{equation}\label{inverse_p2_l2_e4}
\Big\lan\Big\lan\xi,\sum_{r\in[m]}c_rR_{b,h_r}X\Big\ran\Big\ran_b=
4\pi\de_{\cal T}(b)^{-1}\sum_{r\in[m]}c_r
\Big\lan\frac{D}{ds}(\xi_{h_r}\circ q_S)\Big|_{z=0},X\Big\ran.
\end{equation}
Along with Corollary~\ref{reg_crl2}, equations \e_ref{inverse_p2_l2_e4}
gives
\begin{equation}\begin{split}
\label{inverse_p2_l2_e5}
\Big\|\pi_{b,-}\sum_{r\in[m]}c_rR_{b,h_r}X\Big\|_{b,2}
&\ge C(b)|c_{r^*}|
\sup_{\xi_{[M]}\in\ker(D_b),\|\xi_{[M]}\|=1}
\Big\lan\frac{D}{ds}(\xi_{h_{r^*}}\circ q_S)\Big|_{z=0},X\Big\ran_b\\
&\ge C'(b)|c_{r^*}||X|
\ge C''(b)\Big\|\sum_{r\in[m]}c_rR_{b,h_r}X\Big\|_{b,2},
\end{split}\end{equation}
where $r^*\!\in\![m]$ is such that $|c_{r^*}|\!=\!\sup_r|c_r|$.
Since the right-hand side of \e_ref{inverse_p2_l2_e5}
must be a continuous function of $b$, 
the claim follows.\\

\noindent
The statement of Corollary~\ref{inverse_p2_c3} is precisely
Condition (1) of Definition~\ref{gl-obs_setup_dfn1}.
The other two conditions require that the rate of
change of the $(L^2,\ups)$-orthogonal projection onto
$\tilde{\Ga}_-(\ups)$ be controlled by a function of $b_{\ups}$ only.
This is a consequence of the convergence described
in the Corollary~\ref{inverse_p2_c3},
i.e.~we can use the same argument as described 
in the remark following Lemma~\ref{gl-inj_l1} in~\cite{Z},
but with $\Ga_-(b)$ replaced by 
the appropriate space~$\Ga_-(b;c)$ (depending on~$\ups$).

\subsection{Structure Theorem, $S=\Si$}
\label{si_str_sec}

\noindent
If ${\cal T}\!=\!(\Si,[N],I;j,d)$ 
is a simple bubble type and $\mu$ is an 
$N$-tuple of complex submanifolds of~$\P$ such that
the evaluation map,
$$\ev_{[N]}\equiv\ev_1\times\ldots\times\ev_N\!:
{\cal M}_{\cal T}\lra(\P)^N,$$
is transversal to $\mu_1\!\times\!\ldots\!\times\!\mu_N$,
${\cal M}_{\cal T}(\mu)$ is a complex submanifold of~${\cal M}_{\cal T}$.
Let ${\cal N}^{\mu}{\cal T}$ be its normal bundle.
If ${\cal S}$ is a complex submanifold of ${\cal M}$,
denote its normal bundle by ${\cal NS}$
and an identification of small neighborhoods of
${\cal S}$  in ${\cal N}{\cal S}$ and
in ${\cal M}_{\cal T}$ by~$\phi_{\cal S}$.
For any complex vector bundle $V\!\!\lra\!{\cal M}_{\cal T}$, 
we denote by $\Phi_{\cal S}$ an identification 
of $\phi_{\cal S}^*V$ and 
$\pi_{{\cal NS}}^*V$ such that its restriction
to the fibers over~${\cal S}$ is the identity.
We assume that $\Phi_{\cal S}$ preserves 
\hbox{$F^{\eset}{\cal T}\subset F{\cal T}$}.
Let
$$F^{\eset}{\cal S}=
\{(b,\vec{n},\ups)\!\in\! {\cal NS}\oplus F{\cal S}\!:
(b,\ups)\!\in\! F^{\eset}{\cal T}\big\}.$$
If $\ev_{[N]}|{\cal S}$ is transversal to 
$\mu_1\!\times\!\ldots\!\times\!\mu_N$,
${\cal S}(\mu)\!\equiv\!{\cal S}\cap{\cal M}_{\cal T}(\mu)$
is a complex submanifold of ${\cal S}$ with normal
bundle~${\cal N}^{\mu}{\cal T}$.
Let $\phi_{\cal S}^{\mu}$ and $\Phi_{\cal S}^{\mu}$
be the analogues of $\phi_{\cal S}$ and $\Phi_{\cal S}$
for the bundle 
\hbox{${\cal N}^{\mu}{\cal T}\!\lra\!{\cal S}(\mu)$}.
We assume the bundle ${\cal N}^{\mu}{\cal T}$ is normed.
We call the pair $(\Phi_{\cal S},\Phi_{\cal S}^{\mu})$
a {\it regularization} of ${\cal S}(\mu)$ if it satisfies
a certain minor compatibility condition.
For the purposes of this paper, it suffices to say
that once $\Phi_{\cal S}$ is chosen, 
it is a condition on~$\Phi_{\cal S}^{\mu}|F{\cal T}$;
see Subsection~\ref{gl-orient2_sec} in~\cite{Z} for details.
However, the exact nature of $\Phi_{\cal S}^{\mu}|F{\cal T}$
is irrelevant for our computational purposes.
Finally, we denote by $\bar{C}^{\i}_{(d;N)}(\Si;\mu)$
the space of all bubble maps 
$\big(\Si,[N],I;x,(j,y),u\big)$
such that $\sum_{i\in I}u_{i*}[\Si_{b,i}]\!=\!d\la$, 
where \hbox{$\la\!\in\! H_2(\P;Z)$} is the class of a line,
and $u_{j_l}(y_l)\!\in\!\mu_l$ for \hbox{all $l\!\in\![N]$}.

\begin{thm}
\label{si_str}
Suppose $d$ is a positive integer,
${\cal T}\!=\!(\Si,[N],I;j,\under{d})$ is a simple bubble type
with $d_{\hat{0}}\!=\!0$ and $\sum\limits_{i\in I}d_i\!=\!d$,
${\cal S}\subset{\cal M}_{\cal T}$ is a complex submanifold,
and
$$\nu\!\in\!\Ga^{0,1}\big(\Si\times\P;\La^{0,1}_{J,j}
                          \pi_{\Si}^*T^*\Si\otimes\pi_{\P}^*T\P\big)$$
is a generic section.
Let $\mu$ be an $N$-tuple of complex submanifolds of~$\P$ 
in general position of total codimension
$$\codim_{\Bbb{C}}\mu=d(n+1)-n(g-1)+N.$$
and $(\Phi_{\cal S},\Phi_{\cal S}^{\mu})$
a regularization of~${\cal S}(\mu)$.
Then for every precompact open subset $K$ of ${\cal S}(\mu)$,
there exist a neighborhood $U_K$ of $K$ in 
\hbox{$\bar{C}^{\i}_{(d;N)}(\Si;\mu)$} and $\de,\ep,C>0$ 
with the following property.
For every $t\!\in\!(0,\ep)$, there exist  a section 
$$\varphi_{{\cal S},t\nu}^{\mu}
\in\Ga\big(F^{\eset}{\cal S}_{\de}|K;
\pi_{F{\cal S}}^*{\cal N}^{\mu}{\cal S}\big),
\quad\hbox{with}\quad
\big|\varphi_{{\cal S},t\nu}^{\mu}(\ups)\big|_{b_{\ups}}\le
C\big(t+|\ups|^{\frac{1}{p}}\big),$$
and a sign-preserving bijection between
${\cal M}_{\Si,t\nu,d}(\mu)\cap U_K$ and
the zero set of the~section $\psi_{{\cal S},t\nu}^{\mu}$
defined~by
\begin{alignat*}{2}
\psi_{{\cal S},t\nu}^{\mu}\!\in\!\Ga\big(F^{\eset}{\cal S}_{\de}|K;
\pi_{F{\cal S}}^*({\cal H}_{\Si}^{0,1}\otimes\ev^*T\P)\big),&&~~&
\Pi_{b_{\ups},\phi_{\cal S}^{\mu}\varphi_{{\cal S},t\nu}^{\mu}(\ups)}
\psi_{{\cal S},t\nu}^{\mu}(\ups)\!=\!
\psi_{{\cal S},t\nu}\big(\Phi_{\cal S}^{\mu}
(\varphi_{{\cal S},t\nu}^{\mu}(\ups))\big);\\
\psi_{{\cal S},t\nu}\!\in\!\Ga
\big(F^{\eset}{\cal S}_{\de}\big|({\cal S}\cap U_K);
\pi_{F{\cal S}}^*({\cal H}_{\Si}^{0,1}\otimes\ev^*T\P)\big),&&~~&
\Pi_{b_{\ups},\phi_{\cal S}(\ups)}\psi_{{\cal S},t\nu}(\ups)\!=\!
\psi_{{\cal T},t\nu}\big(\Phi_{\cal S}(\ups)\big);\\
\psi_{{\cal T},t\nu}\!\in\!
     \Ga\big(F^{\eset}{\cal T}_{\de}\big|({\cal M}_{\cal T}\cap U_K);
\pi_{F{\cal T}}^*({\cal H}_{\Si}^{0,1}\otimes\ev^*T\P)\big),&&~~&
\tilde{R}_{\ups}\psi_{{\cal T},t\nu}(\ups)\!=\! \pi_{\ups,-}^{0,1}\!\big(
t\nu_{\ups,t}\!\!-\!\bar{\partial}u_{\ups}\!-\!
              D_{\ups}\xi_{\ups,t\nu}\!\big),
\end{alignat*}
where $\Pi_{b,b'}$ denotes the $g_{\P,b}$-parallel transport
along the $g_{\P,b}$-geodesics from $\ev(b)$ to $\ev(b')$
whenever $d_{\P}\big(\ev(b),\ev(b')\big)<r_{\P}$, 
$\xi_{\ups,t\nu}\!\in\!\tilde{\Ga}_+(\ups)$,
$$\big\|\nu_{\ups,t}-\nu\big\|_{\ups,2}\le
C\big(t+|\ups|^{\frac{1}{p}}\big),\quad\hbox{and}\quad
\big\|\xi_{\ups,t\nu}\big\|_{\ups,p,1}\le 
C\big(t+|\ups|^{\frac{1}{p}}\big).$$
\end{thm}

\noindent
{\it Proof:} This theorem follows immediately 
from Theorem~\ref{gl-si_str} in~\cite{Z} applied to
the obstruction bundle setup of 
Subsections~\ref{obs_sec} and~\ref{tangent_sec}.
The only refinement is that we drop the term $\tilde{\eta}_{\ups,t\nu}$
from the definition of~$\psi_{{\cal T},t\nu}$.
This is because it vanishes on the support of 
the $(0,1)$-forms in $\Ga^{0,1}_-(\ups)$, 
provided $\de$ is sufficiently small.
\hbox{Thus, $\pi^{0,1}_{\ups,-}\tilde{\eta}_{\ups,t\nu}\!=\!0$}.

\subsection{Structure Theorem, $S=S^2$}

\noindent
In this subsection, we define sections 
${\cal D}_{\lr{\cal T},k}^{(m)}$, where $k\!\in\! I\!-\!\hat{I}$, 
of the bubble $L_{{\cal T},k}^{*\otimes m}\otimes\ev^*T\P$ 
over $\bar{\cal U}_{\lr{\cal T}}(\mu)$,
and describe their behavior with respect to 
the gluing~maps near each space ${\cal U}_{\cal T}(\mu)$.
In Section~\ref{resolvent_sec}, 
the number of elements of ${\cal M}_{\Si,t\nu,d}(\mu)$ lying near each 
space ${\cal M}_{\cal T}(\mu)$ will be expressed
as the number of zeros of affine maps between certain bundles.
These affine maps will involve 
the sections~${\cal D}_{\bar{\cal T},k}^{(m)}$.
Their behavior near various boundary  is the foundations for
the local computations of Section~\ref{comp_sect}.\\

\noindent
If $b\!=\!\big(S^2,M,I;x,(j,y),u\big)\!\in\!
{\cal B}_{\cal T}$, $m\!\ge\!1$, and $k\!\in\! I$, let
$${\cal D}_{{\cal T},k}^{(m)}b=
\frac{2}{(m-1)!}
\frac{D^{m-1}}{ds^{m-1}}\frac{d}{ds}(u_k\circ q_S)
\Big|_{(s,t)=0},$$
where the covariant derivatives are taken with respect to the metric
$g_{\P,b}$ and $s+it\in\Bbb{C}$.
If ${\cal T}^*$ is a basic bubble type, 
the maps ${\cal D}_{{\cal T},k}^{(m)}$ with 
${\cal T}\!<\!{\cal T}^*$ and $k\!\in\!I\!-\!\hat{I}$
induce a continuous section of $\ev^*T\P$ over 
$\bar{\cal U}_{{\cal T}^*}^{(0)}$ and a continuous section of the bundle 
$L_{{\cal T}^*,k}^{*\otimes m}\otimes\ev^*T\P$ over 
$\bar{\cal U}_{{\cal T}^*}$, described~by 
$${\cal D}_{{\cal T}^*,k}^{(m)}[b,c_k]
=c_i^m{\cal D}_{{\cal T},k}^{(m)}b,
\quad\hbox{if}~~b\in{\cal U}_{\cal T}^{(0)},~c_k\in\Bbb{C}.$$ 
We will often write 
${\cal D}_{{\cal T},k}^{(1)}$ instead of~${\cal D}_{{\cal T},k}$.
If~${\cal T}$ is simple, we will abbreviate 
${\cal D}_{{\cal T},k}^{(m)}$ as~${\cal D}^{(m)}$.
If ${\cal T}\!=\!(\Si,[N],I;j,d)$ is a simple bubble type
and $k\!\in\!\hat{I}$, 
let ${\cal D}_{{\cal T},k}^{(m)}$ denote 
 the section~${\cal D}_{\bar{\cal T},k}^{(m)}$.

\begin{thm}
\label{str_global}
If ${\cal T}^*\!=\!(S^2,M,I^*;j,d)$ is a basic bubble type and
$\mu$ is an $M$-tuple of constraints in general position,
the spaces $\bar{\cal U}_{{\cal T}^*}^{(0)}(\mu)$ and 
$\bar{\cal U}_{{\cal T}^*}(\mu)$ are oriented topological orbifolds.
If ${\cal T}\!\le\!{\cal T}^*$, 
there exist \hbox{$G_{{\cal T}^*}$-invariant} functions
$\de,C\in C^{\i}\big({\cal U}_{\cal T}^{(0)}(\mu);\Bbb{R}^+\big)$
and $G_{{\cal T}^*}$-equivarent continuous map 
$$\tilde{\ga}_{\cal T}^{\mu}\!: 
     F{\cal T}_{\de}\Big|{\cal U}_{\cal T}^{(0)}(\mu)
                \lra \bar{\cal U}_{{\cal T}^*}^{(0)}(\mu),$$
which is an orientation-preserving homeomorphism onto
an open neighborhood of ${\cal U}_{\cal T}^{(0)}(\mu)$
in $\bar{\cal U}_{{\cal T}^*}^{(0)}(\mu)$ and 
is identity on ${\cal U}_{\cal T}^{(0)}(\mu)$.
This map is smooth on~$F^{\eset}{\cal T}_{\de}$.
Furthermore, for any
$$\ups=\big[(b,v_h)_{h\in\hat{I}}\big]=
\big[\big(S^2,M,I;x,(j,y),u\big),(v_h)_{v\in\hat{I}}\big]
\in F{\cal T}_{\de}\Big|{\cal U}_{\cal T}^{(0)}(\mu),$$
\begin{equation*}\begin{split}
\Big|\Pi_{b_{\ups},\ev(\tilde{\ga}_{\cal T}^{\mu}(\ups))}^{-1}
\big( {\cal D}_{{\cal T}^*,k}\!
\tilde{\ga}_{\cal T}^{\mu}(\ups)\!\big)-
2\sum_{h\in I_k,\chi_{\cal T}h=1}
\Big(\prod_{i\in\hat{I},h\in\bar{D}_i{\cal T}}\!\!\!\!v_i\!\Big) 
        \big(du_h|_{\i}e_{\i}\big)\Big|&\\
\le 
C(b_{\ups})|\ups|^{\frac{1}{p}}&\sum_{h\in I_k,\chi_{\cal T}h=1}\!\!
\Big(\prod_{i\in\hat{I},h\in\bar{D}_i{\cal T}}\!\!\!\!|v_i|\!\Big),
\end{split}\end{equation*}
where $I_k\!\subset\! I$ is the rooted tree containing $k$.
\end{thm}

\noindent
{\it Remark:} This theorem states that there exists
an identification 
$\ga_{\cal T}^{\mu}\!:{\cal FT}_{\de}\!\lra\!\bar{\cal U}_{{\cal T}^*}(\mu)$
of neighborhoods  of ${\cal U}_{\cal T}$ in ${\cal FT}$
and in~$\bar{\cal U}_{{\cal T}^*}(\mu)$.
Furthermore, with appropriate identifications,
\begin{gather}\label{str_global_e}
\Big|{\cal D}_{{\cal T}^*,k}\ga_{\cal T}^{\mu}(\ups)
-\al_{\cal T}\big(\rho_{\cal T}(\ups)\big)\Big|
\le C(b_{\ups})|\ups|^{\frac{1}{p}} \big|\rho_{\cal T}(\ups)\big|,
\qquad\hbox{where}\\
\rho_{\cal T}(\ups)=
\big(b,(\tilde{v}_h)_{\chi_{\cal T}h=1}\big)
\in\tilde{\cal F}{\cal T}\equiv
\bigoplus_{\chi_{\cal T}h=1}
L_h{\cal T}\otimes L_{\tilde{\io}_h}^*{\cal T};~~
\tilde{v}_h=\prod_{i\in\hat{I},h\in\bar{D}_i{\cal T}}\!\!\!\!v_i;~~
\tilde{\io}_h\!\in\!I-\hat{I},~ 
h\!\in\!\bar{D}_{\tilde{\io}_h}{\cal T};\notag\\
\al_{\cal T}\big(b,(\tilde{v}_h)_{\chi_{\cal T}h=1}\big)
=\sum_{h\in I_k,\chi_{\cal T}h=1}
{\cal D}_{{\cal T},h}\tilde{v}_h.\notag
\end{gather}
This estimate is used frequently in Section~\ref{comp_sect}.
Note that if ${\cal T}$ is a semiprimitive bubble type,
the bundle ${\cal FT}$ is defined over~$\bar{\cal U}_{\cal T}(\mu)$.
However, ${\cal FT}$ is {\it not} the normal bundle
of $\bar{\cal U}_{\cal T}(\mu)$ in $\bar{\cal U}_{\lr{\cal T}}(\mu)$
unless $M_{\hat{0}}{\cal T}\!+\!H_{\hat{0}}{\cal T}$ is a two-element set;
see~\cite{P2}.
The theorem implies only that the restrictions of the normal
bundle of $\bar{\cal U}_{\cal T}(\mu)$ in $\bar{\cal U}_{\lr{\cal T}}(\mu)$
and of $F{\cal T}$ to ${\cal U}_{\cal T}(\mu)$ are isomorphic.
\\

\noindent
{\it Proof:} (1) 
All statements of this theorem, except for the analytic estimate,
follow immediately from Theorem~\ref{gl-str_global} in~\cite{Z}.
We deduce the analytic estimate from (2) of Theorem~\ref{gl-str_global}.
Let 
$$\ga_{\cal T}^{\mu}(\ups)=
\big(S^2,M,I(\ups);x(\ups),(j(\ups),y(\ups)),\tilde{u}_{\ups}\big).$$
By Theorem~\ref{gl-str_global}, 
there exist a holomorphic bubble map 
$$b'\!=\!\big[S^2,M,I;x',(j,y'),u'\big]$$
such that $d_{C^k}(b,b')\le C(b_{\ups})|\ups|^{\frac{1}{p}}$ and
with appropriate identifications,
$\tilde{u}_{\ups}\!=\!\exp_{b',u_{b'}\circ q_{\ups}}\xi$
for some $\xi\!\in\!\Ga(u_{b'}\circ q_{\ups})$ with 
$\|\xi\|_{b,C^0}\!\le\! C(b_{\ups})|\ups|^{\frac{1}{p}}$.
Thus, for the purposes of proving the analytic estimate,
we can assume that
$u_{\ups}\!=\!\exp_{b,u_b\circ q_{\ups}}\xi_{\ups}$ 
for $\xi\!\in\!\Ga(u\circ q_{\ups})$  with 
$\|\xi_{\ups}\|_{b,C^0}\!\le\! C(b_{\ups})|\ups|^{\frac{1}{p}}$,
i.e.~it is enough to prove 
the estimate for the map $\tilde{\ga}_{\cal T}$
as defined in~\cite{Z} with ${\cal T}$ a simple bubble type.
If $d_k\!\neq\!0$, the claim is immediate from 
the usual Sobolev and elliptic estimates near~$(k,\i)$.
Thus, we assume that $d_{\hat{0}}\!=\!0$.
For future use, we obtain equations describing
the behavior of $\big({\cal D}^{(m)}\tilde{\ga}_{\cal T}(\ups)\big)$
for all~$m\!\ge1\!$.\\
(2)  We identify 
$B_{g_{\P,b}}\big(\ev(b),\frac{1}{2}r_{\P}\big)$
with an open subset of $\Bbb{C}^n$
via the $g_{\P,\ev(b)}$-parallel transport along the geodesics from~$\ev(b)$.
We assume that 
$\de\!\in\! C^{\i}({\cal U}_{\cal T}^{(0)};\Bbb{R}^+)$ satisfies
$$C(b)\de(b)^{\frac{1}{2p}}+\de(b)^{\frac{1}{2}}
      \Big(\sum_{i\in M}\|du_i\|_{b,C^0}\Big)
<\frac{1}{2}r_{\P}.$$
Let $q\!: B_1(0;\Bbb{C})\lra S^2$
be the local stretching map as in Subsection~\ref{gl-basic_gluing} of~\cite{Z}
with $\!v=\!1$
defined with respect to the standard metric on~$\Bbb{C}$.
Let $f_{\ups}=u_{\ups}\circ q$ and 
$\tilde{f}_{\ups}=\tilde{u}_{\ups}\circ q$.
We denote the usual complex coordinate on $\Bbb{C}$ by $z$.
For any $z\!\in\! B_1(0;\Bbb{C})$, 
let $i_{\ups}(z)$ be such that  $q_{\ups}(q(z))\!\in\!\Si_{b,i_{\ups}(z)}$. 
If $X\!\in\! T_{\ev(b)}\P$ and $m\!\ge\!1$,
define $R_{\ups}X\psi^{(m)}\!\in\!\Ga^{0,1}(\tilde{f}_{\ups})$~by 
$$R_{\ups}X\psi^{(m)}\big|_z=
\begin{cases}
X\bar{z}^{m-1}d\bar{z},& \hbox{if~}\chi_{\cal T}i_{\ups}(z)=0;\\
\be\big(\de(b_{\ups})|q_{\ups}(q(z))|\big)X\bar{z}^{m-1}d\bar{z},&  
                   \hbox{if~}\chi_{\cal T}i_{\ups}(z)=1;\\
0,&\hbox{if~}\chi_{\cal T}i_{\ups}(z)=2.
\end{cases}$$
Note that if $\chi_{\cal T}i_{\ups}(z)\!=\!0$, or 
$\chi_{\cal T}i_{\ups}(z)\!=\!1$
and $\be\big(\de(b_{\ups})|q_{\ups}(q(z))|\big)\neq0$,
$\tilde{f}_{\ups}(z)$ lies in 
$B_{g_{\P,b}}\big(\ev(b),\frac{1}{2}r_{\P}\big)$.
Thus, $R_{\ups}X\psi^{(m)}$ is well-defined.
We now compute $\llan\bar{\partial}\tilde{f}_{\ups},R_{\ups}X\psi^{(m)}\rran$
in two ways and compare the results.
First, note that the map
$\tilde{f}_{\ups}$ is holomorphic outside of the annulus
$$A_{\hat{0}}(\ups)\equiv
B_1(0;\Bbb{C})-B_{\frac{1}{2}}(0;\Bbb{C}).$$
Thus, by the same computation as in the proof of Lemma~\ref{dbar_gen2},
we see that
\begin{equation}\label{limit_lmm_e1}
\llan\bar{\partial}\tilde{f}_{\ups},R_{\ups}X\psi^{(m)}\rran
=-\frac{\pi}{m}
\Big\lan \big({\cal D}^{(m)}\tilde{\ga}_{\cal T}(\ups)\big),X\Big\ran.
\end{equation}
(3) Since $\tilde{f}_{\ups}=\exp_{\ev(b),f_{\ups}}(\xi_{\ups}\circ q)$
and $f_{\ups}$ is constant on $A_{\hat{0}}(\ups)$,
\begin{equation}\label{limit_lmm_e2}
2i\llan\bar{\partial}\tilde{f}_{\ups},R_{\ups}X\psi^{(m)}\rran=
\int_{A_{\hat{0}}(\ups)}\Big\lan
\frac{\bar{\partial}}{\partial\bar{z}}(\xi_{\ups}\circ q),X\Big\ran
z^{m-1}d\bar{z}\w dz
\end{equation}
Denote by $A_{\hat{0}}^+(\ups)$ and  $A_{\hat{0}}^-(\ups)$ the outer and 
inner boundary of $A_{\hat{0}}(\ups)$, respectively.
For every $h\!\in\!\hat{I}$ with $\chi_{\cal T}h\!=\!1$, let
$$A_h(\ups)=q_{\ups,\io_h}^{-1}\Big(\big\{
 z\!\in\!\Si_{b_{\ups},\io_h}\!: 
4\de(b_{\ups})^{-1}|v_h|\!\le\!|\phi_{b,h}^{-1}z|\!\le\!|v_h|^{\frac{1}{2}}
\big\}\Big)\subset\Si_{b_{\ups},\hat{0}}.$$
Denote by $A_h^{\pm}(\ups)$ the outer and inner boundary of $A_h(\ups)$.
Let $w$ be the complex coordinate on 
\hbox{$\Bbb{C}\!\subset\!\Si_{b_{\ups},\hat{0}}\!=\!S^2$.} 
Note that $q$ is holomorphic inside of $A_{\hat{0}}^-(\ups)$ and 
outside of $q^{-1}(A_h^-(\ups))$.
Furthermore, since $u_{b}$ and $\tilde{u}_{\ups}$ are both holomorphic,
on the image of this set under $q$
$$\frac{\bar{\partial}}{\partial\bar{w}}\xi_{\ups}=
-\frac{\bar{\partial}}{\partial\bar{w}}u_{\ups}.$$
The last quantity vanishes outside of the annuli $A_h(\ups)$.
Thus by integration by parts,
\begin{equation}\label{limit_lmm_e3}\begin{split}
&\int_{A_{\hat{0}}(\ups)} \Big\lan
\frac{\bar{\partial}}{\partial\bar{z}}(\xi_{\ups}\circ q),X\Big\ran
z^{m-1}d\bar{z}\w dz\\ 
&~=\sum_{\chi_{\cal T}h=1} 
          \bigg(  \int\limits_{~q^{-1}(A_h(\ups))}\!\!\! \Big\lan
 \Big( \frac{\bar{\partial}u_{\ups}}{\partial\bar{w}}\Big)
           \ov{\Big(\frac{\partial q}{\partial z}\Big)}
 ,X\Big\ran z^{m-1}d\bar{z}\w dz
+\int\limits_{q^{-1} 
     (A_h^-(\ups))}\!\!\! \big\lan \xi_{\ups}\circ q,X\big\ran 
      z^{m-1}dz\bigg)\\
&~=\sum_{\chi_{\cal T}h=1}\bigg(
 \int_{A_h(\ups)}\Big\lan\frac{\bar{\partial}u_{\ups}}{\partial\bar{w}}
 ,X\Big\ran gd\bar{w}\w dw
+\int_{A_h^-(\ups)} \big\lan \xi_{\ups},X\big\ran gdw\bigg),
\end{split}\end{equation}
where $g(w)\!=\!w^{m-1}$.
Since $\xi_{\ups}\circ q$ is constant on $A_{\hat{0}}^+(\ups)$,
the second boundary term is zero.
Note that the radius of $A_h^-(\ups)$ in $\Bbb{C}\!\subset\! S^2$
is bounded by $C(b_{\ups})|\tilde{v}_h|$.
Furthermore, $|g|\!\le\! C_m(b_{\ups})$ on~$A_h^-(\ups)$.
It follows~that 
\begin{equation}\label{limit_lmm_e4}
\Big|\int_{A_h^-(\ups)} \big\lan \xi_{\ups},X\big\ran gdw\Big|\le
C_m(b_{\ups})|\ups|^{\frac{1}{p}}|\tilde{v}_h|.
\end{equation}
On the other hand, by the same computation as in
the proof of Lemma~\ref{dbar_gen2},
\begin{equation}\label{limit_lmm_e5}
\int_{A_h(\ups)}\Big\lan\frac{\bar{\partial}u_{\ups}}{\partial\bar{w}}
 ,X\Big\ran  gd\bar{w}\w dw
=-2i\sum_{m'=1}^{m'=m}\frac{\pi a_{m',h}(\ups)}{m'}
 \tilde{v}_h^{m'}\big({\cal D}_{{\cal T},h}^{(m')}b\big),
\end{equation}
where the numbers $a_{m',h}(\ups)\!\in\!\Bbb{C}$ are defined by
$$g\big(q_{\ups,\io_h}^{-1}(\io_h,y)\big)=
\sum_{m'=1}^{m'=m} a_{m',h}(\ups)
\cdot\big(\phi_{b_{\ups},h}y\big)^{m'-1}$$
whenever $y\!\in\! S^2$ is close to $x_h$.
Combining equations \e_ref{limit_lmm_e1}-\e_ref{limit_lmm_e5},
we see that 
\begin{equation}\label{limit_lmm_e6}
\Big|
\Big\lan {\cal D}^{(m)}\tilde{\ga}_{\cal T}(\ups),X\Big\ran
-2m\sum_{\chi_{\cal T}h=1}a_{1,h}(\ups)
\tilde{v}_h \big( du_h\big|_{\i}e_{\i}\big)\Big|
\le C(b_{\ups})|\ups|^{\frac{1}{p}}
\Big(\sum_{\chi_{\cal T}h=1}|\tilde{v}_h|\Big).
\end{equation}
Since $a_{1,h}(\ups)=\big(q_{\ups,\io_h}(\io_h,y)\big)^{m-1}$,
equation~\e_ref{limit_lmm_e6} gives the analytic estimate
of the theorem.\\

\section{Topology}
\label{top_sect}

\subsection{Maps Between Vector Bundles}
\label{top_sec1}

\noindent
In Section~\ref{resolvent_sec}, we express the number of zeros of 
the maps $\psi_{{\cal T},t\nu}^{\mu}$ 
(and $\psi_{{\cal S},t\nu}^{\mu}$ for certain submanifolds
${\cal S}$ of ${\cal M}_{\cal T}$) of Theorem~\ref{si_str}
in terms of the number of zeros of affine maps between
the same vector bundles.
The topological justification for this reduction are discussed
in this subsection.
Subsections~\ref{top_sec2} and~\ref{zeros_sec} are used in the explicit
computations of Section~\ref{comp_sect}.
For simplicity, we state all the results for smooth vector bundles
over smooth manifolds, but similar statements apply in the 
orbifold category.
However, in the cases of $g\!=\!2$, $n\!=\!2,3,4$, the spaces
involved are actually manifolds.
\\

\noindent
Let ${\cal I}$ denote the unit interval $[0,1]$.
If ${\cal Z}$ is a compact oriented zero-dimensional manifold,
we denote the signed cardinality of~${\cal Z}$
by~$^{\pm}|{\cal Z}|$.

\begin{dfn}
\label{top_dfn1}
Suppose ${\cal M}$ is a smooth manifold and
$F\!=\!\bigoplus\limits_{i=1}^{i=k}F_i,{\cal O}\!\lra\!{\cal M}$
are smooth normed complex vector bundles.\\
(1) A bundle map $\al\!:F\!\lra\!{\cal O}$ 
is a \under{polynomial of degree~$d_{[k]}$} if
for each $i\!\in\![k]$ there  exists
$$p_i\!\in\!\Ga({\cal M};F_i^{*\otimes d_i}\!\otimes\!{\cal O})
\hbox{~~for~}i\!\in\![k]
\qquad\hbox{s.t.}\qquad
\al(\ups)=\sum_{i=1}^{i=k}p_i\big(\ups_i^{d_i}\big)
\quad\forall \ups=(v_i)_{i\in[k]}\in\bigoplus_{i=1}^{i=k}F_i.$$
(2) If $\al\!:F\!\lra\!{\cal O}$ is a polynomial, 
the \under{rank of~$\al$} is the number
$$\rk~\al\equiv\max\{\rk_b\al\!: b\!\in\!{\cal M}\},
\quad\hbox{where}\quad
\rk_b\al=\dim_{\Bbb{C}}\big(\hbox{Im}~\al_b\big).$$
Polynomial $\al\!:F\!\lra\!{\cal O}$ is \under{of constant rank}
if $\rk_b\al=\rk~\al$ for all $b\!\in\!{\cal M}$;
$\al$ is \under{nondegenerate}
if \hbox{$\rk_b\al=rk~F$} for all $b\!\in\!{\cal M}$.\\
(3) If $\Om$ is an open subset of ${\cal I}\!\times\! F$ and 
$$\{\phi_t\}= \big\{\phi_t\!:
\{\ups\!\in\!F :(t,\ups)\!\in\!\Om\}\lra{\cal O}\big\}$$ 
is a family of smooth bundle maps, bundle map $\al\!:F\!\lra\!{\cal O}$ is
a \under{dominant term of $\{\phi_t\}$} 
if there exist \hbox{$\de\!\in\! C^{\i}({\cal M};\Bbb{R}^+)$}
and $\ve\!\in\! C^0({\cal I}\!\times\! F; \Bbb{R})$ such that
$$\big|\phi_t(\ups)-\al(\ups)\big|\le
\ve(t,\ups)\big(t+|\al(\ups)|\big)
\quad\forall(t,\ups)\in\Om_{\de}
\quad\hbox{and}\quad
\lim_{(t,\ups)\lra0}\ve(t,\ups)=0.$$
Dominant term $\al\!:F\!\lra\!{\cal O}$ of $\{\phi_t\}$ is 
the \under{resolvent} of $\{\phi_t\}$ if $\al$ is a polynomial
of constant rank.
\end{dfn}

\noindent 
In (2) above, by $\dim_{\Bbb{C}}\big(\hbox{Im}~\al_b\big)$
we mean the dimension of the image of $\al_b$ as 
an analytic subvariety of the fiber~${\cal O}_b$.
Note that if $\bar{\Om}\!\subset\! {\cal I}\!\times\! F$ contains 
$\{0\}\!\times\!{\cal M}$, 
the resolvent of $\{\phi_t\}$ is unique (if it exists).

\begin{lmm}
\label{top_l1}
Suppose ${\cal M}$ is a smooth manifold,\\
(1) $F\equiv F^-\!\oplus F^+$
and ${\cal O}\equiv{\cal O}^-\!\oplus {\cal O}^+$ are
smooth normed complex vector bundles over~${\cal M}$;\\
(2) $\Om$ is an open subset of ${\cal I}\!\times\! F$ and
$\big\{\phi_t\!: \{\ups\!\in\! F\!:(t,\ups)\!\in\!\Om\}\!\lra\!{\cal O}\big\}$
is a family of smooth maps;\\
(3) $\al\!:F\!\lra\!{\cal O}$ is a dominant term of $\{\phi_t\}$
such that $\al(F^+)\!\subset\!{\cal O}^+$,
$\al^-\!\equiv\!\pi^-\!\circ(\al|F^-)$ is a constant-rank polynomial,
where $\pi^-\!:{\cal O}^-\!\oplus {\cal O}^+\!\!\lra\!{\cal O}^-$
is the projection map, and
\hbox{$(\dim{\cal M}\!+\!2\rk~\al^-)\!<\!2\rk~{\cal O}^-$};\\
(4) $\bar{\nu}\!=\!(\bar{\nu}^-,\bar{\nu}^+)
       \!\in\!\Ga({\cal M};{\cal O}^-\!\oplus{\cal O}^+)$
is generic with respect to~$\al^-$.\\
Then for every compact subset $K$ of ${\cal M}$,
there exists $\de_K>0$ such that the map
$$\psi_t\!: \{\ups\!\in\! F\!: (t,\ups)\!\in\!\Om\}\lra{\cal O},
\quad\psi_t(\ups)=t\bar{\nu}_{\ups}+\phi_t(\ups),$$
has no zeros on $\{\ups\!\in\! F_{\de_K}|K\!: (t,\ups)\!\in\!\Om\}$
for all $t\!\in\!(0,\de_K)$.
\end{lmm}

\noindent
{\it Proof:} (1) 
Suppose $\tilde{\ups}\!\in\!\Om_{\de_K}|K$ and $\psi_t(\tilde{\ups})=0$.
Then by our assumptions on $\psi_t$,
$$|\al(\tilde{\ups})|\le C_K
\big(t+\bar{\ve}_K(\de_K)\big|\al(\tilde{\ups})|\big),$$
where $C_K>0$ depends only on $K$ (and $\bar{\nu}$)
and $\bar{\ve}_K$ is a continuous function vanishing at~zero.
Thus, if $\de_K>0$ is sufficiently small,
\begin{equation}
\label{top_l1e1}
|\al(\tilde{\ups})|\le 2C_Kt
\quad\forall t<\de_K,
~\tilde{\ups}\!\in\! F_{\de_K}|K\hbox{~s.t.~}\psi_t(\tilde{\ups})=0.
\end{equation}
(2) Let $F^-\!=\!\bigoplus\limits_{i=1}^{i=k}F_i\!\lra\!{\cal M}$
be the bundles and
\hbox{$p_i\!\in\!\Ga({\cal M};F_i^{*\otimes d_i}\!\otimes\!{\cal O}^-)$}
the sections as in (1) of Definition~\ref{top_dfn1}
corresponding to~$\al^-$.
Define
$$\varphi_t\in\Ga\big({\cal M};\hbox{End}(F^-,F^-)\big)
\quad\hbox{by}\quad
\varphi_t(\ups_i)=t^{-1/d_i}\ups_i
\hbox{~~if~}\ups
_i\!\in\!F_i.$$
Then by our assumption on $\psi_t$ and equation~\e_ref{top_l1e1},
\begin{equation}
\label{top_l1e2}
\big|\bar{\nu}^-+
\al^-\big(\varphi_t(\tilde{\ups}^-)\big)\big|\le 
  \tilde{C}_K\bar{\ve}_K(\de_K)
\quad\forall t<\de_K,
~\tilde{\ups}\!\in\! F_{\de_K}|K\hbox{~s.t.~}\psi_t(\tilde{\ups})=0,
\end{equation}
where $\tilde{C}_K$ is determined by $K$.
Since $\al^-$ has constant rank, the image of $\al^-$ is closed
and is the total space of a bundle of affine analytic varieties
of complex dimension
\hbox{$\rk~\al^-\!\!<\rk~{\cal O}^-\!\!-\!\frac{1}{2}\dim{\cal M}$}.
Thus, by  assumption (4) of the lemma, 
$\bar{\nu}^-$ does not intersect the image of~$\al^-$, and
there exists \hbox{$\ep_K>0$} such~that
\begin{equation}
\label{top_l1e3}
\big|\bar{\nu}^-+\al^-(\ups^-)\big|\ge\ep_K
\quad\forall \ups\!\in\! F^-|K.
\end{equation}
If $\ep_K>\tilde{C}_K\bar{\ve}_K(\de_K)$, by \e_ref{top_l1e2}
and \e_ref{top_l1e3}, $\pi^-\circ \psi_t$ (and thus $\psi_t$)
has no zeros on $F_{\de_K}|K$.\\

\noindent
We will call family
$\big\{\phi_t\!: \{\ups\!\in\! F\!:
(t,\ups)\!\in\!\Om\}\!\lra\!{\cal O}\big\}$
of smooth maps {\it hollow}
if it admits a dominant term~$\al$ that satisfies 
hypothesis (3) of Lemma~\ref{top_l1}.

\begin{dfn}
\label{top_dfn2}
Suppose ${\cal M}$ is a smooth manifold and
$F\!\lra\!{\cal M}$ is a smooth vector bundle.\\
(1) Subset $Y$ of $F$ is \under{small} if $Y$ contains no fiber of $F$ and 
there exists a smooth manifold~$Z$ of dimension~$(\dim F\!-\!1)$
and a smooth map $f\!:Z\lra F$ such that the image of $f$
is closed in $F$ and contains~$Y$.\\	
(2) If $F\!,\tilde{F}\!\!\lra\!{\cal M}$ are smooth complex vector bundles,
$\rho\!\in\!\Ga({\cal M};F^{*\otimes d}\otimes \tilde{F})$
\under{induces a $\tilde{d}$-to-$1$ cover}
\under{$F\!\!\lra\!\tilde{F}$}
if the map 
$$F_b\lra\tilde{F}_b,\quad\ups\lra\rho(\ups)\equiv
\rho\big(\ups^d\big),$$
is $\tilde{d}$-to-$1$ on a  dense open subset of
every fiber~$F_b$ of~$F$.
\end{dfn}

\begin{lmm}
\label{top_l2l}
Suppose ${\cal M}$ is a smooth manifold,
$F\!=\!\bigoplus\limits_{i=1}\limits^{i=k}F_i$ and ${\cal O}$
are smooth complex vector bundles over~${\cal M}$,
and 
$$\al=\sum_{i=1}^{i=k}p_i:F\lra{\cal O},
\quad\hbox{where}\quad
p_i\!\in\!\Ga\big({\cal M};F_i^{*\otimes d_i}\otimes{\cal O}\big),$$
is a nondegenerate polynomial.
Then there exists a small subset $Y_{\al}$
of $F=\bigoplus\limits_{i=1}\limits^{i=k}F_i$,
which is invariant under scalar multiplication in each component
separately, with the following property.
If $K$ is a compact subset of ${\cal O}-\al(Y_{\al})$,
there exists $C_K>0$ such that 
$$|\ups|\le C_K|\al(\ups)|
\quad\forall~ \ups\!\in\! F
\hbox{~~s.t.~}\al(\ups)\!\in\! K.$$
\end{lmm}

\noindent
{\it Proof:} (1) 
Let $Y_{\al}\!\subset\! F$ be the closed subset on which 
the differential of the fiberwise map $\ups\!\lra\!\al(\ups)$ 
does not have full rank, i.e.~its rank is less than~$\rk~F$.
Since~$\al$ is nondegenerate, $Y_{\al}$ contains no fiber of~$F$.
By our assumptions on~$\al$,
$$D(\al|_{F_b})\big|_{\ups}
=\big( D(p_1|_{F_{1,b}})\big|_{\ups_1},\ldots,
D(p_k|_{F_{k,b}})\big|_{\ups_k}\big)\!:
F_1\oplus\ldots\oplus F_k\lra{\cal O},
\quad\forall~b\!\in\!{\cal M},~
\ups\!=\!\ups_{[k]}\!\in\! \bigoplus\limits_{i=1}^{i=k}F_i.$$
Since $p_i|_{F_{i,b}}$ is a homogeneous polynomial of 
degree~$d_i$, its derivative is a homogeneous polynomial
of degree~$(d_i\!-\!1)$.
Thus, $Y_{\al}$ is preserved
under scalar multiplication in each component separately.
It also clearly satisfies the second condition of (1)
of Definition~\ref{top_dfn2}.\\
(2) On $F-Y_{\al}$, $\al$ is a covering map onto its image
with the number of leaves bounded by some number~$N_{\al}$.
Thus, if $K$ is any compact subset of ${\cal O}-\al(Y_{\al})$,
$\al^{-1}(K)$ is a compact subset of~$F$.
Therefore, there exists $C_K$ such~that 
$$|\ups|\le C_K|\al(\ups)| \quad\forall
\ups\!\in\! F\hbox{~~s.t.~}\al(\ups)\!\in\! K.$$
Note that if $0\!\not\in\!\al(Y_{\al})$,
then $\al$ is a linear injection on every fiber,
and the above inequality holds on all of~$F$.

\begin{lmm}
\label{top_l2}
Suppose ${\cal M}$ is a smooth manifold,\\
(1) $F\!=\!\bigoplus\limits_{i=1}^{i=k}F_i$ and ${\cal O}$ are
smooth normed complex vector bundles over~${\cal M}$
with \hbox{$\rk~F\!+\!\frac{1}{2}\dim{\cal M}\!=\!\rk~{\cal O}$};\\
(2) $Y$ is a small subset of $F\!=\!\bigoplus\limits_{i=1}^{i=k}F_i$,
which is invariant under the scalar multiplication in each component 
separately;\\
(3) $\Om$ is an open subset of ${\cal I}\!\times\! F$  such that 
$\Om\cup (\{0\}\!\times\!X)$  is a neighborhood
of $\{0\}\!\times\! X$ \hbox{in~${\cal I}\!\times\! \big(F\!-\!(Y\!-\!X)\big)$};\\
(4)
$\big\{\phi_t\!: \{\ups\!\in\! F\!:(t,\ups)\!\in\!\Om\}\!\lra\!{\cal O}\big\}$
is a family of smooth maps;\\
(5) nondegenerate polynomial
$\al\!:F\!\lra\!{\cal O}$ is the resolvent of $\{\phi_t\}$;\\
(6) $\bar{\nu}\!\in\!\Ga({\cal M};{\cal O})$
is generic with respect to~$(Y,\al)$, and the map
\begin{equation}\label{top_l2_e}
F\lra {\cal O},\quad
\ups\lra\bar{\nu}_{\ups}+\al(\ups),
\end{equation}
has a finite number of (transverse) zeros.\\
If $\psi_t$ is transversal to zero for all $t$,
there exists a compact subset $K_{\al,\bar{\nu}}$ of~${\cal M}$
with the following property.
If~$K$ is a precompact open subset of~${\cal M}$ containing 
$K_{\al,\bar{\nu}}$, 
there exist $\de_K,\ep_K>0$ such that for all
$t\!\in\!(0,\ep_K)$, 
$$^{\pm}\big|\big\{\ups\!\in\! F_{\de_K}|K\!: 
(t,\ups)\!\in\!\Om, \psi_t(\ups)=0\big\}\big|=
^{\pm}\big|\big\{\ups\!\in\! F\!:\bar{\nu}_{\ups}+\al(\ups)=0\big\}\big|,$$
where $\psi_t(\ups)=t\bar{\nu}_{\ups}+\phi_t(\ups)$ as before.
Furthermore, all the zeros of $\psi_t\big|(F_{\de_K}|K)$
lie over~$K_{\al,\bar{\nu}}.$
\end{lmm}

\noindent
{\it Proof:}
(1) Since the map in \e_ref{top_l2_e} has a finite number
of zeros, all of them lie in 
the interior of $F_{C_{\al,\bar{\nu}}}|K_{\al,\bar{\nu}}$
for some compact subset $K_{\al,\bar{\nu}}$ of~${\cal M}$
and number~$C_{\al,\bar{\nu}}\!>\!0$.
Suppose $K\!\subset\!{\cal M}$ is a precompact open subset containing
$K_{\al,\bar{\nu}}$, $\de_K\!>\!0$ is such that
\hbox{$\big(F_{\de_K}|K\!\!-\!Y\big)\!\subset\!\Om$,} and
$\tilde{\ups}\!\in\!\Om_{\de_K}|K$ is such 
\hbox{that $\psi_t(\tilde{\ups})=0$.}
By the same argument as in the proof of 
Lemma~\ref{top_l1}, if \hbox{$\de_K\!>\!0$} is sufficiently small,
\begin{equation}\label{top_l2_e2}
\big|\al(\tilde{\ups})\big|\le C_Kt
\quad\hbox{and}\quad
\big|t\bar{\nu}_{\tilde{\ups}}+\al(\tilde{\ups})\big|\le \bar{\ve}_K(\de_K)t
\quad\forall t<\de_K,
~\tilde{\ups}\!\in\! F_{\de_K}|K\hbox{~s.t.~}\psi_t(\tilde{\ups})=0,
\end{equation}
where $C_K$ and $\bar{\ve}_K=\bar{\ve}_K(\de_K)$ depend only on $K$,
and $\bar{\ve}_K(\de_K)$ tends to zero with~$\de_K$.
Let \hbox{$\phi_t\!:F\!\lra\! F$} be the map defined in (2) of 
the proof of Lemma~\ref{top_l1}, with $F^-$ replaced by~$F$.
By~\e_ref{top_l2_e2},
\begin{equation}\label{top_l2_e5}\begin{split}
\al\big(\phi_t(\tilde{\ups})\big) \!\in\!
{\cal K}_{\bar{\nu}}\big(K;C_K,\bar{\ve}_K(\de_K)\big)\equiv
\big\{\vp\!\in\!{\cal O}_{C_K}\!\!: 
|\bar{\nu}_{\vp}+\vp|\!\le\!\bar{\ve}_K(\de_K)\big\}\qquad\quad&\\
\forall t<\de_K,
~\tilde{\ups}\!\in\! F_{\de_K}|K\hbox{~s.t.~}\psi_t(\tilde{\ups})=0.&
\end{split}\end{equation}
(2) If $\bar{\nu}$ is generic, the map in~\e_ref{top_l2_e} does not vanish
on~$Y_{\al}$,  where $Y_{\al}$ is as in Lemma~\ref{top_l2l}.
Since $\al(Y_{\al})$ is a closed subset of~${\cal O}$,
there exists $\ep_K>0$ such~that
$$\big|\bar{\nu}_{\ups}+\al(\ups)\big|>\ep_K
\quad\forall~\ups\!\in\!Y_{\al}|K.$$
Thus, if $\bar{\ep}_K(\de_K)\!<\!\ep_K$, 
${\cal K}_{\bar{\nu}}\big(K;C_K,\bar{\ve}_K(\de_K)\big)$
is a compact subset of ${\cal O}$ disjoint from~$\al(Y_{\al})$.
Then by \e_ref{top_l2_e5} and Lemma~\ref{top_l2l},
\begin{equation}\label{top_l2_e6}
\big|\phi_t(\tilde{\ups})\big|\le C_K^*
\quad\forall t<\de_K,
~\tilde{\ups}\!\in\! F_{\de_K}|K\hbox{~s.t.~}\psi_t(\tilde{\ups})=0,
\end{equation}
where $C_K^*$ depends only on $K$.\\
(3) There is a one-to-one sign-preserving correspondence between
the  zeros of $\psi_t$ on $\Om_{\de_K}|K$ and the zeros of
$$\tilde{\psi}_t\!:
\Om_{\de_K}(K,t)\equiv\big\{ \ups\!\in\! F\!:
\big(t,\phi_t^{-1}(\ups)\big)\!\in\!
          \Om_{\de_K}|K\Big\}\lra{\cal O},\quad
\tilde{\psi}_t(\ups)=
t^{-1}\psi_t\big(\phi_t^{-1}(\ups)\big).$$
By~\e_ref{top_l2_e6}, all the zeros of $\tilde{\psi}_t$ on 
$\Om_{\de_K}(K,t)$ are in fact contained in~$F_{C_K^*}|K$.
We can assume that $C_K^*>C_{\al,\bar{\nu}}$.
By our assumptions on~$\phi_t$,
\begin{equation}\label{top_l2_e9}
\big|\tilde{\psi}_t(\ups)-\big(\bar{\nu}_{\ups}+\al(\ups)\big)\big|
\le C_K\bar{\ep}_K(\de_K) \quad
\forall\ups\in\Om_{\de_K}(K,t)\cap (F_{C_K^*}|K),  
\end{equation}
where $C_K>0$ depends only on~$K$.
We define a cobordism between the zeros of $\tilde{\psi}_t$
and the zeros of $\bar{\nu}+\al$ on 
$\Om_{\de_K}(K,t)\cap (F_{C_K^*}|K)$ by
$$\Psi\!: {\cal I}\!\times\! \Om_{\de_K}(K,t)\cap (F_{C_K^*}|K)
\lra{\cal O},\quad
\Psi_{\tau}(\ups)=\tau\tilde{\psi}_t(\ups)+
(1-\tau)  
 \big(\bar{\nu}_{\ups}+\al(\ups)\big)
+\eta_{\tau}(\ups),$$  
where $\eta\!: {\cal I}\!\times\! \Om_{\de_K}(K,t)\!\lra\!{\cal O}$ 
is any smooth function with very small $C^0$-norm such that
\hbox{$\eta_0\!=\!\eta_1\!=\!0$} and 
$\Psi$ is transversal to zero.
It remains to see that $\Psi^{-1}(0)$ is compact.
Suppose $\Psi_{\tau_r}(\ups_r)\!=\!0$ and $(\tau_r,\ups_r)$
converges $(\tilde{\tau},\tilde{\ups})\!\in \!
                    {\cal I}\!\times\! F_{2C_K^*}|\bar{K}$;
we need to show that 
\hbox{$\tilde{\ups}\!\in\!\Om_{\de_K}(K,t)\cap (F_{C_K^*}|K)$.}
By equation~\e_ref{top_l2_e9},
\begin{equation}\label{top_l2_e11}
\big|\bar{\nu}_{\ups_r}+\al(\ups_r)\big|
\le C_K\bar{\ep}_K(\de_K)+\|\eta\|_{C^0}~~\forall r\Lra
\big|\bar{\nu}_{\tilde{\ups}}+\al(\tilde{\ups})\big|
\le C_K\bar{\ep}_K(\de_K)+\|\eta\|_{C^0}.
\end{equation}
On the other hand, since $\bar{\nu}$ is generic, 
the map in~\e_ref{top_l2_e} does not vanish on~$Y$.
Furthermore, all the zeros of this map are contained 
in the interior of~$F_{C_{\al,\bar{\nu}}}|K_{\al,\bar{\nu}}$.
Thus, by compactness,
\begin{equation}\label{top_l2_e12}
\tilde{\ep}_K\!\equiv\!\inf
\Big\{\big|\bar{\nu}_{\ups}+\al(\ups)\big|\!:
\ups\!\in\! \big(Y\cap F_{2C_K^*}\big)\cup
\big(F_{C_K^*}|K\!-\!F_{C_{\al,\bar{\nu}}}|K_{\al,\bar{\nu}}\big)
\Big\}\!>\!0,
\end{equation}
where $\tilde{\ep}_K$ depends only on~$K$.
If $\tilde{\ep}_K>C_K\bar{\ep}_K(\de_K)+\|\eta\|_{C^0}$,
by \e_ref{top_l2_e11} and~\e_ref{top_l2_e12},
$$\tilde{\ups}\in F_{C_{\al,\bar{\nu}}}|K_{\al,\bar{\nu}}\subset
F_{C_K^*}|K-Y\subset\Om_{\de_K}(K,t).$$
The last inclusion follows from the very first assumption on
$\de_K$ above. We conclude that
$\Psi^{-1}(0)$ is compact.

\begin{crl}
\label{top_l2c}
Suppose ${\cal M}$ is a smooth oriented manifold,\\
(1) $F\!\equiv\! F^-\!\oplus F^+$, $\tilde{F}^-$, and
 ${\cal O}\!\equiv\!{\cal O}^-\!\oplus {\cal O}^+$ are
smooth normed complex vector bundles over ${\cal M}$ with
$$rk~F^-=\rk~\tilde{F}^-= \rk~{\cal O}^- -\frac{1}{2}\dim~{\cal M}
\quad\hbox{and}\quad
\hbox{rk}~F^+=\hbox{rk}~{\cal O}^+;$$
(2) $\rho\!\in\!\Ga({\cal M};F^{-*\otimes k}\otimes\tilde{F}^-)$
induces a $\tilde{d}$-to-$1$ cover $F\!\lra\!\tilde{F}$, and
$\al^-\!\in\!\Ga({\cal M};\tilde{F}^{-*}\otimes{\cal O}^-)$;\\
(3) $\al\!:F\!\!\lra\!{\cal O}$ is a nondegenerate polynomial such that
$\al^+\!\equiv\!\al|_{F^+}\!:F^+\!\!\lra\!{\cal O}^+$ 
is linear and
$\pi^-\!\circ\!\al\!=\!\al^-\!\circ\!\rho$;\\
(4) $Y$ is a small subset of $F\!=\!F^-\!\oplus\! F^+$,
which is invariant under the scalar multiplication in each component 
separately;\\
(5) $\Om$ is an open subset of ${\cal I}\!\times\! F$ such that 
$\Om\cup X$  is a neighborhood
of $\{0\}\!\times\! X$ 
\hbox{in~${\cal I}\!\times\! \big(F\!-\!(Y\!-\!X)\big)$};\\
(6) $\big\{\phi_t\!: 
 \{\ups\!\in\! F\!:(t,\ups)\!\in\!\Om\}\!\lra\!{\cal O}\big\}$
is a family of smooth maps with resolvent~$\al$;\\
(7) $\bar{\nu}\!=\!(\bar{\nu}^-,\bar{\nu}^+)
       \!\in\!\Ga({\cal M};{\cal O}^-\!\oplus{\cal O}^+)$
is generic with respect to $(\al^+,\al^-,\rho,Y)$,
and the map
\begin{equation}\label{top_l2c_e}
\tilde{F}^-\lra {\cal O}^-,\quad
\vp\lra\bar{\nu}_{\vp}^-+\al^-(\vp),
\end{equation}
has a finite number of (transverse) zeros.\\
If $\psi_t$ is transversal to zero for all $t$,
there exists a compact subset $K_{\al,\bar{\nu}}$ of ${\cal M}$
with the following property.
If $K$ is  precompact open subset of~$K{\cal M}$ containing
$K_{\al,\bar{\nu}}$, there exist $\de_K,\ep_K>0$ such that for all
$t\!\in\!(0,\ep_K)$, 
$$^{\pm}\big|\big\{\ups\!\in\! F_{\de_K}|K\!: 
(t,\ups)\!\in\!\Om, \psi_t(\ups)=0\big\}\big|=\tilde{d}\cdot
{^{\pm}\!\big|}\big\{\vp\!\in\!\tilde{F}^-\!:
\bar{\nu}_{\vp}^- +\al^-(\vp)=0\big\}\big|,$$
where $\psi_t(\ups)=t\bar{\nu}_{\ups}+\phi_t(\ups)$.
Furthermore, all the zeros of $\psi_t\big|(F_{\de_K}|K)$
lie over~$K_{\al,\bar{\nu}}.$
\end{crl}

\noindent
{\it Proof:} Let $K_{\al,\bar{\nu}}$ and $\de_K>0$
be as in Lemma~\ref{top_l2}.
Then if $K$ is a precompact open subset of ${\cal M}$,
for all $t\!\in\!(0,\ep_K)$
the signed number of zeros of  $\psi_t|(\Om_{\de_K}|K)$
is the same as the signed number of solutions of
\begin{equation}\label{top_l2_e7}
F|K\lra{\cal O},\quad
\begin{cases}
\bar{\nu}_{\ups}^- +\al^-\big(\rho(\ups^-)\big)=0\in{\cal O}^-;\\ 
\bar{\nu}_{\ups}^+ +\al^+(\ups^+)+\pi^+\big(\al(\ups^-)\big)=0\in{\cal O}^+.
\end{cases}
\end{equation}
For every solution of the first equation, there is a unique 
solution of the second equation.
Since $\al^+$ is complex-linear on the fibers, 
the signed number of solutions of~\e_ref{top_l2_e7} is 
the same as the signed number of solutions of the first equation.
Since the first equation has no solutions on $Y_{\al^-}$
if $\bar{\nu}$ is generic
and $\rho$~is $\tilde{d}$-to-$1$ outside of~$Y_{\al^-}$,
$\rho$~induces a $\tilde{d}$-to-$1$ 
sign-preserving map from the set of zeros of~\e_ref{top_l2_e} to
the set of solutions of the first equation.

\subsection{Contributions to the Euler Class}
\label{top_sec2}

\noindent
If $\bar{\cal M}$ is a smooth oriented compact manifold of dimension~$n$,
and $V\!\lra\!\bar{\cal M}$ is an oriented vector bundle of rank~$n$,
then the Euler class of $V$ is the number of zeros
of any section \hbox{$s\!:\bar{\cal M}\!\lra\! V$}
which is transverse to the zero set.
In this subsection, under slightly more topological assumptions on 
$\bar{\cal M}$ and $V$, 
we discuss a relationship between subsets of 
the zero set of a non-transverse section and the Euler class of~$V$.

\begin{dfn}
\label{euler_dfn1}
(1) Compact oriented topological manifold 
$\bar{\cal M}\!=\!{\cal M}_n\sqcup
\!\bigsqcup\limits_{i=0}\limits^{i=n-2}\!\!{\cal M}_i$
of dimension~$n$ is \under{mostly smooth}, or~\under{ms}, if\\
(1a) each ${\cal M}_i$ is a smooth oriented manifold
of dimension~$i$;\\
(1b) ${\cal M}\!\equiv\!{\cal M}_n$ is a dense open subset 
of $\bar{\cal M}$;\\
(1c) for each $i\!\in\![n\!-\!2]$, there exists
an ms-manifold $\tilde{\cal M}_i\!=\!
{\cal M}_i\sqcup
\!\bigsqcup\limits_{j=0}\limits^{j=i-2}\!\!\tilde{\cal M}_{i,j}$
and a continuous map 
\hbox{$\varphi_i\!:\tilde{\cal M}_i\!\lra\!\bar{\cal M}$}
such that 
$\varphi_i\big|_{\tilde{\cal M}_{i,j}}\!\!:
\tilde{\cal M}_{i,j}\!\lra\!{\cal M}_j$
is an orientation-preserving diffeomorphism onto
an open subset of~${\cal M}_j$ and 
\hbox{$\varphi_i\big|_{{\cal M}_i}\!=\!Id$}.\\
(2) If $\bar{\cal M}$ is an ms-manifold as above,
topological normed complex vector bundle $V\!\lra\!\bar{\cal M}$ 
is an \under{ms-bundle}
if $V|{\cal M}$ is a smooth complex normed vector bundle
and $\varphi_i^*V\lra\tilde{\cal M}_i$ 
is an ms-bundle for \hbox{all $i\!\in\![n\!-\!2]$}.\\
(3) If $V\!\!\lra\!\bar{\cal M}$  is an ms-bundle,
continuous section $s\!:{\cal M}\!\lra\! V$ is  an \under{ms-section} 
if $s|_{\cal M}$ is $C^2$-smooth and
$\varphi_i^*s_i\!:\tilde{\cal M}_i\!\lra\!\varphi_i^*V$
is an ms-section.
\end{dfn}

\noindent
Since any negative-dimensional topological manifold 
is empty, by (1a) of Definition~\ref{euler_dfn1},
a zero- or one-dimensional ms-manifold 
is simply a compact oriented smooth manifold
of the same dimension.
Thus, Definition~\ref{euler_dfn1} inductively describes
ms-manifolds, ms-bundles, and ms-sections in all dimensions.
We denote the space of ms-sections 
of $V\!\lra\!\bar{\cal M}$ by~$\Ga(\bar{\cal M};V)$.
Using (3) of Definition~\ref{euler_dfn1},
we define a ms-polynomial map between two ms-bundles
analogously to (1) of Definition~\ref{top_dfn1}.\\

\noindent
We topologize $\Ga(\bar{\cal M};V)$ as follows.
If $s_k,s\!\in\!\Ga(\bar{\cal M};V)$, 
the sequence $\{s_k\}$ {\it converges to} $s$
if $s_k$ converges to $s$ in the $C^0$-norm on all of $\bar{\cal M}$
and in the $C^2$-norm on compact subsets of~${\cal M}_i$
for all $i\!\in\![n\!-\!2]$ \hbox{and $i\!=\!n$}.
The $C^0$-norm is defined with respect to the norm 
on $V\!\!\lra\!\bar{\cal M}$.
In order to define the $C^2$-norm on compact subsets of~${\cal M}_i$,
we fix a connection 
in each smooth bundle \hbox{in $V\!\!\lra\!{\cal M}_i$}.

\begin{dfn}
\label{euler_dfn2}
Let  $\bar{\cal M}$ be an ms-manifold of dimension~$2n$
as in Definition~\ref{euler_dfn1}.\\
(1) If ${\cal Z}\!\subset\!{\cal M}_i$ is a smooth oriented submanifold, 
$F\!\lra\!{\cal Z}$ is a smooth complex normed vector bundle,
$Y$ is a small subset of~$F$,
\hbox{$\de\!\in\! C^{\i}({\cal Z};\Bbb{R}^+)$}, and
$\vt\!:\big(F_{\de}-(Y-{\cal Z})\big)\!\lra\!{\cal M}\cup{\cal Z}$ 
is a continuous map,
$(F,Y,\de,\vt)$ is a \under{normal-bundle model for ${\cal Z}$}~if\\
(1a) 
$\vt\!:\big(F_{\de}-(Y-{\cal Z})\big)\!\lra\!{\cal M}\cup{\cal Z}$
is an open map and $\vt|_{\cal Z}$ is the identity;\\
(1b) $\vt|_{F_{\de}-Y}\!: F_{\de}\!-\!Y\lra{\cal M}$
is an orientation-preserving diffeomorphism onto an open
subset of ${\cal M}$.\\
(2) Suppose $V\!\!\lra\!\bar{\cal M}$ is an ms-bundle of rank~$n$
and $s\!\in\!\Ga(\bar{\cal M};V)$.
Smooth submanifold \hbox{${\cal Z}\!\subset\!{\cal M}_i$}
of dimension~$2k$, such that $s|{\cal Z}\!=\!0$, is \under{$s$-regular} if 
there exist a normal-bundle model $(F,Y,\de,\vt)$ for ${\cal Z}$
and a bundle isomorphism $\vt_V\!\!:\vt^*V\!\lra\!\pi_F^*V$,
covering the identity on $\big(F_{\de}-(Y-X)\big)$,
such that $\vt_V|_{F_{\de}-Y}$ is smooth,
$\vt_V|_{\cal Z}$ is the identity, and\\
(2a) 
$\{\phi_0\equiv\vt_V\circ\vt^*s\!: 
          \big(F_{\de}-Y)\big)\!\lra\! V\}$
is a hollow family, OR\\
(2b) there exist an ms-manifold 
$\tilde{\cal Z}={\cal Z}\sqcup\!
\bigsqcup\limits_{i=0}\limits^{i=2k-2}\!\!\tilde{\cal Z}_i$,
continuous map 
$\varphi\!:\tilde{\cal Z}\!\lra\!\bar{\cal M}$,
ms-bundle $E\!\lra\!\tilde{\cal Z}$
and ms-polynomial \hbox{$\al\!:E\!\lra\!\varphi^*V$}
such~that\\
(2b-i) $\varphi|_{\cal Z}=Id$
and $\varphi|_{\tilde{\cal Z}_i}\!:\tilde{\cal Z}_i\!\lra\!{\cal Z}_i$
is a local diffeomorphism onto a smooth submanifold
${\cal Z}_i$ of ${\cal M}_{j(i)}$ for some $j(i)\!\in\![2n]$;\\
(2b-ii) $E|{\cal Z}=\varphi^*F$,
$\al|_{\cal Z}$ is nondegenerate and
is the resolvent for
$\{\phi_0\equiv\vt_V\circ\vt^*s\!: 
          \big(F_{\de}-Y)\big)\!\lra\! V\}$, and
$Y$ is preserved under  scalar multiplication
in each of the components of $F$ for the splitting 
corresponding to $\al$ as in (1) of Definition~\ref{top_dfn1}.
\end{dfn}

\begin{lmm}
\label{euler_l1}
If $(\bar{M},V,s)$ and $({\cal Z},F,Y,\de,\vt)$
are as in (2) of Definition~\ref{euler_dfn2},
there exist a number ${\cal C}_{\cal Z}(s)\!\in\!\Bbb{Z}$,
which equals zero in the Case (2a),
and a dense open subset $\Ga_{\cal Z}(s)\!\subset\!\Ga(\bar{\cal M};V)$
with the following properties.
For every $\nu\!\in\!\Ga_{\cal Z}(s)$,\\ 
(1) there exists $\ep_{\nu}\!>\!0$ such that for all $t\!\in\!(0,\ep_{\nu})$,
all the zeros of $t\nu\!+\!s$ are contained in ${\cal M}$
and $(t\nu\!+\!s)\big|_{\cal M}$ is transversal to 
the zero set in~$V$;\\
(2) there exist a compact subset $K_{\nu}\!\subset\!{\cal Z}$,
open neighborhood $U_{\nu}(K)$ of $K$ in $\bar{\cal M}$
for each compact subset $K\!\subset\!{\cal Z}$,
and $\ep_{\nu}(U)\in(0,\ep_{\nu})$ for each open subset $U$ of $\bar{\cal M}$
such that
$$^{\pm}\big|\{b\!\in\! U\!: t\nu(b)\!+\!s(b)=0\}\big|
={\cal C}_{\cal Z}(s)
\qquad\hbox{if~}
t<\ep_{\nu}(U),~K_{\nu}\subset K\subset U\subset U_{\nu}(K).$$
\end{lmm}

\noindent
{\it Proof:}
It is clear that we can choose 
a dense open subset $\Ga_{\cal Z}'(s)\!\subset\!\Ga(\bar{\cal M};V)$
such that every $\nu\!\in\!\Ga_{\cal Z}'(s)$
satisfies requirement~(1) of the lemma.
If we are in the Case (2a) of Definition~\ref{euler_dfn2},
we also need that $\bar{\nu}\equiv\nu|{\cal Z}$  
is generic with respect to the corresponding polynomial $\al^-$ 
in the sense of the proof of Lemma~\ref{top_l1}.
We can then take $K_{\nu}\!=\!\eset$.
In the Case (2b) of Definition~\ref{euler_dfn2},  
let $\bar{\nu}\!=\!\varphi^*\nu\!\in\!\Ga(\tilde{\cal Z};\varphi^*V)$.
The second part of (6) of 
Lemma~\ref{top_l2} is satisfied if the~map
$$E\lra V,\quad\ups\lra\bar{\nu}_{\ups}+\al(\ups),$$
has no zeros over $\tilde{\cal Z}\!-\!{\cal Z}$.
This is easily insured by transversality on each smooth strata.
The other requirements on $\bar{\nu}$ in 
Lemma~\ref{top_l2} are finitely many transversality properties.
We then take 
$${\cal C}_{\cal Z}(s)=
{^\pm}\!\big|\big\{\ups\!\in\!F\!:
\bar{\nu}_{\ups}\!+\!\al(\ups)\!=\!0\big\}\big|.$$
A cobordism argument similar to (3) of the proof
of Lemma~\ref{top_l2} shows that ${\cal C}_{\cal Z}(s)$
is independent of a generic choice of 
\hbox{$\bar{\nu}\!\in\!\Ga(\tilde{\cal Z};\varphi^*V)$}.
\\

\noindent
The total number of zeros of a section $t\nu\!+\!s$
satisfying condition (1) of Lemma~\ref{euler_l1}
is precisely the Euler class~$e(V)$ of the bundle 
\hbox{$V\!\lra\!\bar{\cal M}$}.
Thus, due to (2) of Lemma~\ref{euler_l1},
we call ${\cal C}_{\cal Z}(s)$ the $s$-{\it{contribution}}
(or simply contribution) of ${\cal Z}$ to~$e(V)$.
If ${\cal Z}$ is any smooth submanifold of ${\cal M}_i$
such that ${\cal Z}\cap s^{-1}(0)$
satisfies the requirements of (2) of Definition~\ref{euler_dfn2},
let ${\cal C}_{\cal Z}(s)={\cal C}_{{\cal Z}(s)\cap s^{-1}(0)}(s)$.
In addition, if $Z$ is a closed subset of $\bar{\cal M}$
such that \hbox{$s^{-1}(0)-Z$} is also closed,
we can easily define ${\cal C}_{\cal Z}(s)$ by Lemma~\ref{euler_l1}.

\begin{crl}
\label{euler_crl}
Let $V\!\lra\!\bar{\cal M}$ be an ms-bundle of rank~$n$
over an ms-manifold of dimension~$2n$.
Suppose ${\cal U}$ is an open subset of ${\cal M}$
and $s\!\in\!\Ga(\bar{\cal M};V)$ is such that
$s|_{\cal U}$ is transversal to the zero~set.\\
(1) If $s^{-1}(0)\cap{\cal U}$ is a finite set,
$^{\pm}|s^{-1}(0)\cap{\cal U}|
=\lan e(V),[\bar{\cal M}]\ran-{\cal C}_{\bar{\cal M}-{\cal U}}(s)$.\\
(2) If $\bar{\cal M}-{\cal U}=
  \bigsqcup\limits_{i=1}\limits^{i=k}\!{\cal Z}_i$,
where each ${\cal Z}_i$ is $s$-regular, then 
$s^{-1}(0)\cap{\cal U}$ is finite, and 
$$^{\pm}|s^{-1}(0)\cap{\cal U}|
=\lan e(V),[\bar{\cal M}]\ran-{\cal C}_{\bar{\cal M}-{\cal U}}(s)
=\lan e(V),[\bar{\cal M}]\ran-\sum_{i=1}^{i=k}{\cal C}_{{\cal Z}_i}(s).$$
If ${\cal Z}_i$ satisfies (2a) of Definition~\ref{euler_dfn2},
${\cal C}_{{\cal Z}_i}(s)\!=\!0$.
If ${\cal Z}_i$ satisfies (2b) of Definition~\ref{euler_dfn2}
and \hbox{$\al_i\!:E_i\!\lra\!V$} is the corresponding polynomial,
$${\cal C}_{{\cal Z}_i}(s)=
^{\pm}\big|\{\ups\!\in\!E_i\!:\bar{\nu}_{\ups}\!+\!\al_i(\ups)=0\}\big|,$$
where $\nu\!\in\!\Ga(\tilde{\cal Z}_i;V)$ is a generic section.
Finally, if $\al_i\!\in\!
\Ga(\tilde{\cal Z}_i;E_i^{*\otimes k}\otimes V)$  has constant rank 
over~$\tilde{\cal Z}_i$, 
$${\cal C}_{{\cal Z}_i}(s)=
k\big\lan e\big(V/\al_i(E_i)\big),[\tilde{\cal Z}_i]\big\ran.$$
\end{crl}

\noindent
All statements of this corollary have already been proved.
A splitting of the zero set as in (2) of Corollary~\ref{euler_crl}
always exists in the complex-analytic category.
If the $s$-regularity condition is restated topologically,
i.e.~in terms of orientation classes, 
such a splitting should exist
for a very large class (if not all) of sections over
any compact oriented topological manifold.
For the cases that we encounter in Section~\ref{comp_sect},
the analytic version of $s$-regularity of Definition~\ref{euler_dfn2}
suffices.

\subsection{Zeros of Polynomial Maps}
\label{zeros_sec}

\noindent
We now present a procedure for computing  the number of zeros 
of a polynomial map between two complex vector bundles 
over a compact oriented manifold.
All the polynomials we encounter in Section~\ref{comp_sect}
are of degree-one.
Thus, we focus in this subsection on the degree-one case,
but discuss the general case at the end for the sake of completeness.\\

\noindent
Suppose  $\bar{\cal M}$ is  an ms-manifold,
$E,{\cal O}\!\lra\!\bar{\cal M}$ are ms-bundles such that
$\rk~E\!+\!\frac{1}{2}\dim\bar{\cal M}\!=\!\rk~{\cal O}$,
and $\al\!\in\!\Ga(\bar{\cal M};E^*\otimes{\cal O})$
is an ms-section.
Let $\bar{\nu}\!\in\!\Ga(\bar{\cal M};{\cal O})$
be such that  $\bar{\nu}$ has no zeros, the map
$$\psi\!\equiv\!\bar{\nu}\!+\!\al\!: E\!\lra\!{\cal O}$$
is transversal to the zero set in ${\cal O}$ 
on~$E|{\cal M}$, and all its zeros
are contained in~$E|{\cal M}$.
The first step in our procedure of determining the number of
zeros of~$\psi$ reduces this issue to the case $E$ is a line bundle.
Let $\Bbb{P}E$ be the projectivization of~$E$
\hbox{(over $\Bbb{C}$)}
and $\ga_E\!\lra\!\Bbb{P}E$ the tautological line bundle.
Then $\al$ induces an ms-section
$\al_E\!\in\!\Ga(\Bbb{P}E;\ga_E^*\otimes\pi_E^*{\cal O})$,
where $\pi_E\!:\Bbb{P}E\lra\bar{\cal M}$ is the bundle projection map. 
The number of zeros of~$\psi$ is the same
as the number of zeros of the induced map
$$\psi_E\!\equiv\! \pi_E^*\bar{\nu}+\al_E\!:
\ga_E\!\lra\!\pi_E^*{\cal O}.$$
Thus, we can always reduce the computation
to the case $E$ is a line bundle.\\

\noindent
The second step describes the number of zeros of $\psi$
topologically in the case $E$ is a line bundle.
Since $\bar{\nu}$ has no zeros, it spans a trivial subbundle
$\Bbb{C}\bar{\nu}$ of~${\cal O}$.
Let ${\cal O}^{\perp}$ be 
the quotient of~${\cal O}$ by this trivial subbundle.
Denote the $\Bbb{C}\bar{\nu}$-
and ${\cal O}^{\perp}$-components of $\al$ by
$\al^t$ and $\al^{\perp}$, respectively.
Then the zeros of~$\psi$ are described by
\begin{equation}\label{zeros_e1}
\begin{cases}
\bar{\nu}_b+\al_b^t(v)=0\in\Bbb{C}\bar{\nu};\\ 
\al_b^{\perp}(v)=0\in{\cal O}^{\perp};
\end{cases}
\quad b\!\in\! X,~v\!\in\! E_b.
\end{equation}
Since $\bar{\nu}$ does not vanish, 
all solutions of the first equations~\e_ref{zeros_e1} are nonzero.
The solution of the second equation with nonzero~$v$ is 
$$\big(E-\bar{\cal M}\big)\big|\al^{\perp-1}(0).$$
Furthermore, if $b\!\in\!\al^{\perp-1}(0)$ and $\al(b)\!\neq\!0$,
$\al^t\!\!: E\!\lra\! (\Bbb{C}\bar{\nu})_b$ is an isomorphism.
Thus, for every 
\hbox{$b\!\in\!\al^{\perp-1}(0)\!-\!\al^{-1}(0)$}, there exists a unique 
$v\!\in\! E_b$ solving the first equation in~\e_ref{zeros_e1},
and the sign of~$(b,v)$ as a zero of $\psi$ agrees with
the sign of~$b$ as a zero of~$\al^{\perp}$.
On the other hand, \e_ref{zeros_e1} has no solutions on~$E|\al^{-1}(0)$.
It follows that the number of zeros of~$\psi$ is 
the number of zeros of $\al^{\perp}$ 
\hbox{on $\bar{\cal M}-\al^{-1}(0)$},~i.e.
\begin{equation}
\label{zeros_sec_e}
^{\pm}|\psi^{-1}(0)|=
\big\lan e(E^*\otimes{\cal O}^{\perp}),[\bar{\cal M}]\big\ran
-{\cal C}_{\al^{-1}(0)}(\al^{\perp});
\end{equation}
see Corollary~\ref{euler_crl}.\\

\noindent
As discussed in the previous subsection,
computing ${\cal C}_Z(s)$ in reasonably good cases
reduces to counting the number of zeros
of polynomial maps between vector bundles over ms-manifolds,
but with the rank of the target bundle one less than 
the rank of the bundle~${\cal O}$ we started with.
Thus, this process will eventually terminate.
The lemma below summarizes the last two paragraphs.
Let $\la_E=c_1(\ga_E^*)$.

\begin{lmm}
\label{zeros_main}
Suppose  $\bar{\cal M}$ is  an ms-manifold and
$E,{\cal O}\lra\bar{\cal M}$ are ms-bundles such that
$$\rk~E+\frac{1}{2}\dim\bar{\cal M}= \rk~{\cal O}\equiv m.$$
If
$\al\!\in\!\Ga(\bar{\cal M};E^*\otimes{\cal O})$ and
$\bar{\nu}\!\in\!\Ga(\bar{\cal M};{\cal O})$
are such that $\bar{\nu}$ has no zeros, the map
$$\psi\!\equiv\!\bar{\nu}\!+\!\al\!: E\lra{\cal O}$$
is transversal to the zero set on~$E|{\cal M}$, and all its zeros
are contained in~$E|{\cal M}$,~then
$$^{\pm}|\psi^{-1}(0)|=
\sum_{k=0}^{k=m-1}
\big\lan c_k({\cal O})\la_E^{m-1-k},[\Bbb{P}E]\big\ran
-{\cal C}_{\al_E^{-1}(0)}(\al_E^{\perp}).$$
Furthermore, if the rank of $E$ is $n$,
\begin{equation}\label{zeros_main_e}
\la_E^n+\sum_{k=1}^{k=n}c_k(E)\la_E^{n-k}=0\in H^{2n}(\Bbb{P}E)
\quad\hbox{and}\quad
\big\lan \mu\la_E^{n-1},[\Bbb{P}E]\big\ran=
\big\lan\mu,[\bar{\cal M}]\big\ran~~
\forall\mu\!\in\!H^{2m-2n}(\bar{\cal M}).
\end{equation}
Finally, if $\al$ has constant rank, 
$$^{\pm}|\psi^{-1}(0)|=
\big\lan e\big({\cal O}/(Im~\al)\big),[\bar{\cal M}]\big\ran.$$
\end{lmm}

\noindent
{\it Proof}: The first statement follows from 
equation~\e_ref{zeros_sec_e} by observing that
\hbox{$c({\cal O}^{\perp})=c({\cal O})$}.
The relations~\e_ref{zeros_main_e} are well-known; see~\cite{BT}.
We state them here for easy reference.
The last claim is clear.\\

\noindent{\it Remark:}
If $\al\!:E\!\lra\!{\cal O}$ is a polynomial, and 
not just a linear map, the first step in computing
the number of zeros of the map $\psi\!=\!\bar{\nu}\!+\!\al$
would be to reduce to the case $\al$ is a linear map
via a projectivization construction similar to
the one in the second paragraph of this subsection.
For example, suppose $\al\!=\!p_1\!+\!p_2$,
where $p_i\!\in\!\Ga(\bar{\cal M};E_i^{*\otimes d_i}\otimes{\cal O})$
and $E\!=\!E_1\!\oplus\! E_2$.
Then the number of zeros of $\psi$ is the same as 
the number of zeros of
$$\psi_{E_1}\!\equiv\!
\pi_{E_1}^*\bar{\nu}\!+\!p_{1,E_1}\!+\!\pi_{E_1}^*p_2\!:
\ga_{E_1}\oplus\pi_{E_1}^*E_2\lra\pi_{E_1}^*{\cal O}$$
over $\Bbb{P}E_1$,
where $p_{1,E_1}\!\in\!\Ga(\Bbb{P}E_1;\ga_{E_1}^{*\otimes d_1})$
is the section induced by~$p_1$.
If $\bar{\nu}$ is generic, this number is $d_1$-times
the number of zeros of the map
$$\tilde{\psi}_{E_1}\!\equiv\!
\pi_{E_1}^*\bar{\nu}\!+\!p_{1,E_1}\!+\!\pi_{E_1}^*p_2\!:
\ga_{E_1}^{\otimes d_1}\oplus\pi_{E_1}^*E_2\lra\pi_{E_1}^*{\cal O}.$$
Note that $p_{1,E_1}$ is {\it linear} on~$\ga_{E_1}^{\otimes d_1}$.
Taking the projection of $\pi_{E_1}^*E_2$
over~$\Bbb{P}E_1$ and repeating the above procedure,
we obtain a linear map 
$$\al_{E_1,E_2}\!:
\pi_{E_2}^*\ga_{E_1}^{\otimes d_1}\oplus
\ga_{\pi_{E_1}^*E_2}^{\otimes d_2}\lra
\pi_{\pi_{E_1}^*E_2}^*\pi_{E_1}^*{\cal O}.$$

\section{Resolvents for $\{\psi_{{\cal T},t\nu}^{\mu}\}$ and 
$\{\psi_{{\cal S},t\nu}^{\mu}\}$}
\label{resolvent_sec}

\subsection{A Power Series Expansion
 for $\pi^{0,1}_{\ups,-}\bar{\partial}u_{\ups}$}
\label{dbar_section}

\noindent
Throughout this section, we assume that 
${\cal T}\!=\!(\Si,[N],I;j,d)$
is a simple bubble type, with $d_{\hat{0}}\!=\!0$
and $\sum\limits_{i\in I}d_i\!=\!d$,
and $\mu$  is an $N$-tuple of constraints in general position
of total codimension 
$$\codim_{\Bbb{C}}\mu=d(n+1)-n(g-1)+N.$$
Our goal is to extract leading-order terms from  
the bundle map $\psi_{{\cal T},t\nu}^{\mu}$ of Theorem~\ref{si_str}
and to describe the zero set of 
$\psi_{{\cal T},t\nu}^{\mu}$ as the union of the zero sets
of affine maps between finite-rank vector bundles.
The main topological tool is Subsection~\ref{top_sec1}.
\\

\noindent
Nearly all of this subsection is devoted to obtaining
the power series expansion for $\pi^{0,1}_{\ups,-}\bar{\partial}u_{\ups}$
of Proposition~\ref{dbar_gen2prp}.
However, we first state an estimate 
for $\pi^{0,1}_{\ups,-}\nu_{\ups,t}$,
which is immediate from Theorem~\ref{si_str}.\\

\noindent
Let $\{\psi_j\}$ denote an orthonormal basis for ${\cal H}^{0,1}_{\Si}$.
Given $q\!\in\!\P$ and an orthonormal basis $\{X_i\}$ for $T_q\P$, put
$$\bar{\nu}_q=\sum_{i=1,j=1}^{i=n,j=g}
\bigg(\int_{z\in\Si}\big\lan\nu(z,q),X_i\psi_j\big\ran_z\bigg)X_i\psi_j
\equiv\pi_{{\cal H}^{0,1}_{\Si}}\nu(\cdot,q)\in
 {\cal H}^{0,1}_{\Si}\otimes T_q\P.$$
Note that $\bar{\nu}$ is well-defined.

\begin{lmm}
\label{nuterm}
There exist $\de,C\!\in\! C^{\i}({\cal M}_{\cal T}^{(0)};\Bbb{R}^+)$ 
such that for all 
$\ups\!\in\! F^{(\eset)}{\cal T}_{\de}$ and $t\!\in\!\big(0,\de(t)\big)$,
$$\big\|\pi^{0,1}_{\ups,-}\nu_{\ups,t}-
 \tilde{R}_{\ups}\bar{\nu}_{\ev(b_{\ups})}\big\|_{\ups,2}\le 
C(b_{\ups})\big(t+|\ups|^{\frac{1}{p}}\big).$$\\
\end{lmm}

\noindent
Suppose 
$\ups\!=\!\big((\Si,[N],I;x,(j,y),u),(v_h)_{h\in\hat{I}}\big)
\!\in\! F^{(\eset)}{\cal T}$ is such that $q_{\ups}$ is defined.
For any $h\!\in\!\hat{I}$, define $\tilde{h}({\cal T})\!\in\!\hat{I}$ by 
$\io_{\tilde{h}({\cal T})}\!=\!\hat{0}$ and 
$h\!\in\!\bar{D}_{\tilde{h}({\cal T})}{\cal T}$.
By the constriction in~\cite{Z}, see Subsection~\ref{gl-basic_gluing},
$$\tilde{v}_h=
d\phi_{b_{\ups},\tilde{h}({\cal T})}\big|_{\tilde{x}_h(\ups)}
\big(dq_{\ups,\io_h}^{-1}\big|_{\tilde{x}_h(\ups)}
d\phi_{b_{\ups},h}^{-1}\big|_0v_h\big)
=\prod_{i\in\hat{I},h\in\bar{D}_i{\cal T}}v_i
\in T_{x_{\tilde{h}({\cal T})}}\Si,$$
where $\phi_{b_{\ups},h}$ is a holomorphic identification
of neighborhoods of $x_h$ in $\Si_{b_{\ups},\io_h}$
and in~$F^{(0)}_{b_{\ups},h}\equiv T_{x_h}\Si_{b_{\ups},h}$.
If $\Si_{b_{\ups},h}\!=\!S^2$, we also identify 
$T_{x_h}\Si_{b_{\ups},h}$ with $\Bbb{C}$ with the map~$q_N$.

\begin{lmm}
\label{dbar_gen1}
For all $\ups\!\in\! F^{(\eset)}{\cal T}$ such that $q_{\ups}$ is defined,
$\bar{\partial}u_{\ups}$ vanishes outside of the annuli
$A_{\ups,h}^-$ with $\chi_{\cal T}h\!=\!1$  and $A_{\ups,h}^{\pm}$
with $\chi_{\cal T}h\!=\!2$.
Furthermore, there exists $\de\!\in\! C^{\i}({\cal M}_{\cal T};\Bbb{R}^+)$
such that for all $\ups\!\in\! F^{(\eset)}{\cal T}_{\de}$
and $h\!\in\!\hat{I}$ with $\chi_{\cal T}h\!=\!1$, on 
$\tilde{A}_{\ups,h}^-\equiv \big\{z\!\in\! F_{h,b_{\ups}}^{(0)}\!: 
\frac{1}{2}|v_h|^{\frac{1}{2}}\!\le\! |z|_{b_{\ups}}
\!\le\!|v_h|^{\frac{1}{2}}\big\}$,
\begin{equation*}\begin{split}
& \Pi_{b_{\ups},\bar{u}_{\ups}(z)}^{-1}
\bar{\partial}\big(u_{\ups}\circ q_{\ups,\io_h}^{-1}\big)
\circ d\phi_{b_{\ups},h}^{-1}\Big|_z\\
&\qquad\qquad=-|v_h|^{-\frac{1}{2}}\bigg(\sum\limits_{m\ge 1}
\big(1-\be_{|v_h|}(2|z|)\big)^{(m-1)}
{\cal D}_{{\cal T},h}^{(m)}
\Big( \Big[ b_{\ups}, \Big(\frac{v_h}{z}\Big) \Big] \Big) \bigg)
\bar{\partial}\be\big|_{2|v_h|^{-\frac{1}{2}}z},
\end{split}\end{equation*}
where $\bar{u}_{\ups}(z)\!\in\! T_{\ev(b_{\ups})}\P$ is given by
$$\exp_{b_{\ups},\ev(b_{\ups})}\bar{u}_{\ups}(z)
=u_{\ups}\big(q_{\ups,\io_h}^{-1}\phi_{b_{\ups},h}^{-1}(z)\big)
=u_h\big(q_{h,(x_h,v_h)}\phi_{b_{\ups},h}^{-1}(z)\big)
,\quad
|\bar{u}_{\ups}(z)|_{b_{\ups}}<r_{\P}.$$
This sum converges uniformly on $\tilde{A}_{\ups,h}^-$.\\
\end{lmm}

\noindent
{\it Remark:}
By construction, $q_{\ups}=
q_{\ups,(x_h,v_h)}\circ q_{\ups,\io_h}$ on $A_{\ups,-}$,
and on $q_{\ups,\io_h}(A_{\ups,-})$
$$q_{\ups,(x_h,v_h)}(z)=\big(h,q_Sp_{h,(x_h,v_h)}(z)\big),
\quad{where}\quad
p_{h,(x_h,v_h)}(z)=\big(1-\be_{|v_h|}(2|\phi_{b_{\ups},h}z|)\big)
\ov{\Big(\frac{v_h}{\phi_{b_{\ups},h}z}\Big)}.$$
\\

\noindent
{\it Proof:} The first claim follows from (G3); 
see~Subsection~\ref{basic_setup}.
If $y\!\in\!\Si_{b_{\ups},h}$ and 
\hbox{$|q_S^{-1}(y)|\!\le\!2\de_{\cal T}(b_{\ups})$},
define \hbox{$\bar{u}_h(y)\!\in\! T_{\ev(b_{\ups})}\P$ by}
$$\exp_{b_{\ups},\ev(b_{\ups})}\bar{u}_h(y)=u_h(y),   \quad
|\bar{u}_h(y)|_{b_{\ups}}<r_{\P}.$$
By construction, 
$u_{\ups}\circ q_{\ups,\io_h}^{-1}=u_h\circ q_{\ups}\circ q_{\ups,\io_h}^{-1}$
on $q_{\ups,\io_h}(A_{\ups,h}^-)$.
Since $\Pi_{b_{\ups},\bar{u}_{\ups}}^{-1}\circ du_h$ is $\Bbb{C}$-linear
on $q_{\ups}(A_{\ups,h}^-)$, for any $z\in\tilde{A}_{\ups,h}^-$
\begin{equation}\label{dbar_gen1_e1}\begin{split}
&\Pi_{b_{\ups},\bar{u}_{\ups}(\cdot)}^{-1}
 \bar{\partial}(u_{\ups}\circ q_{\ups,\io_h}^{-1})\circ 
            d\phi_{b_{\ups},h}^{-1} \Big|_z
=\Pi_{b_{\ups},\bar{u}_{\ups}(\cdot)}^{-1}du_h\circ
\bar{\partial}(q_{\ups}\circ q_{\ups,\io_h}^{-1})\circ d\phi_{b_{\ups},h}^{-1}
\Big|_z\\
&\qquad\qquad\qquad\qquad=-2|v_h|^{-\frac{1}{2}}\l(\frac{v_h}{z}\r)
\Pi_{b_{\ups},\bar{u}_{\ups}(\cdot)}^{-1}
(du_h\circ dq_S)\Big|_{p_{h,(x_h,v_h)}\phi_{\ups,h}^{-1}(z)}
\circ\partial\be\Big|_{2|v_h|^{-\frac{1}{2}}z};
\end{split}\end{equation}
see  Lemma~\ref{gl-basic_gluing_lmm} in~\cite{Z}.
Since $g_{\P,b_{\ups}}$ is flat on $u_{\ups}(A_{\ups,h}^-)$
by our choice of metrics,
\begin{equation}\label{dbar_gen1_e2}
\Pi_{b_{\ups},\bar{u}_{\ups}}^{-1}(du_h\circ dq_S)
=d(\bar{u}_h\circ q_S)
\end{equation}
on $q_S^{-1}q_{\ups}(A_{\ups,h}^-)$.
Since $\bar{u}_h\circ q_S$ is antiholomorphic and 
the metric $g_{\P,b_{\ups}}$ is flat near $\ev(b_{\ups})$,
\begin{equation}\begin{split}\label{dbar_gen1_e3}
\d(\bar{u}_h\circ q_S)\big|_x\Big(\frac{\partial}{\partial s}\Big)
&=d(\bar{u}_h\circ q_S)\big|_x
                     \Big(\frac{\partial}{\partial\bar{y}}\Big)
=\sum_{m\ge1} \frac{\ov{x}^{m-1}}{(m-1)!}
 \frac{d^m}{d\bar{y}^m} \big(\bar{u}_h\circ q_S\big)\big|_{(s,t)=0}\\
&=\sum_{m\ge 1}\frac{\ov{x}^{m-1}}{(m-1)!}
\frac{D^{m-1}}{ds^{m-1}}\frac{d}{ds}(u_h\circ q_S)\big|_{(s,t)=0},
\end{split}\end{equation}
for any $x\!\in\! q_{\ups}(A_{\ups,h}^-)$, where 
$y\!=\!s\!+\!it\!\in\!\Bbb{C}$ 
is the complex coordinate.
The second claim follows from equations 
\e_ref{dbar_gen1_e1}-\e_ref{dbar_gen1_e3}.
For the last claim, note that the sum converges uniformly
on $\tilde{A}_{\ups,h}^-$ as long as
$q_{\ups}(A_{\ups,h}^-)$ is contained
in the ball of convergence for the power series expansion
for $\bar{u}_h$ at~$0$.\\

\noindent
If $\psi\!\in\!{\cal H}_{\Si}^{0,1}$,
$b\!\in\!{\cal M}^{(0)}_{\cal T}$, $m\!\ge\! 1$, and
the metric $g_{b,\hat{0}}$ is flat near~$x$, 
we define  $D_{b,x}^{(m)}\psi\!\in\! T_x^{0,1}\Si^{\otimes m}$
as follows. 
If $(s,t)$ are conformal coordinates centered at $x$
such that $s^2\!+\!t^2$  is the square of the $g_{b,\hat{0}}$-distance 
to $x$, let  
$$\{D_{b,x}^{(m)}\psi\}\Big(\frac{\partial}{\partial s}\Big)\equiv
\{D_{b,x}^{(m)}\psi\}\bigg(
\underset{m}{\underbrace{\frac{\partial}{\partial s},\ldots,
\frac{\partial}{\partial s}}}\bigg)=\frac{\pi}{m!}
\Big\{\frac{D^{m-1}}{ds^{m-1}}\psi_j\Big|_{(s,t)=0}\Big\}
\Big(\frac{\partial}{\partial s}\Big),$$
where the covariant derivatives are taken with respect to the 
metric~$g_{b,\hat{0}}$. 
Since $\psi_j\!\in\! {\cal H}_{\Si}^{0,1}$, 
\hbox{$\psi_j\!=\!f(ds\!-\!idt)$}
for some anti-holomorphic function~$f$.
Since $g_{b,\hat{0}}$ is flat near $x$, it follows that 
$D_{b,x}^{(m)}\psi\!\in\! T_x^{0,1}\Si^{\otimes m}$.
If $\{\psi_j\}$ is an orthonormal basis for ${\cal H}_{\Si}^{0,1}$,
let $s_{b,x}^{(m)}\!\in\! T_x^*\Si^{\otimes m}\otimes{\cal H}_{\Si}^{0,1}$
be given~by
$$s_{b,x}^{(m)}(v)\equiv
s_{b,x}^{(m)}(\underset{m}{\underbrace{v,\ldots,v}})=
\sum_{j\in[g]}\ov{\Big\{D_{b,x}^{(m)}\psi_j\Big\}(v)}\psi_j.$$
The section $s_{b,x}^{(m)}$ is always independent of the choice
of  a basis for ${\cal H}_{\Si}^{0,1}$, but is dependent
on the choice of the metric $g_{b,\hat{0}}$ \hbox{if $m\!>\!1$.}
However, $s_{b,x}^{(1)}$ depends only on~$(\Si,j)$;
we denote this section by~$s_{\Si,x}$.
By \cite[p246]{GH}, $s_{\Si,x}$  does not vanish and thus spans
a subbundle of $\Si\!\times\!{\cal H}_{\Si}^{0,1}\!\lra\!\Si$.
We denote this subbundle by ${\cal H}_{\Si}^+$ and its orthogonal
complement by ${\cal H}_{\Si}^-$.
A slightly different description of these bundles is given in
Subsection~\ref{obs_sec}.
Let 
$$\pi^+,\pi^-\in\Ga\big(\Si;(\Si\times{\cal H}_{\Si}^{0,1})^*\otimes
{\cal H}_{\Si}^{\pm}\big)$$
be the corresponding orthogonal projection maps.
Denote by $s_{b,x}^{(m,\pm)}$ the composition 
$\pi_x^{\pm}\circ s_{b,x}^{(m,\pm)}$.

\begin{lmm}
\label{dbar_gen2}
There exists $\de\!\in\! C^{\i}({\cal M}_{\cal T};\Bbb{R}^+)$
such that for all $\ups\!\in\! F^{(\eset)}{\cal T}_{\de}$,
$X\!\in\! T_{\ev(b_{\ups})}\P$, and $\psi\!\in\!{\cal H}_{\Si}^{0,1}$,
$$\big\lan\big\lan\pi_{\ups,-}^{0,1}\bar{\partial}u_{\ups},
R_{\ups}X\psi\big\ran\big\ran_{\ups,2}=
-\sum_{m\ge 1}\sum_{\chi_{\cal T}h=1}
\big\lan{\cal D}_{{\cal T},h}^{(m)}b_{\ups},X\big\ran
\ov{\Big(      \big\{D_{b_{\ups},\tilde{x}_h(\ups)}^{(m)}\psi\big\}
\big( (d\phi_{b_{\ups},x_{\tilde{h}({\cal T})}}|_{\tilde{x}_h(\ups)})^{-1}
           \tilde{v}_h\big)    \Big)}.$$
Furthermore, the sum is absolutely convergent.
\end{lmm}

\noindent
{\it Proof:} Since $\big\lan\bar{\partial}u_{\ups},R_{\ups}X\psi\big\ran=0$
outside of the annuli $A_{\ups,h}^-$ with $\chi_{\cal T}h=1$,
\begin{equation}\label{dbar_gen2_e1}
\big\lan\big\lan\pi^{0,1}_{\ups,-}\bar{\partial}u_{\ups},
    R_{\ups}X\psi\big\ran\big\ran=
\big\lan\big\lan\bar{\partial}u_{\ups},R_{\ups}X\psi\big\ran\big\ran
=\sum_{\chi_{\cal T}h=1}
\int_{A_{\ups,h}^-}\big\lan\bar{\partial}u_{\ups},R_{\ups}X\psi\big\ran.
\end{equation}
Since $q_{\ups,\io_h}^{-1}\circ\phi_{b_{\ups},h}^{-1}$ is holomorphic
on $\tilde{A}_{\ups,h}^-$, $\Pi_{b_{\ups},\bar{u}_{\ups}}^{-1}$
is unitary on $u_{\ups}(A_{\ups,h}^-)$, and 
the inner-product of one-forms is conformally invariant,
\begin{equation}\label{dbar_gen_e2}\begin{split}
\int_{A_{\ups,h}^-}\big\lan\bar{\partial}u_{\ups},R_{\ups}X\psi\big\ran
&=\int_{\tilde{A}_{\ups,h}^-}\big\lan
\bar{\partial}(u_{\ups}\circ q_{\ups,\io_h}^{-1})
            \circ d\phi_{b_{\ups},h}^{-1},
R_{\ups}X\psi\circ dq_{\ups,\io_h}^{-1}\circ d\phi_{b_{\ups},h}^{-1}
          \big\ran\\
&=\int_{\tilde{A}_{\ups,h}^-}\big\lan\Pi_{b_{\ups},\bar{u}_{\ups}}^{-1}
\bar{\partial}(u_{\ups}\circ q_{\ups,\io_h}^{-1})
                                \circ d\phi_{b_{\ups},h}^{-1},
       X\psi\circ dq_{\ups,i_h}^{-1}\circ \phi_{b_{\ups},h}^{-1}\big\ran,
\end{split}\end{equation}
since  $\Pi_{b_{\ups},\bar{u}_{\ups}}^{-1}R_{\ups}X\psi\!=\!X\psi$
on $A_{\ups,h}^-$.
If $\io_h\!=\!\hat{0}$, we identify $F^{(0)}_{h,b_{\ups}}\!=\!T_{x_h}\Si$
with $\Bbb{C}$ in a $g_{b_{\ups},\hat{0}}$-unitary way. 
In all cases, we can then write 
$$\psi\circ dq_{\ups,\io_h}^{-1}\circ d\phi_{b_{\ups},h}^{-1}=fd\bar{z}.$$
Since $\psi$ is harmonic and $q_{\ups,\io_h}^{-1}\circ \phi_{b_{\ups},h}^{-1}$
is holomorphic on $\tilde{A}_{\ups,h}$, $f$ is anti-holomorphic.
Using the change of variables $2|v_h|^{-\frac{1}{2}}z=re^{i\th}$,
we obtain
\begin{alignat}{1}
&\int_{\tilde{A}_{\ups,h}}\Big\lan |v_h|^{-\frac{1}{2}}
                         \big( 1-\be_{|v_h|}(2|z|) \big)^{m-1}
{\cal D}_{{\cal T},h}^{(m)} \Big(\Big[ b_{\ups},\frac{v_h}{z} \Big]\Big)
  \bar{\partial}\be\big|_{2|v_h|^{-\frac{1}{2}}z},
X\psi\circ dq_{\ups,\io_h}^{-1}\circ d\phi_{b_{\ups},h}^{-1}\Big\ran\notag\\
\label{dbar_gen2_e3}
&\quad
=\Big\lan {\cal D}_{{\cal T},h}^{(m)}b_{\ups},X \Big\ran v_h^m
\int_{\tilde{A}_{\ups,h}}
\Big\{\big(1-\be(2|v_h|^{-\frac{1}{2}}|z|)\big)^{m-1}
     \be'\big|_{2|v_h|^{-\frac{1}{2}}|z|}\Big\}
               |v_h|^{-\frac{1}{2}}z^{-m}\frac{z}{|z|}\bar{f}\\
&\quad
=\Big\lan {\cal D}_{{\cal T},h}^{(m)}b_{\ups},X \Big\ran v_h^m
|v_h|^{-\frac{m-1}{2}}2^{m-2}
\frac{1}{m}\int_1^2\int_0^{2\pi}
\Big\{\big(1-\be(r)\big)^m\Big\}'(re^{i\th})^{-(m-1)}
  \bar{f}\big(\frac{1}{2}|v_h|^{\frac{1}{2}}re^{i\th}\big)d\th dr\notag
\end{alignat}
Since $\bar{f}$ is holomorphic, for any $r>0$,
\begin{equation}\label{dbar_gen2_e4}\begin{split}
\int_0^{2\pi}(re^{i\th})^{-(m-1)}&
\bar{f}\big(\frac{1}{2}|v_h|^{\frac{1}{2}}re^{i\th}\big)d\th
=-i\int\limits_{|z|=r}z^{-m}\bar{f}
                \big(\frac{1}{2}|v_h|^{\frac{1}{2}}z\big)dz\\
&=\frac{2\pi}{(m-1)!}\frac{d^{(m-1)}}{dz^{(m-1)}}
\bar{f}\big(\frac{1}{2}|v_h|^{\frac{1}{2}}z\big)\Big|_{z=0}
=\frac{2\pi}{(m-1)!}2^{-(m-1)}|v_h|^{\frac{m-1}{2}}\bar{f}^{\lr{m-1}}(0).
\end{split}\end{equation}
Since the metric $g_{b,\hat{0}}$ is flat near $\tilde{x}_h$,
\begin{equation}\label{dbar_gen2_e5}\begin{split}
\frac{\pi}{m!} v_h^{m}\bar{f}^{\lr{m-1}}(0)
&=\ov{\big\{  D_{b_{\ups},\tilde{x}_h(\ups)}^{(m)}\psi\big\} 
\big(dq_{\ups,\io_h}^{-1}\big|_{x_h}d\phi_{b_{\ups},h}^{-1}\big|_0v_h\big)}\\
&=\ov{\big\{D_{b_{\ups},\tilde{x}_h(\ups)}^{(m)}\psi\big\}
\big((d\phi_{b_{\ups},x_{\tilde{h}({\cal T})}}|_{\tilde{x}_h(\ups)})^{-1}
                  \tilde{v}_h\big)}.
\end{split}\end{equation}
The claim follows from equations \e_ref{dbar_gen2_e1}-\e_ref{dbar_gen2_e5}
and Lemma~\ref{dbar_gen1}.

\begin{prp}
\label{dbar_gen2prp}
If ${\cal T}\!=\!(\Si,[N],I;j,d)$ is a simple bubble type
with $d_{\hat{0}}\!=\!0$,
there exists \hbox{$\de\!\in\! C^{\i}({\cal M}_{\cal T};\Bbb{R}^+)$}
such that 
$$\pi_{\ups,-}^{0,1}\bar{\partial}u_{\ups}=
-\tilde{R}_{\ups}\sum_{m\ge 1}\sum_{\chi_{\cal T}h=1}
\big({\cal D}_{{\cal T},h}^{(m)}b\big)
\Big(s_{b,\tilde{x}_h(\ups)}^{(m)}
\big(d\phi_{b,x_{\tilde{h}({\cal T})}}|_{\tilde{x}_h(\ups)}^{-1}
           \tilde{v}_h\big)\Big)
\quad\forall~\ups\!=\!\big[b,(v_h)_{h\in\hat{I}}\big]
\!\in\! F^{\eset}{\cal T}_{\de}.$$
Furthermore, the sum is absolutely convergent.
\end{prp}

\noindent
{\it Proof:} This proposition follows from Lemma~\ref{dbar_gen2}
and equation~\e_ref{new_isom}.

\subsection{First-Order Estimate for $\psi_{{\cal T},t\nu}^{\mu}$}
\label{order1}

\noindent
For any $\ups\!=\!
 \big[(\Si,[N],I;x,(j,y),u),(v_h)_{h\in\hat{I}}\big]
   \!\in\! F{\cal T}$ and $h\!\in\!\hat{I}$ such that 
\hbox{$\chi_{\cal T}h\!=\!1$},
let
$$\al_{{\cal T},h}^{(k)}(\ups)=
\big({\cal D}_{{\cal T},h}^{(k)}b_{\ups}\big)
s_{b_{\ups},x_{\tilde{h}({\cal T})}}^{(k)}(\tilde{v}_h),\qquad
\al_{{\cal T},h}^{(k)}(\ups)=\sum_{\chi_{\cal T}h=1}
\al_{{\cal T},h}^{(k)}(\ups).$$
We denote $\al_{{\cal T},h}^{(1)}$ and $\al_{\cal T}^{(1)}$
by $\al_{{\cal T},h}$ and $\al_{\cal T}$, respectively.

\begin{lmm}
\label{dbar_1}
There exist $\de,C\!\in\! C^{\i}({\cal M}_{\cal T};\Bbb{R}^+)$ 
such that for all $\ups\!\in\! F^{\eset}{\cal T}_{\de}$,
$$\big\|\pi_{\ups,-}^{0,1}\bar{\partial}u_{\ups}+
\tilde{R}_{\ups}\al_{\cal T}(\ups)\big\|_2\le
C(b_{\ups})|\ups|\sum_{\chi_{\cal T}h=1}|\ups|_h.$$
\end{lmm}

\noindent
{\it Proof:} This is immediate from Proposition~\ref{dbar_gen2prp}, since  
\begin{gather*}
\big\|s_{\tilde{x}_h(\ups)}
\big(d\phi_{b_{\ups},x_{\tilde{h}({\cal T})}}|_{\tilde{x}_h(\ups)}^{-1}
           \tilde{v}_h\big)-
s_{x_{\tilde{h}({\cal T})}}\big(\tilde{v}_h\big)\big\|_2\le
C(b_{\ups})\big|\phi_{b_{\ups},x_{\tilde{h}({\cal T})}}
\tilde{x}_h(\ups)\big|_{b_{\ups}} \big|\tilde{v}_h\big|
\le C'(b_{\ups})|\ups|\big|\tilde{v}_h\big|_b;\\
\sum_{m\ge 2}
\big|{\cal D}_{{\cal T},h}^{(m)}b_{\ups}\big||\tilde{v}_h|^m
\le C(b_{\ups})|\tilde{v}_h|^2,
\end{gather*}
for all $h\!\in\!\hat{I}$ with $\chi_{\cal T}h\!=\!1$ and 
$\ups\!\in\! F{\cal T}_{\de}$
with $\de\!\in\! C^{\i}({\cal M}_{\cal T};\Bbb{R}^+)$
sufficiently small.

\begin{lmm}
\label{psit_11}
There exist $\de,C\!\in\! C^{\i}({\cal M}_{\cal T}(\mu);\Bbb{R}^+)$ 
such that for all $\ups\!\in\! F^{\eset}{\cal T}_{\de}$,
$$\Big\|\psi_{{\cal T},t\nu}^{\mu}(\ups)-
\big(t\bar{\nu}_{\ev(b_{\ups})}+\al_{\cal T}(\ups)\big)\Big\|_2
\le C(b_{\ups})\big(t+|\ups|^{\frac{1}{p}}\big)
\Big(t+\sum_{\chi_{\cal T}h=1}|\ups|_h\Big),$$
where $\psi_{{\cal T},t\nu}^{\mu}$ denotes  
$\psi_{{\cal M}_{\cal T},t\nu}^{\mu}$.
\end{lmm}

\noindent
{\it Proof:} By Lemma~\ref{approx_coker} and Theorem~\ref{si_str},
$$\big\|\pi_{\ups,-}^{0,1}D_{\ups}\xi_{\ups,t\nu}\big\|_2
\le C(b_{\ups})
\Big(\sum_{\chi_{\cal T}h=1}|\ups|_h\Big)
\|D_{\ups}\xi_{\ups,t\nu}\|_{\ups,p,1}
\le C'(b_{\ups})\big(t+|\ups|^{\frac{1}{p}}\big)
\sum_{\chi_{\cal T}h=1}|\ups|_h.$$
Combining this estimate with Lemmas~\ref{nuterm} and~\ref{dbar_1},
we obtain
\begin{equation}\label{psit_11_e1}
\Big\|\psi_{{\cal T},t\nu}(\ups)-
\big(t\bar{\nu}_{\ev(b_{\ups})}+\al_{\cal T}(\ups)\big)\Big\|_2
\le C(b_{\ups})\big(t+|\ups|^{\frac{1}{p}}\big)
\Big(t+\sum_{\chi_{\cal T}h=1}|\ups|_h\Big)
\end{equation}
for all $\ups\!\in\! F^{\eset}{\cal T}_{\de}$, provided
$\de\!\in\! C^{\i}({\cal M}_{\cal T};\Bbb{R}^+)$ is
sufficiently small.
On the other hand, if $b_{\ups}\!\in\!{\cal M}_{\cal T}$,
\begin{gather}
\big\|\varphi_{{\cal T},t\nu}^{\mu}(\ups)\big\|_{b_{\ups}}
\le C(b_{\ups})\big(t+|\ups|^{\frac{1}{p}}\big)\Lra\notag\\
\label{psit_11_e2}
\begin{split}
\Big\|\big(
t\bar{\nu}_{\ev(\phi_{\cal T}^{\mu}\varphi_{{\cal T},t\nu}^{\mu}(\ups))}
\!+\!\al_{\cal T}(\Phi_{\cal T}^{\mu}\varphi_{{\cal T},t\nu}^{\mu}(\ups))\big)
-\Pi_{b_{\ups},\phi_{\cal T}^{\mu}\varphi_{{\cal T},t\nu}^{\mu}(\ups)}
\big(t\bar{\nu}_{\ev(b_{\ups})}\!+\!\al_{\cal T}(\ups)\big)\Big\|_2
\qquad\qquad&\\
\le C(b_{\ups})\big(t+|\ups|^{\frac{1}{p}}\big)
\Big(t\!+\!\!\sum_{\chi_{\cal T}h=1}\!|\ups|_h\Big),&
\end{split}
\end{gather}
where $\varphi_{{\cal T},t\nu}^{\mu}=\varphi_{{\cal M}_{\cal T},t\nu}^{\mu}$
is the section of Theorem~\ref{si_str} for any fixed regularization
$\big(\Phi_{\cal T}\!\equiv\! Id,\Phi_{\cal T}^{\mu}
\big)$ of~${\cal M}_{\cal T}(\mu)$.
The claim follows from~\e_ref{psit_11_e1}
and~\e_ref{psit_11_e2}.\\

\noindent
Our next step is to apply Lemma~\ref{top_l1} or Corollary~\ref{top_l2c}
to the map $\psi_{{\cal T},t\nu}^{\mu}$ whenever possible.
In terms of notation of Subsection~\ref{top_sec1},
we~take
\begin{gather*}
F^+\!={\cal O}^+\!=\{0\},~\quad 
F^-\!=F{\cal T},\quad
{\cal O}^-\!={\cal H}_{\Si}^{0,1}\otimes\ev^*T\P,\quad
\tilde{F}^-\!=\!\bigoplus_{\chi_{\cal T}h=1}
\bigotimes_{i\in\hat{I},h\in\bar{D}_i{\cal T}}F_i{\cal T};\\
\phi_h\big([b,v_{\hat{I}}]\big)=
\big[b,\bigotimes_{i\in\hat{I},h\in\bar{D}_i{\cal T}}v_i\big]
=\big[b,\tilde{v}_h\big],\quad
\al^-(\phi(\ups))\equiv \al_{\cal T}(\ups),
\end{gather*}
where $\phi_h$ denotes the $h$th component of 
$\phi\!: F^-\!\!\lra\! \tilde{F}^-$.
Note that $\al^-\!\in\!\Ga({\cal M}_{\cal T};
\tilde{F}^{-*}\otimes {\cal O}^-)$ is well-defined.
A priori, $\al^-$ may not have full rank on every fiber 
over~${\cal M}_{\cal T}(\mu)$.
We will call a subset $K\!\subset\!{\cal M}_{\cal T}(\mu)$ 
${\cal T}$-{\it{regular}} if $\al^-$ has full rank over~$K$.
From Theorem~\ref{si_str}, Lemma~\ref{top_l1}, 
and Corollary~\ref{top_l2c}, we then obtain

\begin{crl}
\label{order1_contr}
Suppose $d$ is a positive integer,
${\cal T}\!=\!(\Si,[N],I;j_{[N]},\under{d})$ is
a simple bubble type, with $d_{\hat{0}}\!=\!0$ and
$\sum\limits_{i\in I}d_i\!=\!d$, and 
$\mu$ is an $N$-tuple of constraints in general position 
such that
$$\codim_{\Bbb{C}}\mu=d(n+1)-n(g-1)+N.$$
Let $\nu\!\in\!\Ga(\Si\times\P;
          \La^{0,1}\pi_{\Si}^*T^*\Si\otimes\pi_{\P}^*T\P)$
be a generic section.
If $\io_h\!\neq\!\hat{0}$ for some $h\!\in\!\hat{I}$, for every 
regular compact subset~$K$ of~${\cal M}_{\cal T}(\mu)$,
there exist a neighborhood $U_K$ of $K$ in $\bar{C}^{\i}_{(d;[N])}(\Si;\mu)$
and  $\ep_K\!>\!0$ such that for any $t\!\in\!(0,\ep_K)$,
\hbox{$U_K\cap{\cal M}_{\Si,d,t\nu}(\mu)\!=\!\eset$.}
If $\io_h\!=\!\hat{0}$ for all $h\!\in\!\hat{I}$,
there exists a compact regular subset $K_{\cal T}$ 
of~${\cal M}_{\cal T}(\mu)$ with the following property. 
If $K$ is a compact subset of ${\cal M}_{\cal T}(\mu)$
containing~$K_{\cal T}$,
there exist a neighborhood $U_K$ of $K$ in $\bar{C}^{\i}_{(d;[N])}(\Si;\mu)$
and $\ep_K\!>\!0$ such that for all $t\!\in\!(0,\ep_K)$,
the signed cardinality of $U_K\cap{\cal M}_{\Si,d,t\nu}(\mu)$ 
equals to the signed number of zeros of the map
\begin{equation}\label{order1_contr_e}
F{\cal T}\big|{\cal M}_{\cal T}(\mu)\lra 
{\cal H}_{\Si}^{0,1}\otimes\ev^*T\P,\quad
\ups\lra \bar{\nu}_{\ev(b_{\ups})}+\al_{\cal T}(\ups).
\end{equation}
\end{crl}

\noindent
{\it Proof:} In either case, by Theorem~\ref{si_str},
there exist a neighborhood $U_K$ of $K$ in $\bar{C}^{\i}_{(d;[N])}(\Si;\mu)$
and $\de_K,\ep_K\!>\!0$ such that for any $t\!\in\!(0,\ep_K)$,
there exists a sign-preserving bijection between  
\hbox{$U_K\cap{\cal M}_{\Si,d,t\nu}(\mu)$} and the zeros of
$\psi_{{\cal T},t\nu}^{\mu}$ on 
$F^{\eset}{\cal T}_{\de_K}\big|\big(U_K\cap {\cal M}_{\cal T}(\mu)\big)$,
provided  $U_K\cap {\cal M}_{\cal T}(\mu)$ is precompact 
in~${\cal M}_{\cal T}(\mu)$.
Furthermore, $\de_K$ can be required to be arbitrarily small.
If $K$ is regular, $U_K$ can be chosen so that the closure of
$U_K\cap {\cal M}_{\cal T}(\mu)$ in ${\cal M}_{\cal T}(\mu)$
is also regular.
Then by Lemma~\ref{psit_11},
$$\Big\|\psi_{{\cal T},t\nu}^{\mu}(\ups)-
\big(t\bar{\nu}_{\ev(b_{\ups})}+\al_{\cal T}(\ups)\big)\Big\|_2
\le C_K\big(t+|\ups|^{\frac{1}{p}}\big)
\big(t+|\al_{\cal T}(\ups)|\big)\quad\forall\ups\!\in\!
F^{\eset}{\cal T}_{\de_K}\big|K,$$
where $C_K\!>\!0$ depends only on $K$.
Thus, the first claim follows from Lemma~\ref{top_l1}.
The second follows from Corollary~\ref{top_l2c}, provided that for
a generic $\nu$ the set of zeros of the map in~\e_ref{order1_contr_e}
is ${\cal T}$-regular and finite; see 
Lemma~\ref{order1_contr_l1}~below.

\begin{lmm}
\label{order1_contr_l1}
If $\io_h\!=\!\hat{0}$ for all $h\!\in\!\hat{I}$, 
the set of zeros of the map in~\e_ref{order1_contr_e}
is ${\cal T}$-regular and~finite.
\end{lmm}

\noindent
{\it Proof:} We first show that the set of zeros is finite.
Since any degree-zero holomorphic map is constant,
\begin{equation}\label{cart_split1}
{\cal M}_{\cal T}(\mu)=
\tilde{\cal M}_{\Si,\hat{I}}\times
{\cal U}_{\bar{\cal T}}(\mu),
\quad\hbox{where}\quad
\tilde{\cal M}_{\Si,\hat{I}}=
\big\{x_{\hat{I}}\!\in\!\Si^{\hat{I}}\!:
x_{i_1}\!\neq\! x_{i_2}~\forall i_1\!\neq\! i_2\big\}.
\end{equation} 
Via this decomposition, the map in~\e_ref{order1_contr_e} 
extends continuously to the total space of a bundle 
over~$\Si^{\hat{I}}\!\times\!\bar{\cal U}_{\bar{\cal T}}(\mu)$.
The complement of ${\cal M}_{\cal T}(\mu)$ in 
$\Si^{\hat{I}}\!\times\!\bar{\cal U}_{\bar{\cal T}}(\mu)$
is the union of smooth manifolds of dimension less than
the dimension of ${\cal M}_{\cal T}(\mu)$, as long
as the constraints are in general position.
Thus, if $\nu$ is generic, the extension
of the map in~\e_ref{order1_contr_e} does not vanish over
the complement of ${\cal M}_{\cal T}(\mu)$ in 
$\Si^{\hat{I}}\!\times\!\bar{\cal U}_{\bar{\cal T}}(\mu)$.
Since \hbox{$\Si^{\hat{I}}\!\times\!\bar{\cal U}_{\bar{\cal T}}(\mu)$}
is compact, the set of zeros of (the extension of) the map 
in~\e_ref{order1_contr_e} is finite.
Furthermore, by (2) of the proof of Lemma~\ref{order1_contr_l3},
the subset of ${\cal M}_{\cal T}(\mu)$ on which
$\al^-\!\!=\!\al_{\cal T}$ is 
a proper submanifold of~${\cal M}_{\cal T}(\mu)$.
Thus, if $\nu$ is generic,  the map in~\e_ref{order1_contr_e}
does not vanish of this subset.

\subsection{Consequences of the First-Order Estimate for 
$\psi_{{\cal T},t\nu}^{\mu}$}

\noindent
In this subsection, we show that ${\cal M}_{\cal T}(\mu)$ 
is ${\cal T}$-regular for most bubble types~${\cal T}$ 
under consideration,
and nearly all of them fall under the first 
case of Corollary~\ref{order1_contr}.
We call ${\cal T}$ {\it effective},
if for some generic choice of $\nu$ and of 
the constraints $\mu_1,\ldots,\mu_N$,
$\ov{\bigcup\limits_{t<1}{\cal M}_{\Si,t\nu,d}(\mu)}$ 
intersects~$\bar{\cal M}_{\cal T}(\mu)$.
If $K$ is a compact subset of $\bar{\cal M}_{\cal T}(\mu)$, 
we call $K$ {\it effective} if
$\ov{\bigcup\limits_{t<1}{\cal M}_{\Si,t\nu,d}(\mu)}$ intersects~$K$.

\begin{lmm}
\label{order1_contr_l2}
Let ${\cal T}\!=\!(\Si,[N],I;j,d)$ be a simple bubble type.
If $j_l\!=\!\hat{0}$ for some $l\!\in\![N]$
and $K$ is a ${\cal T}$-regular subset of ${\cal M}_{\cal T}(\mu)$,
then $K$ is not effective.
\end{lmm}

\noindent
{\it Proof:} By Corollary~\ref{order1_contr},
it is sufficient to show that the map
$$\bar{\nu}\!+\!\al_{\cal T}\!: F{\cal T}
\lra{\cal H}^{0,1}_{\Si}\otimes\ev^*T\P$$
has no zeros for a generic $\nu$.
For a generic $\nu$, the zero set of this section is zero-dimensional.
However, if $j_l\!=\!\hat{0}$ for some $l\!\in\![N]$, 
we can move $y_l\!\in\!\Si$ freely,
without changing the value \hbox{of $\bar{\nu}\!+\!\al_{\cal T}$.}
Thus, if the zero-set of the section is nonempty, it must
be at least one-dimensional, which is not the case for a generic $\nu$.

\begin{lmm}
\label{order1_contr_l3}
Let ${\cal T}\!=\!(\Si,[N],I;j,d)$
be a bubble type with $d_{\hat{0}}\!=\!0$.
If
$$ng-|\hat{I}|-\Big(|H_{\hat{0}}{\cal T}|+|M_{\hat{0}}{\cal T}|\!\!+\!\!
 \sum_{i\in\hat{I},d_i=0}\big(|H_i{\cal T}|+|M_i{\cal T}|-2\big)\Big)
\le n-\big|\{h\!:\chi_{\cal T}h\!=\!1\}\big|,$$
${\cal M}_{\cal T}(\mu)$ is ${\cal T}$-regular.
Furthermore, if the number on the left-hand side above is negative,
then ${\cal M}_{\cal T}(\mu)$ is empty.
\end{lmm}

\noindent
{\it Proof:} (1) The dimension of ${\cal M}(\mu)$ is given by 
$$\dim {\cal M}_{\cal T}(\mu)=
\big(d(n+1)+n+N-|\hat{I}|\big)-\big(d(n+1)-n(g-1)+N\big)
=ng-|\hat{I}|.$$
However, given $b\!=\!\big(\Si,[N],I;x,(j,y),u\big)
\!\in\!{\cal M}_{\cal T}(\mu)$, we are free
to vary $x_h$ if $\io_h\!=\!\hat{0}$ \hbox{(i.e. $x_h\!\in\!\Si$)} and 
$y_l$ \hbox{if $j_l\!=\!\hat{0}$.}
Similarly, if $i\!\in\!\hat{I}$, $d_i\!=\!0$, and 
\hbox{$|H_i{\cal T}|\!+\!|M_i{\cal T}|\!>\!2$},
we can vary \hbox{$|H_i{\cal T}|\!+\!|M_i{\cal T}|\!-\!2$} 
marked and singular
points on $\Si_{b,i}$.
Thus, the space ${\cal M}_{\cal T}(\mu)$ must have dimension at least
$$d_{\min}({\cal T})\equiv 
|H_{\hat{0}}{\cal T}|+|M_{\hat{0}}{\cal T}|+
\sum_{i\in\hat{I},d_i=0}\big(|H_i{\cal T}|+|M_i{\cal T}|-2\big),$$
if ${\cal M}_{\cal T}(\mu)$ is nonempty.
Therefore, we can assume $\big|\{h\!:\chi_{\cal T}h\!=\!1\}\big|\!\le\! n$.\\
(2) Let $h_1,\ldots,h_{|\{h:\chi_{\cal T}h=1\}|}$
be the elements of $\{h\!:\chi_{\cal T}h\!=\!1\}$.
The section 
$s_{\Si}\!\in\!\Ga(\Si;T^*\Si\otimes {\cal H}^{0,1}_{\Si})$
does not vanish; see \cite[p246]{GH}.
Thus, the section $\al^-$ defined above has rank at least $k$ if 
the~section 
$$\bar{\cal D}_{{\cal T};k}\!\in\! \Ga\Big({\cal M}_{\cal T}(\mu);
 \Big(\bigoplus_{m\le k}L_{{\cal T},h_m}^*\Big)\otimes\ev^*T\P\Big),
\quad
\bar{\cal D}_{{\cal T};k}\l(\l[b,c_{\{h_m: m\le k\}}\r]\r)=
 \sum_{m\le k}{\cal D}_{{\cal T},h_m}\l([b,c_{h_m}]\r),$$
has rank $k$.
We prove inductively that under the assumptions of the lemma
this is the case for all 
\hbox{$k\!\le\!\big|\{h\!:\chi_{\cal T}h\!=\!1\}\big|$.}
If $k\!=\!0$, there is nothing to prove. 
So we can assume that $k\!>\!0$ and that the statement has been 
shown to be true for $k\!-\!1$. 
The $k\!-\!1$ statement shows that the image of $\bar{\cal D}_{{\cal T};k-1}$
is a rank $k\!-\!1$ subbundle of $\ev^*T\P$.
Let $\pi_{k-1}^{\perp}$ denote the orthogonal projection onto the orthogonal
complement of this rank $(k-1)$-subbundle in~$\ev^*T\P$
with respect to the standard metric in $\P$.
We need to show that the section
$$\pi_{k-1}^{\perp}\circ{\cal D}_{{\cal T};k}\in 
\Ga\big({\cal M}_{\cal T}(\mu);
 L_{{\cal T},h_k}^*\otimes\pi_{k-1}^{\perp}(\ev_{\cal T}^*T\P)\big)$$
does not vanish.
By Corollary~\ref{reg_crl2}, 
$\pi_{k-1}^{\perp}\!\circ\!{\cal D}_{{\cal T};k}$
is transverse to zero for a generic choice of the constraints
$\mu_1,\ldots,\mu_N$.
Its zero set must have dimension at least $d_{\min}({\cal T})$, if nonempty,
since the movements of points described in (1) do not effect
$\pi_{k-1}^{\perp}\circ{\cal D}_{{\cal T};k}$.
Thus, $\pi_{k-1}^{\perp}\circ{\cal D}_{{\cal T};k}$ does not vanish~if
$$\dim({\cal M}_{\cal T}(\mu))-d_{\min}({\cal T})<n-(k-1).$$
By the assumption of the lemma, this is the case as long as 
$k\!\le\!\big|\{h\!:\chi_{\cal T}h\!=\!1\}\big|$.

\begin{crl}
\label{order1_contr_n2}
Let ${\cal T}\!=\!(\Si,[N],I;j,\under{d})$
be an effective bubble type with $d_{\hat{0}}\!=\!0$.
If $g\!=\!2$ and $n\!=\!2$, then either\\
(1) $|\hat{I}|\!=\!1$  and $j_l\!\neq\!\hat{0}$ for all $l\!\in\![N]$, or\\
(2) $|\hat{I}|\!=\!2$, $H_{\hat{0}}{\cal T}\!=\!\hat{I}$, and 
$j_l\!\neq\!\hat{0}$ for all $l\!\in\![N]$.\\
Furthermore, in  Case (2) $\al_{\cal T}$ has full rank over all of 
${\cal M}_{\cal T}(\mu)$.
\end{crl}

\noindent
{\it Proof:} (1) By  Lemma~\ref{order1_contr_l3},
${\cal M}_{\cal T}(\mu)$ is empty, unless $ng-|\hat{I}|\!\ge\! 1$, 
i.e.~$|\hat{I}|\!\le\! 3$. Suppose $|\hat{I}|\!=\!3$. 
If $|H_{\hat{0}}{\cal T}|\!\ge\! 2$,
$$ng-|\hat{I}|-|H_{\hat{0}}{\cal T}|\le 4-3-2<0,$$
and thus ${\cal M}_{\cal T}(\mu)$ is empty by  Lemma~\ref{order1_contr_l3}.
If $|H_{\hat{0}}{\cal T}|\!=\!1$, 
$\chi_{\cal T}h\!\neq\! 1$ for some $h\!\in\!\hat{I}$,
and~thus
$$n-\big|\{h\!:\chi_{\cal T}h\!=\!1\}\big|\ge 2-(|\hat{I}|-1)= 0= 
 ng -|\hat{I}|-|H_{\hat{0}}{\cal T}|,$$
and by Lemma~\ref{order1_contr_l3}
every compact subset of ${\cal M}_{\cal T}(\mu)$ is ${\cal T}$-semiregular.
The space ${\cal M}_{\cal T}(\mu)$ is compact, since by the above
${\cal M}_{{\cal T}'}(\mu)=\emptyset$ if ${\cal T}'<{\cal T}$.
Corollary~\ref{order1_contr} then implies that  ${\cal M}_{\cal T}(\mu)$ 
is not effective, i.e. ${\cal T}$ is not effective.\\
(2) Suppose $|\hat{I}|\!=\!2$.
If $|H_{\hat{0}}{\cal T}|\!=\!2$ and $j_l\!=\!\hat{0}$ 
for some $l\!\in\![N]$, 
$$ng-|\hat{I}|-|H_{\hat{0}}{\cal T}|-|M_{\hat{0}}{\cal T}|\le 4-2-2-1<0,$$
and thus ${\cal M}_{\cal T}(\mu)$ is empty by Lemma~\ref{order1_contr_l3}.
If $|H_{\hat{0}}{\cal T}|\!=\!1$, then 
$\chi_{\cal T}h\!\neq\! 1$ for some $h\!\in\!\hat{I}$, and~thus
$$n-\big|\{h\!:\chi_{\cal T}h\!=\!1\big\}|=2-1
= ng-|\hat{I}|-|H_{\hat{0}}{\cal T}|,$$
and it follows from 
Lemma~\ref{order1_contr_l3} and Corollary~\ref{order1_contr},
that 
every compact subset of ${\cal M}_{\cal T}(\mu)$ is not effective. 
Furthermore, 
$\bar{\cal M}_{\cal T}(\mu)\!-\!{\cal M}_{\cal T}(\mu)$
consists of three-bubble strata, all of which are not effective 
by (1) above. Thus, ${\cal T}$ is not effective, 
unless $\io_h\!=\!\hat{0}$ for all $h\!\in\!\hat{I}$ and  
$j_l\!\neq\!\hat{0}$ for all $l\!\in\![N]$.
The second statement about the $|\hat{I}|\!=\!2$ case is immediate from 
Lemma~\ref{order1_contr_l3}.\\
(3) Finally, suppose $|\hat{I}|\!=\!1$ and $j_l\!=\!\hat{0}$ 
for some $l\!\in\![N]$.     Then,
$$n-\big|\{h\!:\chi_{\cal T}h\!=\!1\}\big|=2-1
\ge ng-|\hat{I}|-|H_{\hat{0}}{\cal T}|-|M_{\hat{0}}{\cal T}|,$$
and thus by Lemmas~\ref{order1_contr_l2} and~\ref{order1_contr_l3},
every compact subset of ${\cal M}_{\cal T}(\mu)$ is not effective.
Furthermore, $\bar{\cal M}_{\cal T}(\mu)\!-\!{\cal M}_{\cal T}(\mu)$
consists of two- and three-bubble strata that by (1) and (2) 
are not effective.
It~follows that ${\cal T}$ is not effective.

\begin{crl}
\label{order1_contr_n3}
Let ${\cal T}\!=\!(\Si,[N],I;j,\under{d})$
be an effective bubble type with $d_{\hat{0}}\!=\!0$.
If $g\!=\!2$ and $n\!=\!3$, then either\\
(1) $|\hat{I}|\!=\!1$, or\\
(2a) $|\hat{I}|\!=\!2$, $H_{\hat{0}}{\cal T}\!=\!\hat{I}$, 
 and $j_l\!\neq\!\hat{0}$ for all $l\!\in\![N]$,  or\\
(2b) $|\hat{I}|\!=\!2$, $H_{\hat{0}}{\cal T}\!\neq\!\hat{I}$,
 and $j_l\!\neq\!\hat{0}$ for all $l\!\in\![N]$,  or\\
(3a) $|\hat{I}|\!=\!3$, $H_{\hat{0}}{\cal T}\!=\!\hat{I}$,
and $j_l\!\neq\!\hat{0}$ for all $l\!\in\![N]$, or\\
(3b) $|\hat{I}|\!=\!3$,  $\io_h\!=\!\hat{1}$ for some $\hat{1}\!\in\!\hat{I}$
and all $h\!\in\!\hat{I}\!-\!\{\hat{1}\}$,
$d_{\hat{1}}\!=\!0$, and $j_l\!\neq\!\hat{0},\hat{1}$ for 
all $l\!\in\![N]$.\\
Furthermore, in Case (3a) $\al_{\cal T}$ has full rank on all of
${\cal M}_{\cal T}(\mu)$.
\end{crl}

\noindent
{\it Proof:} (1) Similarly to the proof of Corollary~\ref{order1_contr_n2},
${\cal M}_{\cal T}(\mu)$ is empty unless $|\hat{I}|\!\le\! 5$.
If $|\hat{I}|\!=\!5$, ${\cal M}_{\cal T}(\mu)$ is compact and
$|H_{\hat{0}}{\cal T}|\!=\!1$.
Let $\hat{1}\!\in\!\hat{I}$ be such that $\io_{\hat{1}}\!=\!\hat{0}$.
If $d_{\hat{1}}\!>\!0$,
$$n-\big|\{h\!:\chi_{\cal T}h\!=\!1\}\big|=3-1>0=
ng-|\hat{I}|-|H_{\hat{0}}{\cal T}|,$$
and ${\cal M}_{\cal T}(\mu)$ is not effective by Lemma~\ref{order1_contr_l3}
and Corollary~\ref{order1_contr}.
Suppose $d_{\hat{1}}\!=\!0$.
Then, $|H_{\hat{1}}{\cal T}|\!=\!2$;
otherwise ${\cal M}_{\cal T}(\mu)$ is empty by Lemma~\ref{order1_contr_l3}.
It follows that
$$n-\big|\{h\!:\chi_{\cal T}h\!=\!1\}\big|\ge
3-(|\hat{I}|-2)= 0=ng-|\hat{I}|-|H_{\hat{0}}{\cal T}|.$$
Thus, by Lemma~\ref{order1_contr_l3} and Corollary~\ref{order1_contr},
${\cal T}$ is not effective.\\
(2) Suppose $|\hat{I}|\!=\!4$. If $|H_{\hat{0}}{\cal T}|\!\ge\! 3$,
${\cal M}_{\cal T}(\mu)$ is empty by Lemma~\ref{order1_contr_l3}.
Let $\hat{1}\!\in\!\hat{I}$ be as above.
If $|H_{\hat{0}}{\cal T}|\!=\!2$,
$$n-\big|\{h\!:\chi_{\cal T}h\!=\!1\}\big|\ge
3-(|\hat{I}|-1)= 0=ng-|\hat{I}|-|H_{\hat{0}}{\cal T}|.$$
If $|H_{\hat{0}}{\cal T}|\!=\!1$ and $d_{\hat{1}}\!>\!0$,
$$n-\big|\{h\!:\chi_{\cal T}h\!=\!1\}\big|=
3-1>1=ng-|\hat{I}|-|H_{\hat{0}}{\cal T}|.$$
If $|H_{\hat{0}}{\cal T}|\!=\!1$, $d_{\hat{1}}\!=\!0$, and
$|H_{\hat{1}}{\cal T}|\!=\!3$,
$$n-\big|\{h\!:\chi_{\cal T}h\!=\!1\}\big|\ge
3-(|\hat{I}|-1)=0=
ng-|\hat{I}|-|H_{\hat{0}}{\cal T}|-(|H_{\hat{1}}{\cal T}|-2).$$
Finally, if $|H_{\hat{0}}{\cal T}|\!=\!1$, $d_{\hat{1}}\!=\!0$, and 
$|H_{\hat{1}}{\cal T}|\!=\!2$,
$$n-\big|\{h\!:\chi_{\cal T}h\!=\!1\}\big|\ge
3-(|\hat{I}|-2)=1=ng-|\hat{I}|-|H_{\hat{0}}{\cal T}|.$$
Thus, by Corollary~\ref{order1_contr} and Lemma~\ref{order1_contr_l3},
in all four cases, no compact subset of ${\cal M}_{\cal T}(\mu)$ is effective.
Since $\bar{\cal M}_{\cal T}(\mu)\!-\!{\cal M}_{\cal T}(\mu)$
consists of five-bubble strata that are not effective by (1) above,
it follows that ${\cal T}$ is not effective.\\
(3) Suppose $|\hat{I}|\!=\!3$. 
If  $H_{\hat{0}}{\cal T}\!=\!\hat{I}$ and $j_l\!=\!\hat{0}$ for some 
$l\!\in\![N]$,
$$ng-|\hat{I}|-|H_{\hat{0}}{\cal T}|-|M_{\hat{0}}{\cal T}|=6-3-3-1<0,$$
and thus ${\cal M}_{\cal T}(\mu)$ is empty by Lemma \ref{order1_contr_l3}.
If $|H_{\hat{0}}{\cal T}|\!=\!2$, 
$$n-\big|\{h\!:\chi_{\cal T}h\!=\!1\}\big|\ge 3-(|\hat{I}|-1)=1
\ge ng-|\hat{I}|-|H_{\hat{0}}{\cal T}|.$$
If $|H_{\hat{0}}{\cal T}|\!=\!1$ and $d_{\hat{1}}\!>\!0$,
$$n-\big|\{h\!:\chi_{\cal T}h\!=\!1\}\big|=2
=ng-|\hat{I}|-|H_{\hat{0}}{\cal T}|.$$
If $|H_{\hat{0}}{\cal T}|\!=\!1$ and $|H_{\hat{1}}{\cal T}|\!=\!1$,
$$n-\big|\{h\!:\chi_{\cal T}h\!=\!1\}\big|=2
=ng-|\hat{I}|-|H_{\hat{0}}{\cal T}|.$$
Thus, in all three cases, by 
Lemma~\ref{order1_contr_l3} and Corollary~\ref{order1_contr},
no compact subset of ${\cal M}_{\cal T}(\mu)$ is effective.
Since 
$\bar{\cal M}_{\cal T}(\mu)\!-\!{\cal M}_{\cal T}(\mu)$
consists of four- and five-bubble strata that are not effective by 
(1) and (2) above, ${\cal T}$ is not effective
in these three cases.
On the other hand, if $|H_{\hat{0}}{\cal T}|\!=\!2$,
$j_l\!=\!\hat{0}$ or $j_l\!=\!\hat{1}$ for some $l\!\in\![N]$,
and $d_{\hat{1}}\!=\!0$,
$$n-\big|\{h\!:\chi_{\cal T}h\!=\!1\}\big|\ge 1
\ge ng-|\hat{I}|-|H_{\hat{0}}{\cal T}|-|M_{\hat{0}}{\cal T}|-
\big(|H_{\hat{1}}{\cal T}|+|M_{\hat{1}}{\cal T}|-2\big).$$
Thus, by Lemmas~\ref{order1_contr_l2} and \ref{order1_contr_l3},
no compact subset of ${\cal M}_{\cal T}(\mu)$ is effective.
Similarly to the above, it follows that ${\cal T}$ is not effective.\\
(4) Suppose $|\hat{I}|\!=\!2$ and $j_l\!=\!\hat{0}$ for some $l\!\in\![N]$.
If $|H_{\hat{0}}{\cal T}|\!=\!2$,
$$n-\big|\{h\!:\chi_{\cal T}h\!=\!1\}\big|\ge 1 \ge 
ng-|\hat{I}|-|H_{\hat{0}}{\cal T}|-|M_{\hat{0}}{\cal T}|.$$
If $|H_{\hat{0}}{\cal T}|\!=\!1$,
$$n-\big|\{h\!:\chi_{\cal T}h\!=\!1\}\big|=2\ge 
ng-|\hat{I}|-|H_{\hat{0}}{\cal T}|-|M_{\hat{0}}{\cal T}|.$$
Thus, in either case,  no compact subset of ${\cal M}_{\cal T}(\mu)$
is effective by Lemmas~\ref{order1_contr_l2} and \ref{order1_contr_l3}.
Furthermore,
$$\bar{\cal M}_{\cal T}(\mu)-{\cal M}_{\cal T}(\mu)=
\bigcup_{{\cal T}'<{\cal T}}{\cal M}_{{\cal T}'}(\mu),$$
where ${\cal T}'$ is either a four- or five-bubble strata,
or a three bubble-strata 
${\cal T}'=(\Si,[N],I';j',d')$
such that either  $|H_{\hat{0}}{\cal T}|\!=\!1$, or
$d_{\hat{1}'}'\!=\!0$ and $j_l'\!=\!\hat{0}$ or $\hat{1}'$.
By (1)-(3) above, none of such bubble types is effective,
and thus ${\cal T}$ is not effective.

\subsection{Second-Order Estimate for $\psi_{{\cal T},t\nu}^{\mu}$, Case 1}
\label{order2_case1}

\noindent
We now refine the first-order estimate for $\psi_{{\cal T},t\nu}^{\mu}$ 
along the sets on which  the section $\al^-$ defined above
does not have full rank.
These are precisely the sets on which the section
$\bar{\cal D}_{{\cal T},|\{h:\chi_{\cal T}h=1\}|}$
defined in the proof of Lemma~\ref{order1_contr_l3}
does not have full rank.\\

\noindent
One set on which  $\bar{\cal D}_{{\cal T},|\{h:\chi_{\cal T}h=1\}|}$
fails to have full rank is the zero set of ${\cal D}_{{\cal T},h_1}$.
If $n\!=\!2,3$, by Lemma~\ref{order1_contr_l3}, ${\cal D}_{{\cal T},h_1}$
does not vanish unless $h_1$ is the only element of the set
$\{h\!:\chi_{\cal T}h\!=\!1\}$.
Thus, we assume that this is the case.
We denote the zero-locus of ${\cal D}_{{\cal T},h_1}$
by ${\cal S}_{{\cal T},1}\!\subset\!{\cal M}_{\cal T}$,
which will be abbreviated as ${\cal S}$ in this subsection.
Since ${\cal D}_{{\cal T},h_1}$ is transversal to zero
by Corollary~\ref{reg_crl2}, ${\cal S}$
is a complex submanifold of ${\cal M}_{\cal T}$ of codimension~$n$.
Its normal bundle ${\cal NS}$ in ${\cal M}_{\cal T}$ 
is the restriction of $L_{{\cal T},h_1}^*\otimes\ev^*T\P$ 
to~${\cal S}_{{\cal T},1}$.
Let $(\Phi_{\cal S},\Phi_{\cal S}^{\mu})$
be a  regularization 
of ${\cal S}_{{\cal T},1}(\mu)
\!\equiv\!{\cal S}\cap{\cal M}_{\cal T}(\mu)$.
This regularization can be chosen so~that
\begin{equation}\label{order2_e}
{\cal D}_{{\cal T},h_1}\tilde{\phi}_{\cal S}(b,X)=
\Pi_{b,\tilde{\phi}_{\cal S}(b,X)}X
\quad\forall~ (b,X)\in {\cal N}\tilde{\cal S}=\ev^*T\P,
\end{equation}
where $\tilde{\phi}_{\cal S}$ is the lift of $\phi_{\cal S}$
to the preimage $\tilde{\cal S}$ of ${\cal S}$
and its normal bundle ${\cal N}\tilde{\cal S}$
in ${\cal M}_{\cal T}^{(0)}$; see Subsection~\ref{gl-orient2_sec}
in~\cite{Z}.
The bundle ${\cal NS}$ carries a natural norm induced
by the $g_{\P,\ev}$-metric on~$\P$.
Denote by $F{\cal S}$ and $F^{\eset}{\cal S}$
the bundles described in Subsection~\ref{si_str_sec}
corresponding to the submanifold~${\cal S}_{{\cal T},1}$.
Let \hbox{$\hat{1}\!\in\! H_{\hat{0}}{\cal T}$} be the
unique element such that \hbox{$h_1\!\in\!\bar{D}_{\hat{1}}{\cal T}$}.
If $\big[b;X,\ups\big]\!\in\! F{\cal S}={\cal NS}\oplus F{\cal T}$, put
$${^{(2)}\al_{{\cal T};1}}(X,\ups)=
X(b_{\ups})s_{\Si,x_{\hat{1}}}\tilde{v}_{h_1}+
\al_{{\cal T},h_1}^{(2)}(\ups).$$

\begin{lmm}
\label{dbar_2}
There exist $\de,C\!\in\! C^{\i}({\cal S};\Bbb{R}^+)$ such that for all
$\vp\!=\![(b;X,\ups)]\!\in\! F^{\eset}{\cal S}_{\de}$,
$$\Big\|\pi_{\Phi_{\cal S}(\vp),-}^{0,1}\bar{\partial}u_{\Phi_{\cal S}(\vp)}
+\tilde{R}_{\Phi_{\cal S}(\vp)}\Pi_{b,\phi_{\cal S}(X)}
{^{(2)}\al_{{\cal T};1}}(X,\ups)\Big\|_2
\le C(b)|\ups|\big(|\ups|_{h_1}^2+|X||\ups|_{h_1}\big).$$
\end{lmm}

\noindent
{\it Proof:} The proof is almost identical to the proof
of Lemma~\ref{dbar_1}.
The only difference is that we use two terms of 
the power series of Proposition~\ref{dbar_gen2prp}.
We then make use of the assumption~\e_ref{order2_e} on $\phi_{\cal S}$
and smooth dependence of ${\cal D}_{{\cal T},h_1}^{(2)}$
on~$X$.

\begin{lmm}
\label{psit_21}
There exist $\de,C\!\in\! C^{\i}({\cal S}_{{\cal T},1}(\mu);\Bbb{R}^+)$
such that for all
$\vp\!=\![(b;X,\ups)]\!\in\! F^{\eset}{\cal S}_{\de}$,
$$\Big\| \psi_{{\cal S},t\nu}^{\mu}(\vp)-
\big(t\bar{\nu}_{\ev(b)}+{^{(2)}\al_{{\cal T};1}}(X,\ups)\big)\Big\|_2 
\le C(b)\big(t+|\ups|^{\frac{1}{p}}\big) 
            \big(t+|\ups|_{h_1}^2+|X||\ups|_{h_1}\big).$$
\end{lmm}

\noindent
{\it Proof:} This claim follows from Lemmas~\ref{nuterm} 
and~\ref{dbar_2} in a way analogous to the proof of Lemma~\ref{psit_11}.
The only difference is that we need to improve
the estimate on $\pi_{\ups,-}^{0,1}D_{\ups}\xi_{\ups,t\nu}$
made in the proof of Lemma~\ref{psit_11}.
Let $\{\psi_j\}$ be an orthonormal basis for ${\cal H}_{\Si}^{0,1}$,
such that $\psi_1\!\in\!{\cal H}_{\Si}^+\big(\tilde{x}_{h_1}(\ups)\big)$,
and $\{X_i\}$ an orthonormal basis for~$T_{\ev(\phi_{\cal S}(X))}\P$.
By Theorem~\ref{si_str}, with $\ups(X)=\Phi_{\cal S}(\vp)$,
\begin{equation}\label{psit_21_e1}\begin{split}
\Big|\big\lan\big\lan \pi_{\ups(X),-}^{0,1}
D_{\ups(X)}\xi_{\ups(X),t\nu},
 R_{\ups(X)}X_i\psi_j\big\ran\big\ran\Big|
&=\Big|\big\lan\big\lan \xi_{\ups(X),t\nu},
      D_{\ups(X)}^*R_{\ups(X)}X_i\psi_j\big\ran\big\ran\Big|\\
&\le C(b)(t+|\ups|^{\frac{1}{p}})
\|D_{\ups(X)}^*R_{\ups(X)}X_j\psi_j\|_{C^0}.
\end{split}\end{equation}
Since $\xi\!\in\!\tilde{\Ga}_+(\ups)$,
by construction in Subsection~\ref{tangent_sec}, 
\begin{equation}\label{psit_21_e2}
\big\lan\big\lan \xi_{\ups(X),t\nu},
          D_{\ups(X)}^*R_{\ups(X)}X_i\psi_1\big\ran\big\ran=0.
\end{equation}
On the other hand, since $\psi_2|_{\tilde{x}_{h_1}(\ups)}\!=\!0$
and $\|\na\psi_2\|_{g_{\phi_{\cal S}(X),\hat{0}},C^0}\!\le\! C(b)$, 
by equation~\e_ref{coker_e4}
\begin{equation}\label{psit_21_e4}
\big\|D_{\ups(X)}^*R_{\ups(X)}X_i\psi_2
     \big\|_{C^0(\tilde{A}_{\ups(X),h_1})}
\le C(b)|\ups|_{h_1}^2,
\end{equation}
where $\tilde{A}_{\ups(X),h_1}$ is the annulus defined in 
Lemma~\ref{approx_coker}. 
By equations~\e_ref{psit_21_e1}-\e_ref{psit_21_e4},
$$\Big| \pi_{\ups(X),-}^{0,1}D_{\ups}\xi_{\ups(X),t\nu}
\Big|\le C(b)(t+|\ups|^{\frac{1}{p}})|\ups|_{h_1}^2.$$
\\

\noindent
The next step is to apply Lemma~\ref{top_l1} or Corollary~\ref{top_l2c}
whenever possible.     Let 
\begin{gather*}
F^+\!=\ev^*T\P\otimes\!
\bigotimes_{i\in\hat{I},h_1\in D_i{\cal T}}\!\! F_i{\cal T},\quad
F^-\!=F{\cal T},\quad
\tilde{F}^-\!=\!
\Big(\!\bigotimes_{i\in\hat{I},h_1\in\bar{D}_i{\cal T}}\!\! F_i{\cal T}
\Big)^{\otimes2}\!,\quad
{\cal O}^{\pm}\!={\cal H}_{\Si}^{\pm}\otimes\ev^*T\P;\\
\al^+\big([X,\ups]\big)=Xs_{\Si,x_{\hat{1}}}\tilde{v}_{h_1},\quad
\phi\big([b,v_{\hat{I}}]\big)
=\big[b,\tilde{v}_{h_1}\otimes\tilde{v}_{h_1}\big],~~
\al^-\big(\phi(\ups)\big)\equiv 
\pi_{x_{\hat{1}}(b_{\ups})}^-\al_{\cal T}^{(2)}(\ups).
\end{gather*}
Note that $\al^+\!\!\in\!\Ga({\cal S};F^{+*}\otimes{\cal O}^+)$,
since $\pi^-\!\circ\! Xs_{\Si}\!=\!0$.
Since the map \hbox{$(X,\ups)\!\lra\! (X\!\otimes\!\tilde{v}_{h_1},\ups)$} 
is injective
on $F^{\eset}{\cal T}$, we can view $\psi_{{\cal S},t\nu}^{\mu}$
as a map on an open subset \hbox{of $F^-\oplus F^+$.}
Analogously to the first-order case of Subsection~\ref{order1},
subset $K\subset{\cal S}_{{\cal T},1}(\mu)$ will be called
{\it{second-order regular}} if 
$\al^-$ has full rank over~$K$.

\begin{crl}
\label{order2_1_contr}
Suppose $d$ is a positive integer, 
${\cal T}\!=\!(\Si,[N],I;j,\under{d})$ is
a simple bubble type, with $d_{\hat{0}}\!=\!0$ and
$\sum\limits_{i\in I}d_i\!=\!d$, and 
$\mu$ is an $N$-tuple of constraints in general position 
such that
$$\codim_{\Bbb{C}}\mu=d(n+1)-n(g-1)+N.$$
Let $\nu\!\in\!\Ga(\Si\times\P;
          \La^{0,1}\pi_{\Si}^*T^*\Si\otimes\pi_{\P}^*T\P)$
be a generic section.
If $|\hat{I}|\!>\!1$, for every 
second-order regular compact subset~$K$ of~${\cal S}_{{\cal T},1}(\mu)$,
there exist a neighborhood $U_K$ of $K$ in $\bar{C}^{\i}_{(d;[N])}(\Si;\mu)$
and $\ep_K\!>\!0$ such that for any $t\!\in\!(0,\ep_K)$,
\hbox{$U_K\cap{\cal M}_{\Si,d,t\nu}(\mu)\!=\!\eset$.}
If $|\hat{I}|\!=\!1$,
there exists a compact regular subset $K_{{\cal T},1}$ 
of~${\cal S}_{{\cal T},1}(\mu)$ with the following property.
If $K$ is a compact subset of ${\cal S}_{{\cal T},1}(\mu)$
containing~$K_{{\cal T},1}$,
there exist a neighborhood $U_K$ of $K$ in $\bar{C}^{\i}_{(d;[N])}(\Si;\mu)$
and $\ep_K\!>\!0$ such that for all $t\!\in\!(0,\ep_K)$,
the signed cardinality of $U_K\cap{\cal M}_{\Si,d,t\nu}(\mu)$ 
equals to twice the signed number of zeros of the~map
\begin{equation}\label{order2_1contr_e}
T{\Si}^{\otimes2}\otimes L_{{\cal T},\hat{1}}^{\otimes 2}
    \big|{\cal S}_{{\cal T},1}(\mu)\lra 
{\cal H}_{\Si}^-\otimes\ev^*T\P,\quad
[b,v]\lra \bar{\nu}_b^-+\al^{(2,-)}\big([b,v]\big).
\end{equation}
\end{crl}

\noindent
{\it Proof:} In either case, by Theorem~\ref{si_str},
there exist a neighborhood $U_K$ of $K$ in $\bar{C}^{\i}_{(d;[N])}(\Si;\mu)$
and $\de_K,\ep_K\!>\!0$ such that for any $t\!\in\!(0,\ep_K)$,
there exists a sign-preserving bijection between  
\hbox{$U_K\cap{\cal M}_{\Si,d,t\nu}(\mu)$} and the zeros of
$\psi_{{\cal S},t\nu}^{\mu}$ on 
$F^{\eset}{\cal S}_{\de_K}\big|\big(U_K\cap {\cal S}_{{\cal T},1}(\mu)\big)$,
provided  $U_K\cap {\cal S}_{{\cal T},1}(\mu)$ is precompact 
in~${\cal S}_{{\cal T},1}(\mu)$.
If $K$ is second-order regular, 
$U_K$ can be chosen so that the closure of
$U_K\cap {\cal S}_{{\cal T},1}(\mu)$ in ${\cal S}_{{\cal T},1}(\mu)$
is also second-order regular.
Since $K$ is regular and $\al^+$ is injective on all fibers,
$$|\ups|_{h_1}^2=|\phi(\ups)|\le C_K\big|\al^-(\phi(\ups))\big|
\Lra 
|\ups|_{h_1}^2+|X||\ups|_{h_1}\le C_K'
\big|{^{(2)}\al_{{\cal T};1}}(X,\ups)\big|
\quad\forall(X,\ups)\!\in\!  F^{\eset}{\cal S}_{\de_K}\big|K,$$
where $C_K,C_K'\!>\!0$ depend only on $K$.
Thus, by Lemma~\ref{psit_21},
$$\Big\|\psi_{{\cal S},t\nu}^{\mu}(\vp)-
\big(t\bar{\nu}_{\ev(b_{\vp})}+
{^{(2)}\al_{{\cal T};1}}(\vp)\big)\Big\|_2
\le C_K\big(t+|\vp|^{\frac{1}{p}}\big)
\big(t+\big|{^{(2)}\al_{{\cal T};1}}(\vp)\big|\big)
\quad\forall\vp\!\in\! F^{\eset}{\cal S}_{\de_K}\big|K,$$
where $C_K\!>\!0$ depends only on $K$.
The first claim now follows from Lemma~\ref{top_l1}.
The second follows from Corollary~\ref{top_l2c}, provided that for
a generic $\nu$ the set of zeros of the map in~\e_ref{order2_1contr_e}
is second-order regular and finite.
This last fact is proved by the same argument as
Lemma~\ref{order1_contr_l1}.

\subsection{Third-Order Estimate for $\psi_{{\cal T},t\nu}^{\mu}$, Case 1}
\label{order3_case1}

\noindent
We continue with the case of Section~\ref{order2_case1}.
Then
$$\al^-([b,\tilde{v}_{h_1}])=
({\cal D}_{{\cal T},h_1}^{(2)}b)s_{b,x_{\hat{1}}}^{(2,-)}(\tilde{v}_{h_1}).$$
By Corollary~\ref{reg_crl2}, for a generic choice of the constraints
$\mu_1,\ldots,\mu_N$, ${\cal D}_{{\cal T},h_1}^{(2)}$ is transversal
to zero along ${\cal S}_{{\cal T},1}(\mu)$ \hbox{if $d_{h_1}\!\ge\! 2$.}
Since the zero set of ${\cal D}_{{\cal T},h_1}^{(2)}$ must
have dimension at least $d_{\min}({\cal T})\!\ge\! 1$
by the same argument as in the proof of Lemma~\ref{order1_contr_l3},
${\cal D}_{{\cal T},h_1}^{(2)}$ does not vanish along 
${\cal S}_{{\cal T},1}(\mu)$ \hbox{if $d_{h_1}\!\ge\! 2$.}
On the other hand, if $d_{h_1}\!=\!1$, ${\cal S}_{{\cal T},1}\!=\!\emptyset$,
since the differential of any degree-one holomorphic map from
$S^2$ to $\P$ is nowhere~zero.
In fact, ${\cal S}_{{\cal T},1}(\mu)\!=\!\emptyset$ even for $d_{h_1}\!=\!2$,
since the image of any degree-two map with a somewhere vanishing 
differential is a line, and no line intersects $\mu_1,\ldots,\mu_N$
\hbox{if $n\!=\!2,3$.}
Thus, we can \hbox{assume $d_{h_1}\!\ge\! 3$}.
It follows that  the only way 
the above homomorphism $\al^-$ 
can fail to have full rank on $\tilde{F}^-$ is 
\hbox{if  $s_{b,x_{\hat{1}}}^{(2,-)}\!=\!0$.}
While $s_{b,x_{\hat{1}}}^{(2)}$ depends on the choice of the 
metric $g_{b,{\hat{0}}}$ on $\Si$, the~section 
$s^{(2,-)}\!\in\!
\Ga\big(\Si;T^*\Si^{\otimes2}\otimes{\cal H}_{\Si}^-\big)$
is independent of the metric and is globally defined on $\Si$.
This can be seen by a direct computation.
It has transverse zeros at the six branch points 
of the double cover $\Si\!\lra\!\Bbb{P}^1$ induced by $s_{\Si}$;
see \cite[p246]{GH}.
Denote by $z_1,\ldots,z_6$ these six points.
Then the set on which $\al^-$ fails to have full rank 
is $\bigcup\limits_{m\in[6]}{\cal S}_{{\cal T},1}^{(m)}(\mu)$,
where
$${\cal S}_{{\cal T},1}^{(m)}=\big\{b\!\in\!{\cal S}_{{\cal T},1}\!: 
 x_{\hat{1}}(b)\!=\!z_m\},\quad
{\cal S}_{{\cal T},1}^{(m)}(\mu)={\cal S}_{{\cal T},1}^{(m)}\cap
{\cal M}_{\cal T}(\mu).$$
The sets ${\cal S}_{{\cal T},1}^{(m)}$ are obviously disjoint.\\

\noindent
Since the normal bundle of ${\cal S}_{{\cal T},1}^{(m)}$ in 
${\cal S}_{{\cal T},1}$  is $T_{z_m}\Si$,
the normal bundle ${\cal NS}$ of ${\cal S}_{{\cal T},1}^{(m)}$
in ${\cal M}_{\cal T}(\mu)$ is \hbox{$T_{z_m}\Si\oplus {\cal NS}_1$},
where ${\cal NS}_1$ is the normal bundle of ${\cal S}_{{\cal T},1}$  
in ${\cal M}_{\cal T}(\mu)$, as described in the previous subsection.
Let 
$\big(\Phi_{\cal S},\Phi_{\cal S}^{\mu}\big)$
be a regularization of ${\cal S}_{{\cal T},1}^{(m)}(\mu)$
induced by the regularization of ${\cal S}_{{\cal T},1}(\mu)$
described in Subsection~\ref{order2_case1}.
In particular,
\begin{equation}\label{order3_e}
{\cal D}_{{\cal T},h_1}\tilde{\phi}_{\cal S}(b,w,X)=
\Pi_{b,\tilde{\phi}_{\cal S}(b,w,X)}X
\quad\forall~ 
(b,w,X)\in T_{z_m}\Si\oplus {\cal N}\tilde{\cal S}_1=T_{z_m}\Si\oplus\ev^*T\P,
\end{equation}
where $\tilde{\phi}_{\cal S}$ is the lift of $\phi_{\cal S}$
to~${\cal M}_{\cal T}^{(0)}$.
We can also assume that $\Phi_{\cal S}^{\mu}$ is given
by the $g_{\P,b}$-parallel transport on~${\cal N}_b{\cal S}_1$.
The bundle ${\cal NS}$ carries a natural norm induced
by the $g_{\P,\ev}$-metric on~$\P$ and $g_{\cdot,\hat{0}}$-metric on~$\Si$.
Denote by $F{\cal S}$ and $F^{\eset}{\cal S}$
the bundles described in Subsection~\ref{si_str_sec}
corresponding to the submanifold~${\cal S}_{{\cal T},1}^{(m)}$.
If $(b,w,X,\ups)\!\in\! F^{\eset}{\cal S}$ is sufficiently small, let
$$\tilde{x}_{\hat{1}}(w,\ups)
=\tilde{x}_{\hat{1}}\big(\phi_{\cal S}(w,X,\ups)\big)
=\tilde{x}_{\hat{1}}\big(\phi_{\cal S}(w,0,\ups)\big)\in\Si.$$
We identify a small neighborhood of $z_m$ in $\Si$ with a neighborhood 
of $0$ in $T_{z_m}\Si$ via the $g_{b,\hat{0}}$-exponential map.
Put
$$\tilde{\al}(w,X,\ups)
\!=\! (Xb) s_{\Si,\tilde{x}_{\hat{1}}(w,\ups)}(\tilde{v}_{h_1})
\!+\Pi_{b,\phi_{\cal S}(b,X)}^{-1}
\!\big({\cal D}_{{\cal T},h_1}^{(2)}\phi_{\cal S}(b,X)\big)
                  s_{b,\tilde{x}_{\hat{1}}(w,\ups)}^{(2)}(\tilde{v}_{h_1})
\!+\!\big({\cal D}_{{\cal T},h_1}^{(3)}b\big)
                  s_{b,z_m}^{(3)}(\tilde{v}_{h_1}).$$
If $(b,w,X,\ups)\!\in\! F^{\eset}{\cal S}|{\cal S}_{{\cal T},1}^{(m)}(\mu)$ 
is sufficiently small, let
$$\tilde{\al}^{\mu}(w,X,\mu)
= (Xb) s_{\Si,\tilde{x}_{\hat{1}}(w,\ups)}(\tilde{v}_{h_1})
+\big({\cal D}_{{\cal S},t\nu}^{\mu,(2)}(w,X,\ups)\big)
                  s_{b,\tilde{x}_{\hat{1}}(w,\ups)}^{(2)}(\tilde{v}_{h_1})
+\big({\cal D}_{{\cal T},h_1}^{(3)}b\big)s_{b,z_m}^{(3)}(\tilde{v}_{h_1}),$$
where, with $\varphi_{{\cal S},t\nu}^{\mu}$ as in 
Theorem~\ref{si_str},
$${\cal D}_{{\cal S},t\nu}^{\mu,(2)}(w,X,\ups)=
\Pi_{\phi_{\cal S}^{\mu}\varphi_{{\cal S},t\nu}^{\mu}(w,X,\ups),
\phi_{\cal S}\Phi_{\cal S}^{\mu}\varphi_{{\cal S},t\nu}^{\mu}(w,X,\ups)}^{-1}
\Pi_{b,\phi_{\cal S}^{\mu}\varphi_{{\cal S},t\nu}^{\mu}(w,X,\ups)}^{-1}
\big({\cal D}_{{\cal T},h_1}^{(2)}\phi_{\cal S}\Phi_{\cal S}^{\mu}
     \varphi_{{\cal S},t\nu}^{\mu}(w,X,\ups)\big).$$

\begin{lmm}
\label{dbar_3}
There exist $\de,C\!\in\! C^{\i}({\cal S}_{{\cal T},1}^{(m)};\Bbb{R}^+)$
such that for all
 $\vp\!=\![(b,w,X,\ups)]\!\in\! F^{\eset}{\cal S}_{\de}$,
$$ \Big\| 
\pi^{0,1}_{\Phi_{\cal S}(\vp),-}\bar{\partial}u_{\Phi_{\cal S}(\vp)}
+\tilde{R}_{\Phi_{\cal S}(\vp)}\Pi_{b,\phi_{\cal S}(X)}
\tilde{\al}(w,X,\ups)\Big\|_2
\le C(b)|\vp||\ups|_{h_1}^3.$$
\end{lmm}

\noindent
{\it Proof:} The proof is the same as that of Lemma~\ref{dbar_2},
except here we use the first three terms of the expansion
of Proposition~\ref{dbar_gen2prp}. Note that 
$|\tilde{x}_{\hat{1}}(w,\ups)|\le C(b)(|w|+|\ups|)$.

\begin{lmm}
\label{psit3_1}
There exist $\de,C\!\in\! C^{\i}({\cal S}_{{\cal T},1}^{(m)}(\mu);\Bbb{R}^+)$
such that for all 
$\vp\!=\!(b,w,X,\ups)\!\in\! F^{\eset}{\cal S}_{\de}$
$$\Big\|\psi_{{\cal S},t\nu}^{\mu}(\vp)-
\big(t\bar{\nu}_{\ev(b)}+
\tilde{\al}^{\mu}(w,X,\ups)\big)\Big\|_2\le
C(b)(t+|\vp|^{\frac{1}{p}})
\big(t+|\ups|_{h_1}^3+|\tilde{x}_1(w,\ups)||\ups|_{h_1}^2\big).$$
\end{lmm}

\noindent
{\it Proof:} 
The proof is similar to the proofs of Lemmas~\ref{psit_11} and~\ref{psit_21}, 
but we need to obtain an even stronger bound on
$$\big\|\pi_{\Phi_{\cal S}(\vp),-}^{0,1}
D_{\Phi_{\cal S}(\vp)}\xi_{\Phi_{\cal S}(\vp),t\nu}\big\|_2.$$
Let $\{\psi_j\}$ be an orthonormal basis for ${\cal H}_{\Si}^{0,1}$
such that $\psi_1\!\in\!{\cal H}_{\Si}^+\big(\tilde{x}_{h_1}(w,\ups)\big)$,
and $\{X_i\}$ an orthonormal basis 
\hbox{for $T_{\ev(\phi_{\cal S}(X,\ups))}\P$.}
Then, as in the proof of Lemma~\ref{psit_21}, 
with $\ups(\vp)=\Phi_{\cal S}(\vp)$,
\begin{gather}\label{psit3_1e1}
\big\lan\big\lan
D_{\Phi_{\cal S}(\vp)}\xi_{\ups(\vp),t\nu},
R_{\ups(\vp)}X_i\psi_1\big\ran\big\ran=0;\\
\label{psit3_1e2}
\begin{split}
&\Big|\big\lan\big\lan\pi_{\ups(\vp),-}^{0,1}
D_{\Phi_{\cal S}(\vp)}\xi_{\ups(\vp),t\nu},
R_{\ups(\vp)}X\psi_2\big\ran\big\ran\Big|
\le C(b)\big(t+|\ups|^{\frac{1}{p}}\big)
\big\|D_{\ups(\vp)}^*R_{\ups(\vp)}X_i\psi_2\big\|_{\ups(\vp),1}.
\end{split}
\end{gather}
The one-form $\psi_2$ vanishes at $\tilde{x}_{h_1}(w,\ups)$ by definition
and $\|\na\psi_2\|_{g_{b,\hat{0}},C^0}\le C|\tilde{x}_{h_1}(w,\ups)|$,
since the derivative of the corresponding one-form for $z_m$
vanishes.
Thus, by equation~\e_ref{coker_e4}
\begin{equation}\label{psit3_1e3}
\big\|D_{\ups(\vp)}^*R_{\ups(\vp)}X_i\psi_2
                     \big\|_{g_{\ups(\vp)},L^1(\tilde{A}_{\ups(\vp),h})}
\le C(b)(|\tilde{x}_{h_1}(w,\ups)||\ups|_{h_1}+|\ups|_{h_1}^2)|\ups|_{h_1},
\end{equation}
as needed for our bound.
Finally, we use our assumption that $\Phi_{\cal S}^{\mu}$
is given by the $g_{b,\hat{0}}$-parallel transport 
on~${\cal N}_b{\cal S}_1$.\\

\noindent
For any 
$(w,X,\ups)\!\in\! F_b^{\eset}{\cal S}|{\cal S}_{{\cal T},1}^{(m)}(\mu)$
sufficiently small, let
\begin{gather*}
Y(w,X,\ups)=(Xb)s_{\Si,\tilde{x}_1(w,\ups)}(\tilde{v}_{h_1})+
\Big({\cal D}_{{\cal S},t\nu}^{\mu,(2)}(w,X,\ups)\Big)
s_{b,\tilde{x}_{\hat{1}}(w,\ups)}^{(2,+)}(\tilde{v}_{h_1},\tilde{v}_{h_1});\\
{^{(3)}\al_{{\cal T};1}^{(m),-}}(w,\ups)=
\big({\cal D}_{{\cal T},h_1}^{(2)}b\big)
s_{b,z_m}^{(3,-)}\big(\tilde{x}_{\hat{1}}(w,\ups),\tilde{v}_{h_1},
                                       \tilde{v}_{h_1}\big)
+\big({\cal D}_{{\cal T},h_1}^{(3)}b\big)
 s_{b,z_m}^{(3,-)}(\tilde{v}_{h_1});\\
r_{{\cal T};1}^+(\ups)=\big({\cal D}_{{\cal T},h_1}^{(3)}b\big)
                  s_{b,z_m}^{(3,+)}(\tilde{v}_{h_1}),\quad
\bar{\nu}^{\pm}_b=\pi_{z_m}\bar{\nu}_b.
\end{gather*}

\begin{crl}
\label{psit3_1c}
There exist $\de,C\!\in\! C^{\i}({\cal S}_{{\cal T},1}^{(m)}(\mu);\Bbb{R}^+)$
such that for all 
$\vp\!=\![(b,w,X,\ups)]\!\in\! F^{\eset}{\cal S}_{\de}$
\begin{equation*}\begin{split}
\Big\|\pi_{x_{\hat{1}}(w,\ups)}^+\psi_{{\cal S},t\nu}^{\mu}(w,X,\ups)-
\big(t\pi_{x_{\hat{1}}(w,\ups)}^+\bar{\nu}_b+
Y(w,X,\ups)+r_{{\cal T};1}^+(\ups)\big)\Big\|_2\qquad\qquad\qquad\qquad&\\
\le C(b)(t+|\vp|^{\frac{1}{p}})
\big(t+|\ups|_{h_1}^3+|\tilde{x}_1(w,\ups)||\ups|_{h_1}^2\big);&\\
\Big\|\pi_{x_{\hat{1}}(w,\ups)}^-\psi_{{\cal S},t\nu}^{\mu}(w,X,\ups)
-\big(t\pi_{z_m}^-\bar{\nu}_b+
{^{(3)}\al_{{\cal T};1}^{(m),-}}(w,\ups)\big)\Big\|_2
\qquad\qquad\qquad\qquad&\\
\le C(b)(t+|\vp|^{\frac{1}{p}})
\big(t+|\ups|_{h_1}^3+|\tilde{x}_1(w,\ups)||\ups|_{h_1}^2\big).&
\end{split}\end{equation*}
\end{crl}

\noindent
{\it Proof:} The first estimate is clear from Lemma~\ref{psit3_1}.
For the second, note that since $s^{(2,-)}_{b,z_m}=0$, 
$|\pi_{\tilde{x}_{\hat{1}}(w,\ups)}^- -\pi_{z_m}^-|\le 
                            C|\tilde{x}_1(w,\ups)|^2$,
and thus
\begin{gather*}
\Big|s_{b,\tilde{x}_{\hat{1}}(w,\ups)}^{(2,-)}(\tilde{v}_{h_1})-
s_{b,z_m}^{(3,-)}(\tilde{x}_{\hat{1}}(w,\ups),\tilde{v}_{h_1},
                                   \tilde{v}_{h_1})\Big|
 \le C|\tilde{x}_{\hat{1}}(w,\ups)|^2|\tilde{v}_{h_1}|^2\Lra\\
\Big|\pi_{x_{\hat{1}}(w,\ups)}^-
\tilde{\al}^{\mu}\big(w,X,\ups)-
{^{(3)}\al_{{\cal T};1}^-}(w,\ups)\Big|
\le C(b)|(t,w,X,\ups)|^{\frac{1}{p}}
\big(|x_{\hat{1}}(w,\ups)||\tilde{v}_{h_1}|^2+
|\tilde{v}_{h_1}|^3\big)
\end{gather*}
Furthermore, $|\varphi_{{\cal S},t\nu}^{\mu}(w,X,\ups)|_b\le 
 C(b)(t+|\vp|^{\frac{1}{p}})|$.\\

\noindent
The next step is to apply Lemma~\ref{top_l1} and Corollary~\ref{top_l2c}.
Let 
\begin{gather*}
F^+\!={\cal H}_{\Si}^+\otimes\ev^*T\P,\quad
F^-\!=T_{z_m}\Si\oplus F{\cal T},\quad
{\cal O}^{\pm}={\cal H}_{\Si}^{\pm}\otimes\ev^*T\P;\\
\tilde{F}^-=T_{z_m}\Si\otimes
\Big(\!\bigotimes_{i\in\hat{I},h_1\in\bar{D}_i{\cal T}}\!\! F_i{\cal T}
\Big)^{\otimes2}\oplus
\Big(\!\bigotimes_{i\in\hat{I},h_1\in\bar{D}_i{\cal T}}\!\! F_i{\cal T}
\Big)^{\otimes3};\\
\phi\big([b;w,v_{\hat{I}}]\big)
=\big[b,x_{\hat{1}(w,\ups)}\otimes\tilde{v}_{h_1}\otimes\tilde{v}_{h_1},
\tilde{v}_{h_1}\otimes\tilde{v}_{h_1}\otimes\tilde{v}_{h_1}\big];\\
\pi^+\al(w,\ups)=r_{{\cal T};1}^+(\ups),\quad
\al^+(Y)=\pi_{z_m}^+Y,\quad
\al^-\big(\phi(w,\ups)\big)\equiv {^{(3)}\al_{{\cal T},1}^{(m),-}}(w,\ups).
\end{gather*}
Note that $\al^-\!\in\!\Ga({\cal S};\tilde{F}^{-*}\otimes{\cal O}^-)$
is well-defined.
Since the map 
$$(w,X,\ups)\!\lra\! \big(Y(w,X,\ups),w,\ups\big)$$
is injective on $F^{\eset}{\cal S}$, we can view $\psi_{{\cal S},t\nu}^{\mu}$
as a map on an open subset \hbox{of $F^-\!\oplus F^+$.}

\begin{crl}
\label{order3_1_contr}
Suppose $d$ is a positive integer, ${\cal T}\!=\!(\Si,[N],I;j,\under{d})$ is
a simple bubble type, with $d_{\hat{0}}\!=\!0$ and
$\sum\limits_{i\in I}d_i\!=\!d$, and 
$\mu$ is an $N$-tuple of constraints in general position 
such that
$$\codim_{\Bbb{C}}\mu=d(n+1)-n(g-1)+N.$$
Let $\nu\!\in\!\Ga(\Si\times\P;
          \La^{0,1}\pi_{\Si}^*T^*\Si\otimes\pi_{\P}^*T\P)$
be a generic section.
If $|\hat{I}|\!>\!1$, 
for every compact subset~$K$ of~${\cal S}_{{\cal T},1}^{(m)}(\mu)$,
there exist a neighborhood $U_K$ of $K$ in $\bar{C}^{\i}_{(d;[N])}(\Si;\mu)$
and $\ep_K\!>\!0$ such that for any $t\!\in\!(0,\ep_K)$,
\hbox{$U_K\cap{\cal M}_{\Si,d,t\nu}(\mu)=\eset$.}
If $|\hat{I}|\!=\!1$,
there exists a compact subset $\tilde{K}_{{\cal T},1}^{(m)}$ 
of~${\cal S}_{{\cal T},1}(\mu)$ with the following property.
If $K$ is a compact subset of ${\cal S}_{{\cal T},1}^{(m)}(\mu)$
containing~$\tilde{K}_{{\cal T},1}^{(m)}$,
there exist a neighborhood $U_K$ of $K$ in $\bar{C}^{\i}_{(d;[N])}(\Si;\mu)$
and $\ep_K\!>\!0$ such that for all $t\!\in\!(0,\ep_K)$,
the signed cardinality of \hbox{$U_K\cap{\cal M}_{\Si,d,t\nu}(\mu)$}
equals to three times the signed number of zeros of the map
\begin{gather}
T_{z_m}\Si^{\otimes3}\otimes 
\big( L_{{\cal T},\hat{1}}^{\otimes2} 
\oplus L_{{\cal T},\hat{1}}^{\otimes3}\big)\big|
{\cal U}_{\bar{\cal T}}(\mu) \lra 
{\cal H}_{\Si}^-\otimes\ev^*T\P,\notag\\
\label{order3_1_contr_e}
[b,w,v_{\hat{1}}]\lra
\bar{\nu}_b^- 
+\big({\cal D}_{{\cal T},\hat{1}}^{(2)}b\big)
s_{b,z_m}^{(3,-)}(w)
+\big({\cal D}_{{\cal T},\hat{1}}^{(3)}b\big)s_{b,z_m}^{(3,-)}(v).
\end{gather}
\end{crl}

\noindent
{\it Proof:}
The proof is similar to the proofs of Corollaries~\ref{order1_contr}
and~\ref{order2_1_contr}, but two modifications are 
needed to be mentioned.
First, we need to show that $\al^-$ always has full rank.
Since we are assuming that $d_{h_1}\ge3$,
the sections ${\cal D}_{{\cal T},h_1}^{(1)}$,
${\cal D}_{{\cal T},h_1}^{(2)}$, and ${\cal D}_{{\cal T},h_1}^{(3)}$ 
over ${\cal M}_{\cal T}$ have transverse images in~$T\P$.
Thus, the sections of 
\hbox{$\Bbb{P}(\ev^*T\P)\!\lra\! S_{{\cal T},1}^{(m)}(\mu)$}
induced by 
${\cal D}_{{\cal T},h_1}^{(2)}$ and ${\cal D}_{{\cal T},h_1}^{(3)}$ 
are mutually transversal.
However, the fiber dimension of $\Bbb{P}(\ev^*T\P)$ is~$n\!-\!1$,
while the dimension of $S_{{\cal T},1}^{(m)}(\mu)$ \hbox{is $n-2$.}
Thus, the two sections do not intersect and $\al^-$
has full rank on all fibers over~$S_{{\cal T},1}^{(m)}(\mu)$.
The second difference with the proofs of Corollaries~\ref{order1_contr}
and~\ref{order2_1_contr} is that we replace the section 
$\psi_{{\cal S},t\nu}^{\mu}$ by the map
$$(w,\ups,X)\lra 
\pi_{z_m}^+\pi_{x_{\hat{1}}(w,\ups)}^+\psi_{{\cal S},t\nu}^{\mu}(w,\ups,X)+
\pi_{z_m}^-\pi_{x_{\hat{1}}(w,\ups)}^-\psi_{{\cal S},t\nu}^{\mu}(w,\ups,X),$$
which has exactly the same zeros provided $w$ and $\ups$ 
are sufficiently small (depending only on~$\Si$).

\subsection{Second-Order Estimate for $\psi_{{\cal T},t\nu}^{\mu}$, Case 2a}
\label{order2_case2a}

\noindent
We now understand all cases except for (2a) and (3b) of
Corollary~\ref{order1_contr_n3}.
Let $\{h_1,h_2\}\!=\!\{\hat{1},\hat{2}\}$ in Case (2a) and 
$\{\hat{2},\hat{3}\}$ in~(3b).
By dimension count as in the proof of Lemma~\ref{order1_contr_l3},
${\cal D}_{{\cal T},h_1}$ and ${\cal D}_{{\cal T},h_1}$ 
do not vanish on  ${\cal M}_{\cal T}(\mu)$ in these two cases.
By Corollary~\ref{reg_crl2}, 
$\pi_b^{\perp}\circ{\cal D}_{{\cal T},h_2}$
is transversal to zero, where $\pi_b^{\perp}$ denotes the projection onto
the orthogonal complement~$E_1$ of the image of ${\cal D}_{{\cal T},h_1}$
in $\ev^*T\P$. 
Since
$$\al_{\cal T}(\ups)=
\big({\cal D}_{{\cal T},h_1}b_{\ups}\big)s_{\Si,x_{\tilde{h}_1({\cal T})}}
      (\tilde{v}_{h_1})+
\big({\cal D}_{{\cal T},h_2}b_{\ups}\big)s_{\Si,x_{\tilde{h}_2({\cal T})}}
(\tilde{v}_{h_2}),$$
$\al_{\cal T}$ can fail to have the full rank only on the zero set of
$\pi_b^{\perp}\circ {\cal D}_{{\cal T},h_2}$.
Furthermore, $s_{\Si,x_{\tilde{h}_1}}$ and $s_{\Si,x_{\tilde{h}_2}}$
must have the same image in ${\cal H}_{\Si}^{0,1}$.
This is automatic in Case (3b), since
 \hbox{$\tilde{h}_1({\cal T})\!=\!\tilde{h}_2({\cal T})\!=\!\hat{1}$},
but in Case~(2a), this means that 
$x_{\hat{1}}$ and $x_{\hat{2}}$ differ by
the nontrivial holomorphic automorphism of~$\Si$;
\hbox{see~\cite[p254]{GH}}.\\

\noindent
We first treat Case (2a); so we can assume 
$h_1\!=\!\hat{1}$, $h_2\!=\!\hat{2}$.
Let ${\cal S}\!\equiv\!{\cal S}_{{\cal T},2}$ 
denote the subset of~${\cal M}_{\cal T}$
on which the section $\al_{\cal T}$ has rank one.
By Corollary~\ref{reg_crl2},
this is a complex submanifold of~${\cal M}_{\cal T}$.
Furthermore,
${\cal S}\!=\!{\cal S}_0\!\times\! {\cal S}_1$,
where ${\cal S}_1$ is the subspace of ${\cal U}_{\bar{\cal T}}$
on which the operator $\bar{\cal D}_{{\cal T},2}$, 
defined as in the proof of  Lemma~\ref{order1_contr_l3}, 
has rank one,
$${\cal S}_0=\big\{(x_{\hat{1}},-x_{\hat{1}})\!: 
x_{\hat{1}}\!\in\!\Si^*\big\},$$
$-x_{\hat{1}}\!\in\!\Si$ denotes the image of $x_{\hat{1}}$ under
the nontrivial automorphism of~$\Si$, and $\Si^*$ is the subset
of $\Si$ which is not fixed by this automorphism,
i.e. the complement of the points $z_1,\ldots,z_6$
described in Subsection~\ref{order3_case1}.
By Corollary~\ref{reg_crl2},
${\cal S}_1$ is a complex submanifold of ${\cal U}_{\bar{\cal T}}$.
The normal bundle of ${\cal S}$ in ${\cal M}_{\cal T}$~is 
$${\cal NS}={\cal NS}_0\oplus{\cal NS}_1,
\quad\hbox{where}\quad
{\cal NS}_0=\pi_{\Si,\hat{2}}^*T\Si,~~
{\cal NS}_1=L_{{\cal T},\hat{2}}^*\otimes E_1,$$
and $\pi_{\Si,h}\!:{\cal S}_0\!\subset\!\Si\!\times\!\Si\!\lra\!\Si$ 
is the projection on the $h$th component.
Let $(\Phi_{\cal S},\Phi_{\cal S}^{\mu})$
be a regularization  of 
${\cal S}_{{\cal T},2}(\mu)
\!\equiv\!{\cal S}\cap{\cal M}_{\cal T}(\mu)$.
This regularization can be chosen so that
\begin{equation}\label{order22_e}
\pi_{\phi_{\cal S}(b,X)}^{\perp}
{\cal D}_{{\cal T},\hat{2}}\tilde{\phi}_{\cal S}(b,X)=
\Pi_{b,\tilde{\phi}_{\cal S}(b,X)}X
\quad\forall (b,X)\!\in\! {\cal N}\tilde{\cal S}_1=E_1,
\end{equation}
where $\tilde{\phi}_{\cal S}$ is the lift of $\phi_{\cal S}$
to ${\cal M}_{\cal T}^{(0)}$.
We also assume that $\Phi_{\cal S}^{\mu}$ is given
by the $g_{\P,b}$-parallel transport on~${\cal N}_b{\cal S}_1$.
Since the section $s$ is invariant under the automorphism group of~$\Si$,
we identify $\pi_{\Si,\hat{2}}^*T\Si|{\cal S}_0$ with
$\pi_{\Si,\hat{1}}^*T\Si|{\cal S}_0$.
If $(b;w)\!\in\!{\cal NS}_0$ is sufficiently small, let 
$$x_{\hat{2}}(w)=\exp_{b,x_{\hat{2}}}w.$$
The bundle ${\cal NS}$ carries a natural norm induced
by the $g_{\P,\ev}$-metric on~$\P$ and $g_{\cdot,\hat{0}}$-metric
on~$\Si$.
Denote by $F{\cal S}$ and $F^{\eset}{\cal S}$
the bundles described in Subsection~\ref{si_str_sec}
corresponding to the submanifold~${\cal S}_{{\cal T},2}$.
If $(w,X,\ups)\!\in\! F{\cal S}\!=\!{\cal NS}\oplus F{\cal T}$, put
\begin{equation*}\begin{split}
\tilde{\al}(w,X,\ups)=
\Pi_{b,\phi_{\cal S}(b,X)}^{-1}
\Big( \big({\cal D}_{{\cal T},\hat{1}}\phi_{\cal S}(b,X)\big)
 s_{\Si,x_{\hat{1}}}(v_{\hat{1}})+
\big({\cal D}_{{\cal T},\hat{2}}\phi_{\cal S}(b,X)\big)
 s_{\Si,x_{\hat{2}}(w)}(v_{\hat{2}})\Big)\qquad\qquad&\\
+\Big(\big({\cal D}_{{\cal T},\hat{1}}^{(2)}b\big) 
      s_{b,x_{\hat{1}}}^{(2)}(v_{\hat{1}})+
\big({\cal D}_{{\cal T},\hat{2}}^{(2)}b\big) 
      s_{b,x_{\hat{1}}}^{(2)}(v_{\hat{2}})\Big)&.
\end{split}\end{equation*}
If $(w,X,\ups)\!\in\! F^{\eset}{\cal S}|{\cal S}_{{\cal T},2}(\mu)$ 
is sufficiently small, let
\begin{equation*}\begin{split}
\tilde{\al}^{\mu}(w,X,\ups)=
\Big( \big({\cal D}_{{\cal S},t\nu,\hat{1}}^{\mu}(w,X,\ups)\big)
s_{\Si,x_{\hat{1}}}(v_{\hat{1}})+
\big({\cal D}_{{\cal S},t\nu,\hat{2}}^{\mu}(w,X,\ups)\big)
s_{\Si,x_{\hat{2}}(w)}(v_{\hat{2}})\Big)\quad\qquad&\\
+\Big(\big({\cal D}_{{\cal T},\hat{1}}^{(2)}b\big) 
      s_{b,x_{\hat{1}}}^{(2)}(v_{\hat{1}})+
\big({\cal D}_{{\cal T},\hat{2}}^{(2)}b\big) 
      s_{b,x_{\hat{1}}}^{(2)}(v_{\hat{2}})\Big)&,
\end{split}\end{equation*}
where, with $\varphi_{{\cal S},t\nu}^{\mu}$ as in Theorem~\ref{si_str},
$${\cal D}_{{\cal S},t\nu,h}^{\mu}(w,X,\ups)=
\Pi_{\phi_{\cal S}^{\mu}\varphi_{{\cal S},t\nu}^{\mu}(w,X,\ups),
\phi_{\cal S}\Phi_{\cal S}^{\mu}\varphi_{{\cal S},t\nu}^{\mu}(w,X,\ups)}^{-1}
\Pi_{b,\phi_{\cal S}^{\mu}\varphi_{{\cal S},t\nu}^{\mu}(w,X,\ups)}^{-1}
\big({\cal D}_{{\cal T},h}\phi_{\cal S}\Phi_{\cal S}^{\mu}
     \varphi_{{\cal S},t\nu}^{\mu}(w,X,\ups)\big).$$

\begin{lmm}
\label{dbar_2a2}
There exist $\de,C\!\in\! C^{\i}({\cal S};\Bbb{R}^+)$
such that for all 
$\vp\!=\![(b,w,X,\ups)]\!\in\! F^{\eset}{\cal S}_{\de}$,
\begin{equation*}
\Big\|  \pi^{0,1}_{\Phi_{\cal S}(\vp),-}
           \bar{\partial}u_{\Phi_{\cal S}(\vp)}+
\tilde{R}_{\Phi_{\cal S}(\vp)}
\Pi_{b,\phi_{\cal S}(X)}\tilde{\al}(w,X,\ups)\Big\|_2
\le C(b)|\vp||\ups|^2.
\end{equation*}
\end{lmm}

\noindent
{\it Proof:} The proof is analogous to the proof of Lemma~\ref{dbar_3};
here we use Proposition~\ref{dbar_gen2prp} with 
two terms for $h\!=\!\hat{1}$ and two terms for $h\!=\!\hat{2}$.\\

\begin{lmm}
\label{psi_2a2}
There exist $\de,C\!\in\! C^{\i}({\cal S}_{{\cal T},2}(\mu);\Bbb{R}^+)$
such that for all
$\vp\!=\![(b,w,X,\ups)]\!\in\! F^{\eset}{\cal S}_{\de}$,
\begin{equation*}
\Big\| \psi_{{\cal S},t\nu}^{\mu} (\vp)-
\big(t\bar{\nu}_{\ev(b)}+\tilde{\al}^{\mu}(w,X,\ups)\big)
\Big\|_2
\le C(b)(t+|\vp|^{\frac{1}{p}})
\big(t+|\ups|^2+|w||v_{\hat{2}}|\big).
\end{equation*} 
\end{lmm}

\noindent
{\it Proof:} As in the proof of Lemmas~\ref{psit_21} and~\ref{psit3_1},
we need to obtain an appropriate estimate on
$$\big\|D_{\Phi_{\cal S}(\vp)}^*R_{\Phi_{\cal S}(\vp)}X_i\psi_2
\big\|_{L^1},$$
where $\psi_2$ is a $(0,1)$-form vanishing at $x_{\hat{1}}$ and 
with norm~$1$.
From equation~\e_ref{coker_e3}, we see that the $L^1$-norm
over the small annulus centered at $x_{\hat{1}}$ is bounded by 
$C(b)|v_{\hat{1}}|^2$; see also the proof of Lemma~\ref{psit_21}.
Furthermore, since $x_{\hat{2}}$ is ``dual'' to $x_{\hat{1}}$, 
$\psi_2$ also vanishes at~$x_{\hat{2}}$.
Thus, the $L^1$-norm over the small annulus centered at~$x_{\hat{2}}(w)$
is bounded by $C(b)(|w|+|v_{\hat{2}}|)|v_{\hat{2}}|$ as can be seen from
equation~\e_ref{coker_e3}.\\

\noindent 
Let $\tilde{s}_{b,m}^{(2,+)}\!\in\! T_{z_m}^*\Si$ be given by
$s_{b,m}^{(2,+)}(v,v)\!=\!\tilde{s}_{b,m}^{(2,+)}(v)s_{\Si,z_m}(v)$.
For any $b\!\in\! {\cal S}_{{\cal T},2}(\mu)$, define
\begin{gather*}
\ka(b)\!\in\!(L_{{\cal T},\hat{2}}^*\otimes L_{{\cal T},\hat{1}}
      \!-\!\{0\}) 
\quad\hbox{and}\quad 
\mu(b)\!\in\! L_{{\cal T},\hat{2}}^*\otimes L_{{\cal T},\hat{1}}
\quad\hbox{by}\\
\big({\cal D}_{{\cal T},\hat{2}}b\big)=
\ka(b)\big({\cal D}_{{\cal T},\hat{1}}b\big),\quad
\pi_b\big({\cal D}_{{\cal T},\hat{1}}^{(2)} b\big)=
\mu(b)\big({\cal D}_{{\cal T},\hat{1}}^{(1)} b\big),
\end{gather*}
where $\pi_b\!:\ev^*T\P\!\lra\!\hbox{Im}({\cal D}_{{\cal T},\hat{1}})$
is the orthogonal projection~map.
If $(w,X,\ups)\!\in\! F^{\eset}{\cal S}|{\cal S}_{{\cal T},2}(\mu)$ 
is sufficiently small, let
\hbox{$\tilde{\ka}(w,X,\ups)\!\in\!\Bbb{C}^*$} be given by
$$\pi_{\phi_{\cal S}\Phi_{\cal S}^{\mu}
     \varphi_{{\cal S},t\nu}^{\mu}(w,X,\ups)}
\big({\cal D}_{{\cal T},\hat{2}}\phi_{\cal S}\Phi_{\cal S}^{\mu}
     \varphi_{{\cal S},t\nu}^{\mu}(w,X,\ups)\big)=
\tilde{\ka}(w,X,\ups)
\big({\cal D}_{{\cal T},\hat{1}}\phi_{\cal S}\Phi_{\cal S}^{\mu}
     \varphi_{{\cal S},t\nu}^{\mu}(w,X,\ups)\big).$$
Note that by Theorem~\ref{si_str}, 
$|\tilde{\ka}(w,X,\ups)-\ka(b)|\le 
C(b)(t+|\vp|^{\frac{1}{p}})$.
Let
\begin{gather*}
Y^t(w,X,\ups)\!=\!
\big({\cal D}_{{\cal S},t\nu,\hat{1}}^{\mu}(w,X,\ups)\big)
s_{\Si,x_{\hat{1}}}\big(v_{\hat{1}}\!+\!\tilde{\ka}(w,X,\ups)v_{\hat{2}}
\!+\!\mu(b)\tilde{s}^{(2,+)}_{\Si,x_{\hat{1}}}(v_{\hat{1}})v_{\hat{1}}\big),~~
Y^{\perp}\big(X,v_{\hat{2}}\big)\!=\!
Xs_{\Si,x_{\hat{1}}}\big(v_{\hat{2}}\big);\\
{^{(2)}\al_{{\cal T};2}^-}\big(w,v_{\hat{2}}\big)
=\big({\cal D}_{{\cal T},\hat{1}}b\big)
        s_{b,x_{\hat{1}}}^{(2,-)}\big(w,v_{\hat{2}}\big)
+\big({\cal D}_{{\cal T},\hat{1}}^{(2)}b\big)
  s_{b,x_{\hat{1}}}^{(2,-)}\big(\ka(b)v_{\hat{2}}\big)
+\big({\cal D}_{{\cal T},\hat{2}}^{(2)}b\big)
       s_{b,x_{\hat{1}}}^{(2,-)}\big(v_{\hat{2}}\big);\\
r_{{\cal T};2}^+(w,\ups)=
\big({\cal D}_{{\cal T},\hat{1}}^{(1)}(b)\big)
s_{b,x_{\hat{1}}}^{(2,+)}\big(w,v_{\hat{2}}\big)+
\pi_b^{\perp}\big({\cal D}_{{\cal T},\hat{1}}^{(2)}(b)\big)
s_{b,x_{\hat{1}}}^{(2,+)}\big(\ka(b)v_{\hat{2}}\big)+
\big({\cal D}_{{\cal T},\hat{2}}^{(2)}(b)\big)
s_{b,x_{\hat{1}}}^{(2,+)}\big(v_{\hat{2}}\big).
\end{gather*}
Let $Y=Y^t+Y^{\perp}$ and 
$\bar{\nu}_b^{\pm}=\pi_{x_{\hat{1}}}^{\pm}\bar{\nu}_b$.

\begin{crl}
\label{psit_2a2c}
There exist $\de,C\!\in\! C^{\i}({\cal S}_{{\cal T},2}(\mu);\Bbb{R}^+)$
such that for all
$\vp\!=\![(b,w,X,\ups)]\!\in\! F^{\eset}{\cal S}_{\de}$,
\begin{equation*}\begin{split}
\Big\| \pi_{x_{\hat{1}}}^+\psi_{{\cal S},t\nu}^{\mu}(\vp)\!- \!
 \big( t\bar{\nu}_b^+\!+\! Y(w,X,\ups)\!+\!r_{{\cal T};2}^+(w,\ups)\big)
\Big\|_2
&\!\le\! C(b)(t+|\vp|^{\frac{1}{p}})
\big(|\ups|^{2}\!+\!|w||v_{\hat{2}}|\!+\!|Y|\big);\\
\Big\|\pi_{x_{\hat{1}}}^-\psi_{{\cal S},t\nu}^{\mu}(\vp)\!-\!
\big( t\bar{\nu}_b^- \!+\!
{^{(2)}\al_{{\cal T};2}^-}(w,v_{\hat{2}})\big)\Big\|_2
&\!\le\! C(b)(t+|\vp|^{\frac{1}{p}})
\big(|\ups|^{2}\!+|w||v_{\hat{2}}|\!+|Y|\big).\\
\end{split}\end{equation*} 
\end{crl}

\noindent
{\it Proof:} The proof is similar to that of Corollary~\ref{psit3_1c}, 
but we use 
$$\Big| s_{\Si,x_{\hat{2}}(w)}(v_{\hat{2}})-
 \big(s_{\Si,x_{\hat{1}}}(v_{\hat{2}})+
            s_{b,x_{\hat{1}}}^{(2)}(w,v_{\hat{2}})\big)\Big|
\le C(b)|w|^2|v_{\hat{2}}|.$$
We also use $\big|\big({\cal D}_{{\cal S},t\nu,\hat{1}}^{\mu}
 (w,X,\ups)\big)\big|
\ge C(b)^{-1}$.\\

\noindent
The next step is to apply Corollary~\ref{top_l2c}.  Let 
\begin{gather*}
F^+\!={\cal H}_{\Si}^+\otimes\ev^*T\P,~~
F^-\!\!=\pi_{\Si,\hat{1}}^* T\Si\oplus F_{\hat{2}}{\cal T},~~
{\cal O}^{\pm}={\cal H}_{\Si}^{\pm}\otimes\ev^*T\P,~~
\tilde{F}^-\!\!= \pi_{\Si,\hat{1}}^* T\Si \otimes F_{\hat{2}}{\cal T}
\oplus F_{\hat{2}}{\cal T}^{\otimes2};\\
\phi\big([b;w,v_{\hat{2}}]\big)
=\big[b,w\otimes v_{\hat{2}},
v_{\hat{2}}\otimes v_{\hat{2}}\big],~~
\al^-\big(\phi(w,v_{\hat{2}})\big)\equiv 
{^{(2)}\al_{{\cal T},2}^-}(w,v_{\hat{2}}),~~
\pi^+r(w,\ups)=r_{{\cal T};1}^+(w,\ups).
\end{gather*}
Note that $\al^-\!\in\!\Ga({\cal S};\tilde{F}^{-*}\otimes{\cal O}^-)$
is well-defined.
Since the map 
$$(w,X,\ups)\lra \big(Y(w,X,\ups),w,v_{\hat{2}}\big)$$
is injective on $F^{\eset}{\cal S}$ as long as 
$\de\!\in\! C^{\i}({\cal S}_{{\cal T},2}(\mu);\Bbb{R}^+)$
is sufficiently small, we can view $\psi_{{\cal S},t\nu}^{\mu}$
as a map on an open subset \hbox{of $F^-\oplus F^+$.}

\begin{crl}
\label{order2_2a_contr}
Suppose $d$ is a positive integer, 
${\cal T}\!=\!(\Si,[N],I;j,\under{d})$ is
a simple bubble type, with 
$\hat{I}\!=\!\{\hat{1},\hat{2}\}$, $M_{\hat{0}}{\cal T}\!=\!\eset$,
$d_{\hat{0}}\!=\!0$, and
$\sum\limits_{i\in I}d_i\!=\!d$, and 
$\mu$ is an $N$-tuple of constraints in general position 
such that
$$\codim_{\Bbb{C}}\mu=d(n+1)-n(g-1)+N.$$
Let $\nu\!\in\!\Ga(\Si\times\P;
          \La^{0,1}\pi_{\Si}^*T^*\Si\otimes\pi_{\P}^*T\P)$
be a generic section.
Then there exists a compact subset~$\tilde{K}_{{\cal T},2}$ 
of~${\cal S}_{{\cal T},2}(\mu)$ with the following property.
If $K$ is a compact subset  of ${\cal S}_{{\cal T},2}(\mu)$
containing~$\tilde{K}_{{\cal T},1}$,
there exist a neighborhood $U_K$ of $K$ in $\bar{C}^{\i}_{(d;[N])}(\Si;\mu)$
and $\ep_K\!>\!0$ such that for any $t\!\in\!(0,\ep_K)$,
the signed cardinality of $U_K\cap{\cal M}_{\Si,d,t\nu}(\mu)$ 
equals to twice the signed number of zeros of the map
\begin{gather}
\pi_{\Si}^*T\Si^{\otimes2}\otimes
\big( L_{\bar{\cal T},\hat{2}}^{\otimes2}\oplus
L_{\bar{\cal T},\hat{2}}^{\otimes2}\big)\big|
\Si^*\otimes{\cal U}_{\bar{\cal T}}(\mu)\lra 
 {\cal H}_{\Si}^-\otimes\ev^*T\P\notag,\\
\label{order2_2a_contr_e}
\big[(x,b);(w,v)\big]\lra
\bar{\nu}_b^-
+\big({\cal D}_{\bar{\cal T},\hat{2}}b\big) s_x^{(2,-)}(w,v)
+\big({\cal D}_{\bar{\cal T},\hat{1}}^{(2)}b\big) s_x^{(2,-)}\big(\ka(b)v\big)
+\big({\cal D}_{\bar{\cal T},\hat{2}}^{(2)}b\big)s_x^{(2,-)}(v).
\end{gather}
\end{crl}

\noindent
{\it Proof:} The proof is similar to that of 
Corollary~\ref{order2_1_contr}.
We only need to see that the section~$\al^-$ defined above
has rank~two.
If $d_{\hat{1}}\!=\!d_{\hat{2}}\!=\!1$, 
the space ${\cal S}_{{\cal T},2}(\mu)\!=\!\emptyset$,
since any two tangent lines in $\P$ agree, and no line passes
through all of the constraints $\mu_1,\ldots,\mu_N$ \hbox{if $n\!=\!3$}.
Thus, it can be assumed \hbox{that $d_{\hat{1}}\ge 2$}.
Note that $S_{{\cal T},2}(\mu)$ is one-dimensional, 
with the only dimension coming from the singular 
\hbox{point $x_{\hat{1}}\!\in\!\Si$}.
Thus, by Corollary~\ref{reg_crl2}, if the constraints 
$\mu_1,\ldots,\mu_N$ are in general position, the image of 
${\cal D}_{{\cal T};\hat{1}}^{(2)}$ does not lie in the linear
span of $({\cal D}_{{\cal T},\hat{2}}b)$ 
and~$({\cal D}_{{\cal T};\hat{2}}^{(2)}b)$.
\hbox{Furthermore, $({\cal D}_{{\cal T},\hat{2}}b)\!\neq\!0$}.

\subsection{Second-Order Estimate for $\psi_{{\cal T},t\nu}^{\mu}$, Case 2b}
\label{order2_case2b}

\noindent
We now treat Case (3b) of Corollary~\ref{order1_contr_n3};
we can assume $h_1\!=\!\hat{2}$, $h_2\!=\!\hat{3}$.
Let ${\cal S}\!\equiv\!{\cal S}_{{\cal T},2}$ denote 
the subset of ${\cal M}_{\cal T}$ on which the operator
$\bar{\cal D}_{{\cal T},2}$ of Lemma~\ref{order1_contr_l3} has rank~one.
Similarly to the case of Subsection~\ref{order2_case2a},
${\cal S}$ is a regular submanifold of ${\cal M}_{\cal T}$
with normal \hbox{bundle ${\cal NS}\!=\!L_{{\cal T},\hat{3}}^*\otimes E_1$}. 
As before, we can choose a regularization 
$\big(\Phi_{\cal S},\Phi_{\cal S}^{\mu}\big)$
of \hbox{${\cal S}_{{\cal T},2}(\mu)\!
\equiv\!{\cal S}\cap{\cal M}_{\cal T}(\mu)$} such that
\begin{equation}
\pi_{\phi_{\cal S}(b,X)}^{\perp}
{\cal D}_{{\cal T},\hat{3}}\tilde{\phi}_{\cal S}(b,X)=
\Pi_{b,\tilde{\phi}_{\cal S}(b,X)}X
\quad\forall (b,X)\in {\cal N}\tilde{\cal S}_1=E_1,
\end{equation}
where $\tilde{\phi}_{\cal S}$ is the lift of $\phi_{\cal S}$
to ${\cal M}_{\cal T}^{(0)}$, and $\Phi_{\cal S}^{\mu}$ is given
by the $g_{\P,b}$-parallel transport on~${\cal N}_b{\cal S}$.
Denote by $F{\cal S}$ and $F^{\eset}{\cal S}$
the bundles described in Subsection~\ref{si_str_sec}
corresponding to the submanifold~${\cal S}_{{\cal T},2}$.
If $(X,\ups)\in F{\cal S}={\cal NS}\oplus F{\cal T}$ is sufficiently
small, let
\begin{gather*}
\tilde{\al}(X,\ups)=
\big({\cal D}_{{\cal T},\hat{2}}\phi_{\cal S}(b,X)\big)
 s_{\Si,\tilde{x}_{\hat{2}}(\ups)}(\tilde{v}_{\hat{2}})+
\big({\cal D}_{{\cal T},\hat{3}}\phi_{\cal S}(b,X)\big)
 s_{\Si,\tilde{x}_{\hat{3}}(\ups)}(\tilde{v}_{\hat{3}});\\
\tilde{\al}^{\mu}(X,\ups)=
\big({\cal D}_{{\cal S},t\nu,\hat{2}}^{\mu}(X,\ups)\big)
s_{\Si,\tilde{x}_{\hat{2}}(\ups)}(\tilde{v}_{\hat{2}})+
\big({\cal D}_{{\cal S},t\nu,\hat{3}}^{\mu}(X,\ups)\big)
s_{\Si,x_{\hat{3}}(\ups)}(v_{\hat{3}}),
\end{gather*}
where, with $\varphi_{{\cal S},t\nu}^{\mu}$ as in Theorem~\ref{si_str},
$${\cal D}_{{\cal S},t\nu,h}^{\mu}(X,\ups)=
\Pi_{\phi_{\cal S}^{\mu}\varphi_{{\cal S},t\nu}^{\mu}(X,\ups),
\phi_{\cal S}\Phi_{\cal S}^{\mu}\varphi_{{\cal S},t\nu}^{\mu}(X,\ups)}^{-1}
\Pi_{b,\phi_{\cal S}^{\mu}\varphi_{{\cal S},t\nu}^{\mu}(X,\ups)}^{-1}
\big({\cal D}_{{\cal T},h}\phi_{\cal S}\Phi_{\cal S}^{\mu}
     \varphi_{{\cal S},t\nu}^{\mu}(X,\ups)\big).$$

\begin{lmm}
\label{dbar_22}
There exist $\de,C\!\in\! C^{\i}({\cal S}_{{\cal T},2};\Bbb{R}^+)$
such that for all 
$\vp\!=\![(b,X,\ups)]\!\in\! F^{\eset}{\cal S}_{\de}$,
$$\Big\| \pi^{0,1}_{\Phi_{\cal S}(\vp),-}
 \bar{\partial}u_{\Phi_{\cal S}(\vp)}
+\tilde{R}_{\Phi_{\cal S}(\vp)}\tilde{\al}(X,\ups)\Big\|_2
\le C(b)\big(|\tilde{v}_{\hat{2}}|^2+|\tilde{v}_{\hat{3}}|^2\big).$$
\end{lmm}

\noindent
{\it Proof:} This lemma is immediate from Proposition~\ref{dbar_gen2prp} 
applied with one term for each $h\!=\!\hat{2},\hat{3}$.\\

\begin{lmm}
\label{psit_22}
There exist $\de,C\!\in\! C^{\i}({\cal S}_{{\cal T},2}(\mu);\Bbb{R}^+)$
such that for all 
$\vp\!=\![(b,X,\ups)]\!\in\! F^{\eset}{\cal S}_{\de}$,
$$\Big\|\psi_{{\cal S},t\nu}^{\mu}(\vp)
-\big(t\bar{\nu}_b+\tilde{\al}^{\mu}(X,\ups)\big)\Big\|_2
\le C(b)(t+|\vp|^{\frac{1}{p}})
\big(t+|v_{\hat{1}}|(|\tilde{v}_{\hat{2}}|+|\tilde{v}_{\hat{3}}|)\big).$$
\end{lmm}

\noindent
{\it Proof:} As usually, we only need to obtain a good bound on
$$\big\|D_{\Phi_{\cal S}(\vp)}^*R_{\Phi_{\cal S}(\vp)}
 X_i\psi_2\|_{L^1},$$
where the notation is as in the proof of Lemma~\ref{psi_2a2}.
By equation~\e_ref{coker_e3}, the $L^1$-norm on the small annulus 
centered at $\tilde{x}_{\hat{2}}(\ups)$ is bounded 
by~$|\tilde{v}_{\hat{2}}|^2$.
Since $g_{b,\hat{0}}$-distance between $\tilde{x}_{\hat{2}}(\ups)$ and 
$\tilde{x}_{\hat{3}}(\ups)$ is bounded by~$C(b)|v_{\hat{1}}|$,
the $L_1$-norm over the annulus centered at
$\tilde{x}_{\hat{3}}(\ups)$ is bounded 
by~$|v_{\hat{1}}||\tilde{v}_{\hat{3}}|$.\\

\noindent
For any $b\!\in\!{\cal S}_{{\cal T},2}(\mu)$, let 
$\ka(b)\!\in\!(L_{{\cal T},\hat{3}}^*\otimes L_{{\cal T},\hat{2}}-\{0\})$
be given by $\big({\cal D}_{{\cal T},\hat{3}}b\big)\!=\!
\ka(b)\big({\cal D}_{{\cal T},\hat{2}}b\big)$.
Define \hbox{$\tilde{\ka}(X,\ups)\!\in\!
(L_{{\cal T},\hat{3}}^*\otimes L_{{\cal T},\hat{2}}-\{0\})$}
for $(X,\ups)\!\in\! F^{\eset}{\cal S}|{\cal S}_{{\cal T},2}(\mu)$ 
sufficiently small by
$$\pi_{\phi_{\cal S}\Phi_{\cal S}^{\mu}
     \varphi_{{\cal S},t\nu}^{\mu}(X,\ups)}
\big({\cal D}_{{\cal T},\hat{3}}\phi_{\cal S}\Phi_{\cal S}^{\mu}
     \varphi_{{\cal S},t\nu}^{\mu}(X,\ups)\big)=
\tilde{\ka}(X,\ups)
\big({\cal D}_{{\cal T},\hat{2}}\phi_{\cal S}\Phi_{\cal S}^{\mu}
     \varphi_{{\cal S},t\nu}^{\mu}(X,\ups)\big).$$
Note that by Theorem~\ref{si_str}, 
$|\tilde{\ka}(X,\ups)-\ka(b)|\le C(b)(t+|\vp|^{\frac{1}{p}})$.
Let
\begin{gather*}
Y^t(X,\ups)\!=\!
\big({\cal D}_{{\cal S},t\nu,\hat{2}}^{\mu}(X,\ups)\big)
\Big( s_{\Si,x_{\hat{1}}}\big(\tilde{v}_{\hat{2}}+
                    \tilde{\ka}(X,\ups)\tilde{v}_{\hat{3}}\big)
+s_{b,x_{\hat{1}}}^{(2,+)}
\big(v_{\hat{1}},x_{\hat{2}}\tilde{v}_{\hat{2}}+
x_{\hat{3}}\tilde{\ka}(b)\tilde{v}_{\hat{3}})\Big);\\
Y^{\perp}(X,\ups)\!=\!
Xs_{\Si,x_{\hat{1}}}\big(\tilde{v}_{\hat{3}}\big),\quad
{^{(2)}\al_{{\cal T};2}^-}(\ups)
=\big({\cal D}_{{\cal T},\hat{2}}b\big)s_{b,x_{\hat{1}}}^{(2,-)}
\big(v_{\hat{1}},x_{\hat{2}}\tilde{v}_{\hat{2}}+
           x_{\hat{3}}\ka(b)\tilde{v}_{\hat{3}}\big).
\end{gather*}
Let $Y=Y^t+Y^{\perp}$ and
$\bar{\nu}_b^{\pm}=\pi_{x_{\hat{1}}}^{\pm}\bar{\nu}_b$.

\begin{crl}
\label{psit2c_2b}
There exist $\de,C\!\in\! C^{\i}({\cal S}_{{\cal T},2}(\mu);\Bbb{R}^+)$
such that for all
$\vp\!=\![(b,X,\ups)]\!\in\! F^{\eset}{\cal S}_{\de}$,
\begin{gather*}
\Big\|\pi_{x_{\hat{1}}}\psi_{{\cal S},t\nu}^{\mu}(X,\vp)-
\big(t\bar{\nu}_b^+ +Y(X,\ups) \big)\Big\|_2
\le C(b)(t+|\vp|^{\frac{1}{p}})\big(
t+|v_{\hat{1}}|(|\tilde{v}_{\hat{2}}|+|\tilde{v}_{\hat{3}}|)+
 |Y^{\perp}(X,\ups)|
\big);\\
\Big\|\pi_{x_{\hat{1}}}^-\psi_{{\cal S},t\nu}^{\mu}(\vp)-
\big(t\bar{\nu}_b^- +{^{(2)}\al_{{\cal T};2}^-}(\ups)\big)\Big\|_2
\le C(b)(t+|\vp|^{\frac{1}{p}})
\big(t+|v_{\hat{1}}|(|\tilde{v}_{\hat{2}}|+|\tilde{v}_{\hat{3}}|)\big).
\end{gather*} 
\end{crl}

\noindent
{\it Proof:} This claim is proved similarly to Corollary~\ref{psit_2a2c}.\\

\noindent
The next step is to apply Lemma~\ref{top_l1}.  Let 
\begin{gather*}
F^+\!={\cal H}_{\Si}^+\otimes E_1,\quad
F^-\!\!=F{\cal T},~~
{\cal O}^{\pm}={\cal H}_{\Si}^{\pm}\otimes\ev^*T\P,~~
\tilde{F}^-\!\!= 
\pi_{\Si}^* T\Si \otimes F_{\hat{2}}{\cal T};\\
\phi\big([b;\ups]\big)
=\big[b,v_{\hat{1}}\otimes 
  (x_{\hat{2}}\tilde{v}_{\hat{2}}+x_{\hat{3}}\ka(b)\tilde{v}_{\hat{3}})\big],
\quad
\al^-\big(\phi(\ups)\big)\equiv 
{^{(2)}\al_{{\cal T},2}^-}(\ups),\quad
\al(X,\ups)=Y(X,\ups)+{^{(2)}\al_{{\cal T},2}^-}(\ups).
\end{gather*}
Note that $\al^-\!\in\!\Ga({\cal S};\tilde{F}^{-*}\otimes{\cal O}^-)$
is well-defined.
Since the map 
$$(X,\ups)\lra \big(Y^{\perp}(X,\ups),\ups\big)$$
is injective on $F^{\eset}{\cal S}$, 
we can view $\psi_{{\cal S},t\nu}^{\mu}$
as a map on an open subset \hbox{of $F^+\oplus F^-$.}

\begin{crl}
\label{order2_2b_contr}
Suppose $d$ is a positive integer, 
${\cal T}\!=\!(\Si,[N],I;j,\under{d})$ is a simple bubble type, with 
$\hat{I}\!=\!\{\hat{1},\hat{2},\hat{3}\}$, 
$H_{\hat{1}}{\cal T}\!=\!\{\hat{2},\hat{3}\}$,
$d_{\hat{0}}\!=\!0$, and
$\sum\limits_{i\in I}d_i\!=\!d$, and 
$\mu$ is an $N$-tuple of constraints in general position 
such that
$$\codim_{\Bbb{C}}\mu=d(n+1)-n(g-1)+N.$$
Let $\nu\!\in\!\Ga(\Si\times\P;
          \La^{0,1}\pi_{\Si}^*T^*\Si\otimes\pi_{\P}^*T\P)$
be a generic section.
For every compact subset~$K$ of ${\cal S}_{{\cal T},2}(\mu)$,
such that $x_{\hat{1}}(b)\in\Si^*$ for all $b\in K$,
there exist a neighborhood $U_K$ of $K$ in $\bar{C}^{\i}_{(d;[N])}(\Si;\mu)$,
where $d=\sum d_h$, and $\ep_K\!>\!0$ such that for any $t\!\in\!(0,\ep_K)$,
\hbox{$U_K\cap{\cal M}_{\Si,d,t\nu}(\mu)\!=\!\eset$.}
\end{crl}

\noindent
{\it Proof:} The set  
${\cal S}_{{\cal T},2}^*(\mu)\!\equiv\!
\{b\!\in\!{\cal S}_{{\cal T},2}(\mu)\!:x_{\hat{1}}\!\in\!\Si^*\}$ 
is an open subset of ${\cal S}_{{\cal T},2}(\mu)$
on which the section~$\al^-$ has full rank, since 
${\cal D}_{{\cal T},\hat{2}}$ does not vanish 
on~${\cal S}_{{\cal T},2}(\mu)$.
Note that the dimension of ${\cal S}_{{\cal T},2}(\mu)$ is~$1$,
the rank of $\tilde{F}^-$ is also~$1$, while
the rank ${\cal O}^-$ is~3.
Thus, the claim follows from Theorem~\ref{si_str},
Lemma~\ref{top_l1}, and Corollary~\ref{psit2c_2b}, provided
$$|v_{\hat{1}}|\big(|\tilde{v}_{\hat{2}}|+|\tilde{v}_{\hat{3}}|\big)
\le C(b)\big(
|v_{\hat{1}}| 
 |x_{\hat{2}}\tilde{v}_{\hat{2}}+x_{\hat{3}}\ka(b)\tilde{v}_{\hat{3}}|
+|Y^t(X,\ups)|\big)$$
for some $C\!\in\! C^{\i}({\cal S}_{{\cal T},2}^*(\mu);\Bbb{R}^+)$.
By definition of $Y^t(X,\ups)$,
$$|\tilde{v}_{\hat{2}}+\ka(b)\tilde{v}_{\hat{3}}|\le
|Y^t(X,\ups)|+C(b)
|x_{\hat{2}}\tilde{v}_{\hat{2}}+x_{\hat{3}}\ka(b)\tilde{v}_{\hat{3}}|.$$
Since $x_{\hat{2}}\!\neq\! x_{\hat{3}}$, 
\begin{equation*}\begin{split}
|v_{\hat{1}}|\big(|\tilde{v}_{\hat{2}}|+|\tilde{v}_{\hat{3}}|\big)
&\le C(b)|v_{\hat{1}}|\big(
|\tilde{v}_{\hat{2}}+\ka(b)\tilde{v}_{\hat{3}}|+
|x_{\hat{2}}\tilde{v}_{\hat{2}}+x_{\hat{3}}\ka(b)\tilde{v}_{\hat{3}}|\big)\\
&\le C'(b)|v_{\hat{1}}|\big(
 |x_{\hat{2}}\tilde{v}_{\hat{2}}+x_{\hat{3}}\ka(b)\tilde{v}_{\hat{3}}|
+|Y^t(X,\ups)|\big).
\end{split}\end{equation*}

\subsection{Third-order Estimate for $\psi_{{\cal T},t\nu}^{\mu}$, Case 2b}
\label{order3_case2b}

\noindent
It remains to consider gluing along the subset
${\cal S}_{{\cal T},2}^{(m)}(\mu)$ of ${\cal S}_{{\cal T},2}(\mu)$
consisting of bubble maps $b$ such that $x_{\hat{1}}(b)\!=\!z_m$,
one of the six distinguished points of~$\Si$.
Let 
$${\cal S}={\cal S}_{{\cal T},2}^{(m)}=\{b\!\in\!{\cal S}_{{\cal T},2}\!:
x_{\hat{1}}(b)\!=\!z_m\}.$$
The normal bundle of ${\cal S}_{{\cal T},2}^{(m)}$ in 
${\cal M}_{\cal T}$
\hbox{is ${\cal NS}=T_{z_m}\Si\oplus{\cal NS}_1$},
where ${\cal NS}_1$ is the normal bundle of ${\cal S}_{{\cal T},2}$
in  ${\cal M}_{\cal T}$
described in the previous subsection.
Let  $\big(\Phi_{\cal S},\Phi_{\cal S}^{\mu}\big)$
be a regularization of ${\cal S}_{{\cal T},2}^{(m)}(\mu)$
induced by the regularization of ${\cal S}_{{\cal T},2}(\mu)$
described in Subsection~\ref{order2_case2b}.
In particular,
$$\pi_{\phi_{\cal S}(b,X)}^{\perp}
   {\cal D}_{{\cal T},\hat{3}}\tilde{\phi}_{\cal S}(b,w,X)=
\Pi_{b,\tilde{\phi}_{\cal S}(b,w,X)}X
\quad\forall 
(b,w,X)\in T_{z_m}\Si\oplus {\cal N}\tilde{\cal S}_1=T_{z_m}\Si\oplus E_1,$$
where $\tilde{\phi}_{\cal S}$ is the lift of $\phi_{\cal S}$
to ${\cal M}_{\cal T}^{(0)}$.
We also assume that $\Phi_{\cal S}^{\mu}$ is given
by the $g_{\P,b}$-parallel transport on~${\cal N}_b{\cal S}_1$.
The bundle ${\cal NS}$ carries a natural norm induced
by the $g_{\P,\ev}$-metric on~$\P$ and $g_{\cdot,\hat{0}}$-metric on~$\Si$.
Denote by $F{\cal S}$ and $F^{\eset}{\cal S}$
the bundles described in Subsection~\ref{si_str_sec}
corresponding to the submanifold~${\cal S}_{{\cal T},2}^{(m)}$.
If $(b,w,X,\ups)\!\in\! F^{\eset}{\cal S}$ is sufficiently small, let
$$\tilde{x}_h(w,\ups)
=\tilde{x}_h(\phi_{\cal S}(w,X,\ups))
=\tilde{x}_h(\phi_{\cal S}(w,0,\ups))\in\Si,\quad h=\hat{2},\hat{3}.$$
We identify a small neighborhood of $z_m$ in $\Si$ with a neighborhood 
of $0$ in $T_{z_m}\Si$ via the $g_{b,\hat{0}}$-exponential map. Put
\begin{gather*}
\begin{split}
\tilde{\al}(w,X,\ups)=
\Pi_{b,\phi_{\cal S}(b,X)}^{-1}
\Big( \big({\cal D}_{{\cal T},\hat{2}}\phi_{\cal S}(b,X)\big)
 s_{\Si,\tilde{x}_{\hat{2}}(w,\ups)}(\tilde{v}_{\hat{2}})+
\big({\cal D}_{{\cal T},\hat{3}}\phi_{\cal S}(b,X)\big)
 s_{\Si,\tilde{x}_{\hat{3}}(w,\ups)}(\tilde{v}_{\hat{3}})\quad\qquad&\\
+\big({\cal D}_{{\cal T},\hat{2}}^{(2)}\phi_{\cal S}(b,X)\big) 
      s_{b,z_m}^{(2)}(\tilde{v}_{\hat{2}})+
\big({\cal D}_{{\cal T},\hat{3}}^{(2)}\phi_{\cal S}(b,X)\big) 
      s_{b,z_m}^{(2)}(\tilde{v}_{\hat{3}})\Big);&
\end{split}\\
\begin{split}
\tilde{\al}^{\mu}(w,X,\ups)=
\Big( \big({\cal D}_{{\cal S},t\nu,\hat{2}}^{\mu}(w,X,\ups)\big)
s_{\Si,\tilde{x}_{\hat{2}}(w,\ups)}(\tilde{v}_{\hat{2}})+
\big({\cal D}_{{\cal S},t\nu,\hat{3}}^{\mu}(w,X,\ups)\big)
s_{\Si,x_{\hat{3}}(w,\ups)}(v_{\hat{3}})\Big)\quad\qquad&\\
+\Big(\big({\cal D}_{{\cal S},t\nu,\hat{2}}^{\mu,(2)}b\big) 
      s_{b,z_m}^{(m)}(\tilde{v}_{\hat{2}})+
\big({\cal D}_{{\cal S},t\nu,\hat{3}}^{\mu,(2)}b\big) 
      s_{b,z_m}(\tilde{v}_{\hat{3}})\Big)&,
\end{split}
\end{gather*}
where, with $\varphi_{{\cal S},t\nu}^{\mu}$ as in Theorem~\ref{si_str},
$${\cal D}_{{\cal S},t\nu,h}^{\mu,(k)}(w,X,\ups)=
\Pi_{\phi_{\cal S}^{\mu}\varphi_{{\cal S},t\nu}^{\mu}(w,X,\ups),
\phi_{\cal S}\Phi_{\cal S}^{\mu}\varphi_{{\cal S},t\nu}^{\mu}(w,X,\ups)}^{-1}
\Pi_{b,\phi_{\cal S}^{\mu}\varphi_{{\cal S},t\nu}^{\mu}(w,X,\ups)}^{-1}
\big({\cal D}_{{\cal T},h}^{(k)}\phi_{\cal S}\Phi_{\cal S}^{\mu}
     \varphi_{{\cal S},t\nu}^{\mu}(w,X,\ups)\big).$$
With $\ka(b)$ as in the previous subsection, let
\begin{gather*}
\al^+(\ups)=\big({\cal D}_{{\cal T},\hat{2}}b\big)
s_{\Si,z_m}\big(\tilde{v}_{\hat{2}}+\ka(b)\tilde{v}_{\hat{3}}\big),\quad
\al^-_{\hat{2}}(w,\ups)
 =\big({\cal D}_{{\cal T},\hat{2}}b\big)s_{b,z_m}^{(3,-)}
\big(\tilde{x}_{\hat{2}}(w,\ups),(x_{\hat{2}}-x_{\hat{3}})v_{\hat{1}},
    \tilde{v}_{\hat{2}}\big);\\
\al^-_{\hat{3}}(w,\ups)=
\big({\cal D}_{{\cal T},\hat{3}}b\big)s_{b,z_m}^{(3,-)}
\big(\tilde{x}_{\hat{3}}(w,\ups),(x_{\hat{3}}-x_{\hat{2}})v_{\hat{1}},
    \tilde{v}_{\hat{3}}\big).
\end{gather*}

\begin{lmm}
\label{dbar3_2b}
There exist $\de,C\!\in\! C^{\i}({\cal S}_{{\cal T},2}^{(m)};\Bbb{R}^+)$ 
such that  for all
$\vp\!=\!(b,w,X,\ups)\!\in\! F^{\eset}{\cal S}_{\de}$,
$$\Big\| \pi^{0,1}_{\Phi_{\cal S}(\vp),-}
  \bar{\partial}u_{\Phi_{\cal S}(\vp)}
-\tilde{R}_{\Phi_{\cal S}(\vp)}\tilde{\al}(w,X,\ups)\Big\|_2
\le C(b)|\vp|
\big(|\tilde{v}_{\hat{2}}|^2+|\tilde{v}_{\hat{3}}|^2\big).$$
\end{lmm}

\noindent
{\it Proof:} This lemma follows  from Proposition~\ref{dbar_gen2prp} 
applied with first- and second-order terms.

\begin{lmm}
\label{psit3_2b}
There exist $\de,C\!>\!0$ such that for all 
$\vp\!=\!(b,w,X,\ups)\!\in\! F^{\eset}{\cal S}_{\de}\big|
{\cal S}^{(m)}_{{\cal T};2}(\mu)$,
$$\Big\|\psi_{{\cal S},t\nu}^{\mu}(\vp)
-\big(t\bar{\nu}_b+
\tilde{\al}^{\mu}(w,X,\ups)\big)\Big\|_2\le 
C(t+|\vp|^{\frac{1}{p}})\big(t+
(|v_{\hat{1}}|^2+|v_{\hat{1}}||w|)
(|\tilde{v}_{\hat{2}}|+|\tilde{v}_{\hat{3}}|)\big).$$
\end{lmm}

\noindent
{\it Proof:} 
Note that the space ${\cal S}_{{\cal T},2}^{(m)}(\mu)$ is zero-dimensional
and compact if $n=3$.
As before,  we need to~bound 
$$\big\|D_{\Phi_{\cal S}(\vp)}^*
    R_{\Phi_{\cal S}(\vp)}X_i\psi_2\big\|_{L^1},$$
where the notation is as in the proof of Lemma~\ref{psi_2a2}.
By equation~\e_ref{coker_e3}, the $L^1$-norm on the annulus centered
at $\tilde{x}_{\hat{2}}=\tilde{x}_{\hat{2}}(w,\ups)$ is bounded by 
$(|\tilde{x}_{\hat{2}}||\tilde{v}_{\hat{2}}|+
 |\tilde{v}_{\hat{2}}|^2)|\tilde{v}_{\hat{2}}|$,
while the norm over the other annulus is bounded by
$(|\tilde{x}_{\hat{2}}||v_{\hat{1}}|+|v_{\hat{1}}|^2)|\tilde{v}_{\hat{3}}|$,
since the $g_{b,{\hat{0}}}$-distance between 
$\tilde{x}_{\hat{2}}$ and $\tilde{x}_{\hat{3}}$
is bounded by~$C|v_{\hat{1}}|$.
See the proof of Lemma~\ref{psit3_1} for more detail.
The claim follows from 
\hbox{$\tilde{x}_{\hat{2}}=w+x_{\hat{2}}v_{\hat{1}}$}.

\begin{lmm}
\label{psit3l_2b}
There exist $\de,C\!>\!0$ such that for all
 $\vp\!=\!(b,w,X,\ups)\!\in\! F^{\eset}{\cal S}_{\de}
\big|{\cal S}_{{\cal T};2}^{(m)}(\mu)$,
\begin{gather*}
\Big\|\tilde{\al}^{\mu}(w,X,\ups)-\al^+(w,\ups)\Big\|_2   
\le C(t+|\vp|^{\frac{1}{p}})
\big( |\tilde{v}_{\hat{2}}|+|\tilde{v}_{\hat{3}}|\big);\\
\Big\|\pi_{\tilde{x}_{\hat{2}}(w,\ups)}^-\tilde{\al}^{\mu}(w,X,\ups)
-\al_{\hat{3}}^-(w,\ups)\Big\|
\le C(t+|\vp|^{\frac{1}{p}})
\big(|v_{\hat{1}}|+|w|\big)  
   |v_{\hat{1}}| \big(|\tilde{v}_{\hat{2}}|+|\tilde{v}_{\hat{3}}|\big);\\
\Big\|\pi_{\tilde{x}_{\hat{3}}(w,\ups)}^-\tilde{\al}^{\mu}(w,X,\ups)
-\al_{\hat{2}}^-(w,\ups)\Big\|
\le C(t+|\vp|^{\frac{1}{p}})
\big(|v_{\hat{1}}|+|w|\big)  
   |v_{\hat{1}}| \big(|\tilde{v}_{\hat{2}}|+|\tilde{v}_{\hat{3}}|\big).
\end{gather*} 
\end{lmm}

\noindent
{\it Proof:} The first bound is clear from the definition of 
$\tilde{\al}^{\mu}$, since 
$$\big({\cal D}_{{\cal T},\hat{3}}b\big)
=\ka(b)\big({\cal D}_{{\cal T},\hat{2}}b\big),
\quad\hbox{\quad}\quad
|\varphi(w,X,\ups)|_b\le C(t+|\vp|^{\frac{1}{p}}).$$
Since $s_{b,z_m}^{(2,-)}\!=\!0$,
\begin{equation}\label{psit3l_2b_e1}
\big|\pi_{\tilde{x}_{\hat{2}}}^-s_{\tilde{x}_h}^{(2)}(\tilde{v}_h)\big|
\le  C\big(|\tilde{x}_{\hat{2}}|+|v_{\hat{1}}|\big)|\tilde{v}_h|^2.
\end{equation}
where $\tilde{x}_h=\tilde{x}_h(w,\ups)$.
Since $\tilde{x}_{\hat{3}}-\tilde{x}_{\hat{2}}
                        =(x_{\hat{3}}-x_{\hat{2}})v_{\hat{1}}$,
$$\Big| s_{b,\tilde{x}_{\hat{3}}}(\tilde{v}_{\hat{3}})\!-\!
\big( s_{b,\tilde{x}_{\hat{2}}}\big(\tilde{v}_{\hat{3}})
\!+\!s_{b,\tilde{x}_{\hat{2}}}^{(2)}\big(
  (x_{\hat{3}}\!-\!x_{\hat{2}})v_{\hat{1}},\tilde{v}_{\hat{3}}\big)
\!+\!s_{b,\tilde{x}_{\hat{2}}}^{(3)}
 \big((x_{\hat{3}}\!-\!x_{\hat{2}})v_{\hat{1}},
      (v_{\hat{3}}\!-\!v_{\hat{2}})v_{\hat{1}},\tilde{v}_{\hat{3}}\big)\big)
\Big|\!\le\! C|v_{\hat{1}}|^3|\tilde{v}_{\hat{3}}|.$$
Since $\pi_{\tilde{x}_{\hat{2}}}^-s_{\Si,\tilde{x}_{\hat{2}}}=0$ and 
$s_{b,z_m}^{(2,-)}=0$, 
\begin{gather*}
\Big|
\pi_{\tilde{x}_{\hat{2}}}^-s_{b,\tilde{x}_{\hat{2}}}^{(2)}
  \big( (x_{\hat{3}}-x_{\hat{2}})v_{\hat{1}},\tilde{v}_{\hat{3}}\big)
     -s_{b,z_m}^{(3,-)}\big(\tilde{x}_{\hat{2}},
            (x_{\hat{3}}-x_{\hat{2}})v_{\hat{1}},\tilde{v}_{\hat{3}}\big)\Big|
 \le C |\tilde{x}_{\hat{2}}|^2|v_{\hat{1}}||\tilde{v}_{\hat{3}}|;\\
\Big| \pi_{\tilde{x}_{\hat{2}}}^- s_{b,\tilde{x}_{\hat{2}}}^{(3)}\big(
(x_{\hat{3}}-x_{\hat{2}})v_{\hat{1}},(x_{\hat{3}}-x_{\hat{2}})v_{\hat{1}},
           \tilde{v}_{\hat{3}}\big)
-s_{b,z_m}^{(3,-)}\big( (x_{\hat{3}}-x_{\hat{2}})v_{\hat{1}},
  (x_{\hat{3}}-x_{\hat{2}})v_{\hat{1}},\tilde{v}_{\hat{3}}\big)\Big|
 \le C |\tilde{x}_{\hat{2}}||v_{\hat{1}}|^2|\tilde{v}_{\hat{3}}|.
\end{gather*}
Putting the last three equations together,
we see that
\begin{equation}\label{psit3l_2b_e4}
\Big| \pi_{\tilde{x}_{\hat{2}}}^-
               s_{b,\tilde{x}_{\hat{3}}}(\tilde{v}_{\hat{3}})-
s_{b,z_m}^{(3,-)}\big(\tilde{x}_{\hat{3}},
 (x_{\hat{3}}-x_{\hat{2}})v_{\hat{1}},\tilde{v}_{\hat{3}}\big)\Big|
\le C\big(|\tilde{x}_{\hat{2}}|+|v_{\hat{1}}|\big)
  \big(|\tilde{x}_{\hat{2}}||v_1|+|v_{\hat{1}}|^2\big)|\tilde{v}_{\hat{3}}|.
\end{equation}
The second bound follows from equations \e_ref{psit3l_2b_e1} and
\e_ref{psit3l_2b_e4}.
The last estimate is proved similarly.

\begin{crl}
\label{psit3c_2b}
There exist $\de,C\!>\!0$ such that for all 
$\vp\!=\!(b,w,X,\ups)\!\in\! F^{\eset}{\cal S}_{\de}\big|
{\cal S}_{{\cal T};2}^{(m)}(\mu)$,
$$\Big\|\psi_{{\cal S},t\nu}^{\mu}(\vp)
-\big(t\bar{\nu}_b+\tilde{\al}^{\mu}(\vp)\big)\Big\|_2\le 
C(t+|\vp|^{\frac{1}{p}})\big(t+|\tilde{\al}^{\mu}(\vp)|\big).$$
\end{crl}

\noindent
{\it Proof:}
In light of Lemma~\ref{psit3_2b}, it is sufficient to show that
\begin{equation}\label{psit3c_2b_e1}
\big(|v_{\hat{1}}|+|w|\big)|v_{\hat{1}}|
\big(|\tilde{v}_{\hat{2}}|+|\tilde{v}_{\hat{3}}|\big)
\le C\big|\tilde{\al}^{\mu}(w,X,\ups)\big|
\end{equation}
for some $C\!>\!0$.
Since $\big({\cal D}_{{\cal T},\hat{2}}b\big)s_{\Si,z_m}$,
$\big({\cal D}_{{\cal T},\hat{2}}b\big)s_{b,z_m}^{(3,-)}$
and $\big({\cal D}_{{\cal T},\hat{3}}b\big)s_{b,z_m}^{(3,-)}$
are nonzero, by Lemma~\ref{psit3l_2b}
\begin{gather*}
\big|\tilde{v}_{\hat{2}}+\ka(b)\tilde{v}_{\hat{3}}\big|
\le C\Big( \big|\tilde{\al}^{\mu}(w,X,\ups)\big|+
(t+|\vp|^{\frac{1}{p}})
      \big(|\tilde{v}_{\hat{2}}|+|\tilde{v}_{\hat{3}}|\big)\Big);\\
|\tilde{x}_h||v_{\hat{1}}||\tilde{v}_h|\le
C\Big( \big|\tilde{\al}^{\mu}(w,X,\ups)\big|+
(t+|\vp|^{\frac{1}{p}})
\big(|v_{\hat{1}}|+|w|\big)|v_{\hat{1}}|
\big(|\tilde{v}_{\hat{2}}|+|\tilde{v}_{\hat{3}}|\big)\Big).
\end{gather*}
Since $\ka(b)\!\neq\!0$, $x_{\hat{2}}\!\neq\! x_{\hat{3}}$, and
$\tilde{x}_h\!=\!w\!+\!x_hv_{\hat{1}}$, we obtain
\begin{equation}\label{psit3c_2b_e5}\begin{split}
&\big(|v_{\hat{1}}|+|w|\big)|v_{\hat{1}}|
\big(|\tilde{v}_{\hat{2}}|+|\tilde{v}_{\hat{3}}|\big)
\le C\big(|\tilde{x}_{\hat{2}}|+|\tilde{x}_{\hat{3}}|\big)
|v_{\hat{1}}|
\big(|\tilde{v}_{\hat{2}}|+|\tilde{v}_{\hat{3}}|\big)\\
&\quad\qquad\qquad\qquad\qquad\le C'\Big(
|\tilde{x}_{\hat{2}}||v_{\hat{1}}|
\big(|\tilde{v}_{\hat{2}}|+
|\tilde{v}_{\hat{2}}+\ka(b)\tilde{v}_{\hat{3}}|\big)+
|\tilde{x}_{\hat{3}}||v_{\hat{1}}|
\big(|\tilde{v}_{\hat{3}}|+
|\tilde{v}_{\hat{3}}+\ka(b)\tilde{v}_{\hat{3}}|\big)\Big)\\
&\quad\qquad\qquad\qquad\qquad\le C''\Big(
\big|\tilde{\al}^{\mu}(w,X,\ups)\big|+
(t+|\vp|^{\frac{1}{p}})
\big(|v_{\hat{1}}|+|w|\big)|v_{\hat{1}}|
\big(|\tilde{v}_{\hat{2}}|+|\tilde{v}_{\hat{3}}|\big)\Big).
\end{split}\end{equation}
If $\de$ is sufficiently small, estimate~\e_ref{psit3c_2b_e1}
follows from~\e_ref{psit3c_2b_e5}.\\

\noindent
The next step is to apply Lemma~\ref{top_l1}.  Let 
\begin{gather*}
F^+\!=L_{{\cal T},\hat{3}}^*\otimes E_1,\quad
F^-\!\!=T_{z_m}\Si\oplus F{\cal T},~~
{\cal O}^{\pm}={\cal H}_{\Si}^{\pm}\otimes\ev^*T\P,~~
\tilde{F}^-\!\!= 
\pi_{\Si}^* T\Si^{\otimes3} \otimes L_{{\cal T},\hat{2}}^{\otimes3};\\
\phi\big([b,w,\ups]\big)
=\big[b,(w+x_{\hat{2}}v_{\hat{1}})\otimes 
((x_{\hat{2}}-x_{\hat{2}})v_{\hat{1}})\otimes 
 \tilde{v}_{\hat{2}}\big];\\
\al^-\big(\phi(w,\ups)\big)\equiv 
\al_{\hat{2}}^-(w,\ups),\quad
\al(X,w,\ups)=\al^{\mu}(X,w,\ups).
\end{gather*}
Note that $\al^-\!\in\!\Ga({\cal S};\tilde{F}^{-*}\otimes{\cal O}^-)$
is well-defined.

\begin{crl}
\label{order3_2b_contr}
Suppose $d$ is a positive integer, 
${\cal T}\!=\!(\Si,[N],I;j,\under{d})$ is
a simple bubble type, with 
$\hat{I}\!=\!\{\hat{1},\hat{2},\hat{3}\}$, 
$H_{\hat{1}}{\cal T}\!=\!\{\hat{2},\hat{3}\}$,
$d_{\hat{0}}\!=\!0$, and
$\sum\limits_{i\in I}d_i\!=\!d$, and 
$\mu$ is an $N$-tuple of constraints in general position 
such that
$$\codim_{\Bbb{C}}\mu=d(n+1)-n(g-1)+N.$$
Let $\nu\!\in\!\Ga(\Si\times\P;
          \La^{0,1}\pi_{\Si}^*T^*\Si\otimes\pi_{\P}^*T\P)$
be a generic section.
There exist a neighborhood $U$ of ${\cal S}_{{\cal T},2}^{(m)}(\mu)$
in $\bar{C}^{\i}_{(d;[N])}(\Si;\mu)$,
and $\ep\!>\!0$ such that for any $t\!\in\!(0,\ep)$,
\hbox{$U\cap{\cal M}_{\Si,d,t\nu}(\mu)\!=\!\eset$.}
\end{crl}

\noindent
{\it Proof:} Analogously to the proof of Corollary~\ref{order3_1_contr},
we apply Lemma~\ref{top_l1} to the map
$$(w,\ups,X)\lra 
\pi_{z_m}^+\pi_{x_{\hat{3}}(w,\ups)}^+\psi_{{\cal S},t\nu}^{\mu}(w,\ups,X)+
\pi_{z_m}^-\pi_{x_{\hat{3}}(w,\ups)}^-\psi_{{\cal S},t\nu}^{\mu}(w,\ups,X)$$
instead of $\psi_{{\cal S},t\nu}^{\mu}$.
The claim then follows from Theorem~\ref{si_str},
Lemma~\ref{top_l1}, and Corollary~\ref{psit3c_2b}.

\subsection{Summary of Section~\ref{resolvent_sec}}
\label{summary4}

\noindent
We conclude Section~\ref{resolvent_sec} 
by reviewing the main results so~far.
Throughout this subsection,
$${\cal T}=(\Si,[N],I;j,d)$$ 
is a simple bubble type, 
with $d\!=\!\sum d_h$ and  $d_{\hat{0}}\!=\!0$, and
$\mu$ is an $N$-tuple of constraints in general position 
such that
\hbox{$\codim_{\Bbb{C}}\mu\!=\!d(n\!+\!1)\!-\!n(g\!-\!1)\!+\!N$}.\\

\noindent 
If $|\hat{I}|\!\ge\! n$, by Corollaries~\ref{order1_contr_n2} 
and~\ref{order1_contr_n3},
there exist a neighborhood $U_{\cal T}$ of $\bar{\cal M}_{\cal T}(\mu)$
in $\bar{C}^{\i}_{(d;[N])}(\Si;\mu)$ and $\ep_{\cal T}\!>\!0$ 
such that for all $t\!\in\!(0,\ep_{\cal T})$,
\hbox{$U_{\cal T}\cap{\cal M}_{\Si,d,t\nu}(\mu)\!=\!\eset$}.
This is also true if  $H_{\hat{0}}{\cal T}\!\neq\!\hat{I}$ 
\hbox{or  $M_{\hat{0}}{\cal T}\!\neq\!\eset$}.
\hbox{If $n\!=\!2$}, this statement is just Corollary~\ref{order1_contr_n2}.
\hbox{If $n\!=\!3$},  we only need to consider Cases~(1), (2b), and (3b)
of Corollary~\ref{order1_contr_n3}.
Case~(3b) follows from
Corollaries~\ref{order1_contr}, \ref{order2_2b_contr},
and~\ref{order3_2b_contr}. 
The claim for Case~(2b) is obtained from
Corollaries~\ref{order1_contr}, \ref{order2_1_contr},
\ref{order3_1_contr} and the same claim for Case~(3b).
Finally, in Case~(1), we use 
Corollaries~\ref{order1_contr}, \ref{order2_2b_contr},
and~\ref{order3_2b_contr}, 
the statement of Corollary~\ref{order1_contr_n3} 
\hbox{for $|\hat{I}|\!\ge\!2$},
and the just stated result for Case~(2b).
\\

\noindent
If $|\hat{I}|\!\le\! n$, $H_{\hat{0}}{\cal T}\!\!=\!\hat{I}$, 
and $M_{\hat{0}}{\cal T}\!=\!\eset$,
i.e.~${\cal T}$ is a primitive bubble type,
by the previous paragraph 
and Corollaries~\ref{order1_contr}, \ref{order2_1_contr},
\ref{order3_1_contr}, and~\ref{order2_2a_contr},
there exist a neighborhood $U_{\cal T}$ of $\bar{\cal M}_{\cal T}(\mu)$
in $\bar{C}^{\i}_{(d;[N])}(\Si;\mu)$ and $\ep_{\cal T}\!>\!0$ 
such that for all $t\!\in\!(0,\ep_{\cal T})$,
the signed cardinality $n_{\cal T}(\mu)$ of
$U_{\cal T}\cap{\cal M}_{\Si,d,t\nu}(\mu)$ is 
the sum of the numbers given by these four corollaries applied
\hbox{to ${\cal T}$}.
\hbox{If $|\hat{I}|\!=\!1$},
\begin{equation}\label{summary_m1}
n_1(\mu)\equiv n_{\cal T}(\mu)
=n_1^{(1)}(\mu)+
 2n_1^{(2)}(\mu)+18n_1^{(3)}(\mu),
\end{equation}
where the numbers $n_1^{(k)}(\mu)$ are described as follows.
The number $n_1^{(1)}(\mu)$ is the signed number of zeros
of the affine map 
\begin{equation}
\psi_1^{(1)}\!:  T\Si\otimes L_{\bar{\cal T},\hat{1}}
\lra{\cal H}_{\Si}^{0,1}\otimes\ev^*T\P,\qquad
\label{summary_m1c1}
\psi_1^{(1)}(x,[b,v_{\hat{1}}])=
\bar{\nu}_b+\big({\cal D}_{{\cal T},\hat{1}}b\big)s_{\Si,x}(v_{\hat{1}}),
\end{equation}
where the bundles are considered over 
$\Si\!\times\!\bar{\cal U}_{\bar{\cal T}}(\mu)\!=\!\bar{\cal M}_{\cal T}(\mu)$
and $\hat{1}$ is the unique element of~$\hat{I}$.
Note that this number is the same as the number of zeros
of the map in~\e_ref{order1_contr_e}, since 
\hbox{$\Si\!\times\!\bar{\cal U}_{\bar{\cal T}}(\mu)\!
   -\!{\cal M}_{\cal T}(\mu)$}
is a finite union of smooth manifolds of dimension less than the dimension
of ${\cal M}_{\cal T}(\mu)$.
Thus, if $\nu$ is generic, $\psi_1^{(1)}$ has no zeros over 
\hbox{$\Si\!\times\!\bar{\cal U}_{\bar{\cal T}}(\mu)
                        \!-\!{\cal M}_{\cal T}(\mu)$}.
The number $n_1^{(2)}(\mu)$ is the signed number of zeros
of the affine~map 
\begin{equation}\label{summary_m1c2}
\psi_1^{(2)}\!:  
T\Si^{\otimes2}\otimes L_{\bar{\cal T},\hat{1}}^{\otimes2}
\lra {\cal H}_{\Si}^-\otimes\ev^*T\P,\qquad
\psi_2^{(2)}(x,[b,v_{\hat{1}}])=
\bar{\nu}_b^-+\big({\cal D}_{{\cal T},\hat{1}}^{(2)}b\big)
                         s_{\Si,x}^{(2,-)}(v_{\hat{1}}),
\end{equation}
where the bundles are considered over
$\Si\!\times\!\bar{\cal S}_1(\mu)$
and $\bar{\cal S}_1(\mu)$ is the closure in 
$\bar{\cal U}_{\bar{\cal T}}(\mu)$ of the space
\begin{equation}\label{summary_m1c2b}
{\cal S}_1(\mu)=
\big\{ b\!\in\!{\cal U}_{\bar{\cal T}}(\mu)\!: 
{\cal D}_{{\cal T},\hat{1}}\big|_b\!=\!0\big\}.
\end{equation}
If $n\!=\!2$, ${\cal S}_1(\mu)$ is a finite set and thus 
\hbox{$\bar{\cal S}_1(\mu)\!=\!{\cal S}_1(\mu)$}.
If $n\!=\!3$,  ${\cal S}_1(\mu)$ is one-dimensional 
over~$\Bbb{C}$.
The boundary 
\hbox{$\bar{\cal S}_1(\mu)\!-\!{\cal S}_1(\mu)$}
is a finite set, as can be seen from the estimate on 
${\cal D}_{\bar{\cal T},\hat{1}}$  of Theorem~\ref{str_global}.
Thus, in either case, the maps in~\e_ref{summary_m1c2}
and~\e_ref{order2_1contr_e} have the same zeros.
Finally, the number $n_1^{(3)}(\mu)$ is the signed number of zeros
of the affine map 
\begin{gather}\label{summary_m1c3}
\psi_1^{(3)}\!:  
T\Si^{\otimes3}\otimes\big( L_{\bar{\cal T},\hat{1}}^{\otimes2}
\oplus L_{\bar{\cal T},\hat{1}}^{\otimes3}\big)
\lra {\cal H}_{\Si}^-\otimes\ev^*T\P,\\
\psi_1^{(3)}(x,[b,v_{\hat{1}},w_{\hat{1}}])=
\bar{\nu}_b^- +\big({\cal D}_{{\cal T},\hat{1}}^{(2)}b\big)
                       s_{b,z_m}^{(3,-)}(v_{\hat{1}})+
\big({\cal D}_{{\cal T},\hat{1}}^{(3)}b\big)
                         s_{b,z_m}^{(3,-)}(w_{\hat{1}}),\notag
\end{gather}
where the bundles are considered over 
$\bar{\cal S}_1(\mu)$
and $z_m$ is one of the six distinguished points of~$\Si$.
By the same argument as above, this number is precisely
the number of zeros of the map in~\e_ref{order3_1_contr_e}.\\

\noindent
If $|\hat{I}|\!=\!2$ and $n\!=\!2$, 
$n_{\cal T}(\mu)\!=\!n_{\cal T}^{(1)}(\mu)$
is the signed number of zeros of the affine map
\begin{gather}\label{summary_m2c1}
\psi_{\cal T}^{(1)}\!:
T\Si_{\hat{1}}\otimes L_{\bar{\cal T},\hat{1}}\oplus
T\Si_{\hat{2}}\otimes L_{\bar{\cal T},\hat{2}}\lra
{\cal H}_{\Si}^{0,1}\otimes\ev^*T\P,\\
\psi_{\cal T}^{(1)}
\big(x_{\hat{1}},x_{\hat{2}},[b,v_{\hat{1}},v_{\hat{2}}]\big)=
\bar{\nu}_b+\big({\cal D}_{{\cal T},\hat{1}}b\big)
    s_{\Si,x_{\hat{1}}}(v_{\hat{1}})
+\big({\cal D}_{{\cal T},\hat{2}}b\big)s_{\Si,x_{\hat{2}}}(v_{\hat{2}}),
\notag
\end{gather}
where the bundles are considered over 
$\Si^2\times\bar{\cal U}_{\bar{\cal T}}(\mu)=
\Si_{\hat{1}}\times\Si_{\hat{2}}\times\bar{\cal U}_{\bar{\cal T}}(\mu)$
and $\hat{1},\hat{2}$ are the two elements of~$\hat{I}$.
By the same argument as before,
the number $n_{\cal T}^{(1)}(\mu)$ is the same as the number of 
zeros of the map~\e_ref{order1_contr}.
If $|\hat{I}|\!=\!2$ and $n\!=\!3$,
\begin{equation}\label{summary_m2}
n_{\cal T}(\mu)=n_{\cal T}^{(1)}(\mu)+2n_{\cal T}^{(2)}(\mu),
\end{equation}
where $n_{\cal T}^{(1)}(\mu)$ is defined the same way as in 
the $n\!=\!2$ case, while $n_{\cal T}^{(2)}(\mu)$ is the signed number
 of zeros of the affine map
\begin{gather}\label{summary_m2c2}
\psi_{\cal T}^{(2)}\!:
T\Si^{\otimes2}\otimes\big( L_{\bar{\cal T},\hat{2}}\oplus
L_{\bar{\cal T},\hat{2}}^{\otimes2}\big)\lra
{\cal H}_{\Si}^-\otimes\ev^*T\P,\\
\psi_{\cal T}^{(2)}
\big(x,[b,v_{\hat{2}},w_{\hat{2}}]\big)\!=\!
\bar{\nu}_b^-\!+
\big({\cal D}_{{\cal T},\hat{2}}b\big)s_{\Si,x}^{(2,-)}(w_{\hat{2}})
+\big({\cal D}_{{\cal T},\hat{1}}^{(2)}b\big)
 s_{\Si,x}^{(2,-)}\big(\ka(b)v_{\hat{2}}\big)+
\big({\cal D}_{{\cal T},\hat{2}}^{(2)}b\big)s_{\Si,x}^{(2,-)}(v_{\hat{2}}),
\notag
\end{gather}
where the bundles are viewed over
$\Si\times{\cal S}_{\bar{\cal T}}(\mu)$,
\begin{equation}\label{summary_m2c2b}
{\cal S}_{\bar{\cal T}}(\mu)=\big\{
b\!\in\!{\cal U}_{\bar{\cal T}}(\mu)\!: 
\pi^{\perp}_{[b]}\!\circ{\cal D}_{{\cal T},\hat{2}}\big|_{[b]}\!=\!0\big\},
\end{equation}
$E_1$~is the quotient of  $\ev^*T\P$
by~$\hbox{Im}({\cal D}_{{\cal T},\hat{1}})$,
$\pi^{\perp}\!\!:\ev^*T\P\!\lra\! E_1$ is the projection map, and
\hbox{$\ka(b)\!\in\! L_{\bar{\cal T},\hat{2}}^*\otimes 
 L_{\bar{\cal T},\hat{1}}$}
is a nonzero homomorphism.
Note that ${\cal S}_{\bar{\cal T}}(\mu)$ is a finite set
with our choice of constraints.
Finally, if $|\hat{I}|\!=\!3$ and $n\!=\!3$, 
$n_{\cal T}(\mu)\!=\!n_{\cal T}^{(1)}(\mu)$
is the signed number of zeros of the affine~map
\begin{gather}\label{summary_m3c1}
\psi_{\cal T}^{(1)}\!:
T\Si_{\hat{1}}\otimes L_{\bar{\cal T},\hat{1}}\oplus
T\Si_{\hat{2}}\otimes L_{\bar{\cal T},\hat{2}}\oplus
T\Si_{\hat{3}}\otimes L_{\bar{\cal T},\hat{3}}\lra
{\cal H}_{\Si}^{0,1}\otimes\ev^*T\P,\\
\psi_{\cal T}^{(1)}
\big(x_{\hat{1}},x_{\hat{2}},x_{\hat{3}},
[b,v_{\hat{1}},v_{\hat{2}},v_{\hat{3}}]\big)=
\bar{\nu}_b+\big({\cal D}_{{\cal T},\hat{1}}b\big)
    s_{\Si,x_{\hat{1}}}(v_{\hat{1}})
+\big({\cal D}_{{\cal T},\hat{2}}b\big)s_{\Si,x_{\hat{2}}}(v_{\hat{2}})
+\big({\cal D}_{{\cal T},\hat{3}}b\big)s_{\Si,x_{\hat{3}}}(v_{\hat{3}}),
\notag
\end{gather}
where the bundles are considered over 
$\Si^3\times\bar{\cal U}_{\bar{\cal T}}(\mu)=
\Si_{\hat{1}}\times\Si_{\hat{2}}\times\Si_{\hat{3}}\times
\bar{\cal U}_{\bar{\cal T}}(\mu)$
and $\hat{1},\hat{2},\hat{3}$ are the three elements of~$\hat{I}$.
As before,
the number $n_{\cal T}^{(1)}(\mu)$ is precisely the number of 
zeros of the map~\e_ref{order1_contr}.
If $m\!\ge\!2$ and $k\!\ge\!1$, we denote by $n_m^{(k)}(\mu)$
the sum of the numbers $n_{\cal T}^{(k)}(\mu)$ 
over all equivalence classes of primitive bubble types ${\cal T}$
with $\!|\hat{I}|\!=\!m$.

\section{Computations}
\label{comp_sect}

\subsection{The Numbers $n_m^{(1)}(\mu)$ with $m=n$}
\label{comp_setup}

\noindent
Our goal now is to compute the numbers $n_{\cal T}^{(k)}(\mu)$ for any 
primitive bubble type 
\hbox{${\cal T}\!=\!(\Si,[N],I;j,d)$},
and thus the genus-two enumerative invariants for $\PP$ and~$\PPP$.
Most of this section is devoted to
expressing the numbers $n_{\cal T}^{(k)}(\mu)$
in terms of intersection numbers of tautological classes
of various spaces of stable rational maps that pass through 
the constraints~$\mu$.
These are shown to be computable in~\cite{P2}.
The procedure for counting the zeros of affine
maps between vector bundles is described in Section~\ref{top_sect}.
We start with the easiest cases.

\begin{lmm}
\label{mequalsn_l}
If ${\cal T}\!=\!(\Si,[N],I;j,d)$ is
a primitive bubble type with $|\hat{I}|\!=\!n$ and
$\mu$ is an $N$-tuple of constraints in general position 
such that
$$\codim_{\Bbb{C}}\mu=(n+1)\sum_{i\in M}d_i-n+N,$$
the set  $\bar{\cal U}_{\bar{\cal T}}(\mu)$ is finite and 
$n_{\cal T}^{(1)}(\mu)=
2^n|\bar{\cal U}_{\bar{\cal T}}(\mu)|$.
\end{lmm}

\noindent
{\it Proof:} The first statement is clear by dimension counting.
By equations~\e_ref{summary_m2c1} and~\e_ref{summary_m3c1},
we need to  apply Lemma~\ref{zeros_main} with
$$\bar{\cal M}=\begin{cases}
\Si_{\hat{1}}\times\Si_{\hat{2}}\times
\bar{\cal U}_{\bar{\cal T}}(\mu),&\!\hbox{if~}n\!=\!2;\\
\Si_{\hat{1}}\times\Si_{\hat{2}}\times\Si_{\hat{3}}\times
\bar{\cal U}_{\bar{\cal T}}(\mu)
,&\!\hbox{if~}n\!=\!3;
\end{cases}\quad
E=
\begin{cases}
T\Si_{\hat{1}}\otimes L_{\bar{\cal T},\hat{1}}\oplus
T\Si_{\hat{2}}\otimes L_{\bar{\cal T},\hat{2}},&\!\hbox{if~}n\!=\!2;\\
T\Si_{\hat{1}}\otimes L_{\bar{\cal T},\hat{1}}\oplus
T\Si_{\hat{2}}\otimes L_{\bar{\cal T},\hat{2}}\oplus
T\Si_{\hat{3}}\otimes L_{\bar{\cal T},\hat{3}},&\!\hbox{if~}n\!=\!3,
\end{cases}$$
${\cal O}\!=\!{\cal H}_{\Si}^{0,1}\otimes\ev^*T\P$,
and $\al$ given by~\e_ref{summary_m2c1} and~\e_ref{summary_m3c1}.
By Lemma~\ref{order1_contr_l3},
\hbox{$\al\!\in\!\Ga(\bar{\cal M};E^*\otimes{\cal O})$} has 
full rank on every fiber of~$E$.
Thus by Lemma~\ref{zeros_main},
\begin{equation}\label{mequalsm_e1}
n_{\cal T}^{(1)}(\mu)
=\big\lan e\big({\cal O}/\al(E)\big),[\bar{\cal M}]\big\ran
=\big\lan c({\cal O})c(E)^{-1},[\bar{\cal M}]\big\ran.
\end{equation}
Since $\bar{\cal U}_{\bar{\cal T}}(\mu)$ is a finite set,
$$E\approx\begin{cases}
T\Si_{\hat{1}}\oplus T\Si_{\hat{2}},&\hbox{if~}n=2;\\
T\Si_{\hat{1}}\oplus T\Si_{\hat{2}}\oplus T\Si_{\hat{3}},&\hbox{if~}n=3;
\end{cases}
\qquad
{\cal O}\approx \bar{\cal M}\times\Bbb{C}^{2n}.$$
Let $y_h\!=\!c_1(T\Si_h)$. 
Thus, if $n\!=\!2$, by~\e_ref{mequalsm_e1}
$$n_{\cal T}^{(1)}(\mu)
=\big\lan \big(1+(y_{\hat{1}}+y_{\hat{2}})+y_{\hat{1}}y_{\hat{2}}\big)^{-1},
[\bar{\cal M}]\big\ran
=\big\lan y_{\hat{1}}y_{\hat{2}},\big[\Si_{\hat{1}}\!\times\!\Si_{\hat{2}}]
\big\ran|\bar{\cal U}_{\bar{\cal T}}(\mu)|=
4|\bar{\cal U}_{\bar{\cal T}}(\mu)|,$$
since $\lan y_h,[\Si_h]\ran\!=\!-2$.
If $n\!=\!3$, we similarly obtain   
$$n_{\cal T}^{(1)}(\mu) 
=\big\lan -y_{\hat{1}}y_{\hat{2}}y_{\hat{3}},
\big[\Si_{\hat{1}}\!\times\!\Si_{\hat{2}}\!\times\!\Si_{\hat{3}}]\big\ran
|\bar{\cal U}_{\bar{\cal T}}(\mu)|=
8|\bar{\cal U}_{\bar{\cal T}}(\mu)|,$$
as claimed.\\

\noindent
Let $\tau_n(\mu)$ denote the sum of 
the numbers of $|\bar{\cal U}_{\bar{\cal T}}(\mu)|$
taken over all equivalence classes of primitive bubble types 
${\cal T}$ with $\hat{I}\!=\!n$.
This is the number of $n$-component connected curves of 
total degree~$d$ passing through 
the constraints $\mu_1,\ldots,\mu_N$ in $\P$
with a choice of a node which belongs to all $n$ components.
From Lemma~\ref{mequalsn_l}, we immediately conclude:

\begin{crl}
\label{mequalsn}
If $n\!=\!2$, $n_2^{(1)}=4\tau_2(\mu)$.
If $n\!=\!3$, $n_3^{(1)}=8\tau_3(\mu)$.
\end{crl}

\subsection{The Numbers 
$n_m^{(2)}(\mu)$ and $n_m^{(3)}(\mu)$  with $m=n-1$}
\label{comp_sec2}

\noindent
In this subsection, we describe the numbers 
$n_{\cal T}^{(2)}(\mu)$ and $n_{\cal T}^{(3)}(\mu)$ with 
$|\hat{I}|\!=\!n\!-\!1$ topologically.
The similarity between these cases is that 
${\cal U}_{\bar{\cal T}}(\mu)$ is two-dimensional 
(over~$\Bbb{C}$), while ${\cal S}_{\bar{\cal T}}(\mu)$ is
a~finite set; see Subsection~\ref{summary4} for notation.\\

\noindent
The numbers $n_{\cal T}^{(2)}(\mu)$ with $|\hat{I}|\!=n\!-\!1\!=\!1$ and
$|\hat{I}|\!=\!n\!-\!1\!=\!2$ are 
the signed cardinalities of the zero sets of the affine maps in
\e_ref{summary_m1c2} and~\e_ref{summary_m2c2}, respectively.
By Subsections~\ref{order2_case1} and~\ref{order2_case2a},
the linear part $\al$ of the affine map $\psi_{\cal T}^{(2)}$
has full rank in these cases, except over the zero set of~$s_{\Si}^{(2,-)}$.
In order to simplify our computations, we replace~$s_{\Si}^{(2,-)}$
by another section that has no zeros on~$\Si$, but so that
the corresponding affine maps have the same number of zeros
as the maps in \e_ref{summary_m1c2} and~\e_ref{summary_m2c2}.
The section
$$s_{\Si}^{(2,-)}\in\Ga(\Si;T^*\Si^{\otimes2}\otimes{\cal H}_{\Si}^-)$$
has transverse zeros at the points $z_1,\ldots,z_6\!\in\!\Si$;
see Subsection~\ref{order3_case1}.
Thus, it induces a nonvanishing section
$$\tilde{s}_{\Si}^{(2,-)}\in
\Ga(\Si;\tilde{T}\Si^*\otimes{\cal H}_{\Si}^-),
\quad\hbox{where}\quad
\tilde{T}\Si=T\Si^{\otimes2}\otimes{\cal O}(z_1)\otimes\ldots
\otimes{\cal O}(z_6)$$
and ${\cal O}(z_m)$ denotes the holomorphic line bundle
corresponding to the divisor~$z_m$ on~$\Si$.
The bundles $\tilde{T}\Si$ and $T\Si^{\otimes2}$
can be identified on~$\Si^*$, the complement of the six points,
in such a way that 
\hbox{$\tilde{s}_{\Si}^{(2,-)}\!=\!\eta s_{\Si}^{(2,-)}$}
on~$\Si^*$ for some \hbox{$\eta\!\in\! C^{\i}(\Si^*;\Bbb{R}^+)$}.
Let $\tilde{\psi}_{\cal T}^{(2)}$ denote the affine maps
obtained by replacing $T\Si^{\otimes2}$ and $s_{\Si}^{(2,-)}$ 
by $\tilde{T}\Si$ and $\tilde{s}_{\Si}^{(2,-)}$,
respectively, in \e_ref{summary_m1c2} and~\e_ref{summary_m2c2}
(depending on~${\cal T}$).
Since $\psi_{\cal T}^{(2)}$ and $\tilde{\psi}_{\cal T}^{(2)}$
have no zeros over~$\{z_m\}$ if $\nu$ is generic
and $s_{\Si}^{(2,-)}$ and $\tilde{s}_{\Si}^{(2,-)}$ differ
by a nonzero multiple on~$\Si^*$, 
there is a sign-preserving bijection between the zeros of
$\psi_{\cal T}^{(2)}$ and of~$\tilde{\psi}_{\cal T}^{(2)}$.
Furthermore, the linear part of~$\tilde{\psi}_{\cal T}^{(2)}$
has full rank on every fiber.\\

\noindent
Denote by ${\cal S}_2(\mu)$ the union of the spaces  
${\cal S}_{\bar{\cal T}}(\mu)$ defined by equation~\e_ref{summary_m2c2b}
taken over all equivalence classes of appropriate bubble types~${\cal T}$.
This set can be identified with the degree-$d$ two-component rational 
curves in $\PPP$ that are connected at a tacnode and pass through 
the constraints~$\mu$.
Similarly, in the $n\!=\!2$,
${\cal S}_1(\mu)$ corresponds to the degree-$d$ cuspidal rational curves
passing through the constraints.

\begin{lmm}
\label{mequalsnm1}
If $n\!=\!2$, $n_1^{(2)}(\mu)=2|{\cal S}_1(\mu)|$ and 
$n_1^{(3)}(\mu)=|{\cal S}_1(\mu)|$.
If $n\!=\!3$, $n_2^{(2)}(\mu)=2|{\cal S}_2(\mu)|$.
\end{lmm}

\noindent
{\it Proof:} Let ${\cal T}\!=\!(\Si,[N],I;j,d)$
be a bubble type that contributes to one of these numbers.
By dimension counting and Corollary~\ref{reg_crl2},
${\cal S}_{\bar{\cal T}}(\mu)$ is zero-dimensional and compact.
Thus, in all cases the bundles $L_{\bar{\cal T},h}$ and $\ev^*T\P$
of equations~\e_ref{summary_m1c2}, \e_ref{summary_m1c3} and 
\e_ref{summary_m2c2} are trivial.
If $n\!=\!2$ and $k\!=\!2$, we are in the case of \e_ref{summary_m1c2}.
By the above, we can apply Lemma~\ref{zeros_main} with
$$E=\tilde{T}\Si,  \quad{\cal O}={\cal H}_{\Si}^-\oplus{\cal H}_{\Si}^-,$$
and 
$\al\!\in\!\Ga(\Si\!\times\!{\cal S}_{\bar{\cal T}}(\mu);E^*\otimes{\cal O})$
that has full rank.
We obtain
$$n_{\cal T}^{(2)}(\mu)
= \big\lan c_1({\cal O})-c_1(E),
  \big[\Si\!\times\!{\cal S}_{\bar{\cal T}}(\mu) ]\big\ran
= \big(4+(4-6)\big)|{\cal S}_{\bar{\cal T}}(\mu)|
= 2|{\cal S}_{\bar{\cal T}}(\mu)|.$$
If $n\!=\!3$ and $k\!=\!2$, we are in the case of \e_ref{summary_m2c2} and
apply Lemma~\ref{zeros_main} with
$$E=\tilde{T}\Si\oplus \tilde{T}\Si,   \quad
{\cal O}={\cal H}_{\Si}^-\oplus{\cal H}_{\Si}^-\oplus{\cal H}_{\Si}^-,$$
and 
$\al\!\in\!\Ga(\Si\!\times\!{\cal S}_{\bar{\cal T}}(\mu);E^*\otimes{\cal O})$
that again has full rank.
Thus, 
$$n_{\cal T}^{(2)}(\mu)=
\big\lan c_1({\cal O})-c_1(E),
  \big[\Si\!\times\!{\cal S}_{\bar{\cal T}}(\mu) ]\big\ran
=(6-4)|{\cal S}_{\bar{\cal T}}(\mu)|=
2|{\cal S}_{\bar{\cal T}}(\mu)|.$$
Finally, 
if $n\!=\!2$ and $k\!=\!3$, we are in the case of \e_ref{summary_m1c3}.
Note that  all the bundles involved are trivial and 
the linear part of $\psi_{\cal T}^{(2)}$ is an isomorphism on every fiber.
Thus, $n_{\cal T}^{(3)}(\mu)\!=\!|{\cal U}_{\bar{\cal T}}(\mu)|$.
\\

\noindent
The next step is to compute the cardinalities of the sets
${\cal S}_{n-1}(\mu)$.
In order to simplify our answers, it is convenient to introduce
cohomology classes $c_1({\cal L}_{{\cal T},k}^*)$
closely related to~$c_1(L_{{\cal T},k}^*)$.
Suppose ${\cal T}\!=\!(S^2,M,I;j,d)$
is a bubble type, 
and $\big\{{\cal T}_k\!=\!(S^2,M_k,I_k;j_k,d_k)\big\}$
are the corresponding simple types; see~\cite{Z}.
For any $k\!\in\! I\!-\!\hat{I}$ and nonempty subset
$M_0$ of $M_k{\cal T}$,
we define bubble types ${\cal T}(M_0)$ and 
 ${\cal T}/M_0$ as follows.~Let
$${\cal T}/M_0=\big(S^2,\hat{I},M-M_0;j|(M-M_0),d|\hat{I}\big).$$
Let ${\cal T}(M_0)\equiv(S^2,M,\hat{I}+_k\hat{1};j',d')$ 
be given~by
$$j'_l=\begin{cases}
k,&\hbox{if~}l\in M_0;\\
\hat{1},&\hbox{if~}l\in M_k{\cal T}-M_0;\\
j_l,&\hbox{otherwise};
\end{cases}\qquad
d'_i=\begin{cases}
0,&\hbox{if~}i=k;\\
d_k,&\hbox{if~}i=\hat{1};\\
d_i,&\hbox{otherwise}.
\end{cases}$$
The tuples ${\cal T}/M_0$ and ${\cal T}(M_0)$ are bubble types as long as 
$d_k\!\neq\!0$ or $M_0\!\neq\! M_{\hat{0}}{\cal T}$.
If all elements of $\hbox{Aut}({\cal T})$ fix~$k$,
\begin{equation}
\label{cart_split2}
\bar{\cal U}_{{\cal T}(M_0)}(\mu)=
 \bar{\cal M}_{0,\{\hat{1}\}+M_0}\times \bar{\cal U}_{{\cal T}/M_0}
\Big(\bigcap_{l\in M_0}\mu_l;\mu\Big),
\end{equation}
where $\bar{\cal M}_{0,\{\hat{1}\}+M_0}$ denotes the
Deligne-Mumford moduli spaces of rational curves with 
$(\{\hat{0},\hat{1}\}+M_0)$-marked points. 
If $k$ is not fixed by $\hbox{Aut}({\cal T})$,
the space $\bar{\cal U}_{{\cal T}/M_0}$ above should be replaced
by a finite cover.
If $l\!\in\! M_k{\cal T}$ for some $k\!\in\! I\!-\!\hat{I}$,
we denote ${\cal T}(\{l\})$ by~${\cal T}(l)$.
If ${\cal T}$ is a basic bubble type,
by Theorem~\ref{str_global} and 
decomposition~\e_ref{cart_split2},
$\bar{\cal U}_{{\cal T}(M_0)}(\mu)$ is an oriented topological 
suborbifold of~$\bar{\cal U}_{\cal T}(\mu)$ of (real) codimension~two.
Thus,
\begin{equation}\label{normal_bundle1}
c_1({\cal L}_{{\cal T},k}^*)\equiv c_1(L_{{\cal T},k}^*)-
\sum_{M_0\subset M_k,M_0\neq\eset}
PD_{\bar{\cal U}_{\cal T}(\mu)}\big[\bar{\cal U}_{{\cal T}(M_0)}(\mu)\big]
\in H^2\big(\bar{\cal U}_{\cal T}(\mu)\big),
\end{equation}
where 
$PD_{\bar{\cal U}_{\cal T}(\mu)}\big[\bar{\cal U}_{{\cal T}(M_0)}(\mu)\big]$
denotes the Poincare Dual of $\big[\bar{\cal U}_{{\cal T}(M_0)}(\mu)\big]$
in~$\bar{\cal U}_{\cal T}(\mu)$,
is a well-defined cohomology class.
Since our constraints~$\mu$ are disjoint, 
\hbox{$\bar{\cal U}_{{\cal T}(M_0)}(\mu)=\eset$} if~$|M_0|\!\ge\!2$.
Furthermore,  it is well-known in algebraic geometry that 
for any $l\!\in\! M_k$
the normal bundle of $\bar{\cal U}_{{\cal T}(l)}(\mu)$ in
$\bar{\cal U}_{\cal T}(\mu)$ is $L_{{\cal T}(l),\hat{1}}$; see~\cite{P2}.
Thus, if $\mu$ is an $M$-tuple of disjoint constraints,
\begin{equation}\label{normal_bundle2}
\big[\bar{\cal U}_{{\cal T}(l)}(\mu)\big] 
                        \cap c_1({\cal L}_{{\cal T},k}^*)=
\big[\bar{\cal U}_{{\cal T}(l)}(\mu)\big]
                        \cap c_1(L_{{\cal T}(l),\hat{1}}^*)=
\big[\bar{\cal U}_{{\cal T}(l)}(\mu)\big]
                    \cap c_1({\cal L}_{{\cal T}(l),\hat{1}}^*),
\end{equation}
since $L_{{\cal T},k}|\bar{\cal U}_{{\cal T}(l)}$ is 
a trivial line bundle.
The above fact from algebraic geometry is only used  to simplify
notation and is not really needed for our computations.
In addition, \e_ref{normal_bundle2} can deduced from 
Subsection~\ref{chern_class}.\\

\noindent
In the $n\!=\!3$ case, we denote by $\bar{\cal V}_2(\mu)$
the disjoint union of the spaces $\bar{\cal U}_{\cal T}(\mu)$
taken over equivalence classes of basic bubble types
${\cal T}\!=\!(S^2,M,I;j,d)$ \hbox{with $|I|\!=\!2$}.
While the components of  $\bar{\cal V}_2(\mu)$ are unordered,
we can still define the chern classes 
$$c_1({\cal L}_1^*)+c_1({\cal L}_2^*),~
c_1^2({\cal L}_1^*)+c_1^2({\cal L}_2^*),~
c_1({\cal L}_1^*)c_1({\cal L}_2^*)
\in H^*\big(\bar{\cal V}_2(\mu)\big).$$
In the notation of the previous paragraph,
$c_1({\cal L}_i^*)$ denotes the cohomology 
class~$c_1({\cal L}_{{\cal T}_{k_i},{k_i}}^*)$,
where we write $I\!=\!\{k_1,k_2\}$.
If ${\cal T}^*\!=\!(S^2,M,I^*;j^*,d^*)$ 
is the unique basic simple type such that $|I|\!=\!1$,
we denote by $\bar{\cal V}_1(\mu)$ the space 
$\bar{\cal U}_{{\cal T}^*}(\mu)$ and by
$c_1({\cal L}^*)\!\in\! H^2\big(\bar{\cal V}_1(\mu)\big)$ 
the cohomology class~$c_1({\cal L}_{{\cal T}^*,\hat{0}}^*)$.

\begin{lmm}
\label{n2cusps}
If $d\!\ge\!1$, the number of rational degree-$d$ cuspidal curves 
passing through a tuple~$\mu$ of $3d\!-\!2$ 
points in general position in~$\PP$ is given by
$$\big|{\cal S}_1(\mu)\big|=
\big\lan 3a^2+3ac_1({\cal L}^*)+c_1^2({\cal L}^*),
 \big[\bar{\cal V}_1(\mu)\big]\big\ran-\tau_2(\mu),$$
where $a=\ev^*\big({\cal O}(1)\big)$.
\end{lmm}

\noindent
{\it Proof:} (1) This result is well-known in algebraic geometry;
see~\cite{V}. 
Nevertheless, for the sake of completeness, we include a proof.
Let~${\cal T}^*$ be as above.
By definition, ${\cal S}_1(\mu)$ is the intersection
of the zero set of the section 
$${\cal D}\equiv
{\cal D}_{{\cal T}^*}\in
\Ga\big(\bar{\cal V}_1(\mu); L^*\otimes\ev^*T\PP\big),
\qquad\hbox{where}\quad L\!=\!L_{{\cal T}^*,\hat{0}},$$
with ${\cal V}_1(\mu)={\cal U}_{{\cal T}^*}(\mu)$.
Thus, by Corollary~\ref{euler_crl},
with $\partial\bar{\cal V}_1(\mu)=
\bar{\cal V}_1(\mu)-{\cal V}_1(\mu)$,
\begin{equation}\label{n2cusps_e1}\begin{split}
|{\cal S}_1(\mu)|
&=\big\lan c_2\big(L^*\otimes\ev^*T\PP),
      \big[\bar{\cal V}_1(\mu)\big]\big\ran-
{\cal C}_{\partial\bar{\cal V}_1(\mu)}({\cal D})\\
&=\big\lan 3a^2+3ac_1(L^*)+
c_1^2(L^*),\big[\bar{\cal V}_1(\mu)\big]
\big\ran-
{\cal C}_{\partial\bar{\cal V}_1(\mu)} ({\cal D}).
\end{split}\end{equation}
(2) Suppose ${\cal T}=(S^2,[N],I;j,d)\!<\!{\cal T}^*$, 
where $N\!=\!3d\!-\!2$, is a bubble type such that
${\cal D}$ vanishes somewhere  on~${\cal U}_{\cal T}(\mu)$.
Since the complex dimension of~${\cal U}_{\cal T}(\mu)$ is at most~one,
by Corollary~\ref{reg_crl2}, $d_{\hat{0}}\!=\!0$.
Let 
$$\rho_{\cal T}\in
\Ga\big({\cal U}_{\cal T}(\mu);
\hbox{Hom}({\cal F}{\cal T};\tilde{\cal F}{\cal T})\big)
\quad\hbox{and}\quad
\al_{\cal T}\in
\Ga\big({\cal U}_{\cal T}(\mu);
\hbox{Hom}(\tilde{\cal F}{\cal T}; L^*\otimes\ev^*T\P)\big)$$
be the sections defined in equation~\e_ref{str_global_e}.
Recall that with appropriate identifications
\begin{equation}\label{n2cusps_e3}
\Big|{\cal D}\big(\ga_{\cal T}^{\mu}(\ups)\big)-
\al_{\cal T}\big(\rho_{\cal T}(\ups)\big)\Big|\le 
C(b_{\ups})|\ups|^{\frac{1}{p}}\big|\rho_{\cal T}(\ups)\big|
\quad\forall\ups\in {\cal FT}_{\de},
\end{equation}
where $\de,C\!\in\!C^{\i}({\cal U}_{\cal T}(\mu);\Bbb{R}^+)$
and $\ga_{\cal T}^{\mu}\!:{\cal FT}_{\de}\lra\bar{\cal V}_1(\mu)$
is an identification of neighborhoods of ${\cal U}_{\cal T}(\mu)$,
which is smooth on the preimage of~${\cal V}_1(\mu)$.
Note that $|\hat{I}|\!\in\!\{1,2\}$
if ${\cal U}_{\cal T}(\mu)$ is nonempty.
By the proof of Lemma~\ref{order1_contr_l3},
$\al_{\cal T}$ has full rank on every 
fiber~$\tilde{\cal F}{\cal T}\!\lra\!{\cal U}_{\cal T}(\mu)$.
Thus, by equation~\e_ref{n2cusps_e3} and Corollary~\ref{euler_crl},
$${\cal C}_{{\cal U}_{\cal T}(\mu)}({\cal D})=0
\qquad\hbox{if~~}H_{\hat{0}}{\cal T}\neq\hat{I}.$$
(3) Suppose $|H_{\hat{0}}{\cal T}|\!=\!|\hat{I}|\!=\!1$.
Then ${\cal T}={\cal T}^*(l)$ for some $l\!\in\![N]$ and
$\tilde{\cal F}{\cal T}\!=\!{\cal F}{\cal T}\!\approx\! 
L_{{\cal T},\hat{1}}$.
Since $\al_{\cal T}\!\circ\!\rho_{\cal T}$ 
has constant rank over $\bar{\cal U}_{\cal T}(\mu)$,
by  Corollary~\ref{euler_crl} and Lemma~\ref{zeros_main},
$${\cal C}_{{\cal U}_{\cal T}(\mu)}({\cal D})
=\big\lan c_1\big(L^*\otimes\ev^*T\PP\big)-
c_1\big(L_{{\cal T},\hat{1}}\big),
\big[\bar{\cal U}_{\cal T}(\mu)\big]\big\ran
=\big\lan 3a+c_1 (L_{{\cal T},\hat{1}}^*),
\big[\bar{\cal U}_{\cal T}(\mu)\big]\big\ran.$$
If $|H_{\hat{0}}{\cal T}|\!=\!|\hat{I}|\!=\!2$,
$\al_{\cal T}\!\circ\!\rho_{\cal T}$ is an isomorphism on every fiber.
Thus, ${\cal C}_{{\cal U}_{\cal T}(\mu)}\!({\cal D})\!=\!
\big|{\cal U}_{\cal T}(\mu)\big|$
by Corollary~\ref{euler_crl}. 
Combining these contributions to the Euler class of
$L^*\otimes\ev^*T\PP$ gives
\begin{equation}\label{n2cusps_e11}\begin{split}
{\cal C}_{\partial\bar{\cal U}}
                       \big({\cal D}_{{\cal T}^*}\big)
&=\sum_{l\in[N]}\big\lan 3a+c_1\big(L_{{\cal T}^*(l),\hat{1}}^*),
\big[{\cal U}_{{\cal T}^*(l)}(\mu)\big]\big\ran+
\sum_{[{\cal T}],|H_{\hat{0}}{\cal T}|=|\hat{I}|=2}
|\bar{\cal U}_{\cal T}(\mu)|\\
&=\sum_{l\in[3d-2]}\big\lan 3a+c_1\big(L_{{\cal T}^*(l),\hat{1}}^*),
\big[{\cal U}_{{\cal T}^*(l)}(\mu)\big]\big\ran
+\tau_2(\mu).
\end{split}\end{equation}
The claim follows by plugging equation~\e_ref{n2cusps_e11}
into~\e_ref{n2cusps_e1}
and using equations~\e_ref{normal_bundle1} and~\e_ref{normal_bundle2}.

\begin{lmm}
\label{n3tangents}
If $d\!\ge\!1$, the number of two-component rational degree-$d$
curves connected at a tacnode and passing through
a tuple $\mu$ of $p$ points and $q$ lines
in general position~$\PPP$, where $2p\!+\!q\!=\!4d\!-\!3$,
is given~by
$$|{\cal S}_2(\mu)|=
\big\lan 6a^2+
4a\big(c_1({\cal L}_1^*)+c_1({\cal L}_2^*)\big)+
\big(c_1^2({\cal L}_1^*)+c_1^2({\cal L}_2^*)\big)+
c_1({\cal L}_1^*)c_1({\cal L}_2^*),
 \big[\bar{\cal V}_2(\mu)\big]\big\ran-3\tau_3(\mu).$$
\end{lmm}

\noindent
{\it Proof:}
(1) Let ${\cal T}^*\!=\!(S^2,[N],I^*;j^*,d^*)$ be a basic 
bubble type such that $I^*\!=\!\{k_1,k_2\}$ is a two-element set,
$d_{k_1}^*,d_{k_2}^*\!>\!0$, 
$d_{k_1}^*\!+\!d_{k_2}^*\!=\!d$, and $N\!=\!p\!+\!q$.
Denote by ${\cal T}_1^*$ and  ${\cal T}_2^*$ the corresponding 
simple types.
The proof is similar to that of Lemma~\ref{n2cusps},
but we pass to the projectivization $\Bbb{P}E$ (over $\Bbb{C}$)
of the bundle
$$E=L_1^*\oplus L_2^*\lra \bar{\cal U}_{{\cal T}^*}(\mu),
\quad\hbox{where}\quad L_i=L_{{\cal T}_i^*,k_i}.$$
The section $\bar{\cal D}_{{\cal T}^*,2}$  of Lemma~\ref{order1_contr_l3}
induces a section ${\cal D}\!\in\!\Ga(\Bbb{P}E;\ga_E^*\otimes\ev^*T\PPP)$
such that ${\cal S}_{{\cal T}^*}(\mu)$ corresponds to the intersection
of the zero set of~${\cal D}$ 
with~$\Bbb{P}E|_{{\cal U}_{{\cal T}^*}(\mu)}$.
If $\Bbb{P}E'$ denotes the restriction of 
$\Bbb{P}E$ to  $\partial\bar{\cal U}\equiv
\bar{\cal U}_{{\cal T}^*}(\mu)-{\cal U}_{{\cal T}^*}(\mu)$,
by Corollary~\ref{euler_crl},
\begin{alignat}{1}\label{n3tangents_e1}
|\bar{\cal S}_{{\cal T}^*}(\mu)|
&=\big\lan c_3\big(\ga_E^*\otimes\ev^*T\PPP\big),
      \big[\Bbb{P}E\big]\big\ran-
{\cal C}_{\Bbb{P}E'}({\cal D})\\
&=\big\lan 6a^2+
4a\big(c_1(L_1^*)+c_1(L_2^*)\big)+
\big(c_1^2(L_1^*)+c_1^2(L_2^*)\big)+c_1(L_1^*)c_1(L_2^*),
\big[\bar{\cal U}_{{\cal T}^*}(\mu)\big]\big\ran-
{\cal C}_{\Bbb{P}E'}({\cal D}).\notag
\end{alignat}
The second equality above is obtained by applying 
\e_ref{zeros_main_e}.\\
(2) Suppose ${\cal T}\!\!=\!(S^2,[N],I;j,d)
\!<\!{\cal T}^*$ is a bubble type such that
${\cal D}$ vanishes somewhere  on~$\Bbb{P}E|_{{\cal U}_{\cal T}(\mu)}$.
Let ${\cal T}_1$ and ${\cal T}_2$ be the corresponding simple types.
Since the constraints are disjoint, up to interchanging the indices,
we must have
$${\cal T}_1={\cal T}_1^*,\qquad
{\cal T}_2=\big(S^2,M_2,I_2;j|M_2,d|I_2\big)<{\cal T}_2^*
\quad\hbox{with}~~d_{k_2}=0.$$
Furthermore, ${\cal D}_{{\cal T}_1^*,k_1}$ does not vanish on 
${{\cal U}_{\cal T}(\mu)}$; see the proof of Lemma~\ref{order1_contr_l3}.
Thus, ${\cal D}$ vanishes only the subspace
$${\cal Z}_{\cal T}\equiv
 \Bbb{P}L_2\big|_{{\cal U}_{\cal T}(\mu)}=
\big\{\big(b,L_2|_b\big)\!: b\!\in\! {\cal U}_{\cal T}(\mu)\big\}.$$
The map $\ga_{\cal T}^{\mu}$ of Theorem~\ref{str_global}
induces an identification of a neighborhood of~$0$~in
$$F{\cal S}\equiv\pi_E^*{\cal FT}\oplus 
     \pi_E^*L_2^*\otimes\pi_E^*L_1\lra{\cal Z}_{\cal T}$$
with a neighborhood of ${\cal Z}_{\cal T}$ in $\Bbb{P}E$.
Similarly to the $n\!=\!2$ case, with appropriate identifications,
\begin{gather}\label{n3tangent_e3}
\Big|{\cal D}\big(\ga_{\cal T}^{\mu}(\ups,u)\big)-
\tilde{\al}_{\cal T}\big(\tilde{\rho}_{\cal T}(\ups,u)\big)\Big|\le 
C(b_{\ups})|\ups|^{\frac{1}{p}}\big|\rho_{\cal T}(\ups)\big|
\quad\forall(\ups,u)\in {\cal FS}_{\de},\\
\hbox{where}\quad
\tilde{\rho}_{\cal T}(\ups,u)=
\big(\rho_{\cal T}(\ups),u\big)\in
\tilde{F}{\cal S}\equiv\pi^*\tilde{\cal F}{\cal T}\oplus 
 \pi_E^*L_2^*\otimes\pi_E^*L_1\lra{\cal Z}_{\cal T},\notag
\end{gather}
and $\tilde{\al}_{\cal T}$ has full rank on every fiber 
by~\e_ref{str_global_e} and Lemma~\ref{order1_contr_l3}.
Thus, similarly to the proof of Lemma~\ref{n2cusps}, 
and ${\cal C}_{\Bbb{P}E'|_{{\cal Z}_{\cal T}}}({\cal D})=0$
if $H_{k_2}{\cal T}\!\neq\!\hat{I}_2$,
and only two cases remain to be considered.\\
(3) If $|H_{k_2}{\cal T}|\!=\!|\hat{I}_2|\!=\!1$,
$\tilde{\al}_{\cal T}\circ\tilde{\rho}_{\cal T}$ 
has full rank over all of~$\bar{\cal Z}_{\cal T}$.
Thus, by Corollary~\ref{euler_crl} and Lemma~\ref{zeros_main},
\begin{equation}\label{n3tangents_e5a}
{\cal C}_{{\cal Z}_{\cal T}}({\cal D})
=\big\lan c_1(\ga_E^*\otimes\ev^*T\PPP)-c_1(F{\cal S}),
\big[\bar{\cal Z}_{\cal T}\big]\big\ran
= \big\lan 4a+ c_1(L_{{\cal T}_2,\hat{1}}^*)+c_1(L_1^*),
\big[\bar{\cal U}_{\cal T}(\mu)\big]\big\ran;
\end{equation}
note that $c_1(\ga_E^*)\!=\!c_1(L_2^*)\!=\!0$
over $\bar{\cal U}_{\cal T}(\mu)$.
If $|H_{k_2}{\cal T}|\!=\!|\hat{I}_2|\!=\!2$,
$\tilde{\al}_{\cal T}\circ\tilde{\rho}_{\cal T}$ is an isomorphism on every
fiber, and thus
\begin{equation}\label{n3tangents_e5b}
{\cal C}_{\Bbb{P}E'|_{{\cal Z}_{\cal T}}}({\cal D})
=|{\cal Z}_{\cal T}|=|{\cal U}_{\cal T}(\mu)|.
\end{equation}
Note that the sum of $|{\cal U}_{\cal T}(\mu)|$ over all equivalence
classes of bubble types ${\cal T}^*$ and~${\cal T}\!<\!{\cal T}^*$
is~$3\tau_3(\mu)$, since one of the three components 
of the image of each bubble map in ${\cal U}_{\cal T}(\mu)$
is distinguished by the bubble type~${\cal T}$.
As before, we now sum up equations~\e_ref{n3tangents_e5a} 
and~\e_ref{n3tangents_e5b} over all equivalence classes of bubble 
types ${\cal T}\!<\!{\cal T}^*$ of the appropriate form,
plug the result back into~\e_ref{n3tangents_e1} and
use equations~\e_ref{normal_bundle1} and~\e_ref{normal_bundle2}.
The claim follows by summing the result over all 
equivalence classes of basic simple bubble types~${\cal T}^*$.

\subsection{The Numbers $n_m^{(1)}$ with $m=n-1$}

\noindent
In this subsection, we give topological formulas for the numbers
$n_{\cal T}^{(1)}$ with $|\hat{I}|\!=\!n\!-\!1$.
As before, the reason these two cases are similar is that
the complex dimension of ${\cal U}_{\cal T}(\mu)$ is two.

\begin{lmm}
\label{n2m1c1}
If $n=2$, $n_1^{(1)}(\mu)=
2\big\lan 6a^2+3ac_1({\cal L}^*),
\big[\bar{\cal V}_1(\mu)\big]\big\ran$.
\end{lmm}

\noindent
{\it Proof:} (1) 
Let $N$, ${\cal T}^*$, $L$, ${\cal D}$ be as in the proof of 
Lemma~\ref{n2cusps}.
Since $s_{\Si}$ does not vanish on~$\Si$,
by equation~\e_ref{summary_m1c1} and Lemma~\ref{zeros_main},
\begin{equation}\label{n2m1c1_e1}\begin{split}
n_1^{(1)}(\mu)
&= \sum_{k=0}^{k=3}\big\lan 
 c_k({\cal O})c_1^{3-k}\big(T^*\Si\otimes L^*\big),
   \big[\Si\times\bar{\cal U}_{{\cal T}^*}(\mu)\big]\big\ran
-{\cal C}_{\Si\times{\cal D}^{-1}(0)}(\al^{\perp}),\\
&=2 \big\lan 15a^2+12ac_1(L^*)+3c_1^2(L^*),
\big[\bar{\cal V}_1(\mu)\big]\big\ran
-{\cal C}_{\Si\times{\cal D}^{-1}(0)}(\al^{\perp}),
\end{split}\end{equation}
where ${\cal O}\!=\!{\cal H}^{0,1}_{\Si}\!\otimes\ev^*T\PP$ and
$\al\!\in\!\Ga(\Si\!\times\!\bar{\cal V}_1(\mu);
T^*\Si\otimes L^*\otimes{\cal O})$ is the linear part
of the affine map~$\psi_1^{(1)}$  of~\e_ref{summary_m1c1}.\\
(2) We first compute 
${\cal C}_{\Si\times{\cal S}_1(\mu)}(\al^{\perp})$.
Since ${\cal V}_1(\mu)$ is a complex manifold
and ${\cal D}$ is transverse to the zero set in $L^*\otimes\ev^*T\PP$
by Corollary~\ref{reg_crl2},
we can identify a neighborhood of~$0$ in
$$F\equiv L^*\otimes\ev^*T\PP\lra {\cal S}_1(\mu)$$
with a neighborhood of ${\cal S}_1(\mu)$ in 
${\cal V}_1(\mu)$ via a map $\ga$ in such a way that
\begin{equation}\label{n2m1c1_e3}
\Pi_{b,\ga(b,X)}^{-1}\big({\cal D}\ga(b,X)\big)=X
\quad\forall (b,X)\in F_{\de}.
\end{equation}
Then with appropriate identifications, 
$$\al^{\perp}\big(\ga(X)\big)=
\pi^{\perp}\circ Xs_{\Si}\equiv\al_{\cal S}(X),$$
where $\pi^{\perp}\!:{\cal O}\!\lra\!{\cal O}^{\perp}$
is the quotient projection map.
In particular, $\al_{\cal S}$ has full rank if 
\hbox{$\bar{\nu}_b\not\in{\cal H}_{\Si}^+\otimes\ev^*T\PP$}
for all $b\!\in\!{\cal S}_1(\mu)$,
i.e.~$\nu$ is generic.
Furthermore, 
$$\big(T^*\Si\otimes L^*\otimes {\cal O}^{\perp}\big)
\big/\big(\hbox{Im}~\al_{\cal S}\big)\approx
T^*\Si\otimes
  \big(\big({\cal H}_{\Si}^-\otimes\Bbb{C}^2\big)/\Bbb{C}\big).$$
Thus, by Corollary~\ref{euler_crl},
\begin{equation}\label{n2m1c1_e6}
{\cal C}_{\Si\times{\cal S}_1(\mu)}(\al^{\perp})
=\big\lan e\big(T^*\Si^{\otimes3}\big),[\Si\!\times\!{\cal S}_1(\mu)]\big\ran
=6|{\cal S}_1(\mu)|.
\end{equation}
(3) It remains to compute the contribution to 
${\cal C}_{\Si\times{\cal D}^{-1}(0)}(\al^{\perp})$ from
$\Si\!\times\!\big(\bar{\cal V}_1(\mu)\!-\!{\cal V}_1(\mu)\big)$.
Suppose 
$${\cal T}=(S^2,[N],I;j,d)<{\cal T}^*$$ 
is a bubble type such that
${\cal D}$ vanishes somewhere  on~${\cal U}_{\cal T}(\mu)$.
As in the proof of Lemma~\ref{n2cusps}, 
$|\hat{I}|\!\in\!\{1,2\}$ and $d_{\hat{0}}\!=\!0$.
Furthermore,
$$\Big| \al^{\perp}\big(x,\ga_{\cal T}^{\mu}(\ups)\big)
-\tilde{\al}_{\cal T}\big(x,b;\rho_{\cal T}(\ups)\big)\Big|\le
C(b_{\ups})|\ups|^{\frac{1}{p}}|\rho_{\cal T}(\ups)|\quad
(x,b;\ups)\!\in\! {\cal FT}_{\de},$$
where 
$\tilde{\al}_{\cal T}\!=\!\pi^{\perp}\circ(s_{\Si}\otimes\al_{\cal T})$.
If $\nu$ is generic, $\tilde{\al}_{\cal T}$ 
has full rank on every fiber,
since $s_{\Si}$ has no zeros.
Thus, by Corollary~\ref{euler_crl},
$${\cal C}_{\Si\times{\cal U}_{\cal T}(\mu)}(\al^{\perp})=0
\qquad\hbox{if~~}H_{\hat{0}}{\cal T}\neq\hat{I}.$$
If $|H_{\hat{0}}{\cal T}|\!=\!|\hat{I}|\!=\!1$,
$\tilde{\al}_{\cal T}\circ\tilde{\rho}_{\cal T}$ has full rank over all 
$\Si\!\times\!\bar{\cal U}_{\cal T}(\mu)$, and thus
by Corollary~\ref{euler_crl} and Lemma~\ref{zeros_main}
\begin{equation}\label{n2m1c1_e9}
{\cal C}_{\Si\times{\cal U}_{\cal T}(\mu)}(\al^{\perp})\!=\!
\big\lan 
c(T^*\Si\otimes {\cal O}^{\perp})c(L_{{\cal T},\hat{1}})^{-1},
\big[\Si\!\times\!\bar{\cal U}_{\cal T}(\mu)\big]\big\ran
\!=\!2\big\lan 12a+3c_1(L_{{\cal T},\hat{1}}^*),   
\big[\bar{\cal U}_{\cal T}(\mu)\big]\big\ran.
\end{equation}
If $|H_{\hat{0}}{\cal T}|\!=\!|\hat{I}|\!=\!2$,
$$\big(T^*\Si\otimes L^*\otimes {\cal O}^{\perp}\big)
\big/\big(\hbox{Im}~\tilde{\al}_{\cal T}\circ\tilde{\rho}_{\cal T}\big)
\approx
T^*\Si\otimes
  \big(\big({\cal H}_{\Si}^-\otimes\Bbb{C}^2\big)/\Bbb{C}\big).$$
Thus, similarly to the computation in (2) above,
\begin{equation}\label{n2m1c1_e7}
{\cal C}_{\Si\times{\cal U}_{\cal T}(\mu)}(\al^{\perp})
=6|{\cal U}_{\cal T}(\mu)|.
\end{equation}
Summing equations~\e_ref{n2m1c1_e7} and~\e_ref{n2m1c1_e9} over all 
equivalence classes of ${\cal T}\!<\!{\cal T}^*$, we obtain
\begin{equation}\label{n2m1c1_e11}
{\cal C}_{\Si\times(\bar{\cal V}_1(\mu)-{\cal V}_1(\mu))}
(\al^{\perp})=
2\sum_{l\in[N]}\big\lan 12a+3c_1(L_{{\cal T},\hat{1}}^*),
\big[{\cal U}_{{\cal T}^*(l)}(\mu)\big]\big\ran
+6\tau_2(\mu).
\end{equation}
The claim follows by plugging \e_ref{n2m1c1_e6} and \e_ref{n2m1c1_e11}
into \e_ref{n2m1c1_e1} and using \e_ref{normal_bundle1},
\e_ref{normal_bundle2}, and Lemma~\ref{n2cusps}.

\begin{lmm}
\label{n3m2c1}
If $n\!=\!3$, 
$n_2^{(1)}(\mu)= 4\big\lan 10a^2
+ 4a\big(c_1({\cal L}_1^*)+c_1({\cal L}_2^*)\big)+
c_1({\cal L}_1^*)c_1({\cal L}_2^*),
 \big[\bar{\cal V}_2(\mu)\big]\big\ran$.
\end{lmm}

\noindent
{\it Proof:} (1) The proof is essentially a mixture of the proofs
of Lemma~\ref{n3tangents} and~\ref{n2m1c1}.
We continue with the notation of the proof of Lemma~\ref{n3tangents},
but take
$$E=T\Si_1\otimes L_1\oplus T\Si_2\otimes L_2.$$
If ${\cal T}'\!\!=\!\!(\Si,[N],I';j',d')$
is the simple bubble type such that $\bar{\cal T}'\!\!=\!\!{\cal T}$,
by equation~\e_ref{summary_m2c1} and Lemma~\ref{zeros_main},
\begin{alignat}{1}\label{n3m2c1_e1}
&n_{{\cal T}'}^{(1)}(\mu)
= \sum_{k=0}^{k=5}\big\lan  c_k({\cal O})\la_E^{5-k},
  \big[\Bbb{P}E\big]\big\ran
-{\cal C}_{\al_E^{-1}(0)}(\al_E^{\perp}),\\
&\quad=4 \big\lan 28a^2\!+\!16a\big(c_1(L_1^*)\!+\!c_1(L_2^*)\big)
\!+\!3\big(c_1^2(L_1^*)\!+\!c_1^2(L_2^*)\big)\!+\!4c_1(L_1^*)c_1(L_2^*),
\big[\bar{\cal V}_1(\mu)\big]\big\ran
-{\cal C}_{\al_E^{-1}(0)}(\al_E^{\perp}),\notag
\end{alignat}
where ${\cal O}\!=\!{\cal H}^{0,1}_{\Si}\!\otimes\ev^*T\PPP$ and
$\al\!\in\!\Ga(\Si_{\hat{1}}\!\times\Si_{\hat{2}}\!\times\!
 \bar{\cal V}_1(\mu);E^*\otimes{\cal O})$ is the linear part
of the affine map~$\psi_{{\cal T}'}^{(1)}$  of~\e_ref{summary_m2c1}.
Let
\begin{gather*}
\Si^{(\pm)}=\big\{ 
 (x_1,x_2)\!\in\!\Si_1^*\!\times\!\Si_2^*\!: 
x_1\!=\!\pm x_2\},\quad
\Si^{(0)}=\big\{(z_m,z_m)\!: m\!\in\![6]\big\};\\
{\cal S}_{{\cal T}^*}^{(\pm)}=
\Si^{(\pm)}\times{\cal S}_{{\cal T}^*}(\mu),\quad
{\cal S}_{{\cal T}^*}^{(0)}=
\Si^{(0)}\times{\cal S}_{{\cal T}^*}(\mu),
\end{gather*}
where $+x_2\equiv x_2$ and 
$-x_2$ is the image of $x_2$
under the nontrivial automorphism of~$\Si$.
The zero set of $\al_E$ is the union of a section of $\Bbb{P}E$
over ${\cal S}_{{\cal T}^*}^{(\pm)}$,
${\cal S}_{{\cal T}^*}^{(0)}$,
and $\Si^2\!\times\!{\cal U}_{\cal T}(\mu)$, where 
${\cal T}$ is as in the proof of Lemma~\ref{n3tangents}.\\
(2) The above section over $\Si^2\!\times\!{\cal U}_{\cal T}(\mu)$ 
is given by
$${\cal Z}_{\cal T}\equiv
\al_E^{-1}(0)\cap\big(\Bbb{P}E|\Si^2\!\times\!{\cal U}_{\cal T}(\mu)\big)=
\big\{ \big(x_1,x_2,b, T_{x_2}\Si_2\otimes L_2|_b\big)\!:
(x_1,x_2,b)\!\in\!\Si_1\!\times\!\Si_2\!\times\!{\cal U}_{\cal T}(\mu)
\big\}.$$
Let ${\cal FS}\!=\!\pi_E^*{\cal FT}\!\oplus\!
T^*\Si_2\otimes L_2^*\otimes T\Si_1\otimes L_1$ and 
$\tilde{\cal F}{\cal S}\!=\!\pi_E^*\tilde{\cal F}{\cal T}\!\oplus\!
T^*\Si_2\otimes L_2^*\otimes T\Si_1\otimes L_1$.
As in the proof of Lemma~\ref{n3tangents},
the map~$\ga_{\cal T}^{\mu}$ of Theorem~\ref{str_global}
induces a homeomorphism
between neighborhoods of ${\cal Z}_{\cal T}$ in ${\cal FS}$
and in $\Bbb{P}E$ such that, with appropriate identifications,
$$\Big|\al_E\big(x_1,x_2,\ga_{\cal T}^{\mu}(\ups)\big)
-\tilde{\al}_{\cal T}(x_1,x_2,b;\tilde{\rho}_{\cal T}(\ups)\big)\Big|\le
C(b_{\ups})|\ups|^{\frac{1}{p}}|\rho_{\cal T}(\ups)|\quad
(x_1,x_2,b;\ups)\!\in\! {\cal FS}_{\de}.$$
The section $\tilde{\al}_{\cal T}\!\in\!\Ga\big({\cal Z}_{\cal T};
\hbox{Hom}(\tilde{\cal F}{\cal S};\ga_E^*\otimes{\cal O}^{\perp})\big)$
has full rank on every fiber as long as 
$\bar{\nu}_b\!\not\in\!{\cal H}_{\Si}^+\otimes\ev^*T\PPP$
for all $b\!\in\!{\cal U}_{\cal T}(\mu)$,
i.e.~if $\nu$ is generic.
This can be seen from Lemma~\ref{order1_contr_l3}
in the same way as in the proof of Lemma~\ref{n2m1c1}.
Thus, as in the proof of Lemma~\ref{n3tangents},
$${\cal C}_{{\cal Z}_{\cal T}}(\al_E^{\perp})=0
\qquad\hbox{if~~} H_{k_2}{\cal T}|\neq\hat{I}.$$
If $|H_{k_2}{\cal T}|\!=\!|\hat{I}|\!=\!1$, 
i.e.~${\cal T}\!=\!{\cal T}^*(l)$ for some $l\!\in\![N]$,
$\tilde{\al}_{\cal T}\!\circ\!\tilde{\rho}_{\cal T}$ has full rank
on $\Si\!\times\!\bar{\cal U}_{\cal T}(\mu)$.
Thus, by Corollary~\ref{euler_crl} and Lemma~\ref{zeros_main},
\begin{equation}\label{n3m2c1_e3}
{\cal C}_{{\cal Z}_{\cal T}}(\al_E^{\perp})
= \big\lan c\big(\ga_E^*\otimes{\cal O}^{\perp}\big)
c\big({\cal FS})^{-1},[\bar{\cal Z}_{\cal T}]\big\ran
=4\big\lan 16a+4c_1(L_1^*)+3c_1(L_{{\cal T},\hat{1}}^*),
\big[\bar{\cal U}_{\cal T}(\mu)\big]\big\ran,
\end{equation}
since $F{\cal S}\approx L_{{\cal T},\hat{1}}\oplus 
T^*\Si_2\otimes T\Si_1\otimes L_1$.
If $|H_{k_2}{\cal T}|\!=\!|\hat{I}|\!=\!2$, we similarly obtain
\begin{equation}\label{n3m2c1_e5}
{\cal C}_{{\cal Z}_{\cal T}}(\al_E^{\perp})
=\big\lan c\big(\ga_E^*\otimes{\cal O}^{\perp}\big)
c\big({\cal FS})^{-1},[\bar{\cal Z}_{\cal T}]\big\ran
=12|{\cal U}_{\cal T}(\mu)|
\end{equation}
Note that $F{\cal S}\approx\Bbb{C}^2\oplus T^*\Si_2\otimes T\Si_1$
in this case
Summing up equations~\e_ref{n3m2c1_e3} and~\e_ref{n3m2c1_e5}
over all equivalence classes of bubble types ${\cal T}\!<\!{\cal T}^*$
of the appropriate form, we obtain
\begin{equation}\label{n3m2c1_e7}
{\cal C}_{\Bbb{P}E|\partial{\cal V}_2(\mu)}(\al_E^{\perp})\!=\!
4\!\!\sum_{l\in M_{k_j},i\neq j} \!\!\!\!\!
 \big\lan 16a+4c_1(L_i^*)+3c_1(L_{{\cal T},\hat{1}}^*),
\big[{\cal U}_{{\cal T}^*(l)}(\mu)\big]\big\ran\! +12\!\!\!
\sum_{{\cal T}<{\cal T}^*,
          |H_{k_i}{\cal T}|=2}\!\!\!\!\!\!|{\cal U}_{\cal T}(\mu)|.
\end{equation}
(3) It remains to compute
${\cal C}_{\Bbb{P}E|{\cal S}_{{\cal T}^*}^{\pm}}(\al_E^{\perp})$
and ${\cal C}_{{\cal S}_{{\cal T}^*}^{(0)}}(\al_E^{\perp})$. 
Let $V_1$ denote the orthogonal complement of the image of
${\cal D}_{{\cal T}^*,k_1}$ in~$\ev^*T\P$. 
Analogously to the $n\!=\!2$ case, we can identify a neighborhood of~$0$ in
$$F\equiv L_2^*\otimes V_1
  \lra {\cal S}_{{\cal T}^*}(\mu)$$
with a neighborhood of ${\cal S}_{{\cal T}^*}(\mu)$ in 
${\cal U}_{{\cal T}^*}(\mu)$ via a map $\ga$ in such a way that
\begin{equation}\label{n3m2c1_e15}
\pi_{V_1}\Pi_{b,\ga(b,X)}^{-1}
\big({\cal D}_{{\cal T}^*,k_2}\ga(b,X)\big)=X
\quad\forall (b,X)\in F_{\de}.
\end{equation}
Let $\ka\!\in\!\Ga({\cal S}_{{\cal T}^*}(\mu); L_2^*\otimes L_1)$
be the nonvanishing section defined in Subsection~\ref{order2_case2a}.
Then 
\begin{gather*}
\al_E^{-1}(0)\cap\big(\Bbb{P}E|{\cal S}_{{\cal T}^*}^{(0)}\big)=
{\cal Z}_{{\cal T}^*}^{(0)}\equiv
\big\{ \big(z_m,z_m,b,[v_1,v_2]\big)\!\in\!
 \Bbb{P}E|_{{\cal S}^{(0)}_{{\cal T}^*}}\!:
                              v_1\!+\!\ka^{\pm}(b)v_2\!=\!0\big\},\\
\al_E^{-1}(0)\cap\big(\Bbb{P}E|{\cal S}_{{\cal T}^*}^{(\pm)}\big)=
{\cal Z}_{{\cal T}^*}^{(\pm)}\equiv
\big\{ \big(x_1,x_2,b,[v_1,v_2]\big)\!\in\!
 \Bbb{P}E|_{{\cal S}^{(\pm)}_{{\cal T}^*}}\!:
                              v_1\!+\!\ka^{(\pm)}(b)v_2\!=\!0\big\},
\end{gather*}
where $\ka^+\!=\!\ka$ is viewed as a section of 
$T^*\Si_2\otimes L_2^*\otimes T^*\Si_1\otimes L_1$ 
along ${\cal S}^{(+)}_{{\cal T}^*}$
and $\ka^-$ is the section of 
$T^*\Si_2\otimes L_2^*\otimes T^*\Si_1\otimes L_1$ 
along ${\cal S}^{(-)}_{{\cal T}^*}$ induced by $\ka$ and 
the differential of the nontrivial automorphism of~$\Si$.
See Subsection~\ref{order2_case1} for more details.
Let 
$$F{\cal Z}_{{\cal T}^*}^{\pm} =\pi_E^*\big({\cal FT}\oplus T\Si_2\big)
\oplus \pi_E^* \big( T^*\Si_2\otimes L_2^*\otimes T\Si_1 \otimes L_1\big)
\approx \Bbb{C}^2\oplus \pi_E^* T\Si_2 \oplus \Bbb{C}
\lra{\cal Z}_{{\cal T}^*}^{\pm}.$$
The above map $\ga$ and the $g_{b,\hat{0}}$-exponential map in $\Si$
induce a diffeomorphism between neighborhoods
of ${\cal Z}_{{\cal T}^*}^{\pm}$ in 
$F{\cal Z}_{{\cal T}^*}^{\pm}$ and in~$\Bbb{P}E$.
With appropriate identifications, the linearization of 
$\al_E\!\in\!\Ga(\Bbb{P}E;\ga_E^*\otimes{\cal O})$
near ${\cal Z}_{{\cal T}^*}^{\pm}$ has the~form
\begin{equation}\label{n3m2c1_e30}
\tilde{\al}_{{\cal T}^*}^{(\pm)}(X,w,u)\!=\!
\big(u({\cal D}_{{\cal T}^*,k_1}b) s_{\Si,x_2}+
({\cal D}_{{\cal T}^*,k_2}b)s_{x_2}^{(2,+)}(w,\cdot)
+Xs_{\Si,x_2}\big)+
  ({\cal D}_{{\cal T}^*,k_2}b) s_{x_2}^{(2,-)}(w,\cdot),
\end{equation}
where $s_{x_2}^{(2,+)}\!\in
 \!\Ga(\Si;T^*\Si\otimes{\cal H}_{\Si}^{0,1})$ 
is defined with respect to the metric~$g_{b,\hat{0}}$; 
see Subsection~\ref{dbar_section}.
Since $s^{(2,-)}$ does not vanish on $\Si^*$,
it follows that  $\tilde{\al}_{{\cal T}^*}^{(\pm)}$ has
full rank on~$F{\cal Z}_{{\cal T}^*}^{\pm}$.
The same is true for the linearization of 
$\al_E^{\perp}$  as long as $\nu$ is generic.
Thus, by Corollary~\ref{top_l2c},
${\cal C}_{{\cal Z}_{{\cal T}^*}^{\pm}}(\al_E^{\perp})$
is given by Lemma~\ref{zeros_main}~with
$$\bar{\cal M}=\bar{\cal Z}_{{\cal T}^*}^{(\pm)}
\approx\Si\times{\cal S}_{{\cal T}^*}(\mu),\quad
E_2=T\Si,\quad
{\cal O}_2=\big(\ga_E^*\otimes{\cal O}^{\perp}\big)/\Bbb{C}^3
\approx \big(T^*\Si\otimes\Bbb{C}^5\big)/\Bbb{C}^3,$$
and $\al_2\!\in\!\Ga(\bar{\cal M};E_2^*\otimes{\cal O}_2)$,
which is the composition of $\tilde{\al}_{{\cal T}^*}^{\pm}$
with the projection 
onto
${\cal O}/\big({\cal H}_{\Si}^+\otimes\Bbb{C}^3\!\oplus\Bbb{C}\bar{\nu}\big)$.
However, similarly to the construction in Subsection~\ref{comp_sec2},
we can replace $E_2$ with 
$$\tilde{E}_2=\tilde{T}'\Si\equiv 
T\Si\otimes{\cal O}(z_1)\otimes\ldots\otimes{\cal O}(z_6)$$
and $s^{(2,-)}$ with 
$\tilde{s}^{(2,-)}\!\in\!\Ga(\Si;\tilde{T}'\Si^*\otimes{\cal O}_2)$
in~\e_ref{n3m2c1_e30} without changing the number of zeros.
Then the corresponding section $\al_2$ has no zero
on $\bar{\cal Z}_{{\cal T}^*}^{\pm}$, 
and thus by Lemma~\ref{zeros_main},
\begin{equation}\label{n3m2c1_e9}
{\cal C}_{{\cal Z}_{{\cal T}^*}^{(+)}\cup{\cal Z}_{{\cal T}^*}^{(-)}} 
(\al_E^{\perp})
=2\big\lan\!-c_1(\tilde{E}_2)+c_1({\cal O}_2), [\Si]\big\ran
=2(-4+10)|{\cal S}_{{\cal T}^*}(\mu)|
=12|{\cal S}_{{\cal T}^*}(\mu)|.
\end{equation}
(4) We next show that 
${\cal C}_{{\cal Z}_{{\cal T}^*}^{(0)}}(\al_E^{\perp})\!=\!0$
The normal bundle of ${\cal Z}_{{\cal T}^*}^{(0)}$
in $\Bbb{P}E$ is given by
$$F{\cal Z}_{\cal T}^{(0)}=\pi_E^*
\big({\cal FT}\oplus T_{z_m}\Si_1\oplus T_{z_m}\Si_2\big)
\oplus\pi_E^*\big(T_{z_m}\Si_2^*\otimes L_2^*
\otimes T_{z_m}\Si_1\otimes L_1\big)
\approx \Bbb{C}^2\oplus\big(\Bbb{C}\oplus\Bbb{C}\big)\oplus\Bbb{C}.$$
With appropriate identifications, 
$\al_E\!\in\!\Ga(\Bbb{P}E;\ga_E^*\otimes{\cal O})$
at $(X,w_1,w_2,u)\!\in\! F{\cal Z}_{{\cal T}^*}^{(0)}$ sufficiently small
is given~by
$$\al_E(X,w_1,w_2,u)=
(u-\ka(b))\big({\cal D}_{{\cal T}^*,k_1}\ga(b,X)\big)s_{\Si,w_1}+
\big({\cal D}_{{\cal T}^*,k_2}\ga(b,X)\big)s_{\Si,w_2}
+Xs_{\Si,w_2}.$$
For any $(X,w_1,w_2,u)\!\in\! F{\cal Z}_{{\cal T}^*}^{(0)}$
with $w_1,w_2\!\in\! T_{z_m}\Si$, let
\begin{gather}
\label{n3m2c1_e21}
\al_{{\cal T}^*}^{(0)}(X,w_1,w_2,u)=
\Big( u\big({\cal D}_{{\cal T}^*,k_1}b\big)s_{\Si,w_2}+
\big({\cal D}_{{\cal T}^*,k_2}b\big)
  \big( s_{\Si,w_2}-s_{\Si,w_1}\big)\Big)+Xs_{\Si,w_2};\\
\al_{{\cal T}^*,1}^{(0,-)}(w_1,w_2)=
\big({\cal D}_{{\cal T}^*,k_2}b\big)s_{z_m}^{(3,-)}(w_1,w_1-w_2,\cdot),
~~
\al_{{\cal T}^*,2}^{(0,-)}(w_1,w_2)=
\big({\cal D}_{{\cal T}^*,k_2}b\big)s_{z_m}^{(3,-)}(w_2,w_2-w_1,\cdot).
\notag
\end{gather}
Since $\big({\cal D}_{{\cal T}^*,k_2}b)=\ka(b)
({\cal D}_{{\cal T}^*,k_1}b)$,
\begin{equation}\label{n3m2c1_e22}
\Big|\al_E(X,w_1,w_2,u)-\al_{{\cal T}^*}^{(0)}(w_1,w_2)\Big|
\le C |(w_1,w_2,u)||X|.
\end{equation}
By the same computation as in the proof of Lemma~\ref{psit3l_2b},
\begin{gather}\label{n3m2c1_e23}
\Big|\pi_{w_1}^-\al_{{\cal T}^*}^{(0)}(X,w_1,w_2,u)
-\al_{{\cal T}^*,2}^{(0,-)}(w_1,w_2)\Big|
\le C\big(|w_1|^2+|w_2|^2\big)|w_1-w_2|;\\
\Big|\pi_{w_2}^-\al_{{\cal T}^*}^{(0)}(X,w_1,w_2,u)
-\al_{{\cal T}^*,1}^{(0,-)}(w_1,w_2)\Big|
\le C\big(|w_1|^2+|w_2|^2\big)|w_1-w_2|.
\end{gather}
Since the two terms in \e_ref{n3m2c1_e21} lie in orthogonal
subspaces, from~\e_ref{n3m2c1_e22} and~\e_ref{n3m2c1_e23},
we obtain
\begin{equation}\label{n3m2c1_e25}\begin{split}
&\Big|\al_E(X,w_1,w_2,u)-\big(\pi_{w_1}^+\al_{{\cal T}^*}^{(0)}(X,w_1,w_2,u)
+\al_{{\cal T}^*,2}^{(0,-)}(w_1,w_2)\big)\Big|\\
&\qquad\qquad\qquad\qquad\qquad\qquad
\le |(w_1,w_2,u)|\big|\pi_{w_1}^+\al_{{\cal T}^*}^{(0)}(X,w_1,w_2,u)
+\al_{{\cal T}^*,2}^{(0,-)}(w_1,w_2)\big|.
\end{split}\end{equation}
If $\nu$ is generic, inequality~\e_ref{n3m2c1_e25} continues
to hold with all the maps replaced by their composition
with the projection to~${\cal O}^{\perp}$.
Note that the rank of the bundle ${\cal O}_2$ defined in~(3) is~two,
while the dimension of ${\cal Z}_{{\cal T}^*}^{(0)}$ is zero.
Thus, applying Lemma~\ref{top_l1}, we conclude that
${\cal C}_{{\cal Z}_{{\cal T}^*}^{(0)}}(\al_E^{\perp})\!=\!0$.
The claim follows by plugging 
\e_ref{n3m2c1_e7} and \e_ref{n3m2c1_e9} into \e_ref{n3m2c1_e1}, 
using \e_ref{normal_bundle1} and \e_ref{normal_bundle2}, 
then summing over all equivalence classes of bubble types
${\cal T}^*$ of appropriate form, and applying 
and Lemma~\ref{n3tangents}.

\subsection{Behavior of ${\cal D}^{(2)}$ and
${\cal D}^{(3)}$ near
$\bar{\cal S}_1(\mu)-{\cal S}_1(\mu)$}
\label{cuspcurves}

\noindent
If $n\!=\!3$, the space ${\cal S}_1(\mu)$ s not compact.
In order to be able to compute the numbers $n_1^{(k)}(\mu)$, 
we thus must understand the structure
of $\bar{\cal S}_1(\mu)$ as well as the behavior of 
${\cal D}_{{\cal T}^*,\hat{0}}^{(2)}$ and
${\cal D}_{{\cal T}^*,\hat{(0)}}^{(3)}$,
where ${\cal T}^*\!=\!(S^2,[N],\{\hat{0}\};\hat{0},d)$,
\hbox{near $\bar{\cal S}_1(\mu)\!-\!{\cal S}_1(\mu)$}.
\\

\noindent
If ${\cal T}\!=\!(S^2,[N],I;j,d)\!<\!{\cal T}^*$,
from Theorem~\ref{str_global} one should expect that
the normal bundle 
of \hbox{${\cal S}_{\cal T}(\mu)\equiv {\cal U}_{\cal T}(\mu)
                            \cap\bar{\cal S}_1(\mu)$}
in~$\bar{\cal S}_1(\mu)$ to be given by
$$F{\cal S}=\Big\{\big[\ups\!=\!\big(b,(v_h)_{h\in\hat{I}}\big)\big]
\!\in\! {\cal FT}\big|{\cal S}_{\cal T}(\mu)\!:
\sum_{\chi_{\cal T}h=1}
\Big(\prod_{i\in\hat{I},h\in\bar{D}_i{\cal T}}v_i\Big)
\big({\cal D}_{{\cal T},h}b\big)=0\Big\}.$$
The next lemma shows that this is indeed the case.
Let ${\cal NS}\!\lra\! F{\cal S}$ denote the normal bundle
of  $F{\cal S}$ in ${\cal FT}\!\lra\! {\cal U}_{\cal T}(\mu)$.
While for the purposes of Lemma~\ref{cuspcurves_l1},
we can use any identification of neighborhoods of $F{\cal S}$ 
in ${\cal NS}$ and in ${\cal FT}\!\lra\!{\cal U}_{\cal T}(\mu)$,
in order to simplify the statement of Lemma~\ref{cuspcurves_l2},
we choose a fairly natural one.
More precisely, denote by $F{\cal S}^{\perp}$ 
a subspace of \hbox{${\cal FT}\lra {\cal S}_{\cal T}(\mu)$} 
complementary to~$F{\cal S}$ 
and by 
\hbox{$\pi_{\cal S}\!:{\cal NS}^{(1)}\!\lra\!  {\cal S}_{\cal T}(\mu)$}
the normal bundle of ${\cal S}_{\cal T}(\mu)$ in~${\cal U}_{\cal T}(\mu)$.
Choose a norm on ${\cal NS}^{(1)}$ and an identification
$\phi_{\cal S}\!: {\cal NS}_{\de}^{(1)}\!\lra\!{\cal S}_{\cal T}(\mu)$
of neighborhoods of ${\cal S}_{\cal T}(\mu)$ in ${\cal NS}^{(1)}$ 
and in~${\cal U}_{\cal T}(\mu)$.
Let 
$\Phi_{\cal S}\!: \pi_{\cal S}^*{\cal FT}\!\lra\! {\cal FT}$
be a lift of $\phi_{\cal S}$ such that 
$\Phi_{\cal S}$ restricts to the identity over 
${\cal S}_{\cal T}(\mu)\!\subset\! {\cal NS}_{\de}^{(1)}$.
Let $\pi\!: {\cal FT}\!\lra\! {\cal S}_{\cal T}(\mu)$
be the bundle projection.
Then
$${\cal NS}=\pi^*{\cal NS}^{(1)}\oplus F{\cal S}^{\perp},
\quad\hbox{and}\quad
\tilde{\phi}_{\cal S}\!: {\cal NS}_{\de}\lra {\cal FT},~~
\tilde{\phi}_{\cal S}\big((b,v),(X,v^{\perp})\big)=
\Phi_{\cal S}\big((b,X),v+v^{\perp}\big),$$
is  an identification of neighborhoods of $F{\cal S}$
in ${\cal NS}$ and \hbox{${\cal FT}\lra {\cal U}_{\cal T}(\mu)$}.

\begin{lmm}
\label{cuspcurves_l1}
For every bubble type 
${\cal T}\!=\!(S^2,[N],I;j,d)\!<\!{\cal T}^*$,
there exist $\de,C\!\in\! C^{\i}({\cal S}_{\cal T}(\mu);\Bbb{R}^+)$
and a section
$\varphi_{\cal S}\!\in\!\Ga( F{\cal S}_{\de};{\cal NS})$
such that 
$$\|\varphi_{\cal S}(\ups)\|\le C(b_{\ups})|\ups|^{\frac{1}{p}},
\qquad
\|\varphi_{F{\cal S}^{\perp}}(\ups)\|\le C(b_{\ups})|\ups|^{1+\frac{1}{p}},$$
where $\varphi_{F{\cal S}^{\perp}}$ denotes 
the $F{\cal S}^{\perp}$-component of $\varphi_{\cal S}$,
and the map
$$\ga_{\cal S}\!: F{\cal S}_{\de}\lra{\cal S}_1(\mu),\qquad
\ga_{\cal S}(\ups)=\ga_{\cal T}^{\mu}
\big(\tilde{\phi}_{\cal S}\varphi_{\cal S}(\ups)\big),$$
is a homeomorphism onto an open neighborhood of ${\cal S}_{\cal T}(\mu)$
in $\bar{\cal S}_1(\mu)$, which is smooth and orientation-preserving
on the preimage of ${\cal U}_{{\cal T}^*}(\mu)$.
\end{lmm}

\noindent
{\it Proof:} 
(1) The proof is similar to that of Lemma~\ref{gl-rat_l1} in~\cite{Z},
so we only describe the differences.
If ${\cal S}_{\cal T}(\mu)\!\neq\!\eset$,
${\cal T}$  must have one of the three forms described
by Lemma~\ref{cuspcurves_l2}.
We restrict the proof to those three cases.
In Case (1), we apply Subsection~\ref{gl-impl} in~\cite{Z},
which contains an application of the Implicit Function Theorem,  
to ${\cal D}_{{\cal T}^*,\hat{0}}$ instead of the evaluation~maps.
By Theorem~\ref{str_global}, 
$$\Big|\Pi_{b,\tilde{\ga}_{\cal T}^{\mu}(\ups)}^{-1}
\big({\cal D}_{{\cal T}^*,\hat{0}}\tilde{\ga}_{\cal T}^{\mu}(\ups)\big)-
\big({\cal D}_{{\cal T}^*,\hat{0}}b_{\ups}\big)\Big|
\le C'(b)|\ups|^{\frac{1}{p}}
\qquad\forall\ups\!\in\! F{\cal S}_{\de}.$$
This estimate suffices for applying an argument 
similar to the proof of~Lemma~\ref{gl-rat_l1}.\\
(2) In Case (2) of Lemma~\ref{cuspcurves_l2}, 
instead of the section ${\cal D}_{{\cal T}^*,\hat{0}}$
of $L_{{\cal T}^*,\hat{0}}^*\otimes\ev^*T\P$, we consider 
the section $\tilde{\cal D}$ of 
$\big(L_{{\cal T}^*,\hat{0}}\otimes {\cal FT}\big)^*\otimes\ev^*T\P$
on a neighborhood of ${\cal U}_{\cal T}(\mu)$
in $\bar{\cal U}_{{\cal T}^*}(\mu)$ defined by
$$\tilde{\cal D}\big|_{\tilde{\ga}_{\cal T}^{\mu}(b,v_{\hat{1}})}
\big(v_{\hat{0}},v_{\hat{1}}\big)=
{\cal D}_{{\cal T}^*,\hat{0}}\big|_{\tilde{\ga}_{\cal T}^{\mu}(b,v_{\hat{1}})}
(v_{\hat{0}})\in\ev^*T\P.$$
This section is well-defined outside of ${\cal U}_{\cal T}(\mu)$
and by Theorem~\ref{str_global} extends over ${\cal U}_{\cal T}(\mu)$~by
$$\tilde{\cal D}\big|_b\big(v_{\hat{0}},v_{\hat{1}}\big)=
v_{\hat{0}}v_{\hat{1}}\big({\cal D}_{{\cal T},\hat{1}}b\big).$$
The restriction of this section to ${\cal U}_{\cal T}(\mu)$
vanishes transversally at ${\cal S}_{\cal T}(\mu)$
by Corollary~\ref{reg_crl2},
while its zero set on ${\cal U}_{{\cal T}^*}(\mu)$
is the same as the zero set of ${\cal D}_{{\cal T}^*,\hat{0}}$.
By Theorem~\ref{str_global}, with appropriate identifications,
$$\Big|\tilde{\cal D}\big|_{\ga_{\cal T}^{\mu}(b,v_{\hat{1}})}
-\tilde{\cal D}\big|_b\Big|
\le C'(b)|\ups|^{\frac{1}{p}}.$$
(3) In the final case of Lemma~\ref{cuspcurves_l2}, 
we replace ${\cal D}_{{\cal T}^*,\hat{0}}$  by a bundle section
over the blowup of ${\cal FT}$ along ${\cal U}_{\cal T}(\mu)$.
Let 
$$\Om_{\cal T}=\big\{(b,v,\ell)\!: (b,v)\!\in\! {\cal FT},~ 
v\!\in\!\ell\!\in\!\Bbb{P}{\cal FT}|_b\big\},\quad
\Om_{\cal T}^*=\big\{(b,v,\ell)\!\in\!\Om_{\cal T}\!: v\!\neq\!0\big\},\quad
{\cal E}_{\cal T}=\Om_{\cal T}\!-\Om_{\cal T}^*.$$
Denote by $\ga\!\lra\!\Om_{\cal T}$ the tautological line bundle.
The normal bundle $\tilde{\cal N}{\cal S}$ of 
$\ga\!\lra\!\Bbb{P}F{\cal S}$ in $\ga\!\lra\!{\cal E}_{\cal T}$
is given~by
$$\tilde{\cal N}{\cal S}=
\pi_{\ga}^*\pi_{{\cal FT}}^*{\cal NS}^{(1)}\oplus
\pi_{\ga}^*\big(\ga^* \otimes\pi_{{\cal FT}}^*F{\cal S}^{\perp}\big),
\quad
\tilde{\phi}_{\tilde{\cal N}{\cal S}}
\big((b,\ell,v),X,\si\big)=
\big(\phi_{\cal S}(b,X),\big[\Phi_{\cal S}(v+\si(v))\big],v+\si(v)\big),$$
where $\pi_{\ga}\!:\ga\!\lra\! \Om_{\cal T}$ is the bundle projection map.
The bundle $L_{{\cal T}^*,\hat{0}}$ pulls back to a bundle~$\tilde{L}$ 
over a neighborhood $\Om_{\de}$ of 
${\cal E}_{\cal T}$ in $\Om_{\cal T}$.
We define a section $\tilde{\cal D}$ of 
$\big(\tilde{L}\otimes\ga)^*\otimes\ev^*T\P$ over $\Om_{\de}$ by
$$\tilde{\cal D}\big|_{(b,v_{\hat{1}},v_{\hat{2}},\ell)}
\big(v_{\hat{0}},v_{\hat{1}},v_{\hat{2}}\big)=
{\cal D}_{{\cal T}^*,\hat{0}}
\big|_{\tilde{\ga}_{\cal T}^{\mu}(b,v_{\hat{1}},v_{\hat{2}})}
(v_{\hat{0}})\in\ev^*T\P.$$
This section is well-defined outside of ${\cal E}_{\cal T}(\mu)$
and by Theorem~\ref{str_global} extends over ${\cal E}_{\cal T}(\mu)$~by
$$\tilde{\cal D}\big|_b\big(v_{\hat{0}},v_{\hat{1}},v_{\hat{2}}\big)=
v_{\hat{0}}
\Big( v_{\hat{1}}\big({\cal D}_{{\cal T},\hat{1}}b\big)
+v_{\hat{2}}\big({\cal D}_{{\cal T},\hat{2}}b\big)\Big).$$
The restriction of this section to ${\cal E}_{\cal T}(\mu)$
vanishes transversally at $\Bbb{P}F{\cal S}\!\lra\!{\cal S}_{\cal T}(\mu)$
by Corollary~\ref{reg_crl2},
while its zero set on $\Om^*$ corresponds to
the zero set of ${\cal D}_{{\cal T}^*,\hat{0}}$ on 
$\ga_{\cal T}^{\mu}\big({\cal FT}_{\de}
                        \!-\!{\cal U}_{\cal T}(\mu)\big)$.
By Theorem~\ref{str_global}, with appropriate identifications,
$$\Big|\tilde{\cal D}\big|_{(b,v_{\hat{1}},v_{\hat{2}},\ell)}
-\tilde{\cal D}\big|_{(b,\ell)}\Big|
\le C'(b)|\ups|^{\frac{1}{p}}.$$
Thus, we can apply the arguments of Lemma~\ref{gl-rat_l1} 
in~\cite{Z} to $\tilde{\cal D}$ to describe 
its zero set near~${\cal E}_{\cal T}$.
We obtain a section 
$\tilde{\varphi}_{\cal S}\!\in\!\Ga( \ga_{\de}|\Bbb{P}F{\cal S};
\tilde{\cal N}{\cal S})$
such that
$\|\tilde{\varphi}_{\cal S}(\ups)\|\!\le\! C(b_{\ups})|\ups|^{\frac{1}{p}}$,
and the~map
$$\tilde{\ga}_{\cal S}\!:\ga_{\de}|\Bbb{P}F{\cal S}\lra\ups_{\cal T},
\quad
\tilde{\ga}_{\cal S}(\ups)=\tilde{\phi}_{\tilde{\cal N}{\cal S}}
\big(\varphi_{\cal S}(\ups)\big),$$
is a homeomorphism onto an open neighborhood of 
$\Bbb{P}F{\cal S}$ in ${\cal D}^{-1}(0)$.
This section $\tilde{\varphi}_{\cal S}$ induces 
the required section $\varphi_{\cal S}$
 with the claimed properties.

\begin{crl}
\label{cuspcurves_l1c}
For every bubble type 
${\cal T}\!=\!(S^2,[N],I;j,d)\!<\!{\cal T}^*$,
there exist $\de\!\in\! C^{\i}({\cal S}_{\cal T}(\mu);\Bbb{R}^+)$
and a map
$$\ga_{\cal S}: \big({\cal NS}^{(1)}\oplus {\cal FT}\big)_{\de}
\big|{\cal S}_{\cal T}(\mu)\lra
\bar{\cal U}_{{\cal T}^*}(\mu)$$
such that $\ga_{\cal S}$ is a homeomorphism onto 
an open neighborhood of ${\cal S}_{\cal T}(\mu)$
in $\bar{\cal U}_{{\cal T}^*}(\mu)$, 
which is smooth and orientation-preserving
on  the preimage of ${\cal U}_{{\cal T}^*}(\mu)$,
and with appropriate identifications,
$$\Pi_{b,\ga_{\cal S}(\ups)}^{-1}
{\cal D}_{{\cal T}^*,\hat{0}}= 
\begin{cases}
X,&\hbox{in Case (1)};\\
v_{\hat{1}}X,&\hbox{in Case (2)};\\
v_{\hat{2}}X+v_{\hat{1}}{\cal D}_{{\cal T},\hat{1}}
+v_{\hat{2}}{\cal D}_{{\cal T},\hat{2}} ,&\hbox{in Case (3)},
\end{cases}$$
where the three cases are the ones described by 
Lemma~\ref{cuspcurves_l2}.
\end{crl}

\noindent
{\it Proof:}
See Subsections~\ref{order2_case1} and~\ref{order3_case1}
for the meaning of $X\!\in\!{\cal NS}^{(1)}$ in the three case.
The proof is just a modification of the proof of 
Lemma~\ref{cuspcurves_l1}.
In all three cases, we work with the same section $\tilde{\cal D}$
as before.
Denote by ${\cal E}_{\cal T}$ the space
${\cal U}_{\cal T}(\mu)$ in the first two cases.
Let $\tilde{\cal S}\!=\!{\cal S}_{\cal T}(\mu)$,
\hbox{$\tilde{F}{\cal S}\!=\!{\cal N}{\cal S}^{(1)}\!\oplus\! {\cal FT}$}
in those two cases
and $\tilde{\cal S}\!=\!\Bbb{P}F{\cal S}$,
$\tilde{F}{\cal S}\!=\!\tilde{\cal N}{\cal S}\!\oplus\!\ga$
in the remaining case.
Choose an identification 
$\ga\!: {\cal NS}^{(1)}_{\de}\!\lra\!{\cal E}_{\cal T}$
 of neighborhoods of $\tilde{\cal S}$
in ${\cal NS}^{(1)}$ and in ${\cal E}_{\cal T}$
as well as of the appropriate line bundle over these
neighborhoods such that
\begin{gather*}
\Pi_{b,\ga_{\cal S}(X)}^{-1}\tilde{\cal D}=X\quad
\forall X\in {\cal NS}^{(1)}_{\de}\qquad\hbox{in Cases (1),(2)};\\
\Pi_{b,\ga_{\cal S}(X)}^{-1}\tilde{\cal D}=
X+\big({\cal D}_{{\cal T},\hat{1}}+{\cal D}_{{\cal T},\hat{2}}\big)
    \circ (1+u)\quad \forall (X,u)\in 
\big({\cal NS}^{(1)}\oplus\ga^*\otimes F{\cal S}^{\perp}\big)_{\de},
\qquad\hbox{in Case (3)}.
\end{gather*}
Note that in Cases (2) and (3), the restriction of 
$L_{{\cal T}^*,\hat{0}}^*$ to ${\cal U}_{\cal T}(\mu)$
is trivial.
By the same argument as in the proof of Lemma~\ref{cuspcurves_l1},
for any $Y\!\in\! \tilde{F}{\cal S}$ sufficiently small,
there exists a unique $Z\!\in\! N{\cal S}^{(1)}$ in Cases 1,2
and $Z\!\in\! \tilde{N}{\cal S}^{(1)}$ in Case~3,
such that 
$$\Pi_{b,\tilde{\ga}_{\cal T}^{\mu}\ga_{\cal S}(Y+Z)}^{-1}\tilde{\cal D}
\big|_{\tilde{\ga}_{\cal T}^{\mu}\ga_{\cal S}(Y+Z)}=
\Pi_{b,\ga_{\cal S}(Y)}^{-1}
\tilde{\cal D}\big|_{\ga_{\cal S}(Y)}.$$
Furthermore, $|Z|\le C(b)|Y|$.

\begin{lmm}
\label{cuspcurves_l2}
If $d\!\ge\!1$, $\mu$ is a tuple of $p$ points and $q$ lines
in general position in $\PPP$ with \hbox{$2p\!+\!q\!=\!4d\!-\!3$},
and $N\!=\!p\!+\!q$, the set $\bar{\cal S}_1(\mu)\!-\!{\cal S}_1(\mu)$,
is finite.
Furthermore, if 
$${\cal T}=(S^2,[N],I;j,d)<{\cal T}^*
\quad\hbox{and}\quad {\cal S}_{\cal T}(\mu)\neq\eset,$$
(1) $\hat{I}\!=\!\{\hat{1}\}$, $d_{\hat{0}}\!>\!0$, and
the images of ${\cal D}_{{\cal T}^*,\hat{0}}^{(2)}$ and 
${\cal D}_{{\cal T}^*,\hat{0}}^{(3)}$ are linearly independent
in every fiber of $\ev^*T\P$ over ${\cal S}_{\cal T}(\mu)$;\\
(2) OR $\hat{I}\!=\!\{\hat{1}\}$,  $d_{\hat{0}}\!=\!0$, $d_{\hat{1}}\!=\!d$, 
and  for all $\ups\!=\!\big[b,v_{\hat{1}}\big]\!\in\! F{\cal S}_{\de}$,
\begin{gather*}
\Big|\Pi_{b,\ga_{\cal S}(\ups)}^{-1}
\big({\cal D}^{(2)}\ga_{\cal S}(\ups)\big)-
v_{\hat{1}}^2\big({\cal D}_{{\cal T},\hat{1}}^{(2)}b\big)\Big|
\le C|v_{\hat{1}}|^{2+\frac{1}{p}};\\
\Big|\Pi_{b,\ga_{\cal S}(\ups)}^{-1}
\Big( \big({\cal D}^{(3)}\ga_{\cal S}(\ups)\big)
-3x_{\hat{1}}\big({\cal D}^{(2)}\ga_{\cal S}(\ups)\big)\Big)
-v_{\hat{1}}^3\big({\cal D}_{{\cal T},\hat{1}}^{(3)}b\big)\Big|
\le C|v_{\hat{1}}|^{3+\frac{1}{p}};
\end{gather*}
(3) OR $\hat{I}\!=\!\{\hat{1},\hat{2}\}$,  $d_{\hat{0}}\!=\!0$, and
for all $\ups\!=\![b,v_{\hat{1}},v_{\hat{2}}]\!\in\! F{\cal S}$
\begin{gather*}
\Big|\Pi_{b,\ga_{\cal S}(\ups)}^{-1}
\big({\cal D}^{(2)}\ga_{\cal S}(\ups)\big)-
2\Big( x_{\hat{1}}v_{\hat{1}}\big({\cal D}_{{\cal T},\hat{1}}b\big)
+x_{\hat{2}}v_{\hat{2}}\big({\cal D}_{{\cal T},\hat{2}}b\big)\Big)
\Big|\le C|\ups|^{1+\frac{1}{p}};\\
\Big|\Pi_{b,\ga_{\cal S}(\ups)}^{-1}
 \Big( 2\big({\cal D}^{(3)}\ga_{\cal S}(\ups)\big)
\!-\!3(x_{\hat{1}}\!+\!x_{\hat{2}})
\big({\cal D}^{(2)}\ga_{\cal S}(\ups)\big)\Big)
\!-\!3\big(x_{\hat{1}}\!-\!x_{\hat{2}})
\Big( v_{\hat{1}}^2\big({\cal D}_{{\cal T},\hat{1}}^{(2)}b\big)
\!+\!v_{\hat{2}}^2\big({\cal D}_{{\cal T},\hat{2}}^{(2)}b\big)\Big)
\Big|\!\le\! C|\ups|^{2+\frac{1}{p}}.
\end{gather*}
\end{lmm}

\noindent
{\it Proof:} (1) The statement about the possible structures of ${\cal T}$
is easily seen from Theorem~\ref{str_global} and dimension count.
The finiteness claim then also follows by dimension count.
In Case (1), if $d_{\hat{0}}\!\ge\!3$, by Corollary~\ref{reg_crl2},
the images of
${\cal D}_{{\cal T}^*,\hat{0}}^{(2)}$ and 
${\cal D}_{{\cal T}^*,\hat{0}}^{(3)}$ are transversal
and thus linearly independent over the finite set~${\cal S}_{\cal T}(\mu)$.
On the other hand, if $d_{\hat{0}}\!<\!3$,
${\cal S}_{\cal T}(\mu)=\eset$; see Subsection~\ref{order3_case1}.\\
(2) The four inequalities in the lemma
will be obtained by refining the proof
of the analytic estimate of Theorem~\ref{str_global}.
We use the same notation.
Combining equations~\e_ref{limit_lmm_e1},
\e_ref{limit_lmm_e2}, \e_ref{limit_lmm_e3}, and~\e_ref{limit_lmm_e5}, 
we~obtain
\begin{equation}\label{cusp_curves_e1}
\big({\cal D}^{(m)}\tilde{\ga}_{\cal T}(\ups)\big)=
m\sum_{\chi_{\cal T}h=1}\sum_{k=1}^{k=m}\frac{a_{k,h}(\ups)}{k}
\tilde{v}^k\big({\cal D}_{{\cal T},h}^{(k)}b\big)
-\frac{m}{2\pi i}\sum_{\chi_{\cal T}h=1}
\int_{A_h^-(\ups)}\xi_{\ups} w^{m-1}dw,
\end{equation}
where the integral is computed by using the same trivializations
as before.
This equality holds for any bubble type.
If $\ga_{\cal T}(\ups)\!\in\!{\cal S}_1(\mu)$ and
${\cal T}$ is as in (2) of the lemma,
\e_ref{cusp_curves_e1} with $m\!=\!1,2,3$ gives
\begin{gather}
\label{cusp_curves_e2a}
0=v_{\hat{1}}\big({\cal D}_{{\cal T},\hat{1}}^{(1)}b\big)
-\frac{1}{2\pi i}\int_{|x_{\hat{1}}-w|=\ep}\xi_{\ups} dw;\\
\label{cusp_curves_e2b}
\big({\cal D}^{(2)}\tilde{\ga}_{\cal T}(\ups)\big)=
2x_{\hat{1}}v_{\hat{1}}\big({\cal D}_{{\cal T},\hat{1}}^{(1)}b\big)
+v_{\hat{1}}^2\big({\cal D}_{{\cal T},\hat{1}}^{(2)}b\big)
-\frac{1}{\pi i}\int_{|x_{\hat{1}}-w|=\ep}\xi_{\ups} wdw;\\
\label{cusp_curves_e2c}
\big({\cal D}^{(3)}\tilde{\ga}_{\cal T}(\ups)\big)=
3x_{\hat{1}}^2v_{\hat{1}}\big({\cal D}_{{\cal T},\hat{1}}^{(1)}b\big)
+3x_{\hat{1}}v_{\hat{1}}^2\big({\cal D}_{{\cal T},\hat{1}}^{(2)}b\big)
+2v_{\hat{1}}^3\big({\cal D}_{{\cal T},\hat{1}}^{(3)}b\big)
-\frac{3}{2\pi i}\int_{|x_{\hat{1}}-w|=\ep}\xi_{\ups} w^2dw.
\end{gather}
where $\ep\!=\!4\de(b_{\ups})^{-1}|v_{\hat{1}}|$.
Subtracting $2x_{\hat{1}}$ times the first equation
from the second, we obtain
\begin{equation}\label{cusp_curves_e3}
\Big|\big({\cal D}^{(2)}\tilde{\ga}_{\cal T}(\ups)\big)-
v_{\hat{1}}^2\big({\cal D}_{{\cal T},\hat{1}}^{(2)}b\big)\Big|
\le C(b)|v_{\hat{1}}|^{2+\frac{1}{p}}.
\end{equation}
Similarly, subtracting $3x_{\hat{1}}$ times~\e_ref{cusp_curves_e2b} from
and adding $3x_{\hat{1}}^2$ times~\e_ref{cusp_curves_e2a} 
to~\e_ref{cusp_curves_e2c} , we obtain
\begin{equation}\label{cusp_curves_e5}
\Big|\Big( \big({\cal D}^{(3)}\tilde{\ga}_{\cal T}(\ups)\big)
-3x_{\hat{1}}\big({\cal D}^{(2)}\tilde{\ga}_{\cal T}(\ups)\big)\Big)
-2v_{\hat{1}}^3\big({\cal D}_{{\cal T},\hat{1}}^{(3)}b\big)\Big|
\le C(b)|v_{\hat{1}}|^{3+\frac{1}{p}}.
\end{equation}
If $\ups\!\in\! F{\cal S}$ is sufficiently small,
the claim in Case (2) follows from 
equations~\e_ref{cusp_curves_e3} and~\e_ref{cusp_curves_e5}
along with Lemma~\ref{cuspcurves_l1} and our choice 
of~$\tilde{\phi}_{\cal S}$.
Note that if $\ups\!\in\! F{\cal S}$, we have to apply
\e_ref{cusp_curves_e3} and \e_ref{cusp_curves_e5} with
$\ups$ replaced by 
$ \Phi_{\cal T}^{\mu} 
 \varphi_{\cal T}^{\mu}\tilde{\phi}_{\cal S}\varphi_{\cal S}(\ups)$,
where $\Phi_{\cal T}^{\mu}$ and
$\varphi_{\cal T}^{\mu}$ are as in 
Subsection~\ref{gl-orient1_sec} of~\cite{Z}.
However, applying the bounds on 
$\varphi_{\cal T}^{\mu}$ and $\tilde{\phi}_{\cal S}$,
we obtain the claimed estimates.\\
(3) In case (3), we proceed similarly.
The analog of equation~\e_ref{cusp_curves_e2b} gives
$$\Big|\big({\cal D}^{(2)}\tilde{\ga}_{\cal T}(\ups)\big)-
2\Big( x_{\hat{1}}v_{\hat{1}}\big({\cal D}_{{\cal T},\hat{1}}b\big)
+x_{\hat{2}}v_{\hat{2}}\big({\cal D}_{{\cal T},\hat{2}}b\big)\Big)
\Big|\le C(b)|\ups|^{1+\frac{1}{p}}.$$
Subtracting $3(x_{\hat{1}}+x_{\hat{2}})$ times
the analog of~\e_ref{cusp_curves_e2b} from
and adding $6x_{\hat{1}}x_{\hat{2}}$ times
the analog of~\e_ref{cusp_curves_e2a} to
twice the analog of~\e_ref{cusp_curves_e2c}, we obtain
$$\Big|\Big( 2\big({\cal D}^{(3)}\tilde{\ga}_{\cal T}(\ups)\big)
-3(x_{\hat{1}}+x_{\hat{2}})
\big({\cal D}^{(2)}\tilde{\ga}_{\cal T}(\ups)\big)\Big)
-3\big(x_{\hat{1}}-x_{\hat{2}})
\Big( v_{\hat{1}}^2\big({\cal D}_{{\cal T},\hat{1}}^{(2)}b\big)
+v_{\hat{2}}^2\big({\cal D}_{{\cal T},\hat{2}}^{(2)}b\big)\Big)
\Big|\le C(b)|\ups|^{2+\frac{1}{p}}.$$
The estimates of Case 3 follow from 
the last two equations and Lemma~\ref{cuspcurves_l1}.
The finer bound on $\varphi_{F{\cal S}^{\perp}}$
of Lemma~\ref{cuspcurves_l1} is essential~here.

\subsection{The Numbers $n_1^{(2)}(\mu)$ and $n_1^{(3)}(\mu)$
in the $n=3$ Case}

\noindent
In this subsection, we express the numbers 
$n_1^{(2)}(\mu)$ and $n_1^{(3)}(\mu)$ in the $n\!=\!3$ case
in terms of intersection numbers on the spaces
$\bar{\cal V}_1(\mu)$, 
$\bar{\cal V}_2(\mu)$, $\bar{\cal V}_3(\mu)$.

\begin{lmm}
\label{n3m1c2}
If $n\!=\!3$,
$n_1^{(2)}(\mu)=
 4\lan 2a+c_1({\cal L}^*),[\bar{\cal S}_1(\mu)]\ran-2|{\cal S}_2(\mu)|$.
\end{lmm}

\noindent
{\it Proof:} (1) We continue with the notation of the previous subsection.
The number $n_1^{(2)}(\mu)$ is the number
of zeros of the affine map in~\e_ref{summary_m1c2}.
As in the proof of Lemma~\ref{mequalsnm1}, 
we can replace $s_{\Si}^{(2,-)}$ by~$\tilde{s}_{\Si}^{(2,-)}$.
Since the linear part of the new affine map does not vanish
on $\Si\!\times\!{\cal S}_1(\mu)$ (see Subsection~\ref{order2_case1}),  
by Lemma~\ref{zeros_main},
\begin{equation}\label{n3m1c2_e1}\begin{split}
n_1^{(2)}(\mu)&=\sum_{k=0}^{k=2}
\big\lan c_1^{2-k}\big(\tilde{T}\Si^*\otimes L^{*\otimes2}\big)c_k({\cal O}),
\big[\Si\times\bar{\cal S}_1(\mu)\big]\big\ran 
-{\cal C}_{\Si\times\partial\bar{\cal S}_1}(\al^{\perp})\\
&=4\big\lan 2a+c_1(L^*),\big[\bar{\cal S}_1(\mu)\big]\big\ran
-{\cal C}_{\Si\times\partial\bar{\cal S}_1}(\al^{\perp}),
\end{split}\end{equation}
where ${\cal O}\!=\!{\cal H}_{\Si}^-\otimes\ev^*T\PPP$,
$\partial\bar{\cal S}_1\!=\!\bar{\cal S}_1(\mu)\!-\!{\cal S}_1(\mu)$,
and $\al$ is the linear part of the affine map in~\e_ref{summary_m1c2},
with $s_{\Si}^{(2,-)}$ replaced by~$\tilde{s}_{\Si}^{(2,-)}$.\\
(2) If ${\cal T}\!\!=\!(S^2,[N],I;j,d)\!<\!\!{\cal T}^*$
and ${\cal S}_{\cal T}(\mu)\!\neq\!\eset$,
${\cal T}$ must have one of the three forms 
given by Lemma~\ref{cuspcurves_l2}.
Since ${\cal D}^{(2)}$ does not vanish on 
${\cal S}_{\cal T}(\mu)$ in Case (1) of Lemma~\ref{cuspcurves_l2},
${\cal C}_{\Si\times{\cal S}_{\cal T}(\mu)}(\al^{\perp})=0$
in this case.
In Case (2), i.e.~${\cal T}\!=\!{\cal T}^*(l)$ for some $l\!\in\![N]$,
${\cal D}_{{\cal T},\hat{1}}^{(2)}$ does not vanish
over ${\cal S}_{\cal T}(\mu)$; see Subsection~\ref{order3_case1}.
Thus, by Corollaries~\ref{euler_crl}, \ref{top_l2c}, and
the first estimate of Lemma~\ref{cuspcurves_l2},  
${\cal C}_{\Si\times{\cal S}_{\cal T}(\mu)}(\al^{\perp})$
is {\it twice} the number of Lemma~\ref{zeros_main} corresponding
$$\bar{\cal M}=\Si\times{\cal S}_{\cal T}(\mu),\quad
E_2={\cal FT}^{\otimes2}\approx\Bbb{C},\quad
{\cal O}_2=\tilde{T}\Si^*\otimes L^*\otimes{\cal O}^{\perp}
\approx \tilde{T}\Si^*\otimes{\cal O}^{\perp},$$
and $\al_2\!\in\!\Ga(\bar{\cal M};E_2^*\otimes{\cal O}_2)$
that has full rank on every fiber.
It follows that
\begin{equation}\label{n3m1c2_e2}
{\cal C}_{\Si\times{\cal S}_{\cal T}(\mu)}(\al^{\perp})
=2\big\lan \!c_1({\cal O}_2) -c_1(E_2),
\big[\Si\!\times\!{\cal S}_{\cal T}(\mu)\big]\big\ran
=4|{\cal S}_{\cal T}(\mu)|.
\end{equation}
(3) Suppose ${\cal T}$ is as in Case (3) of Lemma~\ref{cuspcurves_l2}.
Since $x_{\hat{1}}\!\neq\! x_{\hat{2}}$, ${\cal D}_{{\cal T},\hat{1}}$
and ${\cal D}_{{\cal T},\hat{2}}$ do not vanish on 
${\cal S}_{\cal T}(\mu)\!=\!{\cal U}_{\cal T}(\mu)\cap{\cal S}_2(\mu)$
(see Subsection~\ref{order2_case2a}),
and ${\cal D}_{{\cal T},\hat{1}}\!+\!{\cal D}_{{\cal T},\hat{2}}$
vanishes on $F{\cal S}$,
$x_{\hat{1}}{\cal D}_{{\cal T},\hat{1}}\!+\!
x_{\hat{2}}{\cal D}_{{\cal T},\hat{2}}$ does not 
vanish on~$F{\cal S}$.
Thus, the third estimate of Lemma~\ref{cuspcurves_l2}, 
Corollary~\ref{euler_crl}, and Lemma~\ref{zeros_main},
\begin{equation}\label{n3m1c2_e3}
{\cal C}_{\Si\times{\cal S}_{\cal T}(\mu)}(\al^{\perp})
=\big\lan \!c_1({\cal O}_2)-c_1(E),
\big[\Si\times{\cal S}_{\cal T}(\mu)\big]\big\ran
=2|{\cal S}_{\cal T}(\mu)|.
\end{equation}
Summing up equations~\e_ref{n3m1c2_e2} 
and~\e_ref{n3m1c2_e3} over all appropriate
bubble types ${\cal T}<{\cal T}^*$
and substituting the result into~\e_ref{n3m1c2_e1},
we obtain the claim.

\begin{lmm}
\label{n3m1c3}
If $n\!=\!3$, 
$n_1^{(3)}(\mu)=
  \lan 4a+5c_1({\cal L}^*),[\bar{\cal S}_1(\mu)]\ran
             -3|{\cal S}_2(\mu)|$.
\end{lmm}

\noindent
{\it Proof:} (1) We continue with the notation of Lemma~\ref{n3m1c2}.
The number $n_1^{(3)}(\mu)$ is the number
of zeros of the affine map in~\e_ref{summary_m1c3}.
Let 
$$E=L^{\otimes2}\oplus L^{\otimes3}\lra\bar{\cal S}_1(\mu).$$
Since the linear part $\al$ of the affine map 
has full rank on ${\cal S}_1(\mu)$ (see Subsection~\ref{order3_case1}),
\begin{equation}\label{n3m1c3_e1}
n_1^{(3)}(\mu)
=\sum_{k=0}^{k=2}
\big\lan \la_E^{2-k}c_k({\cal O}),\big[\Bbb{P}E\big]\big\ran 
-{\cal C}_{\Bbb{P}E|\partial\bar{\cal S}_1}(\al_E^{\perp})
=\big\lan 4a+5c_1(L^*),\big[\bar{\cal S}_1(\mu)\big]\big\ran 
-{\cal C}_{\Bbb{P}E|\partial\bar{\cal S}_1}(\al_E^{\perp}),
\end{equation}
where ${\cal O}=\ev^*T\PPP$.\\
(2) As in the proof of Lemma~\ref{n3m1c2},
${\cal C}_{\Bbb{P}E|{\cal S}_{\cal T}(\mu)}(\al_E^{\perp})
\!=\!0$ for bubble types 
${\cal T}$ of Case (1) of Lemma~\ref{cuspcurves_l2}.
Suppose ${\cal T}\!=\!{\cal T}^*(l)$ for some $l\!\in\![N]$,
i.e.~we are in Case (2) of Lemma~\ref{cuspcurves_l2}.
The normal bundle of $\Bbb{P}E|{\cal S}_{\cal T}(\mu)$
in $\Bbb{P}E$ is~$\pi_E^*{\cal FT}\approx\Bbb{C}$.
By the first two estimates of Lemma~\ref{cuspcurves_l2},
with appropriate identifications,
$$\Big|\al_E^{\perp}\big(\ga_{\cal S}(b,v_{\hat{1}})\big)
-\tilde{\al}_{\cal T}(b,v_{\hat{1}})\Big|\le C|v_{\hat{1}}|^{2+\frac{1}{p}}
\qquad\forall (b,v_{\hat{1}})\!\in\! {\cal FT}_{\de},$$
for some 
$\tilde{\al}_{\cal T}\!\in\!\Ga(\Bbb{P}E|{\cal S}_{\cal T}(\mu);
         {\cal FT}^{*\otimes2}\otimes\ga_E^*\otimes{\cal O}^{\perp})$ 
which vanishes only on
$${\cal Z}_{\cal T}\equiv
\big\{(b,[v,w])\!\in\!\Bbb{P}E|{\cal S}_{\cal T}(\mu)\!: 
3x_{\hat{1}}v\!+\!w\!=\!0\big\}.$$
Thus, by Corollaries~\ref{euler_crl} and~\ref{top_l2c}
and Lemma~\ref{zeros_main},
$${\cal C}_{\Bbb{P}E|({\cal S}_{\cal T}(\mu)-{\cal Z}_{\cal T})}
(\al_E^{\perp})
=2\Big( \big\lan \!c_1(\ga_E^*\otimes{\cal O}^{\perp})-c_1({\cal FT}),
\big[\Bbb{P}E|{\cal S}_{\cal T}(\mu)\big]\big\ran-
{\cal C}_{{\cal Z}_{\cal T}}(\tilde{\al}_{\cal T}^{\perp})\Big)
=4|{\cal S}_{\cal T}(\mu)|-
         2{\cal C}_{{\cal Z}_{\cal T}}(\al_2^{\perp}).$$
By the first two estimates of Lemma~\ref{cuspcurves_l2},
${\cal C}_{{\cal Z}_{\cal T}}
(\tilde{\al}_{\cal T}^{\perp})\!=\!|{\cal Z}_{\cal T}|
\!=\!|{\cal S}_{\cal T}(\mu)|$.
Since the images of ${\cal D}_{{\cal T},\hat{1}}^{(2)}$
and ${\cal D}_{{\cal T},\hat{1}}^{(3)}$ are linearly
independent in every fiber of $\ev^*T\P$ over ${\cal S}_{\cal T}(\mu)$,
by the first two estimates of Lemma~\ref{cuspcurves_l2}
and Corollary~\ref{euler_crl},
${\cal C}_{{\cal Z}_{\cal T}}(\al_E^{\perp})\!=\!3|{\cal Z}_{\cal T}|.$
Thus, 
\begin{equation}\label{n3m1c3_e2}
{\cal C}_{\Bbb{P}E|{\cal S}_{\cal T}(\mu)}(\al_E^{\perp})=
\big(4|{\cal S}_{\cal T}(\mu)|-2|{\cal S}_{\cal T}(\mu)|\big)+
3|{\cal S}_{\cal T}(\mu)|=5|{\cal S}_{\cal T}(\mu)|.
\end{equation}
(3) Suppose ${\cal T}$ is as in Case (3).
By the last two estimates of Lemma~\ref{cuspcurves_l2},
with appropriate identifications,
$$\Big|\al_E^{\perp}\big(\ga_{\cal S}(b,v)\big)
-\tilde{\al}_{\cal T}(b,v)\Big|\le C|v_{\hat{1}}|^{1+\frac{1}{p}}
\qquad\forall (b,v)\!\in\! F{\cal S}_{\de},$$
for some 
$\tilde{\al}_{\cal T}\!\in\!\Ga(\Bbb{P}E|{\cal S}_{\cal T}(\mu);
         F{\cal S}^*\otimes\ga_E^*\otimes{\cal O}^{\perp})$ 
which vanishes only on
$${\cal Z}_{\cal T}\equiv
\{(b,[v,w])\!\in\!\Bbb{P}E|{\cal S}_{\cal T}(\mu)\!: 
3(x_{\hat{1}}+x_{\hat{2}})v\!+2\!w\!=\!0\}.$$
Thus, by Corollary~\ref{euler_crl} and Lemma~\ref{zeros_main},
$${\cal C}_{\Bbb{P}E|({\cal S}_{\cal T}(\mu)-{\cal Z}_{\cal T})}
(\al_E^{\perp})
=\big\lan \!c_1(\ga_E^*\otimes{\cal O}^{\perp})-c_1(F{\cal S}),
\big[\Bbb{P}E|{\cal S}_{\cal T}(\mu)\big]\big\ran-
{\cal C}_{{\cal Z}_{\cal T}}(\tilde{\al}_{\cal T}^{\perp})\Big)
=2|{\cal S}_{\cal T}(\mu)|-
         {\cal C}_{{\cal Z}_{\cal T}}(\tilde{\al}_{\cal T}^{\perp}).$$
By the last two estimates of Lemma~\ref{cuspcurves_l2},
${\cal C}_{{\cal Z}_{\cal T}}(\tilde{\al}_{\cal T}^{\perp})
\!=\!|{\cal Z}_{\cal T}|$.
Finally, by Lemma~\ref{cuspcurves_l2} and Corollary~\ref{euler_crl},
${\cal C}_{{\cal Z}_{\cal T}}(\al_E^{\perp})
=2|{\cal Z}_{\cal T}|$.
Thus,
\begin{equation}\label{n3m1c3_e3}
{\cal C}_{\Bbb{P}E|{\cal S}_{\cal T}(\mu)}(\al^{\perp})=
\big(2|{\cal S}_{\cal T}(\mu)|-|{\cal S}_{\cal T}(\mu)|\big)+
2|{\cal S}_{\cal T}(\mu)|=3|{\cal S}_{\cal T}(\mu)|.
\end{equation}
The claim follows by summing up equations~\e_ref{n3m1c3_e2} 
and~\e_ref{n3m1c3_e3} over the appropriate equivalence classes of 
bubble types ${\cal T}\!<\!{\cal T}^*$ and plugging the result 
back into~\e_ref{n3m1c3_e1}.\\

\noindent
The next step is to relate $\lan a,[{\cal S}_1(\mu)]\ran$ and
$\lan c_1({\cal L}^*),[{\cal S}_1(\mu)]\ran$
to intersection numbers on the spaces $\bar{\cal V}_1(\mu)$,
$\bar{\cal V}_2(\mu)$, and $\bar{\cal V}_3(\mu)$.
The approach is similar to the proof of Lemma~\ref{n2cusps},
but first we need to interpret $\lan a,[{\cal S}_1(\mu)]\ran$ and
$\lan c_1({\cal L}^*),[{\cal S}_1(\mu)]\ran$  as 
the zero sets of some bundle sections.
In our case, 
the spaces $\bar{\cal U}_{{\cal T}^*}(\mu)$ and 
$\bar{\cal U}_{{\cal T}^*(l)}(\mu)$ for all $l\!\in\![N]$
are topological manifolds (not just orbifolds).
Thus, $c_1({\cal L}^*)$ represents the first chern class
of some line bundle ${\cal L}^*\!\lra\!\bar{\cal U}_{\cal T}(\mu)$.
It is well-known in algebraic geometry that
a slightly weaker statement is in fact true  for any choice of
constraints, and
$${\cal L}^*=L^*\otimes{\cal O}\Big(-\sum_{l\in[N]}
                             \bar{\cal U}_{{\cal T}^*(l)}\Big).$$
Let $V_1\!=\!\ev^*{\cal O}(1)\!\lra\!\bar{\cal U}_{{\cal T}^*}(\mu)$,
$V_2\!=\!{\cal L}^*\!\lra\!\bar{\cal U}_{{\cal T}^*}(\mu)$, and
$\eta_i\!=\!c_1(V_1)$.
Choose sections \hbox{$s_i\!\in\!\Ga(\bar{\cal U}_{{\cal T}^*}(\mu);V_i)$}
such that $s_i$ is
smooth and transversal to the zero set on all  smooth strata 
${\cal U}_{\cal T}(\mu)\subset\bar{\cal U}_{{\cal T}^*}(\mu)$
and on \hbox{\hbox{${\cal S}_{\cal T}(\mu)\subset\bar{\cal S}_1(\mu)$}}.
The second condition implies that $s_i$ does not vanish
on the finite set~$\partial\bar{\cal S}_1$.

\begin{lmm}
\label{n3cusps}
If $d\!\ge\!1$, $\mu$ is a tuple of $p$ points and $q$ lines
in general position in $\PPP$ with \hbox{$2p\!+\!q\!=\!4d\!-\!3$},
\begin{gather*}
\lan a,[\bar{\cal S}_1(\mu)]\ran=
\big\lan 6a^3c_1({\cal L}^*)+4a^2c_1^2({\cal L}^*)+
ac_1^3({\cal L}^*),[\bar{\cal V}_1(\mu)]\big\ran-
\big\lan 4a^2+a\big(c_1({\cal L}_1^*)+c_2({\cal L}_1^*)\big),
[\bar{\cal V}_2(\mu)]\big\ran;\\
\lan c_1({\cal L}^*),[\bar{\cal S}_1(\mu)]\ran=
\big\lan 4a^3c_1({\cal L}^*)+6a^2c_1^2({\cal L}^*)+
4ac_1^3({\cal L}^*)+c_1^4({\cal L}^*),[\bar{\cal V}_1(\mu)]\big\ran-
\tau_3(\mu).
\end{gather*}
\end{lmm}

\noindent
{\it Proof:}
(1) Similarly to the proof of Lemma~\ref{n2cusps},
\begin{equation}\label{n3cusps_e1}
\lan \eta_i,[\bar{\cal S}_1(\mu)]\ran=
\big\lan \eta_ic_3\big(L^*\otimes\ev^*T\PPP\big),
                        [\bar{\cal V}_1(\mu)]\big\ran-
{\cal C}_{\partial\bar{\cal V}_1(\mu)}\big({\cal D}\oplus s_i\big).
\end{equation}
Suppose  ${\cal T}\!=\!\big(S^2,[N],I;j,d\big)\!<\!{\cal T}^*$
is a bubble type such that ${\cal U}_{\cal T}(\mu)\neq\eset$.
If $d_{\hat{0}}\!\neq\!0$, by our assumptions on $s_i$, 
${\cal D}\!\oplus\! s_i$ does not vanish on ${\cal U}_{\cal T}(\mu)$.
Thus, for the purposes of computing 
\hbox{${\cal C}_{\partial\bar{\cal V}_1(\mu)}\big({\cal D}
                \!\oplus\! s_i\big)$},
we can \hbox{assume $d_{\hat{0}}\!=\!0$}.\\
(2) In order to compute the numbers
${\cal C}_{{\cal U}_{\cal T}(\mu)}\big({\cal D}
                \!\oplus\! s_i\big)$,
we slightly modify the approach of Subsection~\ref{top_sec2},
since we have a great amount of flexibility in choosing the section~$s_i$.
We consider a family $\psi_t\!=\!(t\nu\!+\!{\cal D},s_i)$ of sections of
$L^*\!\otimes\!\ev^*T\PPP\!\oplus\! V_i$, with $\nu$ generic with
respect to~${\cal D}$.
Let $\pi\!: {\cal FT}\!\lra\!{\cal U}_{\cal T}(\mu)$ be the bundle
projection map and fix an identification of 
$\ga_{\cal T}^{\mu*}V_i\!\lra\! {\cal FT}_{\de}$ with~$\pi^*V_i$.
It can be assumed that the section~$s_i$ has been chosen so
that $\ga_{\cal T}^{\mu*}s_i\!\in\!
\Ga({\cal FT}_{\de};\pi^*V_i)$
is constant on the fibers of ${\cal FT}_{\de}$ over 
an open subset $K_{\cal T}$ of ${\cal U}_{\cal T}(\mu)$
that contains all of the finitely many zeros of the affine~map
$${\cal FT}\lra L^*\otimes\ev^*T\PPP\oplus V_i,\quad
(b,v)\lra \big(\nu_b+\al_{{\cal T},b}(v),s_i(b)\big),$$
over $\bar{\cal U}_{\cal T}(\mu)$, 
where $\al_{\cal T}$ is the leading term of 
${\cal D}_{\ga_{\cal T}^{\mu}(b,v)}$ 
given by Theorem~\ref{str_global}.
Note that by our assumptions on~$s_i$, 
the images of $\{{\cal D}_{{\cal T},h}\!:\chi_{\cal T}h\!=\!1\}$
are linearly independent in every fiber of $\ev^*T\P$
over~$s_i^{-1}(0)$.
Thus by Theorem~\ref{str_global}, Corollary~\ref{euler_crl},
and Lemma~\ref{top_l1}, 
${\cal C}_{{\cal U}_{\cal T}(\mu)}\big({\cal D}
                \!\oplus\! s_i\big)\!=\!0$
if $H_{\hat{0}}{\cal T}\!\neq\!\hat{I}$.
Furthermore, if  \hbox{$H_{\hat{0}}{\cal T}\!=\!\hat{I}$},
${\cal C}_{{\cal U}_{\cal T}(\mu)}\big({\cal D}\!\oplus\! s_i\big)$
is the number of zeros of the affine map
\begin{equation}\label{n3cusps_e3}
{\cal FT}\lra L^*\otimes\ev^*T\PPP,\quad
\ups=(b,v)\lra  \bar{\nu}_{{\cal T},b}+\al_{\cal T}(\ups),
\end{equation}
over $s_i^{-1}(0)\cap\bar{\cal U}_{\cal T}(\mu)$,
where $\bar{\nu}_{\cal T}\!\in\!\Ga(\bar{\cal U}_{\cal T}(\mu);
L^*\otimes\ev^*T\PPP)$ is a generic section.
Thus, by Lemma~\ref{zeros_main},
\begin{alignat}{1}
{\cal C}_{{\cal U}_{\cal T}(\mu)}\big({\cal D}\oplus s_i\big)
&=\sum_{k=0}^{k=2}
\big\lan \la_{{\cal FT}}^{2-k}c_k(L^*\!\otimes\!\ev^*T\PPP),
\big[\Bbb{P}{\cal FT}\big|\big(s_i^{-1}(0)\cap
\bar{\cal U}_{\cal T}(\mu)\big)\big]\big\ran
-{\cal C}_{\Bbb{P}{\cal FT}\big|\big(s_i^{-1}(0)\cap\partial
\bar{\cal U}_{\cal T}\big)}(\al_{\cal FT}^{\perp}),\notag\\
\label{n3cusps_e4}
&=\sum_{k=0}^{k=2}
\big\lan \la_{{\cal FT}}^{2-k}c_k(L^*\!\otimes\!\ev^*T\PPP)\eta_i,
\big[\Bbb{P}{\cal FT}]\big\ran
-{\cal C}_{\Bbb{P}{\cal FT}\big|\big(s_i^{-1}(0)\cap\partial
\bar{\cal U}_{\cal T}\big)}(\al_{\cal FT}^{\perp})
\end{alignat}
where $\partial\bar{\cal U}_{\cal T}\!=\!\bar{\cal U}_{\cal T}(\mu)\!-\!
                 {\cal U}_{\cal T}(\mu)$ and 
$\al_{\cal FT}\!\in\!\Ga(\Bbb{P}{\cal FT};
                  \ga_{{\cal FT}}^*\!\otimes\!L^*\!\otimes\!\ev^*T\PPP)$
is the section induced by~$\al_{\cal T}$.
We denote $\ga_{{\cal FT}}^*\!\otimes\!L^*\!\otimes\!\ev^*T\PPP$
by~${\cal O}_{\cal T}$.\\
(3) Suppose $i\!=\!1$, i.e.~$\eta_i\!=\!a$.
If ${\cal T}\!=\!{\cal T}^*(l)$ and
${\cal T}'\!=\!\big(S^2,[N],I';j,d'\big)\!<\!{\cal T}$
is a bubble type such that 
$s_1^{-1}(0)\cap{\cal U}_{{\cal T}'}(\mu)\cap 
\al_{\cal T}^{-1}(0)\neq\eset$,
${\cal T}'$ must have the~form
$$ |I'-I|=2,\quad
H_{\hat{1}}{\cal T}'=\{\hat{2},\hat{3}\},\quad
d_{\hat{1}}'=0,~d_{\hat{2}}'\neq0,~d_{\hat{3}}'\neq0.$$
By Theorem~\ref{str_global} applied to $\bar{\cal T}'\!<\!\bar{\cal T}$,
and Corollary~\ref{euler_crl},
$${\cal C}_{s_1^{-1}(0)\cap{\cal U}_{{\cal T}'}(\mu)}
                         (\al_{\cal T}^{\perp})=
\big|{\cal U}_{{\cal T}'}(\mu)\cap s_1^{-1}(0)\big|
=\lan a,[\bar{\cal U}_{{\cal T}'}(\mu)]\ran.$$
Thus, summing up equation~\e_ref{n3cusps_e4} over
${\cal T}\!=\!{\cal T}^*(l)$ with $l\!\in\![N]$, we obtain
\begin{equation}\label{n3cusps_e7}
\sum_{l\in[N]}
{\cal C}_{{\cal U}_{{\cal T}^*(l)}(\mu)}\big({\cal D}\!\oplus\! s_1\big)
\!=\! \sum_{l\in[N]} 
\big\lan  6a^3\!+\!4a^2c_1(L_{{\cal T}^*(l),\hat{1}}^*)\!+\!
ac_1^2(L_{{\cal T}^*(l),\hat{1}}^*),
\big[\bar{\cal U}_{{\cal T}^*(l)}(\mu)\big]\big\ran  
\!-\!\tau_2^{(1)}(\mu),
\end{equation}
where $\tau_2^{(1)}(\mu)$ is the number of two-component
connected degree-$d$ curves passing through the constraints
with the node at the intersection of one of the constraints
with a generic plane in~$\PPP$.
If $|H_{\hat{0}}{\cal T}|\!=\!|\hat{I}|\!=\!2$, 
$|M_{\hat{0}}{\cal T}|\!=\!0$, and ${\cal T}'$ is as 
above, up to equivalence of bubble types,
$$ |I'-I|=1,\quad\io_{\hat{3}}'=\hat{1},\quad
d_{\hat{1}}'=0,~
d_{\hat{2}}'\neq0,~d_{\hat{3}}'\neq0,$$
i.e.~$\bar{\cal T}'\!=\!\bar{\cal T}(l)$ for some $l\!\in\![N]$.
By  Theorem~\ref{str_global} applied to $\bar{\cal T}'\!<\!\bar{\cal T}$
and Corollary~\ref{euler_crl},
$${\cal C}_{\Bbb{P}{\cal FT}\big|s_1^{-1}(0)\cap{\cal U}_{{\cal T}'}(\mu)}
                         (\al_{\cal FT}^{\perp})=
\big|{\cal U}_{{\cal T}'}(\mu)\cap s_1^{-1}(0)\big|=
\lan a,[\bar{\cal U}_{{\cal T}'}(\mu)]\ran.$$
Thus, summing up equation~\e_ref{n3cusps_e4} over ${\cal T}$ with 
$|H_{\hat{0}}{\cal T}|\!=\!|\hat{I}|\!=\!2$ and 
$|M_{\hat{0}}{\cal T}|\!=\!0$, we thus obtain
\begin{equation}\label{n3cusps_e9}\begin{split}
\sum_{|H_{\hat{0}}{\cal T}|=|\hat{I}|=2,|M_{\hat{0}}{\cal T}|=0}
\!\!\!\!\!\!
\big\lan a\big(4a\!+\!c_1(L_{{\cal T},\hat{1}}^*)\!+\!
c_1(L_{{\cal T},\hat{2}}^*)\big),
\big[\bar{\cal U}_{\cal T}(\mu)\big]\big\ran 
\!-\!2\tau_2^{(1)}(\mu)\qquad\qquad\qquad&\\
=\big\lan 4a^2\!+\!a\big(c_1({\cal L}_1^*)\!+\!c_1({\cal L}_2^*)\big),
\big[\bar{\cal V}_2(\mu)\big]\big\ran.&
\end{split}\end{equation}                      
If $|H_{\hat{0}}{\cal T}|\!=\!|\hat{I}|\!=\!2$ and 
$|M_{\hat{0}}{\cal T}|\!=\!1$,
$\al_{\cal T}$ has full rank on all of~$\bar{\cal U}_{\cal T}(\mu)$.
Thus, by Corollary~\ref{euler_crl},
$${\cal C}_{{\cal U}_{\cal T}(\mu)}({\cal D}\oplus s_i)
=\big\lan c_1(L^*\otimes\ev^*T\PPP)-c_1({\cal FT}),
[\bar{\cal U}_{\cal T}(\mu)\cap s_1^{-1}(0)]\big\ran
=|{\cal U}_{\bar{\cal T}}(\mu)|.$$
Here we used ${\cal FT}\!=\!L\!\otimes
\!\big(L_{\hat{1}}\!\oplus\!L_{\hat{2}}\big)
\!\approx\! L\!\oplus\! L$ and Corollary~\ref{chern_mbar4}.
Thus, summing up equation~\e_ref{n3cusps_e4} over ${\cal T}$ with 
$|H_{\hat{0}}{\cal T}|\!=\!|\hat{I}|\!=\!2$ and 
$|M_{\hat{0}}{\cal T}|\!=\!1$ gives
\begin{equation}\label{n3cusps_e11}
\sum_{|H_{\hat{0}}{\cal T}|=|\hat{I}|=2,|M_{\hat{0}}{\cal T}|=1}\!\!\!
{\cal C}_{{\cal U}_{\cal T}(\mu)}({\cal D}\oplus s_i)
=\tau_2^{(1)}(\mu).
\end{equation}
Finally, if $|H_{\hat{0}}{\cal T}|\!=\!|\hat{I}|\!=\!3$,
\hbox{$\eta_1|\bar{\cal U}_{\cal T}'(\mu)\!=\!0$}.
The first claim follows by plugging the sum of 
equations~\e_ref{n3cusps_e7}-\e_ref{n3cusps_e11} into~\e_ref{n3cusps_e1}.
See also equations~\e_ref{normal_bundle1} and~\e_ref{normal_bundle2}.\\
(4) Suppose $\eta_i\!=\!c_1({\cal L}^*)$.
We continue as in (3) above.
If ${\cal T}\!=\!{\cal T}^*(l)$, 
$\al_{\cal T}$ does not vanish anywhere on 
$s_2^{-1}(0)\cap \bar{\cal U}_{\cal T}(\mu)$.
Thus, by Corollary~\ref{euler_crl},
\begin{equation}\label{n3cusps_e15}\begin{split}
\sum_{l\in[N]}\!
{\cal C}_{{\cal U}_{{\cal T}^*(l)}(\mu)}
              \big({\cal D}\!\oplus\! s_2\big)
&=\!\!\sum_{l\in[N]} \!\!
\big\lan c(L^*\otimes\ev^*T\PPP)c(L_{{\cal T},\hat{1}})^{-1}
\big[\bar{\cal U}_{{\cal T}^*(l)}(\mu)\cap s_2^{-1}(0)\big]\big\ran\\
&=\!\!\sum_{l\in[N]} 
\big\lan c_1({\cal L}^*) 
\big(6a^2\!+\!4ac_1(L_{{\cal T}^*(l),\hat{1}}^*)\!+\!
c_1^2(L_{{\cal T}^*(l),\hat{1}}^*)\big),
\big[\bar{\cal U}_{{\cal T}^*(l)}(\mu)\big]\big\ran.  
\end{split}\end{equation}
If $|H_{\hat{0}}{\cal T}|\!=\!|\hat{I}|\!=\!3$, 
$\al_{\cal T}$ again does not vanish anywhere on 
$s_2^{-1}(0)\cap \bar{\cal U}_{\cal T}(\mu)$, and thus
\begin{equation}\label{n3cusps_e17}
{\cal C}_{{\cal U}_{\cal T}(\mu)}\big({\cal D}\oplus s_2\big)=
\big|\bar{\cal U}_{\cal T}(\mu)\cap s_2^{-1}(0)\big|
=\big\lan c_1({\cal L}^*),[\bar{\cal U}_{\cal T}(\mu)]\big\ran
=\big|\bar{\cal U}_{\bar{\cal T}}(\mu)\big|.
\end{equation}
Here we used Corollary~\ref{chern_mbar4} again.
Note that if $|H_{\hat{0}}{\cal T}|\!=\!|\hat{I}|\!=\!2$, 
$\eta_2|\bar{\cal U}_{\cal T}(\mu)\!=\!0$.
This is immediate in the case $|M_{\hat{0}}{\cal T}|\!=\!0$
and follows from Corollary~\ref{chern_mbar4} and
\e_ref{normal_bundle1} in the case $|M_{\hat{0}}{\cal T}|\!=\!1$.
The second claim of the lemma is obtained by summing
\e_ref{n3cusps_e17} over all equivalence classes
of bubble types ${\cal T}\!<\!{\cal T}^*$ with
$|H_{\hat{0}}{\cal T}|\!=\!|\hat{I}|\!=\!3$, and plugging
the result along with \e_ref{n3cusps_e15} into~\e_ref{n3cusps_e1}.
Note that 
$$a^3\big|\bar{\cal U}_{{\cal T}^*(l)}(\mu)=0
\quad\forall~l\!\in\![N]\Lra
\big\lan 4a^3c_1(L^*),[\bar{\cal V}_1(\mu)]\big\ran
=\big\lan 4a^3c_1({\cal L}^*),[\bar{\cal V}_1(\mu)]\big\ran.$$

\subsection{The Number $n_1^{(1)}(\mu)$ in the $n=3$ Case}

\noindent
We finally compute the remaining number $n_1^{(1)}(\mu)$.
The computation parallels the proof of Lemma~\ref{n2m1c1}.

\begin{lmm}
\label{n3m1c1}
If $n\!=\!3$, 
$$n_1^{(1)}(\mu)=4
\big\lan 10a^3c_1({\cal L}^*)+3a^2c_1^2({\cal L}^*),
\big[\bar{\cal V}_1(\mu)\big]\big\ran-12\tau_2^{(2)}(\mu),$$
where $\tau_2^{(2)}(\mu)$ denotes the number of two-component
connected degree-$d$ curves that pass through the constraints
and with the node on a generic line in~$\PPP$.
\end{lmm}

\noindent
{\it Proof:} (1) We use the same notation as in the proof
of Lemma~\ref{n2m1c1}. 
By equation~\e_ref{summary_m1c1} and Lemma~\ref{zeros_main},
\begin{equation}\label{n3m1c1_e1}\begin{split}
n_1^{(1)}(\mu)
&= \sum_{k=0}^{k=5}\big\lan 
 c_k({\cal O})c_1^{5-k}\big(T^*\Si\otimes L^*\big),
   \big[\Si\times\bar{\cal U}_{{\cal T}^*}(\mu)\big]\big\ran
-{\cal C}_{\Si\times{\cal D}^{-1}(0)}(\al^{\perp}),\\
&=2 \big\lan 112a^3c_1(L^*)+84a^2c_1^2(L^*)+32ac_1^3(L^*)+5c_1^4(L^*),
\big[\bar{\cal U}_{{\cal T}^*}(\mu)\big]\big\ran
-{\cal C}_{\Si\times{\cal D}^{-1}(0)}(\al^{\perp}),
\end{split}\end{equation}
where ${\cal O}\!=\!{\cal H}^{0,1}_{\Si}\!\otimes\!\ev^*T\PPP$ and
$\al\!\in\!\Ga(\Si\times\bar{\cal U}_{{\cal T}^*}(\mu);
T^*\Si\otimes L^*\otimes{\cal O})$ is the linear part
of the affine map~$\psi_1^{(1)}$  of~\e_ref{summary_m1c1}.
Let ${\cal O}_2=T^*\Si\otimes L^*\otimes{\cal O}^{\perp}$.\\
(2) Similarly to (2) of the proof of Lemma~\ref{n2cusps},
\begin{equation}\label{n3m1c1_e3}
{\cal C}_{\Si\times{\cal S}_1(\mu)}(\al^{\perp})\!=\!
\big\lan e\big(T^*\Si\otimes\! L^*\!\otimes
({\cal H}_{\Si}^-\!\otimes\ev^*T\PPP/\Bbb{C})\big),
	\big[\Si\!\times\!\bar{\cal S}_1(\mu)\big]\big\ran
\!=\!2\big\lan 12a\!+\!5c_1(L^*),[\bar{\cal S}_1(\mu)]\big\ran.
\end{equation}
Suppose ${\cal T}\!\!=\!(S^2,[N],I;j,d)\!<\!{\cal T}^*$
is a bubble type such that 
${\cal S}_{\cal T}(\mu)\!\neq\!\eset$.
By Lemma~\ref{cuspcurves_l2},
there are three possibilities for the structure of ${\cal T}$,
but ${\cal C}_{\Si\times{\cal D}^{-1}(0)}(\al^{\perp})\!=\!0$
in all three cases.
This claim follows from Corollaries~\ref{cuspcurves_l1c} and
\ref{euler_crl} and Lemma~\ref{top_l1}.\\
(3) As before, if ${\cal T}\!<\!{\cal T}^*$ and 
${\cal C}_{\Si\times({\cal U}_{\cal T}(\mu)-{\cal S}_{\cal T}(\mu))}
                 (\al^{\perp})\!\neq\!0$,
$d_{\hat{0}}\!=\!0$ and $H_{\hat{0}}{\cal T}\!=\!\hat{I}$.
In such a case, 
\begin{equation}\label{n3m1c1_e5}
{\cal C}_{\Si\times({\cal U}_{\cal T}(\mu)-{\cal S}_{\cal T}(\mu))}
(\al^{\perp})=
\sum_{k=0}^{k=4}\big\lan \la_{\cal FT}^{4-k}c_k({\cal O}_2),
\big[\Bbb{P}{\cal FT}\big]\big\ran-
{\cal C}_{\tilde{\al}_{\cal FT}^{-1}(0)}
(\tilde{\al}_{\cal FT}^{\perp}),
\end{equation}
where 
$\tilde{\al}_{\cal FT}\!\in\!\Ga(\Bbb{P}{\cal FT};
\ga_{\cal FT}^*\otimes{\cal O}_2\big)$ is the section induced by 
the section 
$$\pi^{\perp}\circ(\al_{\cal T}\circ s_{\Si})
\!\in\!\Ga(\Si\!\times\!\bar{\cal U}_{\cal T}(\mu);
{\cal FT}^*\otimes{\cal O}_2\big).$$
(4) If ${\cal T}\!=\!{\cal T}^*(l)$ for some $l\!\in\![N]$,
$F{\cal T}\approx\! L_{{\cal T},\hat{1}}$
over $\bar{\cal U}_{\cal T}(\mu)$, and
\hbox{$\al_{\cal T}^{-1}(0)\!=\!
\Si\!\times\!{\cal D}_{{\cal T},\hat{1}}^{-1}(0)$}.
Thus, by~\e_ref{n3m1c1_e5}, 
$${\cal C}_{\Si\times({\cal U}_{\cal T}(\mu)-{\cal S}_{\cal T}(\mu))}
(\tilde{\al}_{\cal FT}^{\perp})\!=\!
2 \big\lan 112a^3\!+\!84a^2c_1(L_{{\cal T},\hat{1}}^*)\!+\!
32ac_1^2(L_{{\cal T},\hat{1}}^*)\!+\!5c_1^3(L_{{\cal T},\hat{1}}^*),
\big[\bar{\cal U}_{\cal T}(\mu)\big]\big\ran 
\!-\!{\cal C}_{\Si\times{\cal D}_{{\cal T},\hat{1}}^{-1}(0)}
         (\tilde{\al}_{\cal FT}^{\perp}).$$
By Corollaries~\ref{reg_crl2} and \ref{euler_crl},
$${\cal C}_{\Si\times{\cal S}_{\cal T}(\mu)}(\tilde{\al}_{\cal FT}^{\perp})
=\big\lan c\big({\cal FT}^*\otimes{\cal O}_2^{\perp}\big)
c(L_{{\cal T},\hat{1}}^*\otimes\ev^*T\PPP)^{-1},
\big[\Si\!\times\!{\cal S}_{\cal T}(\mu)\big]\big\ran
=10\big|{\cal S}_{\cal T}(\mu)\big|.$$
On the other hand, if 
${\cal T}'\!=\!(S^2,[N],I';j',d')\!<\!{\cal T}$,
we apply Theorem~\ref{si_str} to $\bar{\cal T}'\!<\!\bar{\cal T}$.
Then for the same reason as before,
${\cal C}_{\Si\times{\cal U}_{{\cal T}'}}
         (\tilde{\al}_{\cal FT}^{\perp})\!=\!0$
unless $d_{\hat{1}}\!=\!0$ and 
$H_{\hat{1}}{\cal T}'\!=\!I'\!-\!I$, i.e.
$$\hat{I}'=\{\hat{1},\hat{2},\hat{3}\},\quad
\io_{\hat{2}}'=\hat{1},~\io_{\hat{3}}'=\hat{1},\quad
d_{\hat{1}}'=0,~d_{\hat{2}}'\neq0,~d_{\hat{3}}'\neq0.$$
In such a case,  with $E\!=L_{{\cal T},\hat{2}}\oplus L_{{\cal T},\hat{3}}$,
by Corollary~\ref{euler_crl},
\begin{equation*}\begin{split}
{\cal C}_{\Si\times{\cal U}_{{\cal T}'}(\mu)}
(\tilde{\al}_{\cal FT}^{\perp})
&=\big\lan c\big(F{\cal T}^*\otimes{\cal O}_2\big)
c\big(L_{{\cal T}',\hat{2}}\oplus L_{{\cal T}',\hat{3}}\big)^{-1}
,\big[\Si\times\bar{\cal U}_{{\cal T}'}(\mu)\big]\big\ran\\
&=\big\lan 32a+5\big(c_1(L_{{\cal T}',\hat{2}}^*)+
c_1(L_{{\cal T}',\hat{3}}^*)\big),
\big[\bar{\cal U}_{{\cal T}'}(\mu)\big]\big\ran.
\end{split}\end{equation*}
Summing equation~\e_ref{n3m1c1_e5} over 
${\cal T}\!=\!{\cal T}^*(l)$, we thus obtain
\begin{alignat}{1}\label{n3m1c1_e7}
&\sum_{l\in[N]}\!
{\cal C}_{\Si\times({\cal U}_{{\cal T}^*(l)}(\mu)-
 {\cal S}_{{\cal T}^*(l)}(\mu)}(\al^{\perp})
= -64\tau_2^{(1)}(\mu)
-10\big\lan c_1(L_1^*)+c_1(L_2^*),\big[\bar{\cal V}_{2,1}(\mu)\big]\big\ran\\
&\qquad+ 
2\sum_{l\in[N]}\!\Big(
\big\lan 112a^3\!\!+\!84a^2c_1(L_{\hat{1}}^*)\!+\!
32ac_1^2(L_{\hat{1}}^*)\!+\!5c_1^3(L_{\hat{1}}^*),
\big[\bar{\cal U}_{{\cal T}^*(l)}(\mu)\big]\big\ran
-5\big|\bar{\cal U}_{{\cal T}^*(l)}\cap\bar{\cal S}_1(\mu)\big|
\Big),\notag 
\end{alignat} 
where $L_{\hat{1}}\!=\!L_{{\cal T}^*(l),\hat{1}}$,
$\bar{\cal V}_{2,1}(\mu)\!=
\!\bigcup\limits_{l\in[N]}\bar{\cal V}_{2,1;l}(\mu)$,
and $\bar{\cal V}_{2,1;l}(\mu)$ denotes the union of the spaces
$\bar{\cal U}_{\cal T}(\mu_l;\mu_{[N]-\{l\}})$
taken over all equivalence classes of basic bubble types 
${\cal T}\!=\!(S^2,[N]-\{l\}\{\hat{1},\hat{2\}};j,d)$
with $d_{\hat{1}},d_{\hat{2}}\!>\!0$ and
$d_{\hat{1}}\!+\!d_{\hat{2}}\!=\!d$.\\
(5) If $|H_{\hat{0}}{\cal T}|\!=\!|\hat{I}|\!=\!2$ and
$|M_{\hat{0}}{\cal T}|\!=\!0$,
${\cal FT}\!\approx L_{{\cal T},\hat{1}}\oplus L_{{\cal T},\hat{2}}$ 
over $\bar{\cal U}_{\cal T}(\mu)$ and
$\tilde{\al}_{\cal FT}^{-1}(0)$ consists of a section
${\cal Z}_{\cal T}$ of $\Bbb{P}F{\cal T}$
over $\Si\!\times\!{\cal S}_{\cal T}(\mu)$
and the spaces $\Si\!\times\!\bar{\cal U}_{{\cal T}'}(\mu)$,
with $\bar{\cal T}'$ corresponding to the bubble types~${\cal T}$
described in (2) in the proof of Lemma~\ref{n3tangents}.
By Corollaries~\ref{reg_crl2} and~\ref{euler_crl},
$${\cal C}_{\Si\times{\cal S}_{\cal T}(\mu)}
(\tilde{\al}_{\cal FT}^{\perp})
=\big\lan c\big(\ga_{\cal FT}^*\otimes{\cal O}_2^{\perp}\big)
c(\Bbb{C}^3)^{-1},\big[\Si\!\times\!{\cal S}_{\cal T}(\mu)\big]\big\ran
= 10\big|{\cal S}_{\cal T}(\mu)\big|.$$
On the other hand, if 
${\cal T}'\!=\!(S^2,[N],I';j',d')\!<\!{\cal T}$
and ${\cal C}_{\Bbb{P}{\cal FT}|
    (\Si\times{\cal U}_{\cal T}'(\mu))}(\tilde{\al}_{\cal FT}^{\perp})
\!\neq\!0$,
$$|I'-I|\in\{1,2\},\quad H_{\hat{2}}{\cal T}=I'-I,\quad
d_{\hat{2}}'=0,~d_h'\neq0\hbox{ if }h\!\in\! \hat{I}'\!-\!\{\hat{2}\}.$$
If $|I'\!-\!I|\!=\!1$, by Corollary~\ref{euler_crl},
\begin{equation*}\begin{split}
{\cal C}_{\Bbb{P}{\cal FT}|
    (\Si\times{\cal U}_{\cal T}'(\mu))}(\tilde{\al}_{\cal FT}^{\perp})
&=\big\lan c({\cal O}_2)
c(L_{{\cal T}',\hat{1}}\oplus L_{{\cal T}',\hat{1}}\big)^{-1},
\big[\Si\!\times\!\bar{\cal U}_{{\cal T}'}(\mu)\big]\big\ran\\
&=2\big\lan 32a+5\big(c_1(L_{{\cal T}',\hat{1}}^*)+
                          c_1(L_{{\cal T}',\hat{3}}^*)\big),
\big[\bar{\cal U}_{{\cal T}'}(\mu)\big]\big\ran;
\end{split}\end{equation*}
see the proof of Lemma~\ref{n3cusps} for more details.
If $|I'\!-\!I|\!=\!2$, by Corollary~\ref{euler_crl},
$${\cal C}_{\Bbb{P}{\cal FT}|
    (\Si\times{\cal U}_{\cal T}'(\mu))}(\tilde{\al}_{\cal FT}^{\perp})
=\big\lan c_1({\cal O}_2)-c_1(\Bbb{C}^3),
\big[\Si\!\times\!\bar{\cal U}_{{\cal T}'}(\mu)\big]\big\ran
=10\big|\bar{\cal U}_{{\cal T}'}(\mu)\big|.$$
Thus, summing equation~\e_ref{n3m1c1_e5} over ${\cal T}\!<\!{\cal T}^*$ 
with $|H_{\hat{0}}{\cal T}|\!=\!|\hat{I}|\!=\!2$ and 
$|M_{\hat{0}}{\cal T}|\!=\!0$,
we obtain
\begin{alignat}{1}\label{n3m1c1_e9}
&\sum_{[{\cal T}]}
{\cal C}_{\Si\times({\cal U}_{\cal T}(\mu)-{\cal S}_{\cal T}(\mu))}
(\al^{\perp})
= -30\tau_3(\mu)-10|{\cal S}_2(\mu)|\\
&\qquad\qquad\quad+
2 \big\lan 84a^2+ 32a\big( c_1({\cal L}_1^*)+c_1({\cal L}_2^*)\big)+
5\big(  c_1^2({\cal L}_1^*)+c_1^2({\cal L}_2^*)\big)
+5c_1({\cal L}_1^*)c_1({\cal L}_2^*),
        \big[\bar{\cal V}_2(\mu)\big]\big\ran\notag
\end{alignat} 
(6) If $|H_{\hat{0}}{\cal T}|\!=\!|\hat{I}|\!=\!2$ and 
$|M_{\hat{0}}{\cal T}|\!=\!1$, $\tilde{\al}_{\cal FT}$
does not vanish.
Thus, by Corollary~\ref{euler_crl},
$${\cal C}_{\Si\times{\cal U}_{\cal T}(\mu)}(\al^{\perp})\!= \!
\big\lan c({\cal O}_2)c\big(L^*\otimes
(L_{{\cal T},\hat{1}}\oplus L_{{\cal T},\hat{2}})\big)^{-1}\!,
\big[\Si\!\times\!\bar{\cal U}_{\bar{\cal T}}(\mu)\big]\big\ran
\!=\!2\big\lan 32a\!+\!5\big(c_1(L_{{\cal T},\hat{1}}^*)\!+\!
   c_1(L_{{\cal T},\hat{3}}^*)\big),
\big[\bar{\cal U}_{\bar{\cal T}}(\mu)\big]\big\ran.$$
Here we used the decomposition~\e_ref{cart_split2} and
Corollary~\ref{chern_mbar4}.
Summing equation~\e_ref{n3m1c1_e5} over ${\cal T}\!<\!{\cal T}^*$ 
with $|H_{\hat{0}}{\cal T}|\!=\!|\hat{I}|\!=\!2$ and 
$|M_{\hat{0}}{\cal T}|\!=\!1$, we obtain
\begin{equation}\label{n3m1c1_e11}
\sum_{[{\cal T}]}
{\cal C}_{\Si\times{\cal U}_{\cal T}(\mu)}(\al^{\perp})
= 64\tau_2^{(1)}(\mu)
+10 \big\lan c_1(L_1^*)+c_1(L_2^*),[\bar{\cal V}_{2,1}(\mu)]
\big\ran
\end{equation}
(7) Finally, if  $|H_{\hat{0}}{\cal T}|\!=\!|\hat{I}|\!=\!3$, 
${\cal FT}\approx L^*\oplus L^*\oplus L^*$ 
over $\bar{\cal U}_{\cal T}(\mu)$,
and $\tilde{\al}_{\cal FT}$ again does not vanish 
over $\Si\!\times\!\bar{\cal U}_{\cal T}(\mu)$.
Then by Corollary~\ref{euler_crl},
$${\cal C}_{\Si\times{\cal U}_{\cal T}(\mu)}(\al^{\perp})= 
\big\lan c({\cal O}_2)c\big(L^*\oplus L^*\oplus L^*\big)^{-1},
\big[\Si\!\times\!\bar{\cal U}_{\cal T}(\mu)\big]\big\ran
=10|\bar{\cal U}_{\bar{\cal T}}(\mu)|.$$
Thus, summing equation~\e_ref{n3m1c1_e5} over ${\cal T}\!<\!{\cal T}^*$ 
with $|H_{\hat{0}}{\cal T}|\!=\!|\hat{I}|\!=\!3$,
we obtain
\begin{equation}\label{n3m1c1_e15}
\sum_{|H_{\hat{0}}{\cal T}|=3}
{\cal C}_{\Si\times{\cal U}_{\cal T}(\mu)}(\al^{\perp})
=10\tau_3(\mu).
\end{equation}
From equations~\e_ref{n3m1c1_e1}, \e_ref{n3m1c1_e3}, 
and \e_ref{n3m1c1_e7}-\e_ref{n3m1c1_e15}, we conclude~that
\begin{alignat}{1}
n_1^{(1)}(\mu)=&2
 \big\lan 112a^3c_1({\cal L}^*)+84a^2c_1^2({\cal L}^*)+
32ac_1^3({\cal L}^*)+5c_1^4({\cal L}^*),
\big[\bar{\cal V}_1(\mu)\big]\big\ran\notag\\
&-2 \big\lan 84a^2+ 32a\big( c_1({\cal L}_1^*)+c_1({\cal L}_2^*)\big)+
5\big(  c_1^2({\cal L}_1^*)+c_1^2({\cal L}_2^*)\big)
+5c_1({\cal L}_1^*)c_1({\cal L}_2^*),
        \big[\bar{\cal V}_2(\mu)\big]\big\ran\notag\\
\label{n3m1c1_e19}
&-2\big\lan 12a+5c_1({\cal L}^*),[\bar{\cal S}_1(\mu)]\big\ran
+10\big|{\cal S}_2(\mu)\big|+20\tau_3(\mu).
\end{alignat}
The claim follows by using Lemma~\ref{n3tangents} and~\ref{n3cusps}.

\subsection{Computation of Chern Classes}
\label{chern_class}

\noindent
In this subsection, we  show that all intersection numbers 
of the spaces $\bar{\cal V}_k(\mu)$ involving powers of~$a$ and 
powers of $c_1({\cal L}_i^*)$ are computable.
We can then conclude that the numbers~$n_m^{(k)}(\mu)$
are computable.
The computability of intersection numbers of tautological classes
of $\bar{\cal V}_k(\mu)$,
which include $a$ and $c_1({\cal L}_i^*)$, has been shown 
in~\cite{P2}.
For the sake of completeness, a slightly different approach
is presented~below.\\

\noindent
If $d_{\hat{0}}$ and $d_{\hat{1}}$ are nonnegative integers and 
$\mu$ is an $N$-tuple of any generic constraints in $\P$, 
let $\bar{\cal M}_{(d_{\hat{0}},d_{\hat{1}})}(\mu)$ denote 
the union of the spaces $\bar{\cal U}_{\cal T}(\mu)$, 
where  ${\cal T}$ is a simple bubble type of the~form
$${\cal T}=\big(S^2,[N],\{\hat{0},\hat{1}\};j,
\{d_{\hat{0}},d_{\hat{1}}\}\big).$$
Then $\bar{\cal M}_{{\cal T},(d_{\hat{0}},d_{\hat{1}})}(\mu)$ 
is a complex codimension-one homology class in 
$\bar{\cal U}_{(d_{\hat{0}}+d_{\hat{1}};[N])}$.
\hbox{If $d\!>\!0$, let}
$$\sum_{d_{\hat{0}}+d_{\hat{1}}=d}^{\ge} f(d_{\hat{0}},d_{\hat{1}})
  =\sum_{\stackrel{d_{\hat{0}}+d_{\hat{1}}=d}{d_{\hat{0}},d_{\hat{1}}\ge 0}}
                f(d_{\hat{0}},d_{\hat{1}}),\qquad
\sum_{d_{\hat{0}}+d_{\hat{1}}=d}^{>} f(d_{\hat{0}},d_{\hat{1}})
    =\sum_{\stackrel{d_{\hat{0}}+d_{\hat{1}}=d}{d_{\hat{0}},d_{\hat{1}}>0}}
               f(d_{\hat{0}},d_{\hat{1}}),$$
whenever $f$ is any function defined on an appropriate subset of
$\Bbb{Z}\times\Bbb{Z}$.

\begin{lmm}
\label{chern_l1}
Let ${\cal T}^*\!=\!(S^2,[N],\{\hat{0}\};\hat{0},d)$
 be a bubble type with $d\!>\!0$.
Then in $H^*\big(\bar{\cal U}_{{\cal T}^*}(\mu)\big)$,
$$c_1(L^*)=    \frac{1}{d^2}\Big(  {\cal H}-2da+
\sum_{d_{\hat{0}}+d_{\hat{1}}=d}^{\ge}d_{\hat{1}}^2
      \bar{\cal M}_{(d_{\hat{0}},d_{\hat{1}})}(\mu)\Big),$$
where ${\cal H}$ denotes the subset of elements in 
$\bar{\cal U}_{{\cal T}^*}(\mu)$ that pass 
 through a generic codimension-two linear subspace of $\P$.
\end{lmm}

\noindent
{\it Proof:} (1) 
We restate the proof of \cite{I} in terms of the line
bundle $L^{*\otimes d^2}\!\!\!\lra\!\bar{\cal U}_{{\cal T}^*}(\mu)$,
instead of passing to a cover of $\bar{\cal U}_{{\cal T}^*}(\mu)$.
Define a section 
$\psi\!\in\!\Ga(\bar{\cal U}_{{\cal T}^*}(\mu);L^{*\otimes d^2})$
as follows. 
Let $H_0$ and $H_1$ be two fixed hyperplanes in $\P$,
generic with respect to the constraints $\mu_1,\ldots\mu_N$.
Suppose
$$[b]=\big[(S^2,[N],\hat{0};,({\hat{0}},y),u_{\hat{0}})\big]\in
 \bar{\cal U}_{\cal T}(\mu)$$
is such that $u_{\hat{0}}$ is transversal to $H_0$ and $H_1$.
Then,
$$u_{\hat{0}}^{-1}(H_i)=\big\{[x^{(i)}_1,y^{[i]}_1],\ldots,
        [x^{(i)}_d,y^{[i]}_d]\big\},\quad i=0,1,$$
for some $[x^{(i)}_k,y^{[i]}_k]\in\Bbb{P}^1$.
Define $\psi({[b]})$ by
\begin{equation}\label{dfn_of_sct}
\psi([b,c])=c^{d^2}\prod_{k,l\in[d]}
\Big(\frac{x^{(0)}_k}{y^{(0)}_k}-\frac{x^{(1)}_k}{y^{(1)}_k}\Big).
\end{equation}
While this section could be infinite, it is well-defined,
i.e.~independent of the choice of a representative 
$b\!\in\!{\cal B}_{{\cal T}^*}$ for~$[b]$.
With an appropriate coordinate change on $\Bbb{C}^{n+1}$,
it can be assumed that \hbox{$H_i\!=\!\{X_i\!=\!0\}$}.
The map $u_{\hat{0}}$ corresponds to $(n\!+\!1)$ homogeneous polynomials
of degree~$d$: $p_0,\ldots,p_n$.
Since the right-hand side of~\e_ref{dfn_of_sct} is symmetric 
in the roots of~$p_0$ and separately in the roots of~$p_1$,
$\psi$ is a rational function in the coefficients of $p_0$ and~$p_1$.
Thus, $\psi$ extends over all of~${\cal U}_{{\cal T}^*}(\mu)$.
Furthermore, this section extends by zero over 
\hbox{$\bar{\cal U}_{{\cal T}^*}(\mu)\!-\!{\cal U}_{{\cal T}^*}(\mu)$.}\\
(2) We  now identify the zero set of the section $\psi$
with multiplicities. 
From equation~\e_ref{dfn_of_sct}, it is clear that $\psi$
vanishes with multiplicity one if $p_0$ and $p_1$ have a common
root, i.e. if $u_{\hat{0}}$ passes through $H_0\cap H_1$.
The section~$\psi$ also has a pole of order $d$ along the sets~of~maps 
$$X_0=\{b\!:y^{(0)}_k(b)\!=\!0
\hbox{~for~}!k,~p_1(1,0)\!\neq\!0\},\qquad
X_1=\{b\!:y^{(1)}_k(b)\!=\!0\hbox{~for~}!k,~p_0(1,0)\!\neq\!0\}.$$
Note that $\bar{X}_i\!=\!\ev^{-1}(H_i)$.
Finally, while $\psi$ vanishes outside of ${\cal U}_{{\cal T}^*}(\mu)$,
$\bar{\cal U}_{\cal T}(\mu)$ has (complex) codimension one in 
$\bar{\cal U}_{{\cal T}^*}(\mu)$ if and only if 
${\cal T}\!<\!{\cal T}^*$ is a two-bubble strata, 
i.e.~as described just before the statement of the lemma. 
Let $d_{\hat{0}}$ and $d_{\hat{1}}$ be the corresponding degrees.
It follows from equation~\e_ref{dfn_of_sct} that 
$\psi$ has a zero of order $d_{\hat{1}}^2$ along an open subset 
of~${\cal U}_{\cal T}(\mu)$.
Thus, we obtain
$$c_1\big(L^{*\otimes d^2}\big)=
{\cal H}-2da+
\sum_{d_{\hat{0}}+d_{\hat{1}}=d}^{\ge}d_{\hat{1}}^2
       \bar{\cal M}_{(d_{\hat{0}},d_{\hat{1}})}(\mu).$$

\begin{crl}
\label{chern_c1} 
With notation as in Lemma~\ref{chern_l1},
$$c_1({\cal L}^*)= \frac{1}{d^2}\Big({\cal H}-2da+
\sum_{d_{\hat{0}}+d_{\hat{1}}=d}^{>}d_{\hat{1}}^2
      \bar{\cal M}_{(d_{\hat{0}},d_{\hat{1}})}(\mu)\Big).$$
\end{crl}

\noindent
{\it Proof:} This is immediate from Lemma~\ref{chern_l1}
and equation~\e_ref{normal_bundle1}.\\

\noindent
If  ${\cal T}\!=\!(S^2,[N],I;j,d)$ is
any bubble type, 
let ${\cal T}_{\hat{0}}\!=\!(S^2,M_{\hat{0}}{\cal T}+H_{\hat{0}}{\cal T},
\{\hat{0}\};\hat{0},d_{\hat{0}})$.
Denote by ${\cal T}_k$ for $k\!\in\!H_{\hat{0}}{\cal T}$
the simple bubble types corresponding to~$\bar{\cal T}$.
Then,
\begin{alignat}{1}\label{chern_e3}
\bar{\cal U}_{\cal T}(\mu)
&= \bar{\cal U}_{{\cal T}_{\hat{0}}}(\mu)
 \times_{\prod\limits_{k\in H_{\hat{0}}{\cal T}}(\ev_k\times\ev)}
\prod_{k\in H_{\hat{0}}{\cal T}}\bar{\cal U}_{{\cal T}_k}(\mu)\\
&\equiv\Big\{ \big(b_{\hat{0}},(b_k)_{k\in H_{\hat{0}}{\cal T}}\big)\in
\bar{\cal U}_{{\cal T}_{\hat{0}}}(\mu)\times\!\!
\prod\limits_{k\in H_{\hat{0}}{\cal T}}\bar{\cal U}_{{\cal T}_k}(\mu):
\ev_k(b_{\hat{0}})\!=\!\ev(b_h)
~\forall k\!\in\! H_{\hat{0}}{\cal T}\Big\}.\notag
\end{alignat}

\begin{lmm}
\label{chern_l2} 
With notation as above, if ${\cal T}\!<\!{\cal T}^*$ and
$d_{\hat{0}}\!\neq\!0$,
$$c_1 \big({\cal L}_{{\cal T}^*}^*\big) \Big|\bar{\cal U}_{\cal T}(\mu)
=\Big\{  c_1\big({\cal L}_{{\cal T}_{\hat{0}}}^*\big)
                           \big|\bar{\cal U}_{{\cal T}_{\hat{0}}}(\mu)
+\sum_{\eset\neq M_0\subset M_{\hat{0}}{\cal T}_{\hat{0}},
M_0\cap H_{\hat{0}}{\cal T}\neq\eset}
   \bar{\cal U}_{{\cal T}_{\hat{0}}(M_0)}(\mu)
\Big\}\times_{\prod\limits_{k\in H_{\hat{0}}{\cal T}}(\ev_k\times\ev)}
\prod_{k\in H_{\hat{0}}{\cal T}}
\bar{\cal U}_{{\cal T}_k}(\mu).$$
\end{lmm}

\noindent
{\it Proof:} Since
$L_{{\cal T}^*}\big|\bar{\cal U}_{\cal T}(\mu)\!=\!L_{\cal T}$
and $\bar{\cal U}_{\cal T}(\mu)\cap 
\bar{\cal U}_{{\cal T}^*(M_0)}(\mu)\!=\!\eset$ unless
$M_0\!\subset\! M_{\hat{0}}{\cal T}$,
by equation~\e_ref{normal_bundle1},
$$
c_1\big({\cal L}_{{\cal T}^*}^*\big) \Big|\bar{\cal U}_{\cal T}(\mu)
\!=\!c_1\big(L_{\cal T}^*\big) \Big|\bar{\cal U}_{\cal T}(\mu)
-\!\!\sum_{\eset\neq M_0\subset M_{\hat{0}}{\cal T}}\!\!\!\!
          \bar{\cal U}_{\cal T}(\mu)\cdot\bar{\cal U}_{{\cal T}^*(M_0)(\mu)}
\!=\!c_1\big(L_{{\cal T}_{\hat{0}}}^*\big) \Big|\bar{\cal U}_{\cal T}(\mu)
-\!\!\sum_{\eset\neq M_0\subset M_{\hat{0}}{\cal T}}\!\!\!\!
                  \bar{\cal U}_{{\cal T}(M_0)}(\mu).$$
The claim follows by using equation~\e_ref{normal_bundle1} again.

\begin{crl}
\label{chern_c2} 
All intersection numbers on $\bar{\cal V}_k(\mu)$ involving
only the powers of $a$ and $c_1({\cal L}_k^*)$
are computable.
\end{crl}

\noindent
{\it Proof:} Corollary~\ref{chern_c1} and Lemma~\ref{chern_l2}
reduce the computation of such numbers to understanding the restrictions
$c_1({\cal L}_{{\cal T}^*}^*)\big|
           \bar{\cal U}_{{\cal T}_{\hat{0}}(M_0)}(\mu)$,
where $M_0$ is  a subset of $M_{\hat{0}}{\cal T}_{\hat{0}}$
intersecting~$H_{\hat{0}}{\cal T}$.
By \e_ref{cart_split2},
$$ \bar{\cal U}_{{\cal T}_{\hat{0}}(M_0)}\approx
\bar{\cal M}_{0,\{\hat{0},\hat{1}\}+M_0}\times
      \bar{\cal U}_{{\cal T}_{\hat{0}}/M_0}.$$
We express 
$c_1({\cal L}_{{\cal T}^*}^*)\big| 
            \bar{\cal U}_{{\cal T}_{\hat{0}}(M_0)}$
in terms of cohomology classes of $\bar{\cal M}_{0,\{\hat{0},\hat{1}\}+M_0}$.
By definition, $L_{{\cal T}^*}|\bar{\cal U}_{{\cal T}_{\hat{0}}(M_0)}$ 
comes from a line bundle over $\bar{\cal M}_{0,\{\hat{0},\hat{1}\}+M_0}$.
In fact, 
$$c_1 \big(L_{{\cal T}^*}^*\big)\big| 
\big( \bar{\cal M}_{0,\{\hat{0},\hat{1}\}+M_0}\times
   \bar{\cal U}_{{\cal T}_{\hat{0}}/M_0}\big) =\psi_{\hat{0}}\times1,$$
where $\psi_{\hat{0}}$ is the $\psi$-class of 
$\bar{\cal M}_{0,\{\hat{0},\hat{1}\}+M_0}$
corresponding to the marked point~$\hat{0}$.
Since \hbox{$L_{{\cal T}^*}| 
       \bar{\cal U}_{{\cal T}_{\hat{0}}(M_0)}=L_{{\cal T}_{\hat{0}}(M_0)}$,}
\begin{equation*}\begin{split}
c_1 \big({\cal L}_{{\cal T}^*}^*\big)\big| 
             \bar{\cal U}_{{\cal T}_{\hat{0}}(M_0)}
&=c_1 \big(L_{{\cal T}^*}^*\big)\big| \bar{\cal U}_{{\cal T}_{\hat{0}}(M_0)}
-\sum_{\emptyset\neq M_0'\subset M_{\hat{0}}{\cal T}}
\bar{\cal U}_{{\cal T}_{\hat{0}}(M_0')}\cdot
            \bar{\cal U}_{{\cal T}_{\hat{0}}(M_0)}\\
&=\psi_{\hat{0}}\times 1 
 -\sum_{\emptyset\neq M_0'\subset 
               (M_0-H_{\hat{0}}{\cal T})}
\bar{\cal U}_{{\cal T}_{\hat{0}}(M_0';M_0-M_0')}
=\tilde{\psi}_{M_0- H_{\hat{0}}{\cal T}}\times 1\big|
\bar{\cal M}_{\{\hat{0},\hat{1}\}+M_0}
         \times\bar{\cal U}_{{\cal T}_{\hat{0}}/M_0},
\end{split}\end{equation*}
where ${\cal T}_{\hat{0}}(M_0';M_0-M_0')\!\equiv\!
\{{\cal T}_{\hat{0}}(M_0)\}(M_0')$ and
for any proper subset $\tilde{J}$ of~$J$ we define 
\hbox{$\tilde{\psi}_{\tilde{J}}
 \!\in\! H^2\big(\bar{\cal M}_{0,\{\hat{0},\hat{1}\}+J}\big)$} by
$$\tilde{\psi}_{\tilde{J}}=\psi_{\hat{0}}-
 \sum_{\emptyset\neq J'\subset\tilde{J}}
 \bar{\cal M}_{0,(\{\hat{0}\}+J',\{\hat{1}\}+(J-J'))}.$$
Here 
$\bar{\cal M}_{0,(\{\hat{0}\}+J',\{\hat{1}\}+(J-J'))}$ denotes the closure
in $\bar{\cal M}_{0,\{\hat{0},\hat{1}\}+J}$ of the two-component strata 
such that the marked points on one of the components 
\hbox{are $\{\hat{0}\}\!+\!J'$}.
The numbers 
$$\chi\big(|J|,|\tilde{J}|\big)\equiv
\big\lan \tilde{\psi}_{\tilde{J}}^{|J|-1},
        \big[\bar{\cal M}_{0,\{\hat{0},\hat{1}\}+J}\big]\big\ran$$
are given in Corollary~\ref{chern_l3c},
which is a consequence of the following well-known
lemma; see~\cite{P2} for example.

\begin{lmm}
\label{chern_l3}
(1) For any $j^*\!\in\! J$,  $\tilde{\psi}_{J-\{j^*\}}\!=\!0$ 
in $H^*\big(\bar{\cal M}_{0,\{\hat{0},\hat{1}\}+J}\big)$.\\
(2) If ${\cal N}\bar{\cal M}_{0,(\{\hat{0}\}+J',\{\hat{1}\}+(J-J'))}$ 
is the normal bundle~of 
$$\bar{\cal M}_{0,(\{\hat{0}\}+J',\{\hat{1}\}+(J-J'))}\approx
\bar{\cal M}_{0,\{\hat{0},\hat{1}\}+J'}\times
\bar{\cal M}_{0,\{\hat{0},\hat{1}\}+(J-J')}$$
in $\bar{\cal M}_{0,\{\hat{0},\hat{1}\}+J}$,
$$c_1\big({\cal N}\bar{\cal M}_{0,(\{\hat{0}\}+J',\{\hat{1}\}+(J-J'))}\big)=
-\psi_{\hat{1}}\times 1-1\times\psi_{\hat{0}}.$$
\end{lmm}

\begin{crl}
\label{chern_l3c}
If $m>0$, $\chi(m,0)=1$. If $m>k>0$, $\chi(m,k)=0$.
\end{crl}

\noindent
For our purposes, we can assume that the constraints 
$\mu_1,\ldots,\mu_N$ are disjoint.
In the case of $\PP$, the dimension of the space 
$\bar{\cal V}_k(\mu)$ is at most~$2$.
Thus, by a dimension count, if $\bar{\cal U}_{{\cal T}_{\hat{0}}(M_0)}(\mu)$
is nonempty and appears in the computation of the intersection numbers
of Corollary~\ref{chern_c2} via Lemma~\ref{chern_l2},
then $H_{\hat{0}}{\cal T}$ consists of a single element and 
$M_0\!=\!H_{\hat{0}}{\cal T}$.
The corresponding moduli space $\bar{\cal M}_{0,\{\hat{0},\hat{1}\}+M_0}$
is a single  point and thus
$$\big\lan \tilde{\psi}_{M_0-H_{\hat{0}}{\cal T}}^{|M_0|-1},
        \big[\bar{\cal M}_{\{\hat{0},\hat{1}\}+M_0}\big]\big\ran=1.$$
In the case of $\PPP$, $\bar{\cal V}_1(\mu)$
is four-dimensional,
and we encounter two cases
when $\bar{\cal M}_{0,\{\hat{0},\hat{1}\}+M_0}$ is positive-dimensional.
One possibility is that $H_{\hat{0}}{\cal T}$ is still a single-element set,
but $M_0$ contains one of the $N$ marked points.
In this case, by Corollary~\ref{chern_l3c}
or simply by the first statement of Lemma~\ref{chern_l3},
$$\big\lan \tilde{\psi}_{M_0- H_{\hat{0}}{\cal T}}^{|M_0|-1},
        \big[\bar{\cal M}_{\{\hat{0},\hat{1}\}+M_0}\big]\big\ran=
\chi(2,1)=0.$$
In fact, we can replace the first statement of Lemma~\ref{chern_l3}
with the direct computation of the degree $\psi_{\hat{0}}$ on 
$\bar{\cal M}_{0,4}$ given by Lemma~\ref{chern_l4} below.
The other case when
$\bar{\cal M}_{0,\{\hat{0},\hat{1}\}+M_0}$ is positive-dimensional
is $M_0\!=\!H_{\hat{0}}{\cal T}$ is a two-element set.
Then 
$$\big\lan \tilde{\psi}_{M_0- H_{\hat{0}}{\cal T}}^{|M_0|-1},
        \big[\bar{\cal M}_{\{\hat{0},\hat{1}\}+M_0}\big]\ran=
\chi(2,0)=0.$$

\begin{lmm}
\label{chern_l4}
Let $\bar{\cal M}_{0,4}^{(0)}\!=\!\big\{(y_1,y_2,y_3)\!\in\!\Bbb{C}^3\!:
y_1\!+\!y_2\!+\!y_3\!=\!0,~
\be(|y_1|)\!+\!\be(|y_2|)\!+\!\be(|y_3|)\!=\!\frac{1}{2}\big\}$.
Then the action of $S^1$ on $\bar{\cal M}_{0,4}^{(0)}$ induced 
from the standard action on $\Bbb{C}$ is free,
$$\bar{\cal M}_{0,4}=\bar{\cal M}_{0,4}^{(0)}/S^1\approx\Bbb{P}^1,$$
and the line bundle associated to this quotient is the tautological
line bundle over~$\Bbb{P}^1$.
\end{lmm}

\noindent
{\it Proof:} Identify $\bar{\cal M}_{0,4}^{(0)}$ with 
$S^3\subset\Bbb{C}^2$ $S^1$-equivarently by the map
$$(y_1,y_2,y_3)\lra \frac{(y_1,y_2)}{|y_1|+|y_2|}.$$
Our assumptions on 
$\be$ imply that this map is a diffeomorphism;
see Subsection~\ref{notation_sec}.

\begin{crl}
\label{chern_mbar4}
If ${\cal T}=(S^2,[3],\{\hat{0}\};\hat{0},0)$,
$\lan c_1(L^*),[\bar{\cal U}_{\cal T}]\ran=1$.
\end{crl}

\noindent
{\it Remark:} In \cite{Z}, we extend the definition of 
$\bar{\cal M}_{0,4}^{(0)}$ of Corollary~\ref{chern_l4}
to construct spaces $\bar{\cal M}_{\cal T}^{(0)}$
for all bubble types~${\cal T}$.

\subsection{The Final Formulas}
\label{final_form}

\noindent
We finally put everything together to arrive
at formulas for the numbers $n_{2,d}(\mu)$ in $\PP$ and~$\PPP$.
It~can be assumed that $\mu$ is a tuple of $(3d\!-\!2)$ points 
in the case 
of $\PP$ and of $p$ points and $q$ lines, with \hbox{$2p\!+\!q=4d\!-\!3$},
in the case of~$\PPP$.
In the former case, we write $n_{2,d}$ for $n_{2,d}(\mu)$
and $n_d$ for the number of rational plane degree~$d$ 
curves passing through $3d\!-\!1$ points.\\

\noindent
If $\nu\!\in\!\Ga(\Si\!\times\!\P;
\La^{0,1}\pi_{\Si}^*T^*\Si\otimes\pi_{\P}^*T^*\P)$ is generic,
for all $t\!\in\!(0,1)$, the signed cardinality of the set
${\cal M}_{\Si,t\nu,d}(\mu)$ is the symplectic invariant~$RT_{2,d}(;\mu)$.
If $t\!>\!0$ is sufficiently small, every element of 
${\cal M}_{\Si,t\nu,d}(\mu)$ lies either in a small 
neighborhood~$U$ of the set~${\cal H}_{\Si,d}(\mu)$
or in a small neighborhood~$W$ of the space of all bubble map
with singular domains.
Furthermore,
$$^{\pm}\big|{\cal M}_{\Si,t\nu,d}(\mu)\cap U\big|=
\big|{\cal H}_{\Si,d}(\mu)\big|=2n_{2,d}(\mu).$$
On the other hand, by Subsection~\ref{summary4},
\begin{equation}\label{final_form_e1}
\big|{\cal M}_{\Si,t\nu,d}(\mu)\cap W\big|\!=\!
\begin{cases}
n_1^{(1)}(\mu)\!+\!2n_1^{(2)}(\mu)\!+\!18n_1^{(3)}(\mu)
 \!+\!n_2^{(2)}(\mu),&\!\hbox{if }n\!=\!2;\\
n_1^{(1)}(\mu)\!+\!2n_1^{(2)}(\mu)\!+\!18n_1^{(3)}(\mu)
 \!+\!n_2^{(2)}(\mu)\!+\!2n_2^{(2)}(\mu)+n_3^{(1)}(\mu),&\!\hbox{if }n\!=\!3.
\end{cases}
\end{equation}
Thus, $n_{2,d}(\mu)$ is one-half of the difference between
$RT_{2,d}(;\mu)$ and the number in~\e_ref{final_form_e1}.
We write $CR(\mu)$ for the number given by~\e_ref{final_form_e1}.\\

\noindent
We first restrict to the case $n\!=\!2$.
We abbreviate ${\cal M}_{(d_1,d_2)}(\mu)$ as~${\cal M}_{d_1,d_2}$.
Let 
$${\cal Z}_{2;d}=
\Big(\bigcup_{\stackrel{d_1,d_2>0}{d_1+d_2=0}}{\cal Z}_{d_1,d_2}
\Big)\big/\Bbb{Z}_2,
\quad\hbox{where}\quad
{\cal Z}_{d_1,d_2}=
\bigcup_{j_l=1,2}
\bar{\cal U}_{(S^2,[N],I;j,\{0,d_1,d_2\}}
                            (\mu),$$
where $I\!=\!\{\hat{0},1,2\}$ with the partial ordering $0\!<\!1,2$. 
The set ${\cal Z}_{2;d}$ is the
the zero-dimensional space of three-bubble maps passing
through the (3d-2) points~$\mu$, such that the map is
trivial on the principal component.
Note that
\begin{equation}\label{p2_e1}
\big|{\cal Z}_{d;2}\big|=\tau_2(\mu)=
\frac{1}{2}\sum\limits_{d_1+d_2=d}
\binom{3d-2}{3d_1-1}d_1d_2n_{d_1}n_{d_2}.
\end{equation}
The binomial coefficient counts the number of possible ways
of distributing the constraints between the two nontrivial bubbles.
Without the factor $d_1d_2$, the above number would have been precisely
the number of two-component rational curves passing through 
$(3d\!-\!2)$ generic points in $\P$.
However, we have to account for the image of the evaluation
map at~$\hat{0}$,
which must be one of the $d_1d_2$ points of intersection of two rational
curves of degrees $d_1$ and~$d_2$.

\begin{lmm}
\label{p2_l1}
In the $n\!=\!2$ case, the total correction is given by
$$CR(\mu)=
\big\lan 78a^2+72ac_1({\cal L}^*)+
 22c_1^2({\cal L}^*),\big[\bar{\cal V}_1(\mu)\big]\big\ran-
       18\tau_2(\mu).$$
\end{lmm}

\noindent
{\it Proof:}
The four numbers of \e_ref{final_form_e1} are given by 
Lemmas~\ref{n2m1c1}, \ref{mequalsnm1}, and \ref{mequalsn}.
The cardinality of~${\cal S}_1(\mu)$ is given by Lemma~\ref{n2cusps}.

\begin{lmm}
\label{p2_l2}
With notation as above,
$$\big\lan ac_1({\cal L}^*),[\bar{\cal V}_1(\mu)]\big\ran=
\frac{1}{d}\Big(-n_d+\frac{1}{2}
\sum_{d_1+d_2=d}d_1^2d_2^2\binom{3d-2}{3d_1-1}
n_{d_1}n_{d_2}\Big).$$
\end{lmm}

\noindent
{\it Proof:} By Corollary~\ref{chern_c1},
\begin{equation}\label{p2_l2_e1}
ac_1({\cal L}^*)=\frac{1}{d^2}a\Big({\cal H}-2da+
\sum\limits_{d_1+d_2=d}^{>}d_2^2{\cal M}_{d_1,d_2}\Big).
\end{equation}
Note that
\begin{equation}\label{p2_l2_e2}\begin{split}
\sum\limits_{d_1+d_2=d}d_2^2
\big\lan a,[{\cal M}_{d_1,d_2}]\big\ran
&=\sum\limits_{d_1+d_2=d}d_1(d_1d_2)d_2^2
\binom{3d-2}{3d_1-1}n_{d_1}n_{d_2}\\
&=\frac{1}{2}d\sum\limits_{d_1+d_2=d}d_1^2d_2^2
\binom{3d-2}{3d_1-1}n_{d_1}n_{d_2}.
\end{split}\end{equation}
The reason for the appearance of the factor $d_1d_2$ in
\e_ref{p2_l2_e2} is the same one as in \e_ref{p2_e1}.
On the other hand, the factor $d_1$ appears because
we need to count the number of times the first rational component
intersects a line in $\PP$.
Since 
$$\big\lan a{\cal H},[\bar{\cal V}_1(\mu)]\big\ran= dn_d
\quad\hbox{and}\quad\big\lan a^2,[\bar{\cal V}_1(\mu)]\big\ran= n_d,$$
the claim follows by plugging \e_ref{p2_l2_e2} into \e_ref{p2_l2_e1}.

\begin{lmm}
\label{p2_l3}
With notation as above,
$$\big\lan c_1^2({\cal L}^*),[\bar{\cal V}_1(\mu)]\big\ran=-\frac{1}{2}
\sum\limits_{d_1+d_2=d}\binom{3d-2}{3d_1-1}
d_1d_2n_{d_1}n_{d_2}.$$
\end{lmm}

\noindent
{\it Proof:} By Corollary~\ref{chern_c1},
\begin{equation}\label{p2_l3_e1}
c_1^2({\cal L}^*)=\frac{1}{d^2}c_1({\cal L}^*) \Big({\cal H}-2da+
\sum\limits_{d_1+d_2=d}d_2^2{\cal M}_{d_1,d_2}\Big).
\end{equation}
Since there are no two-component rational curves of total degree
$d$ passing through $(3d\!-\!1)$ generic points in $\P$
and there are no  three-component rational curves of total degree
$d$ passing through $(3d\!-\!2)$ generic points in $\P$,
by Corollary~\ref{chern_c1}
\begin{equation}\label{p2_l3_e2}
\big\lan {\cal H}c_1({\cal L}^*),[\bar{\cal V}_1(\mu)]\big\ran
=\frac{1}{d^2}\big\lan -2da{\cal H},[\bar{\cal V}_1(\mu)]\big\ran=-2n_d.
\end{equation}
Similarly by Corollary~\ref{chern_c1} and Lemma~\ref{chern_l2},
\begin{equation}\label{p2_l3_e3}\begin{split}
\big\lan c_1({\cal L}^*),[{\cal M}_{d_1,d_2}]\big\ran
&=\frac{1}{d_1^2}\big\lan -2d_1a{\cal H},[{\cal M}_{d_1,d_2}]\big\ran
                                        +\big|{\cal Z}_{d_1,d_2}\big|
=-\big|{\cal Z}_{d_1,d_2}\big|\\
&=-d_1d_2\binom{3d-2}{3d_1-1}n_{d_1}n_{d_2}.
\end{split}\end{equation}
Note that by symmetry
\begin{equation}\label{p2_l3_e4}
\sum\limits_{d_1+d_2=d}
d_1d_2^3\binom{3d-2}{3d_1-1}n_{d_1}n_{d_2}=
\frac{1}{2}\sum\limits_{d_1+d_2=d}
d_1d_2\l(d^2-2d_1d_2\r)\binom{3d-2}{3d_1-1}n_{d_1}n_{d_2}.
\end{equation}
The claim now follows from equations \e_ref{p2_l3_e1}-\e_ref{p2_l3_e4}
and Lemma~\ref{p2_l2}.

\begin{crl}
\label{p2_c4}
The total correction term is given by
$$CR(\mu)=78n_d+72\frac{1}{d}\Big(-n_d+\frac{1}{2}
\sum\limits_{d_1+d_2=d}d_1^2d_2^2\binom{3d-2}{3d_1-1}
n_{d_1}n_{d_2}\Big)-
20\sum\limits_{d_1+d_2=d}d_1d_2\binom{3d-2}{3d_1-1}
n_{d_1}n_{d_2}.$$
\end{crl}

\noindent
{\it Proof:} This claim is immediate from 
Lemmas~\ref{p2_l1}-\ref{p2_l3} and equation~\e_ref{p2_e1}.

\begin{lmm}
\label{p2_l5}
The genus-two RT-invariant in $\PP$ is given by
$$RT_{2,d}(;\mu)\equiv RT_{2,d}(;p_{[3d-2]})=6d^2n_d+
\sum\limits_{d_1+d_2=d}d_1^3d_2^3\binom{3d-2}{3d_1-1}
n_{d_1}n_{d_2}.$$
\end{lmm}

\noindent
{\it Proof:} Applying the genus-reducing composition law of \cite{RT} twice,
we obtain
\begin{equation}\label{p2_l5_e1}\begin{split}
&RT_{2,d}(;p_{[3d-2]})=2RT_{1,d}(p,\PP;p_{[3d-2]})
+RT_{1,d}(\ell,\ell;p_{[3d-2]})\\
&\quad=4RT_{0,d}(p,\PP,p,\PP;p_{[3d-2]})+
4RT_{0,d}(p,\PP,\ell,\ell;p_{[3d-2]})
+RT_{0,d}(\ell,\ell,\ell,\ell;p_{[3d-2]})\\
&\quad=0+4RT_{0,d}(p,\ell,\ell;p_{[3d-2]})+
RT_{0,d}(\ell,\ell,\ell,\ell;p_{[3d-2]}).
\end{split}\end{equation}
Since the genus-zero three-point RT-invariant is 
the usual enumerative invariant,
the middle term above is simply $4d^2n_d$.
On the other hand, by the component-splitting  
composition law of~\cite{RT},
\begin{equation}\label{p2_l5_e2}\begin{split}
RT_{0,d}(\ell,\ell,\ell,\ell;p_{[3d-2]})
&=2RT_{0,0}(\ell,\ell,\PP;)
RT_{0,d}(\ell,\ell,p;p_{[3d-2]})\\
&\quad+\!\sum\limits_{d_1+d_2=d} 
\sum\limits_{J_1+J_2=[3d-2]}\!\!\!
RT_{0,d_1}(\ell,\ell,\ell;p_{J_1})
RT_{0,d_2}(\ell,\ell,\ell;p_{J_2})\\
&=2d^2n_d+
\sum\limits_{d_1+d_2=d}d_1^3d_2^3\binom{3d-2}{3d_1-1}
n_{d_1}n_{d_2}.
\end{split}\end{equation}
The lemma follows from equations \e_ref{p2_l5_e1} and \e_ref{p2_l5_e2}.\\

\noindent
Theorem~\ref{p2_final} is nearly proved.
We can simplify the expression in Corollary~\ref{p2_c4}
by using a recursive relation for the numbers $n_d$;
see \cite[p363]{RT}.
The expression of Theorem~\ref{p2_final} is half of
the difference between the quantity of Lemma~\ref{p2_l5}
and Corollary~\ref{p2_c4}.
Note that the numbers $n_d$ with $d\!=\!1,2,3$
have long been known to be zero; see~\cite{ACGH}.
Strictly speaking, our computation does not apply 
to the cases~$d\!=\!1,2$.
However, these two cases do provide a consistency check.\\

\noindent
The case of $\PPP$ is significantly harder than the $n\!=\!2$ case.
An explicit recursive formula as in Theorem~\ref{p2_final}
would be rather long, so we do not provide one.
Instead we express $n_{2,d}(\mu)$ in terms of the corresponding
symplectic invariant and intersection numbers of
the spaces $\bar{\cal V}_1(\mu)$, 
$\bar{\cal V}_2(\mu)$, and~$\bar{\cal V}_3(\mu)$.

\begin{thm}
\label{p3_final}
If $d$ is a positive integer and $\mu$ is a tuple  of $p$ points
and $q$ lines in general position in $\PPP$ with 
\hbox{$2p\!+\!q\!=\!4d\!-\!3$},
\begin{gather*}
2n_{2,d}(\mu)=RT_{2,d}(\cdot;\mu)-CR(\mu),\qquad\hbox{where}\\
\begin{split}
\frac{1}{2}CR(\mu)=&
\l\lan 480a^3c_1({\cal L}^*)+476 a^2c_1^2({\cal L}^*)
+240ac_1^3({\cal L}^*)+49c_1^4({\cal L}^*),\bar{\cal V}_1(\mu)\r\ran\\
&~-\big\lan 144a\big(
c_1({\cal L}_1^*)+c_1({\cal L}_2^*)\big)+
27\big(c_1^2({\cal L}_1^*)+c_1^2({\cal L}_2^*)\big)+
25c_1({\cal L}_1^*)c_1({\cal L}_2^*),
\bar{\cal V}_2(\mu)\big\ran\\
&~-324\tau_2^{(2)}(\mu)+36\tau_3(\mu).
\end{split}
\end{gather*}
Furthermore, $RT(\cdot;\mu)$ and all intersection numbers above
are computable.
\end{thm}

\noindent
{\it Proof:}
The six numbers of \e_ref{final_form_e1} in the $n\!=\!3$ case are given by 
Lemmas~\ref{n3m1c1}, \ref{n3m1c2}, \ref{n3m1c3}, 
\ref{n3m2c1}, \ref{mequalsnm1}, and \ref{mequalsn}, respectively.
The numbers $\lan a,[\bar{\cal S}_1(\mu)]\ran$, 
$\lan c_1({\cal L}^*),[\bar{\cal S}_1(\mu)]\ran$, and
$\big|\bar{\cal S}_2(\mu)\big|$ are given by 
Lemmas~\ref{n3tangents} and~\ref{n3cusps}.
The symplectic invariant $RT_{2,d}(\cdot;\mu)$ is well-known 
to be computable; see~\cite{RT}.
The above intersection numbers are computable by Corollary~\ref{chern_c2}.\\

\noindent
As in the case of $\PP$, we recover the well-known fact 
that all degree-one, -two, and -three numbers are~zero.
The only degree-one number, the number of genus-two degree-one curves 
through a line, is zero because there are no
holomorphic degree-one maps from a positive-genus curve into~$\P$;
see~\cite{ACGH}.
The eight degree-two and -three numbers are zero because
the image of any holomorphic map of degree two or three
from a genus-two curve into $\P$ is a line, see~\cite{ACGH},
while no line passes through the required constraints.
The first three degree-four numbers given below have also been known 
to be zero, since the image of any holomorphic map of degree four
from a genus-two curve into $\P$ must lie in a~plane.
Finally, observe that the fourth degree-four number is
the number~$n_{2,4}$ given by Theorem~\ref{p2_final},
as should be the~case.

\begin{center}
\begin{tabular}{||c|c|c|c|c|c|c||}
\hline\hline
degree&         \multicolumn{5}{c|}{4}& 5\\
\hline
(p,q)&   (6,1)& (5,3)& (4,5)& (3,7)& (0,13)& (5,7)\\
\hline
$\RT_{2,d}(\cdot;\mu)$ &    7,872& 64,960& 548,608& 4,906,304&
                   5,130,826,752& 290,439,680\\
\hline
$CR(\mu)$&     7,872& 64,960& 548,608& 4,877,504& 
          4,998,465,792&258,287,360\\
\hline
$n_{2,d}(\mu)$&    0& 0&0& 14,400& 66,180,480&16,076,160\\
\hline\hline
\end{tabular}
\end{center}

\section{Appendix}

\subsection{A Short Exact Sequence on $\P$}

\noindent 
If $M$ is a Kahler manifold and $E\!\lra\! M$ is a holomorphic vector bundle,
let ${\cal O}(E)$ denote the sheaf of holomorphic sections of $E$.
If $E\!\lra\! M$ is the trivial holomorphic line bundle,
we write  ${\cal O}$ for ${\cal O}(E)$.
Let $H\!\lra\!\P$ be the hyperplane bundle.

\begin{lmm}
\label{ses_lmm}
There is an exact sequence of sheaves over $\P$:
$$0\lra{\cal O}\lra (n+1){\cal O}(H) \lra {\cal O}(T\P)\lra 0.$$
\end{lmm}

\noindent
{\it Proof:} (1) Let $[X_0:\ldots:X_n]$ denote the homogeneous coordinates
on $\P$. Denote by $\bar{X}_i$ the section of the hyperplane bundle 
given by
$$\bar{X}_i|_{[X_0:\ldots:X_n]}(X_0,\ldots,X_n)=X_i\in\Bbb{C}.$$
Then we define a sheaf map ${\cal O}\lra (n+1){\cal O}(H)$ by
$$f\lra (f\bar{X}_0,\ldots,f\bar{X}_n).$$
Let $U_i=\{[X_0:\ldots:X_n]\!: X_i\!\neq\!0\}$. 
On $U_i$, we can use the complex coordinates
$$z_{i,k}=\frac{X_k}{X_i},\qquad k\in\{0,\ldots,n\}-\{i\}.$$
Using these coordinates, we define a sheaf map 
$(n+1){\cal O}(H) \lra {\cal O}(T\P)$ by
\begin{equation}
\label{ses_e1}
(p_0,\ldots,p_n)\lra\sum_{k\neq i}
\big(p_k(z_{i,0},\ldots,z_{i,n})-z_{i,k}p_i(z_{i,0},\ldots,z_{i,n})\big)
\frac{\partial}{\partial z_{i,k}},
\end{equation}
where $z_{i,i}=1$. We need to see that this map is well-defined.
Suppose $j\neq i$. Then,
\begin{equation}
\label{ses_e2}
z_{j,l}=z_{i,j}^{-1}z_{i,l}\Lra
\frac{\partial}{\partial z_{i,k}}=\sum_{l\neq j}
\frac{\partial z_{j,l}}{\partial z_{i,k}}\frac{\partial}{\partial z_{j,l}}
=\begin{cases}
z_{i,j}^{-1}\frac{\partial}{\partial z_{j,k}}, &\hbox{if~}k\neq j;\\
-z_{i,j}^{-2}\Big(
\frac{\partial}{\partial z_{j,i}}+\sum\limits_{l\neq i,j}z_{i,l}
\frac{\partial}{\partial z_{j,l}}\Big), &\hbox{if~}k=j.
\end{cases}
\end{equation}
Since each $p_l$ is a linear functional,  if $k\neq i,j$, 
we can write the $k$th summand in \e_ref{ses_e1} as
\begin{equation}
\label{ses_e3}
\begin{split}
&\big(z_{j,i}^{-1}p_k(z_{j,0},\ldots,z_{j,n})-
z_{j,i}^{-2}z_{j,k}p_i(z_{j,0},\ldots,z_{j,n})\big)
z_{i,j}^{-1}\frac{\partial}{\partial z_{j,k}}\\
&\qquad\qquad\qquad\qquad=\big(p_k(z_{j,0},\ldots,z_{j,n})-
z_{j,i}^{-1}z_{j,k}p_i(z_{j,0},\ldots,z_{j,n})\big)
\frac{\partial}{\partial z_{j,k}}.
\end{split}
\end{equation}
The remaining, $k=j$, summand in \e_ref{ses_e1} equals
\begin{equation}
\label{ses_e4}
\begin{split}
&\big(z_{j,i}^{-1}p_j(z_{j,0},\ldots,z_{j,n})-
z_{j,i}^{-2}p_i(z_{j,0},\ldots,z_{j,n})\big)
(-z_{i,j}^{-2})\Big(
\frac{\partial}{\partial z_{j,i}}+\sum_{k\neq i,j}z_{i,k}
\frac{\partial}{\partial z_{j,k}}\Big)\\
&\qquad\qquad\qquad=
\big(p_i(z_{j,0},\ldots,z_{j,n})-z_{j,i}p_j(z_{j,0},\ldots,z_{j,n})\big)
\Big(\frac{\partial}{\partial z_{j,i}}+\sum_{k\neq i,j}z_{i,k}
\frac{\partial}{\partial z_{j,k}}\Big).
\end{split}
\end{equation}
Since $z_{j,i}z_{i,k}\!=\!z_{j,k}$, collecting similar terms in 
\e_ref{ses_e3} and \e_ref{ses_e4}, we obtain equation 
\e_ref{ses_e1} with $i$ replaced by~$j$.\\
(2) It is clear that the first map is injective, the second is surjective,
and the composite is zero. Finally, if $(p_0,\ldots,p_n)$ is mapped
to zero by the second map, then \e_ref{ses_e1} implies that
$\bar{X}_jp_i=\bar{X}_ip_j$ for all $i$ and $j$. Thus, the function $f$,
given by
$$f([X_0:\ldots:X_n])=\frac{p_i(X_0,\ldots,X_n)}{X_i},$$
is well-defined and holomorphic wherever $(p_0,\ldots,p_n)$ is.

\subsection{On Regularity of Kernel of $D_b$}

\begin{lmm}
\label{les_crl2}
If $u\!: S^2\!\lra\!\P$ is a holomorphic map, 
there is a surjection
$$(n+1)H^1\big(S^2;{\cal O}(u^*H\otimes(-(k+1)p))\big)
\lra H^1\big(S^2;{\cal O}(u^*T\P\otimes(-(k+1)p))\big),$$
where $p$ denotes the divisor corresponding to a point $p\!\in\! S^2$.
If the degree of $u$ is at least $k$, then both cohomology groups are trivial.
\end{lmm}

\noindent
{\it Proof:}
Pulling back the short exact sequence of sheaves of Lemma~\ref{ses_lmm}
by $u$,  tensoring it with $-(k+1)p$, and 
taking the corresponding long exact sequence,
we obtain:
\begin{equation}\label{les_crl2_e1}\begin{split}
&\lra (n+1)H^1\big(S^2;{\cal O}(u^*H\otimes(-(k+1)p))\big)
\lra H^1\big(S^2;{\cal O}(u^*T\P\otimes(-(k+1)p))\big)\\
&\lra H^2\big(S^2;{\cal O}(-(k+1)p)\big)\lra...
\end{split}\end{equation}
Since $S^2$ is a one-dimensional complex manifold, the 
last cohomology group in \e_ref{les_crl2_e1} must vanish,
and the first statement of the lemma follows.
On the other hand, by Kodaira-Serre duality,
\begin{equation}\label{les_crl2_e2}\begin{split}
H^1\big(S^2;{\cal O}(u^*H\otimes(-(k+1)p))\big)
&=H^1\big(S^2;\Om^1(u^*H\otimes(-(k-1)p)\big)\\
&\approx H^0\big(S^2;{\cal O}((u^*H\otimes(-(k-1)p))^*)\big)^*.
\end{split}\end{equation}
The last group in \e_ref{les_crl2_e2} is trivial if 
${\cal O}\big(u^*H\otimes(-(k-1)p)\big)$ is positive, i.e.~if
$$\big\lan c_1\big(u^*H\otimes(-(k-1)p\big),[S^2]\big\ran=d-(k-1)>0,$$
where $d$ is the degree of $u$.

\begin{crl}
\label{reg_crl2}
If $f\!:S^2\!\lra\!\P$ is holomorphic map of degree $d$, 
for any $p\!\in\! S^2$ and nonzero $v\!\in\! T_pS^2$, the~map
\begin{equation*}
\phi_{p,v}^{(k)}\!: \ker D_u\lra \bigoplus_{m\in\lr{k}}T_{u(p)}\P,\quad
\phi_{p,v}^{(k)}\xi=\big(\xi_p,D\xi|_{p,v},\ldots,D^{(k)}\xi|_{p,v}\big),
\end{equation*}
where $D\xi|_{p,v}$ denotes the covariant derivative 
of $\xi$ along $u$ in the direction of $v$,
is surjective provided~$d\!\ge\! k$.
\end{crl}

\noindent
{\it Remark:} If one defines $D^{(k)}\xi$ with respect
to the metric $g_{\P,u(p)}$ on $\P$, 
$D^{(k)}\!\in\! T_{u(p)}\P\otimes T^*S^{2\otimes k}$,
where $T^*S^2$ is viewed as a complex line bundle.
However, the statement is independent of the choice of metric
on~$\P$.\\

\noindent
{\it Proof:} Since $\xi$ is holomorphic, if $\phi_{p,v}^{(k)}\xi$
is zero, $\xi$ has a zero of order $k+1$ at $p$.
Thus, $\phi_{p,v}^{(k)}$ induces 
 a short exact sequence of sheaves on $S^2$:
$$0\lra {\cal O}\big(u^*T\P\otimes(-(k+1)p)\big)\lra {\cal O}(u^*T\P)
\stackrel{\phi_{p,v}^{(k)}}{\lra}
(k+1){\cal O}\big((u^*T\P)_p\big)\lra 0,$$
where we view ${\cal O}\big((u^*T\P)_p\big)$
as a sheaf on $S^2$ via extension by $0$; see \cite[p38]{GH}.
Taking the corresponding long exact sequence in cohomology,
we obtain
\begin{equation}\label{reg_crl2_e1}\begin{split}
&\qquad\qquad\ldots\lra H^0\big(S^2;{\cal O}(u^*T\P)\big)
\stackrel{\phi_{p,v}^{(k)}}{\lra} 
(k+1)H^0\big(S^2;{\cal O}((u^*T\P)_p)\big)\\
&\lra H^1\big(S^2;{\cal O}(u^*T\P\otimes(-(k+1)p))\big)\ldots
\end{split}\end{equation}
By Lemma~\ref{les_crl2}, the last cohomology group in
\e_ref{reg_crl2_e1} is zero if $d\ge k$.
It follows that the map $\phi_{p,v}^{(k)}$ is surjective.

\subsection{Dimension Counts}

\begin{lmm}
\label{les_crl}
If $u\!:M\!\lra\!\P$ is a holomorphic map, there is an isomorphism 
$$H^1\big(M;(n+1){\cal O}(u^*H)\big)\lra 
H^1\big(M;{\cal O}(u^*T\P)\big).$$
\end{lmm}

\noindent
{\it Proof:} The short exact sequence of Lemma~\ref{ses_lmm}
gives a long exact sequence in sheaf cohomology:
\begin{equation}
\label{les_e1}
\ldots H^1(\P;{\cal O})\lra H^1\big(\P;(n+1){\cal O}(H)\big)
\lra H^1\big(\P;{\cal O}(T\P)\big)\lra H^2(\P;{\cal O})\ldots
\end{equation}
Since $\P$ is Kahler, $H^p(\P;{\cal O})={\cal H}^{p,0}(\P)$,
where ${\cal H}^{p,q}$ denotes the harmonic $(p,q)$-forms.
However, all of the cohomology of $\P$ is generated by 
$c_1(H)$ and thus lies in $H^{p,p}(\P;{\cal O})$.
In particular,
\begin{equation}
\label{les_e2}
H^1(\P;{\cal O})=H^2(\P;{\cal O})=0.
\end{equation}
If $u\!: M\!\lra\!\P$ is a homomorphic map, pulling back the long exact
sequence in \e_ref{les_e1} and using \e_ref{les_e2}, we see that
$$H^1\big(M;(n+1){\cal O}(u^*H)\big)\lra 
H^2\big(M;{\cal O}(u^*T\P)\big)$$
is an isomorphism, induced by the sheaf map 
${\cal O}(H)\!\lra\! {\cal O}(T\P)$ defined in 
the proof of Lemma~\ref{ses_lmm}.

\begin{crl}
\label{reg_crl}
Let $\Si$ be a compact Riemann surface.
If $u\!:\Si\!\lra\!\P$ is a holomorphic map, 
the $\bar{\partial}$-operator for the bundle $u^*T\P$,
$$D_u:  \Ga(\Si;u^*T\P)\lra\Ga(\Si;\La^{0,1}T^*\Si\otimes u^*T\P)$$
is surjective, provided $d+\chi(\Si)>0$, where $d$ is the degree of $u$.
\end{crl}

\noindent
{\it Proof:}
The cokernel of $D_u$ is $H_{\bar{\partial}}^1(\Si;u^*T\P)$.
By Dolbeault Theorem and Lemma~\ref{les_crl}, 
\begin{equation}
\label{reg_crl_e1}
H_{\bar{\partial}}^1(\Si;u^*T\P)=
H^1\big(\Si;{\cal O}(u^*T\P)\big)=(n+1)H^1\big(\Si;{\cal O}(u^*H)\big).
\end{equation}
By Kodaira-Serre duality (see \cite[p153]{GH}),
\begin{equation}\label{reg_crl_e2}
\begin{split}
H^1\big(\Si;{\cal O}(u^*H)\big)&
=H^1\big(\Si;\Om^1(T\Si\otimes u^*H)\big)\\ 
&=H^0\big(\Si;{\cal O}((T\Si\otimes u^*H)^*)\big)^*
=H^0_{\bar{\partial}}\big(\Si;(T\Si\otimes u^*H)^*\big)^*.
\end{split}
\end{equation}
The bundle $(T\Si\otimes u^*H)^*$ does not admit any holomorphic
section if it is negative, i.e.~if
$$\big\lan c_1\big((T\Si\otimes u^*H)\big),[\Si]\big\ran=
\big\lan c_1(T\Si)+c_1(u^*H),[\Si]\big\ran=\chi(\Si)+d>0.$$
Thus, the claim follows from equations \e_ref{reg_crl_e1} and 
\e_ref{reg_crl_e2}. 

\begin{prp}
\label{empty_spaces}
Let $\Si$ be a Riemann surface of genus $2$ and let $d$ and $n$ be 
positive integers \hbox{with $n\!\le\! 4$}. 
If $n\!=\!4$, assume that $d\!\neq\! 2$.
Suppose $\mu\!=\!(\mu_1,\ldots,\mu_N)$ is an $N$-tuple of proper 
complex submanifolds of $\P$ of total complex codimension 
$d(n\!+\!1)\!-\!n\!+\!N$ in general position.
If ${\cal T}\!=\!\big(\Si,[N],I;j,d'\big)$ is a bubble type 
such that 
${\cal T}\!<\!\big(\Si,[N],\{\hat{0}\};\hat{0},d\big)$ and 
$d_{\hat{0}}'\!>\!0$, 
then  ${\cal H}_{\cal T}(\mu)\!=\!\emptyset$.
Furthermore,~if
$$b=\big(\Si,[N],\{\hat{0}\};,(\hat{0},y),u\big)
\in{\cal H}_{\big(\Si,[N],\{\hat{0}\};\hat{0},d\big)}(\mu),$$
then the map $u$ is not multiply-covered.
\end{prp}

\noindent
{\it Proof:}
(1) If $d_{\hat{0}}'\!\ge\! 3$, by Corollaries~\ref{reg_crl2} and~\ref{reg_crl} 
and standard arguments such as in~\cite{MS},
the space ${\cal H}_{\cal T}$ is a smooth manifold
and the maps $\ev_l$ are smooth.
If $b\!\in\!{\cal H}_{\cal T}$ a neighborhood
of $b$ in ${\cal H}_{\cal T}$ can be modeled on 
$\ker D_b\oplus\bigoplus\limits_{l=1}\limits^{l=n}T_{y_l}\Si_{b,j_l}$.
In particular, by the Index Theorem,
$$\dim_{\Bbb{C}}{\cal H}_{\cal T}=
\sum_{i\in I}\big( d_i'(n+1)+n(1-g(\Si_{b,i}))\big)-(n-1)|\hat{I}|+N
=d(n+1)-n+|\hat{I}|+N.$$
Thus, if the map
$$\ev_{[N]}\equiv \ev_1\times\ldots\times\ev_N\!: 
{\cal H}_{\cal T}\lra\P\times\ldots\times\P,$$
is smooth and transversal to $\mu_1\times\ldots\times\mu_N$,
${\cal H}_{\cal T}(\mu)$ is a smooth manifold
of (complex) dimension~$|\hat{I}|$.
Since the map $\ev_{[N]}$ is invariant under the action 
of $2|\hat{I}|$-dimensional group 
$${\cal G}_{\cal T}\equiv 
\{g\!\in\! PSL_2\!: g(\i)=\i\}^{\hat{I}},$$
${\cal G}_{\cal T}$ acts smoothly
on~${\cal H}_{\cal T}(\mu)$.
Furthermore, the stabilizer at each point is finite.
Thus, ${\cal H}_{\cal T}(\mu)\!=\!\emptyset$.\\
(2) Suppose $d_{\hat{0}}'\!=\!2$. If 
$b\!=\!\big(\Si,[N],I;x,(j,y),u\big)\!
\in\!{\cal H}_{\cal T}$, 
the map $u_{\hat{0}}$ must factor through a degree-one map 
$\tilde{u}_{\hat{0}}\!: S^2\!\lra\!\P$;
see \cite[p116]{ACGH}.
Thus, it is enough to show that the space 
\hbox{${\cal H}_{{\cal T}'}(\mu)\!=\!\emptyset$,} where
${\cal T}'=\big(S^2,[N],I;j,d''\big)$,
$d''_h\!=\!d_h'$ if $h\!\in\!\hat{I}$ 
and $d''_{\hat{0}}\!=\!0$.
By Corollaries~\ref{reg_crl2} and~\ref{reg_crl} , 
the space ${\cal H}_{{\cal T}'}$ 
is a smooth manifold of dimension 
$$\dim_{\Bbb{C}}{\cal H}_{{\cal T}'}=(d-1)(n+1)+n+|\hat{I}|+N.$$
Similarly to (1) above, it follows that ${\cal H}_{{\cal T}'}(\mu)$
is a smooth manifold of dimension $n\!-\!1\!+\!|\hat{I}|$ 
on which the $(2M\!+\!3)$-dimensional group 
$PSL_2\!\times\!{\cal G}_{\cal T}$ acts with only finite stabilizers.
It follows that ${\cal H}_{{\cal T}'}(\mu)\!=\!\emptyset$ if 
$n\!<\!M\!+\!4$. 
Note that the case $\hat{I}\!=\!\eset$ can occur only if 
$d\!=\!d_{\hat{0}}\!=\!2$.
Finally, if $d_{\hat{0}}\!=\!1$, the entire space ${\cal H}_{\cal T}$
is empty, since there are no holomorphic degree-one maps
from $\Si$ into $\P$; \hbox{see \cite{ACGH}}.\\
(3) Suppose $b=\big(\Si,[N],\{\hat{0}\};,(\hat{0},y),u\big)\!\in\!
{\cal H}_{(\Si,[N],\{\hat{0}\};\hat{0};d)}(\mu)$ and 
$u\!:\Si\!\lra\!\P$ factors through a $k$-fold cover of $S^2$,
where $k\!\ge\! 2$ and $k$ divides $d$. 
Then $b$ arises from the space ${\cal H}_{{\cal T}'}(\mu)$, where
$${\cal T}'=\big(S^2,[N],\{\hat{0}\};\hat{0},d/k\big).$$
Similarly to the above,
this space is a smooth manifold of dimension
$$\big((d/k)(n+1)+n+N\big)-\big(d(n+1)-n+N\big)=-
\frac{k-1}{k}d(n+1)+2n.$$
Thus, ${\cal H}_{{\cal T}'}(\mu)=\emptyset$, provided $d\ge 3$.
In fact, since
${\cal H}_{{\cal T}'}(\mu)$ has a three-dimensional group of symmetry,
${\cal H}_{{\cal T}'}(\mu)\!=\!\emptyset$ unless $d\!=\!2$ 
and $n\!\ge\! 4$.\\
(4) Suppose $b$ is as in (3) and $u$ factors through a $k$-fold
cover of a torus $T$, where $k\!\ge\! 2$ and $k$ divides~$d$. 
Then $b$ arises from the space 
\begin{equation*}\begin{split}
&{\cal H}_{1,d/k}(\mu)\equiv \big\{
 (x,y_{[N]},u)\!: x\!\in\!\Bbb{C}\!-\!\Bbb{R};~ u\!:\Bbb{C}\!\lra\!\P,~ 
 \bar{\partial}u=0,~u_*[T]=\frac{d}{k}\la;\\
&\qquad\qquad\qquad 
u(z\!+\!a\!+\!bx)=u(z)~\forall x\!\in\!\Bbb{C},a,b\!\in\!\Bbb{Z};~
   y_l\!\in\!\Bbb{C},~u(y_l)\in\mu_l~\forall l\in[N]\big\}.
\end{split}\end{equation*}
Similarly to the above, Corollary~\ref{reg_crl} implies that 
${\cal H}_{1,d/k}(\mu)$ is a smooth space of dimension
$$\big((d/k)(n+1)+1+N\big)-
\big(d(n+1)-n+N\big)=
-\Big(\frac{k-1}{k}d-1\Big)(n+1)<1.$$
Since ${\cal H}_{1,d/k}(\mu)$ has a one-dimensional group of symmetries
of $\Bbb{C}$-translations, ${\cal H}_{1,d/k}(\mu)=\emptyset$.

\end{document}